\documentclass[12pt,twoside,reqno,openany]{amsart}

\pdfoutput=1

\usepackage{amsbsy,amscd,amsfonts,amsmath,amssymb,amsthm,color,
fancybox,fancyhdr,footmisc,graphics,graphicx,ifthen,mathrsfs,
multicol,multirow,nccrules,pdfpages,rotating,shuffle,stmaryrd,
textcomp,times,wasysym,wrapfig}

\usepackage[dvipsnames,svgnames,x11names]{xcolor}

\usepackage[all]{xy}
\usepackage[utf8]{inputenc}
\usepackage[T1]{fontenc}
\usepackage{mathtools}
\newtagform{EngelLie}[\scriptstyle]{$}{$}
\makeatletter\newcommand{\leqnomode}{\tagsleft@true}
\newcommand{\reqnomode}{\tagsleft@false}\makeatother

\newtheorem{Theorem}[equation]{Theorem}

\newtheorem{Proposition}[equation]{Proposition}

\newtheorem{Lemma}[equation]{Lemma}

\newtheorem{Assertion}[equation]{Assertion}

\newtheorem{Observation}[equation]{Observation}

\newtheorem{PropositionDefinition}[equation]{Proposition-D\'efinition}

\theoremstyle{definition}

\newtheorem{MainHypothesis}[equation]{Main Hypothesis}
\newtheorem{InductionHypothesis}[equation]{Induction Hypothesis}
\newtheorem{Definition}[equation]{D\'efinition}
\newtheorem{Definition-Notation}[equation]{D\'efinition-Notation}
\newtheorem{Notation}[equation]{Notation}

\newtheorem{Terminology}[equation]{Terminology}

\newtheorem{Example}[equation]{Example}
\newtheorem{Examples}[equation]{Examples}

\newtheorem{Question}[equation]{Question}

\newtheorem{Problem}[equation]{Problem}

\newcommand{\C}{\mathbb{C}}

\newcommand{\N}{\mathbb{N}}

\newcommand{\R}{\mathbb{R}}

\newcommand{\TT}{\text{\sc t}}

\newcommand{\aaux}{{\text{\usefont{T1}{qcs}{m}{sl}a}}}

\newcommand{\caux}{{\text{\usefont{T1}{qcs}{m}{sl}c}}}

\definecolor{blue}{cmyk}{1.,1.,0.,0.63}
\definecolor{red}{cmyk}{0.,1.,1.,0.63}
\definecolor{green}{cmyk}{1.,0.,1.,0.63}
\definecolor{black}{cmyk}{1.,1.,1.,1.}

\newcommand{\green}{\textcolor{green}}
\newcommand{\red}{\textcolor{red}}

\makeatletter
\renewcommand{\@fnsymbol}[1]
{\ensuremath{\ifcase#1\or $*$\or $**$\or $***$\or $****$\or $*****$
\else\@ctrerr\fi}}
\makeatother

\newcommand{\HEAD}[2]{%
\pagestyle{fancy}
\fancyhead[RO]{\tiny\sf\thepage}
\fancyhead[CO]{{\tiny\sf #1}}
\fancyhead[LE]{\tiny\sf\thepage}
\fancyhead[CE]{{\tiny\sf #2}}
\fancyfoot{}}

\numberwithin{equation}{section}

\newcommand{\Section}[1]{
\renewcommand{\thesection}{\bf\arabic{section}}
\section{#1}
\renewcommand{\thesection}{\arabic{section}}}

\newcommand{\SectionHead}[2]{
\Section{\bf #1}
\label{#2}
\HEAD{\ref{#2}.~{\sf 
#1}}{
Jo\"el {\sc Merker}, 
D\'epartement de Math\'ematiques d'Orsay, 
Universit\'e Paris-Saclay, France}}

\newcommand{\style}[1]{{\sf #1}}

\newcommand{\Aff}{\style{Aff}}

\renewcommand{\det}{\style{det}}

\renewcommand{\dim}{\style{dim}}

\renewcommand{\exp}{\style{exp}}

\newcommand{\GL}{\style{GL}}

\newcommand{\Hessian}{\style{Hessian}}

\newcommand{\ind}{\style{ind}}

\newcommand{\Lie}{\style{Lie}}

\renewcommand{\lim}{\style{lim}}

\renewcommand{\min}{\style{min}}

\newcommand{\rank}{\style{rank}}

\newcommand{\Saff}{\style{Saff}}

\newcommand{\Span}{\style{Span}}

\newcommand{\stab}{\style{stab}}

\newcommand{\Sym}{\style{Sym}}

\newcommand{\centersmallbullet}{{}_{{}^{{}^{
\scriptscriptstyle{\bullet\!}}}}}

\newcommand{\Hall}{\Hall}

\newcommand{\smallbullet}{{\scriptscriptstyle{\bullet}}}

\newcommand{\smallsum}[1]{
\underset{#1}{\raisebox{1pt}{$\sum$\,}}
}

\newcommand{\vf}{\vfill\end{document}}

\newcommand{\zero}[1]{\underline{#1}_{\red{\circ}}}

\makeatletter
\def\hlinethick#1{%
\noalign{\ifnum0=`}\fi\hrule \@height #1 \futurelet
\reserved@a\@xhline}
\makeatother

\newlength{\myskip}
\begin{document}

\setcounter{section}{0}

$\:$

\bigskip\bigskip


\begin{center}

{\large\bf Inexistence\footnotemark[1] of Non-Product Hessian Rank 1}
\label{inexistence-dim-5}

\medskip

{\large\bf Affinely Homogeneous Hypersurfaces $H^n \subset \R^{n+1}$}

\medskip

{\large\bf in Dimension $n \geqslant 5$}

\bigskip\bigskip

Joël~{\sc Merker}\footnotemark[2]

\end{center}\bigskip

\footnotetext[1]{\,
This research was supported
in part by the Polish National Science Centre (NCN) 
via the grant number 2018/29/B/ST1/02583,
and by the Norwegian Financial Mechanism
2014--2021 via the project registration number 2019/34/H/ST1/00636.}

\footnotetext[2]{\,\,
D\'epartement de Math\'ematiques d'Orsay,
CNRS, Universit\'e Paris-Saclay, 91405 Orsay Cedex,
France, {\bf joel.merker@universite-paris-saclay.fr}}

\begin{center}
\begin{minipage}[t]{12.5cm}
\parindent 0.53cm
\footnotesize
\noindent
{\sc Abstract}.
Equivalences under the affine group 
$\Aff(\R^3)$ of constant Hessian rank $1$ surfaces $S^2 \subset \R^3$, 
sometimes called {\sl parabolic}, were, among other
objects, studied by
Doubrov, Komrakov, Rabinovich,
Eastwood, Ezhov, Olver, Chen, Merker, Arnaldsson, Valiquette.
Especially, homogeneous models and algebras of differential
invariants in various branches have been fully understood.

{\sl Then what about higher dimensions?} We consider hypersurfaces
$H^n \subset \R^{n+1}$ graphed as $\big\{ u = F(x_1, \dots, x_n) 
\big\}$ whose Hessian matrix $\big( F_{x_i x_j} \big)$, 
a relative affine invariant, is, similarly, of constant rank $1$. 
{\sl Are there homogeneous models?}

Complete explorations were done by the author on a computer in 
dimensions $n = 2, 3, 4, 5, 6, 7$. The first, expected outcome,
was to obtain a complete
classification of homogeneous models in dimensions 
$n = 2, 3, 4$ (forthcoming article, case $n = 2$ already known).
The second, unexpected outcome, was that in dimensions $n = 5, 6, 7$,
{\em there are {\em no} affinely homogenous models!}

(Except those that are affinely equivalent to 
a product of $\R^m$ with a
homogeneous model in dimensions $2, 3, 4$.)

The present article establishes such a non-existence result
in every dimension $n \geqslant 5$, 
based on the production of a normal form for 
$\big\{ u = F(x_1, \dots, x_n) \big\}$,
under $\Aff(\R^{n+1})$ up to order $\leqslant n+5$,
valid in any dimension $n \geqslant 2$.

\end{minipage}
\end{center}

\SectionHead{Introduction}
{introduction-inexistence}

Let $n \geqslant 1$, let $x = (x_1, \dots, x_n) \in \R^n$,
$u \in \R$ and $y = (y_1, \dots, y_n) \in \R^n$, $v \in \R$.
We consider equivalences of local analytic hypersurfaces
$H^n \subset \R_{x,u}^{n+1}$ and $K^n \subset \R_{y,v}^{n+1}$
graphed as:
\[
u
\,=\,
F(x_1,\dots,x_n)
\ \ \ \ \ \ \ \ \ \ \ \ \ \ \ \ \ \ \ \
\text{and}
\ \ \ \ \ \ \ \ \ \ \ \ \ \ \ \ \ \ \ \
v
\,=\,
G(y_1,\dots,y_n),
\]
under {\sl affine transformations} $\Psi \colon \R^{n+1}
\longrightarrow \R^{n+1}$:
\leqnomode\usetagform{default}
\begin{align}
\label{affine-a-b-c-d-tau-sigma}
\aligned
y_1
&
\,=\,
a_{1,1}\,x_1
+\cdots+
a_{1,n}\,x_n
+
b_1\,u
+
\tau_1,
\\
\cdots
&
\cdots\cdots\cdots\cdots\cdots\cdots\cdots\cdots\cdots\cdots\cdots
\cdot\cdot
\\
y_n
&
\,=\,
a_{n,1}\,x_1
+\cdots+
a_{n,n}\,x_n
+
b_n\,u
+
\tau_n,
\\
v
&
\,=\,
c_1\,x_1
+\cdots+
c_n\,x_n
+
d\,u
+
\tau_0,
\endaligned
\end{align}
where the $(n+1) \times (n+1)$ linear-part
$\big( \begin{smallmatrix} a & b \\ c & d \end{smallmatrix} \big)$
matrix is invertible, {\em i.e.} belongs to 
$\GL(n+1,\R)$.
The collection of all these
transformations is the Lie transformation group
$\Aff(\R^{n+1})$. 

The Lie algebra $\mathfrak{aff} (\R^{n+1})$ of $\Aff(\R^{n+1})$
consists of the
vector fields:
\[
\aligned
L
&
\,=\,
T_1\,
\tfrac{\partial}{\partial x_1}
+\cdots+
T_n\,
\tfrac{\partial}{\partial x_n}
+
T_0\,
\tfrac{\partial}{\partial u}
\,+
\\
&
\ \ \ \ \
+
\Big(
A_{1,1}\,x_1
+\cdots+
A_{1,n}\,x_n
+
B_1\,u
\Big)\,
\tfrac{\partial}{\partial x_1}
\,+
\\
&
\ \ \ \ \
+
\cdots\cdots\cdots\cdots\cdots\cdots\cdots\cdots\cdots\cdots\cdots
\cdot
\,+
\\
&
\ \ \ \ \
+
\Big(
A_{n,1}\,x_1
+\cdots+
A_{n,n}\,x_n
+
B_n\,u
\Big)\,
\tfrac{\partial}{\partial x_n}
\,+
\\
&
\ \ \ \ \
+
\Big(
C_1\,x_1
+\cdots+
C_n\,x_n
+
D\,u
\Big)\,
\tfrac{\partial}{\partial u}.
\endaligned
\]

A fixed hypersurface $H^n \subset \R^{n+1}$ has affine
symmetry group the 
{\em local Lie group}\,\,---\,\,for 
background, {\em see}~{\cite[Chap.~3]{Lie-Merker-2015}}\,\,---\,\,:
\[
\Sym(H)
\,:=\,
\big\{
\Psi\in\Aff(\R^{n+1})
\colon\,
\Psi(H)
\subset
H
\big\},
\]
where "$\subset$" is understood up to shrinking $H$,
and where the transformations $\Psi$ are close to
the identity. Then $\Sym\,H$ has Lie algebra:
\[
\Lie\,\Sym(H)
\,=\,
\mathfrak{sym}(H)
\,:=\,
\big\{
L\colon\,
L\big\vert_H\,\,
\text{tangent to}\,\,
H
\big\}.
\]

Since all our considerations will be {\em local}, 
we can assume that everything takes place in some
neighborhood of a fixed point $p_0 \in H$. 
Neighborhood shrinking is allowed
(a finite number of times).

\begin{Definition}
The hypersurface $H$ is said to be 
{\sl (locally) affinely homogeneous} if:
\[
T_{p_0}H
\,=\,
\Span_\R\,
\big\{
L\big\vert_{p_0}
\colon\,
L
\in
\mathfrak{sym}(H)
\big\}.
\]
\end{Definition}

According to Lie theory, the $1$-parameter groups
$p \longmapsto \exp(t\,L)(p)$ stabilize $H$, and 
$\Sym(H)$ is then locally
transitive in a neighborhood of $p_0 \in H$.

\smallskip

The problem of classifying all affinely homogeneous 
$n$-dimensional local analytic 
smooth submanifolds $H^n \subset \R^{n+c}$ is
probably of infinite complexity.
Even for $n = 2 = c$, it is not terminated.

In the hypersurface case $c = 1$
and in dimension $n = 2$,
the classification 
was terminated two decades ago~{\cite{Abdalla-Dillen-Vrancken-1997,
Doubrov-Komrakov-Rabinovich-1996,
Eastwood-Ezhov-1999}} by Doubrov-Komrakov-Rabinovich,
by Eeastwood-Ezhov, 
{\em see} also~{\cite{Arnaldsson-Valiquette-2020,
Chen-Merker-2020}} for
a differential invariants perspective.

In dimension $n = 3$, and codimension $c = 1$, 
there is~{\cite{Doubrov-Komrakov-1998}}
in which all {\sl multiply 
transitive} models were classified, while
for the {\em special} affine subgroup
$\Saff(\R^{3+1}) \subset \Aff(\R^{3+1})$, there are
the (unpublished) Ph.D. thesis of
Marc Wermann~{\cite{Wermann-2001}},
and the works of 
Eastwood-Ezhov~{\cite{Eastwood-Ezhov-2001, Eastwood-Ezhov-2001-2}},
both complete.

Joint with Chen~{\cite{Chen-Merker-2019}}, 
the author has studied the so-called {\sl parabolic 
surfaces} $H^2 \subset \R^3$, those whose Hessian has
constant rank $1$ 
({\em see} also~{\cite{Chen-Merker-2020}}).
Somewhat analogous $5$-dimensional CR structures of dimension $5$
whose Levi form is of constant rank $1$ have
also been studied in~{\cite{Foo-Merker-Ta-2020, 
Foo-Merker-Nurowski-Ta-2021}}.

\begin{Problem}
{\em Study algebras of differential invariants 
and classify homogenous models of
constant Hessian rank $1$ hypersurfaces $H^n \subset \R^{n+1}$.}
\end{Problem}

A similar problem can be formulated in the context of
CR geometry, {\em cf.}~{\cite{Merker-2019}}.

In Winter 2021, using a computer, the author found
all affinely homogenous Hessian rank $1$ hypersurfaces
$H^n \subset \R^{n+1}$ in dimensions $n = 2$, $3$, $4$,
{\em cf.}~{\cite{Merker-2019}}, forthcoming. 

Then exploring dimensions $n = 5, 6, 7$, 
the author was surprised to realize
that there are {\em no} homogenous models,
except the degenerate ones obtained by taking a
product of $\R^m$ with a homogeneous hypersurface
$H^{n-m} \subset \R^{n-m+1}$ so that $2 \leqslant 
n - m \leqslant 4$. 

He then tackled to prove a {\em non-existence} result,
which, incidentally, provides a complete classification
(in the constant Hessian rank $1$ branch). 
But the computational task appeared to be unexpectedly hard,
and it took one year to write a detailed proof in general 
dimension $n \geqslant 5$.

The main result of this paper concerns dimensions 
$n \geqslant 5$, but several results
are true for all $n \geqslant 2$,
and will be useful to~{\cite{Merker-2022}}.

From the corpus of this article, we may emphasize three statements.
The first one appears as Theorem~{\ref{Thm-product-nH}} below.

\begin{Theorem}
\label{Thm-n-H-introduction}
Let $H^n \subset \R^{n+1}$ be a local affinely homogeneous
hypersurface
having constant Hessian rank $1$. Then there exists
an integer $1 \leqslant n_H \leqslant n$ and
affine coordinates $(x_1, \dots, x_n)$ in which:
\[
H^n
\,=\,
H^{n_H}
\times
\R_{x_{n_H+1},\dots,x_n}^{n-n_H-1}
\]
is a product of an affinely homogeneous
hypersurface $H^{n_H} \subset \R^{n_H + 1}$
times a `dumb' $\R^{n-n_H-1}$,
and is graphed as:
\[
\aligned
u
&
\,=\,
\frac{x_1^2}{2}
+
\frac{x_1^2x_2}{2}
+
\sum_{m=3}^{n_H}\,
\Big(
\frac{x_1^mx_m}{m!}
+
x_1^{m-1}
\sum_{i,j\geqslant 2
\atop
i+j=m+1}\,
\tfrac{1}{2}\,
\frac{x_ix_j}{(i-1)!(j-1)!}
+
{\rm O}_{x_2,\dots,x_{m-1}}(3)
\Big)
\\
& 
\ \ \ \ \ \ \ \ \ \ \ \ \ \ \ \ \ \ \ \ \ \ \ \ \ 
+
\sum_{m=n_H+2}^\infty\,
E^m(x_1,\dots,x_{n_H}),
\endaligned
\]
with graphing function $F = F(x_1, \dots, x_{n_H})$ 
independent of $x_{n_H+1}, \dots, x_n$.
\end{Theorem}

Notice that the variables $\big(x_1, \dots, x_{n_H}, u\big)$ 
{\em are present} in such a graphed equation.

Of course, we are interested in the hypersurfaces for which $n_H = 
n$. Such hypersurfaces can be called {\sl nondegenerate}, 
but we will not use such a terminology.

With $n_H = n$, Theorem~{\ref{Thm-n-H-introduction}}
shows the graphing function up to order $n+1$ included.
Up to order $n+3$ included, we prove

\begin{Theorem}
\label{Thm-nf-n-3-intro}
In any dimension $n \geqslant 2$, 
every local hypersurface $H^n \subset
\R^{n+1}$ having constant Hessian rank $1$,
and which is not affinely equivalent to a product 
of $\R^m$ $(1 \leqslant m \leqslant n)$
with a hypersurface $H^{n-m} \subset \R^{n-m+1}$, 
can be affinely normalized up to order $n+3$ as:
\[
\footnotesize
\aligned
u
&
\,=\,
\frac{x_1^2}{2}
+
\frac{x_1^2\,x_2}{2}
+
\sum_{m=3}^n\,
\Big(
\frac{x_1^m\,x_m}{m!}
+
x_1^{m-1}
\sum_{i,j\geqslant 2
\atop
i+j=m+1}\,
\tfrac{1}{2}\,
\frac{x_i\,x_j}{(i-1)!(j-1)!}
\Big)
\\
&
\ \ \ \ \
+
F_{n+1,10\cdots0}\,
\frac{x_1^{n+1}x_2}{(n+1)!}
+
x_1^n
\sum_{i,j\geqslant 2
\atop
i+j=n+2}\,
\tfrac{1}{2}\,
\frac{x_i\,x_j}{(i-1)!(j-1)!}
\\
&
\ \ \ \ \
+
F_{n+3,0\cdots0}\,
\frac{x_1^{n+3}}{(n+3)!}
+
F_{n+2,10\cdots0}\,
\frac{x_1^{n+2}\,x_2}{(n+2)!}
+
F_{n+2,0010\cdots0}\,
\frac{x_1^{n+2}\,x_4}{(n+2)!}
+\cdots+
F_{n+2,0\cdots01}\,
\frac{x_1^{n+2}\,x_n}{(n+2)!}
\\
&
\ \ \ \ \ \ \ \ \ \ \ \ \ \ \ \ \ \ \ \ \ \ \ \ \ \ \ \ \ \
+
F_{n+1,10\cdots0}\,
\frac{x_1^{n+1}\,x_2\,x_2}{n!}
+
x_1^{n+1}
\sum_{i,j\geqslant 2
\atop
i+j=n+3}\,
\tfrac{1}{2}\,
\frac{x_i\,x_j}{(i-1)!(j-1)!}
\\
&
\ \ \ \ \ 
+
{\rm O}_{x_2,\dots,x_n}(3)
+
{\rm O}_{x_1,x_2,\dots,x_n}(n+4).
\endaligned
\]

Furthermore, linear $\GL(n+1, \R)$ 
self-maps (fixing the origin)
$\big( \begin{smallmatrix} y \\ v \end{smallmatrix} \big)
=
\big( \begin{smallmatrix} a & b \\ c & d \end{smallmatrix} \big)\,
\big( \begin{smallmatrix} x \\ u \end{smallmatrix} \big)$
of such a hypersurface 
are necessarily weigted dilations of the form:
\[
y_1
\,=\,
\caux\,x_1,
\ \ \ \ \ \
y_2
\,=\,
0,
\ \ \ \ \ \
y_3
\,=\,
\tfrac{1}{\caux}\,
x_3,
\ \ \ \ \ \
\dots,
\ \ \ \ \ \
y_n
\,=\,
\tfrac{1}{\caux^{n-2}}\,
x_n,
\ \ \ \ \ \
v
\,=\,
\caux^2\,u,
\]
with $\caux \in \R^\ast$.
\end{Theorem}

In other words, the isotropy is at most one-dimensional.
The more advanced Theorem~{\ref{Thm-nf-order-n-5}}
gives terms of orders $n+4$, $n+5$, which are more
complicated, but (unfortunately) necessary in order to establish
our (unexpected) main 

\begin{Theorem}
\label{Thm-inexistence-n-5}
In any dimension $n \geqslant 5$, 
there are {\em no}
affinely homogeneous constant Hessian rank 1
nondegenerate hypersurfaces $H^n \subset \R^{n+1}$.
\end{Theorem}

Here is the key reason why homogeneous models do not exist
when $n \geqslant 5$.
In Sections~{\ref{summary-proof-theorem-orders-n-3-n-4}}
and~{\ref{tangency-equations-orders-n-4-n-5}},
we will obtain the two equations 
\green{\bf I} and \green{\bf II} 
shown in Proposition~{\ref{Prp-equations-I-II}},
that will
contradict homogeneity, at the infinitesimal level.

Readers could admit 
Theorem~{\ref{Thm-nf-n-3-intro}} or 
Theorem~{\ref{Thm-nf-order-n-5}},
and go directly to 
Sections~{\ref{summary-proof-theorem-orders-n-3-n-4}}
and~{\ref{tangency-equations-orders-n-4-n-5}}, 
which are simple to read.
Unfortunately, the core normalizations
done in the previous sections are hard, 
require patience, and indurance. 

All computations were done on 
a computer fully in dimensions $n = 2, 3, 4, 5, 6, 7$,
during $>$ 2 months of exploration, 
from December 2020, to February 2021.
Especially, all technical statements of this article
were {\em constantly checked to be true} on a computer in
dimensions $n = 5, 6, 7$. This acted as guide 
to set up the general dimension by hand. 

In fact, it happened to be unexpectedly hard to write
by hand a detailed
proof in general dimension\,\,---\,\,the 
computer being unable to do that!
It was really necessary to normalize
$\big\{u = F(x_1, \dots, x_n) \big\}$ up to order $n+5$
as stated in Theorem~{\ref{Thm-nf-order-n-5}},
because no contradiction occured
in lower order $\leqslant n+4$ for the 
`concrete' dimensions $n = 5, 6, 7$.

With slightly harder computations, 
it can be shown that in any dimension $n \geqslant 2$:
\[
\dim\,
\Lie\,\Sym(H)
\,\leqslant\,
4,
\]
always. Thus, when $\dim\, H = n \geqslant 5$,
(infinitesimal) homogeneity cannot take place.

All statements hold for $\C$ instead of $\R$,
with the same proofs.

\SectionHead{Hessian Matrix and its Rank Invariancy}
{Hessian-rank}

We may assume that coordinates $(x,u)$ are centered at $p_0 \in M$, so
that $p_0$ is the origin. Denote its image 
by $q_0 := \Psi(p_0)$. We may
also assume that $q_0$ is
the origin in $\R_{y,v}^{n+1}$ too.
Then $\Psi(0) = 0$ forces $0 = \tau_1 = \cdots = \tau_n = \tau_0$
in~({\ref{affine-a-b-c-d-tau-sigma}}), which 
becomes a general $\GL(\R^{n+1})$ transformation:
\leqnomode\usetagform{default}
\begin{align}
\label{linear-a-b-c-d}
\left[
\begin{array}{c}
y_1
\\
\vdots
\\
y_n
\\
v
\end{array}
\right]
\,=\,
\left[
\begin{array}{cccc}
a_{1,1} & \cdots & a_{1,n} & b_1
\\
\vdots & \ddots & \vdots & \vdots
\\
a_{n,1} & \cdots & a_{n,n} & b_1
\\
c_1 & \cdots & c_n & d
\end{array}
\right]\,
\left[
\begin{array}{c}
x_1
\\
\vdots
\\
x_n
\\
u
\end{array}
\right].
\end{align}

Importantly, the basic hypothesis that $\Psi(H) \subset K$, namely
that $0 = -\, v + G(y)$ when $\big( \begin{smallmatrix}
y \\ v \end{smallmatrix} \big) = \big( \begin{smallmatrix}
a & b \\ c & d \end{smallmatrix} \big) \big( \begin{smallmatrix}
x \\ u \end{smallmatrix} \big)$
as written above {\em and when $u$ is replaced by $F(x)$}, 
expresses as the {\sl fundamental equation:}
\leqnomode\usetagform{default}
\begin{align}
\label{eqFG}
\aligned
0
&
\,=\,
-\,
c_1\,x_1
-\cdots-
c_n\,x_n
-
d\,F(x_1,\dots,x_n)
\\
&
\ \ \ \ \
+
G
\Big(
a_{1,1}x_1+\cdots+a_{1,n}x_n+b_1F(x_1,\dots,x_n),\,
\dots\dots,\,
a_{n,1}x_1+\cdots+a_{n,n}x_n+b_nF(x_1,\dots,x_n)
\Big),
\endaligned
\end{align}
which holds identically in $\C\{x_1, \dots, x_n\}$.

One step further, by composing $\Psi$ with two elementary affine
transformations, we may assume that the tangent space $T_0 H^n = \{ u
= 0\}$ is horizontal, and that $T_0 K^n = \{ v = 0\}$ as well. In
other words, $F = {\rm O}_x(2)$ and $G = {\rm O}_y(2)$.

Thus, order $0$ and order $1$ terms are absent in the expansions:
\[
u
\,=\,
F(x)
\,=\,
\sum_{\sigma_1\geqslant 0,\dots,\sigma_n\geqslant 0
\atop
\sigma_1+\cdots+\sigma_n\geqslant 2}\,
\frac{x_1^{\sigma_1}\cdots x_n^{\sigma_n}}{
\sigma_1!\,\cdots\,\sigma_n!}\,
F_{x_1^{\sigma_1}\cdots x_n^{\sigma_n}}(0)
\ \ \ \ \ 
\text{and}
\ \ \ \ \ 
v
\,=\,
G(y)
\,=\,
\sum_{\tau_1\geqslant 0,\dots,\tau_n\geqslant 0
\atop
\tau_1+\cdots+\tau_n\geqslant 2}\,
\frac{y_1^{\tau_1}\cdots y_n^{\tau_n}}{\tau_1!\,\cdots\,\tau_n!}\,
G_{y_1^{\tau_1}\cdots y_n^{\tau_n}}(0).
\]

\begin{Lemma}
\label{Lm-c1-cdots-cn-zero}
The linear transformation $\big( \begin{smallmatrix}
y \\ v \end{smallmatrix} \big) = \big( \begin{smallmatrix}
a & b \\ c & d \end{smallmatrix} \big) \big( \begin{smallmatrix}
x \\ u \end{smallmatrix} \big)$ sends $u = {\rm O}_x(2)$ 
to $v = {\rm O}_y(2)$ if and only if $0 = c_1 = \cdots = 
c_n$.
\end{Lemma}

\proof
Read~({\ref{eqFG}}) modulo ${\rm O}_x(2)$:
\[
0
\,\equiv\,
-\,c_1\,x_1
-\cdots-
c_n\,x_n
-
{\rm O}_x(2)
+
{\rm O}_x(2).
\qedhere
\]
\endproof

This interprets as {\sl group reduction:}
\[
\left[
\begin{array}{cccc}
a_{1,1} & \cdots & a_{1,n} & b_1
\\
\vdots & \ddots & \vdots & \vdots
\\
a_{n,1} & \cdots & a_{n,n} & b_n
\\
c_1 & \cdots & c_n & d
\end{array}
\right]^{\green{\bf 0}}
\,\,\,\leadsto\,\,\,
\left[
\begin{array}{cccc}
a_{1,1} & \cdots & a_{1,n} & b_1
\\
\vdots & \ddots & \vdots & \vdots
\\
a_{n,1} & \cdots & a_{n,n} & b_n
\\
\red{\bf 0} & \cdots & \red{\bf 0} & d
\end{array}
\right]^{\green{\bf 1}}.
\]
Necessarily:
\[
\det\,\big(a_{i,j}\big)
\,\neq\,
0
\,\neq\,
d.
\]

Later, we will need an equivalent
expansion gathering homogeneous terms
of fixed order,
whose proof is left to the reader:
\leqnomode\usetagform{default}
\begin{align}
\label{F-homogeneous-m}
u
\,=\,
\sum_{m\geqslant 2}\,
\tfrac{1}{m!}\,
\sum_{i_1=1}^n\,\cdots\,\sum_{i_m=1}^n\,
x_{i_1}\cdots x_{i_m}\,
F_{x_{i_1}\cdots x_{i_m}}(0)
\,\,=:\,
\sum_{m\geqslant 2}\,
F^m(x).
\end{align}

Next, at order $2$:
\[
u
\,=\,
\tfrac{1}{2}\,
\sum_{i_1=1}^n\,\sum_{i_2=1}^n\,
x_{i_1}x_{i_2}\,f_{i_1,i_2}
+
{\rm O}_x(3)
\ \ \ \ \ \ \ \ \ \ \ \ \ \ \ \ \ \ \ \
\text{and}
\ \ \ \ \ \ \ \ \ \ \ \ \ \ \ \ \ \ \ \
v
\,=\,
\tfrac{1}{2}\,
\sum_{j_1=1}^n\,\sum_{j_2=1}^n\,
y_{j_1}y_{j_2}\,g_{j_1,j_2}
+
{\rm O}_y(3),
\]
where:
\[
f_{i_1,i_2}
\,:=\,
F_{x_{i_1}x_{i_2}}(0)
\,=\,
f_{i_2,i_1}
\ \ \ \ \ \ \ \ \ \ \ \ \ \ \ \ \ \ \ \
\text{and}
\ \ \ \ \ \ \ \ \ \ \ \ \ \ \ \ \ \ \ \
g_{j_1,j_2}
\,:=\,
G_{y_{j_1}y_{j_2}}(0)
\,=\,
g_{j_2,j_1}.
\]
With $x = {}^{\TT} (x_1, \dots, x_n)$ being a column vector,
we may abbreviate:
\[
u
\,=\,
\tfrac{1}{2}\,
{}^{\TT}x\cdot F^{(2)}\cdot x
+
{\rm O}_x(3)
\ \ \ \ \ \ \ \ \ \ \ \ \ \ \ \ \ \ \ \
\text{and}
\ \ \ \ \ \ \ \ \ \ \ \ \ \ \ \ \ \ \ \
v
\,=\,
\tfrac{1}{2}\,
{}^{\TT}y\cdot G^{(2)}\cdot y
+
{\rm O}_y(3).
\]

\begin{Definition}
\label{Def-Hessian-0}
At any point $p \in H^n$, the
{\sl Hessian matrix} is defined, in any system
of coordinates $(x,u)$ centered at $p = (0, 0)$ 
in which $H^n$ is graphed as $u = F(x)$ with
$0 = F(0) = F_{x_1}(0) = \cdots = F_{x_n}(0)$, as:
\[
\Hessian(F)
\,:=\,
\left[
\begin{array}{ccc}
F_{x_1x_1}(0) & \cdots & F_{x_1x_n}(0)
\\
\vdots & \ddots & \vdots
\\
F_{x_nx_1}(0) & \cdots & F_{x_nx_n}(0)
\end{array}
\right]
\,=\,
\left[
\begin{array}{ccc}
f_{1,1} & \cdots & f_{1,n}
\\
\vdots & \ddots & \vdots
\\
f_{n,1} & \cdots & f_{n,n}
\end{array}
\right]
\,=\,
F^{(2)}.
\]
\end{Definition}

At first sight, this definition depends on coordinates,
but a key invariancy lies behind.

\begin{Lemma}
The {\em rank} of the Hessian matrix $\Hessian(F)$ is independent of
the affine coordinates $(x_1, \dots, x_n, u)$ 
centered at $p$
in which $0 = F(0) =
F_{x_1}(0) = \cdots = F_{x_n}(0)$.
\end{Lemma}

\proof
Take another system of coordinates
$(y_1, \dots, y_n, v)$ centered at $p = 0$ in which
the hypersurface is graphed as $v = G(y)$ with
$0 = G(0) = G_{y_1}(0) = \cdots = G_{y_n}(0)$,
namely:
\[
v
\,=\,
{}^{\TT}y
\cdot
G^{(2)}
\cdot
y
+
{\rm O}_y(3).
\]
So there is an invertible linear transformation
of the form~({\ref{linear-a-b-c-d}})
with $0 = c_1 = \cdots = c_n$ by Lemma~{\ref{Lm-c1-cdots-cn-zero}},
namely $y = A\cdot x + B\cdot u$ and $v = d\,u$,
which sends $u = F(x)$ to $v = G(y)$.

Then read the fundamental equation~({\ref{eqFG}})
modulo ${\rm O}_x(3)$, compute,
observe that ${\rm O}_y(3) = 
{\rm O}_x(3)$ when $y = A\cdot x + B\cdot F(x) = 
{\rm O}_x(1)$, factorize:
\[
\aligned
0
&
\,=\,
-\,v
+
G(y)
\\
&
\,\equiv\,
-\,d\,F(x)
+
G\big(
A\cdot x
+
B\cdot F(x)
\big)
\\
&
\,\equiv\,
-\,d\,
\tfrac{1}{2}\,
{}^{\TT}x\cdot F^{(2)}\cdot x
-
{\rm O}_x(3)
+
\tfrac{1}{2}\,
{}^{\TT}
\big[
A\cdot x
+
B\cdot F(x)
\big]
\cdot
G^{(2)}
\cdot
\big[
A\cdot
x
+
B\cdot F(x)
\big]
+
{\rm O}_y(3)
\\
&
\,=\,
\tfrac{1}{2}\,
{}^{\TT}x
\cdot
\Big(
-\,d\,F^{(2)}
+
{}^{\TT}A
\cdot
G^{(2)}
\cdot
A
\Big)
\cdot
x
+
{\rm O}_x(3),
\endaligned
\]
and deduce, since $x \in \R^n$ is arbitrary, that:
\[
d\,F^{(2)}
\,=\,
{}^{\TT}A
\cdot
G^{(2)}
\cdot
A.
\]
Lastly, from $d \neq 0 \neq \det\, A$, conclude
that $\rank\, F^{(2)} = \rank\, G^{(2)}$.
\endproof

Now, suppose that $H = \{u = F(x)\}$ is given
in coordinates $(x,u)$ centered at $p_0 = 0$,
with $0 = F(0) = F_{x_1}(0) = \cdots = F_{x_n}(0)$.
At any other point $p \sim p_0$ close to the origin
with $p = (x_{1,p}, \dots, x_{n,p}, u_p)$, where
$u_p = F(x_p)$, 
use the centered coordinates:
\[
y_1
\,:=\,
x_1-x_{1,p},\,\,
\dots\dots\,\,
y_n
\,:=\,
x_n-x_{n,p},\,\,
v
\,:=\,
u-F(x_p)
-
F_{x_1}(x_p)\big(x_1-x_{1,p}\big)
-\cdots-
F_{x_n}(x_p)\big(x_n-x_{n,p}\big).
\]
The new graphing function:
\[
G(y)
\,:=\,
F(x_p+y)
-
F(x_p)
-
F_{x_1}(x_p)\,y_1
-\cdots-
F_{x_n}(x_p)\,y_n,
\]
satisfies $G(0) = G_{y_1}(0) = \cdots = G_{y_n}(0)$,
hence Definition~{\ref{Def-Hessian-0}}
applies.
But since the correction terms are of degree $\leqslant 1$
in $y_1, \dots, y_n$, they disappear after second order
differentiation:
\[
\left[
\begin{array}{ccc}
G_{y_1y_1}(0) & \cdots & G_{y_1y_n}(0)
\\
\vdots & \ddots & \vdots
\\
G_{y_ny_1}(0) & \cdots & G_{y_ny_n}(0)
\end{array}
\right]
\,=\,
\left[
\begin{array}{ccc}
F_{x_1x_1}(x_p) & \cdots & F_{x_1x_n}(x_p)
\\
\vdots & \ddots & \vdots
\\
F_{x_nx_1}(x_p) & \cdots & F_{x_nx_n}(x_p)
\end{array}
\right].
\]

\begin{PropositionDefinition}
The Hessian of a graphed hypersurface $u = F(x_1, \dots, x_n)$
is a well defined matrix-valued function of $x$:
\[
\Hessian(F)(x)
\,:=\,
\left[
\begin{array}{ccc}
F_{x_1x_1}(x) & \cdots & F_{x_1x_n}(x)
\\
\vdots & \ddots & \vdots
\\
F_{x_nx_1}(x) & \cdots & F_{x_nx_n}(x)
\end{array}
\right]
\]
whose rank is invariant under affine equivalences
(at pairs of points corresponding one to another).\qed
\end{PropositionDefinition}

\SectionHead{Constant Hessian Rank $1$}
{constant-Hessian-rank-1}

Take a hypersurface $H^n \subset \R^{n+1}$ with
$0 \in H^n$ and $T_0 H^n = \{ u = 0\}$, graphed as:
\[
u
\,=\,
F(x)
\,=\,
\tfrac{1}{2}\,
\sum_{i=1}^n\,
\sum_{j=1}^n\,
x_ix_j\,
F_{x_ix_j}(0)
+
{\rm O}_x(3).
\]
Its Hessian at the origin $0$ is represented by the
$n \times n$ matrix $\big( F_{x_i x_j}(0) \big)$.

\begin{Lemma}
Whenever $\big( F_{x_i x_j}(0) \big)$ is not the zero matrix,
there exists an affine change of coordinates $y = A\,x$,
$v = u$, making nonzero (after renaming $y =: x$):
\[
F_{x_1x_1}(0)
\,\neq\,
0.
\]
\end{Lemma}

\proof
If there exists an index $i_\ast$ with $F_{x_{i_\ast} 
x_{i_\ast}}(0) \neq 0$, simply permute affine coordinates
to set $x_1 := x_{i_\ast}$.

Assume therefore that $0 = F_{x_i x_i} (0)$ for all
$i = 1, \dots, n$. By assumption, there exist
$i_\ast$ and $j_\ast \neq i_\ast$ such that
$F_{x_{i_\ast} x_{j_\ast}}(0) \neq 0$. 
Rename $x_1 := x_{i_\ast}$ to have $0 \neq F_{x_1 x_{j_\ast}} (0)$.
Change $j_\ast \geqslant 2$ to be the smallest
satisfying $F_{x_1 x_{j_\ast}}(0) \neq 0$.
Also, abbreviate $f_{i,j} := F_{x_i x_j}(0)$.
So $f_{1,j_\ast} \neq 0$.

Since diagonal terms are all zero, only $\sum_{i < j}$ remains,
and 
we may expand,
modulo ${\rm O}_x(3)$:

\[
\aligned
u
\,\equiv\,
\sum_{i<j}\,
x_ix_j\,f_{i,j}
&
\,=\,
\sum_{
i<j
\atop
j\leqslant j_\ast}\,
x_ix_j\,f_{i,j}
+
\sum_{
i<j
\atop
j_\ast<j}\,
x_ix_j\,f_{i,j}
\\
&
\,=\,
x_1x_{j_\ast}\,f_{1,j_\ast}
+
\sum_{j_\ast<j}\,
x_1x_j\,f_{1,j}
+
\sum_{
2\leqslant i<j
\atop
j\leqslant j_\ast}\,
x_ix_j\,f_{i,j}
+
\sum_{
2\leqslant i<j
\atop
j_\ast<j}\,
x_ix_j\,f_{i,j}.
\endaligned
\]

Now, set $x_{j_\ast} := y_1 + y_{j_\ast}$ while
$x_j := y_j$ for $j \neq j_\ast$ and $u := v$,
whence ${\rm O}_x(3) = {\rm O}_y(3)$,
so that the first monomial becomes:
\[
y_1\,\big(y_1+y_{j_\ast}\big)\,
f_{1,j_\ast},
\]
with nonzero coefficient for $y_1y_1$. 
The three remaining sums cannot
incorporate $y_1 y_1$. Thus, the new graph $v = G(y)$
satisfies $G_{y_1y_1}(0) = 2\, f_{1,j_\ast} \neq 0$.
\endproof

\SectionHead{Independent and Border-Dependent Jets}
{independent-border-dependent-jets}

Up to the end of the article, we will assume:
\[
F_{x_1x_1}(0)
\,\neq\,
0.
\]
Also, our main hypothesis will be that 
the Hessian matrix $\big(F_{x_ix_j}(x) \big)$ has
constant rank $1$ for all $x \sim 0$ in some neighborhood 
of the origin.

It is elementary to verify that, for constants
$\varphi_{i,j} \in \R$ with $\varphi_{1,1} \neq 0$:
\[
1
\,=\,
\rank\,
\left[
\begin{array}{cccc}
\varphi_{1,1} & \varphi_{1,2} & \cdots & \varphi_{1,n}
\\
\varphi_{2,1} & \varphi_{2,2} & \cdots & \varphi_{2,n}
\\
\vdots & \vdots & \ddots & \vdots
\\
\varphi_{n,1} & \varphi_{n,2} & \cdots & \varphi_{n,n}
\end{array}
\right]
\ \ \ \ \
\Longleftrightarrow
\ \ \ \ \
0
\,=\,
\left\vert
\begin{array}{cc}
\varphi_{1,1} & \varphi_{1,j}
\\
\varphi_{i,1} & \varphi_{i,j}
\end{array}
\right\vert
\ \ \ \ \
\forall\,\,2\,\leqslant\,i,j\,\leqslant\,n.
\]

\begin{MainHypothesis}
\label{Main-Hypothesis}
{\sl For all $x \sim 0$
in some neighborhood of the origin:}
\leqnomode\usetagform{default}
\begin{align}
\label{F-Hrank-1}
F_{x_ix_j}(x)
\,\equiv\,
\frac{F_{x_1x_i}(x)\,F_{x_1x_j}(x)}{
F_{x_1x_1}(x)}
\ \ \ \ \ \ \ \ \ \ \ \ \ \ \ \ \ \ \ \
\forall\,\,
i,j=2,\dots,n.
\end{align}
\end{MainHypothesis}

By differentiating this identity with respect to $x_1, x_2, \dots,
x_n$, it is not difficult to prove by induction that
every derivative $F_{x_1^\tau x_2^{i_2} \cdots x_n^{i_n}} (x)$
with $i_2 + \cdots + i_n \geqslant 2$
and $\tau \in \N$ arbitrary,
expresses as a polynomial
in the derivatives
$F_{x_1^{i_1'} x_2^{i_2'} \cdots x_n^{i_n'}} (x)$
with $i_1' \leqslant 1$, divided by
a certain power $\big( F_{x_1x_1} (x) \big)^\ast$.
We will in fact need to know
what formulas hold only for $i_2 + \cdots + i_n = 2$.

\begin{Terminology}
\label{Terminology-border-dependent-jets}
$\bullet$\,
The {\sl independent jets} are the derivatives with $i_1 = 0$ 
or $i_1 = 1$:
\[
F_{x_2^{i_2}\cdots x_n^{i_n}}(x),
\ \ \ \ \ \ \ \ \ \ \ \ \ \ \
F_{x_1x_2^{i_2}\cdots x_n^{i_n}}(x).
\]

\smallskip\noindent$\bullet$\,
The {\sl border-dependent jets} are the derivatives:
\[
F_{x_1^\tau x_ix_j}(x),
\ \ \ \ \ \ \ \ \ 
\tau\,\geqslant\,0, 
\ \ \ \ \
i,j=2,\dots,n.
\]
\end{Terminology}

The way how these
border-dependent jets $F_{x_1^\tau x_i x_j}$ express in terms 
of the independent jets can be seen by just 
differentiating~({\ref{F-Hrank-1}}) $\nu$ times
with respect to $x_1$.

Later, the following abbreviations will be used:
\[
\aligned
x'
&
\,:=\,
(x_2,\dots,x_n),
\\
f_{i_1,i_2,\dots,i_n}
&
\,:=\,
F_{x_1^{i_1}x_2^{i_2}\cdots x_n^{i_n}}(0).
\endaligned
\]

\SectionHead{Normalization at Order $2$}
{normalization-order-2}

Thus with $f_{1,1} \neq 0$, start from:
\[
u
\,=\,
F(x)
\,=\,
\tfrac{1}{2}
\sum_i\,\sum_j\,
f_{i,j}\,x_ix_j
+
{\rm O}_x(3).
\]
Changing $u$ to $-u$ if needed,
we can assume $f_{1,1} > 0$.
The plain dilation $y_1 := 
\sqrt{f_{1,1}}\, x_1$
makes $f_{1,1} = 1$,
namely, using again the letters $x$, $u$:
\[
u
\,=\,
\tfrac{1}{2}\,
x_1^2
+
\smallsum{2\leqslant j\leqslant n}\,
x_1x_j\,f_{1,j}
+
\tfrac{1}{2}\,
\smallsum{2\leqslant i,j\leqslant n}\,
x_ix_j\,f_{i,j}
+
{\rm O}_x(3).
\]
Remind that our Main Hypothesis~{\ref{Main-Hypothesis}}
is (implicitly) assumed in all statements.

\begin{Assertion}
There exists an affine change of coordinates which normalizes
$u = F(x)$ to:
\[
u
\,=\,
\tfrac{1}{2}\,
x_1^2
+
{\rm O}_x(3).
\]
\end{Assertion}

\proof
First, absorb all monomials incorporating $x_1$:
\[
u
\,=\,
\tfrac{1}{2}\,
\Big(
\underbrace{
x_1
+
\smallsum{2\leqslant j\leqslant n}\,
x_j\,f_{1,j}}_{
=:\,\,\text{\tiny\rm new}\,\,x_1}
\Big)^2
-
\tfrac{1}{2}\,
\Big(
\smallsum{2\leqslant j\leqslant n}\,
x_j\,f_{1,j}
\Big)^2
+
\tfrac{1}{2}\,
\smallsum{2\leqslant i,j\leqslant n}\,
x_ix_j\,f_{i,j}
+
{\rm O}_x(3),
\]
to get, with modified $f_{i,j}$:
\[
u
\,=\,
\tfrac{1}{2}\,
x_1^2
+
\tfrac{1}{2}\,
\smallsum{2\leqslant i,j\leqslant n}\,
x_ix_j\,f_{i,j}
+
{\rm O}_x(3).
\]

This means that:
\[
0
\,=\,
F_{x_1x_j}(0)
\eqno
{\scriptstyle{(\forall\,\,2\,\leqslant\,j\,\leqslant\,n)}}.
\]
Lastly, ({\ref{F-Hrank-1}}) yields:
\[
F_{x_ix_j}(0)
\,=\,
0
\ \ \ \ \ \ \ \ \ \ \ \ \ \ \ \ \ \ \ \
{\scriptstyle{(\forall\,\,2\,\leqslant\,i,\,j\,\leqslant\,n)}}.
\qedhere
\]
\endproof

Next, starting from $u = \tfrac{1}{2}\, x_1^2 + {\rm O}_x(3)$, the
goal is to normalize order $3$ terms. But before increasing the order,
we must {\sl conserve} or {\sl stabilize} the normalization at order
$2$. This conducts to {\sl group reduction}:
\[
\footnotesize
\aligned
\overset{\GL(n+1,\R)}{
\left[
\begin{array}{ccccc}
a_{1,1} & a_{1,2} & \cdots & a_{1,n} & b_1
\\
a_{2,1} & a_{2,2} & \cdots & a_{2,n} & b_2
\\
\vdots & \ddots & \vdots & \vdots
\\
a_{n,1} & a_{n,2} & \cdots & a_{n,n} & b_n
\\
c_1 & c_2 & \cdots & c_n & d
\end{array}
\right]^{\green{\bf 0}}}
\,\,\,\leadsto\,\,\,
\overset{G_\stab^1}{
\left[
\begin{array}{ccccc}
a_{1,1} & a_{1,2} & \cdots & a_{1,n} & b_1
\\
a_{2,1} & a_{2,2} & \cdots & a_{2,n} & b_2
\\
\vdots & \ddots & \vdots & \vdots
\\
a_{n,1} & a_{n,2} & \cdots & a_{n,n} & b_n
\\
\red{\bf 0} & \red{\bf 0} & \cdots & \red{\bf 0} & d
\end{array}
\right]^{\green{\bf 1}}}
\,\,\,\leadsto\,\,\,
\overset{G_\stab^2}{
\left[
\begin{array}{ccccc}
a_{1,1} & \red{\bf 0} & \cdots & \red{\bf 0} & b_1
\\
a_{2,1} & a_{2,2} & \cdots & a_{2,n} & b_2
\\
\vdots & \ddots & \vdots & \vdots
\\
a_{n,1} & a_{n,2} & \cdots & a_{n,n} & b_n
\\
\red{\bf 0} & \red{\bf 0} & \cdots & \red{\bf 0} & 
\red{a_{1,1}^2}
\end{array}
\right]^{\green{\bf 2}}}.
\endaligned
\]

\begin{Lemma}
The subgroup $G_\stab^2$ of $G_\stab^1$ which 
sends $u = \frac{1}{2}\, x_1^2 + {\rm O}_x (3)$
to $v = \frac{1}{2}\, y_1^2 + {\rm O}_y (3)$ 
consists of matrices in $G_\stab^1$ with:
\[
0
\,=\,
a_{1,2}
\,=\,\cdots\,=\,
a_{1,n}
\ \ \ \ \ \ \ \ \ \ \ \ \ \ \ \ \ \ \ \
\text{and}
\ \ \ \ \ \ \ \ \ \ \ \ \ \ \ \ \ \ \ \
d
\,=\,
a_{1,1}^2.
\]
\end{Lemma}

\proof
Read~({\ref{eqFG}}) modulo ${\rm O}_x(3)$
and annihilate the coefficients
of $x_1x_2$, \dots, $x_1 x_n$:
\reqnomode\usetagform{EngelLie}
\begin{align}
0
&
\,=\
-\,v
+
\tfrac{1}{2}\,y_1^2
+
{\rm O}_y(3)
\notag
\\
&
\,=\,
-\,d\,u
+
\tfrac{1}{2}\,
\big(
a_{1,1}\,x_1
+
a_{1,2}\,x_2
+\cdots+
a_{1,n}\,x_n
+
b_1\,u
\big)^2
+
{\rm O}_x(3)
\notag
\\
&
\,\equiv\,
-\,d\,\tfrac{1}{2}\,x_1^2
-
{\rm O}_x(3)
+
\tfrac{1}{2}\,
a_{1,1}^2\,x_1^2
+
a_{1,1}\,x_1\,
\big(
a_{1,2}\,x_2
+\cdots+
a_{1,n}\,x_n
+
{\rm O}_x(2)
\big)
\notag
\\
&
\ \ \ \ \ \ \ \ \ \ \ \ \ \ \ \ \ \ \ \ \ \ \ \ \ \ \ \ \ \ \ \ \ \ \ 
\ \ \ \ \ \ \ \ \ \ \ \ \ \ \ \ \ \ 
+
\big(
a_{1,2}\,x_2
+\cdots+
a_{1,n}\,x_n
+
{\rm O}_x(2)
\big)^2
+
{\rm O}_x(3).
\notag
\qedhere
\end{align}
\endproof

\SectionHead{Normalization at Order $3$}
{normalization-order-3}

Now, let order $3$ terms appear:
\[
u
\,=\,
\tfrac{1}{2}\,
x_1^2
+
\sum_{
\sigma_1,\sigma_2,\sigma_3\geqslant 0
\atop
\sigma_1+\sigma_2+\sigma_3=3}\,
\frac{x_1^{\sigma_1}x_2^{\sigma_2}x_3^{\sigma_3}}{
\sigma_1!\,\sigma_2!\,\sigma_3!}\,
f_{\sigma_1,\sigma_2,\sigma_3}
+
{\rm O}_x(4).
\]
In the sum, 
pick the monomial $\frac{1}{6}\, x_1^3\, f_{1,1,1}$. 
Recall $x' = (x_2, \dots, x_n)$.
The remaining cubic terms
are of the form $x_1^2 A(x') + x_1 B(x') + C(x')$.
So they are ${\rm O}_{x'}(1)$. Since they are cubic,
they are of the product form ${\rm O}_{x'}(1)\, {\rm O}_x(2)$.
Thus:
\[
\aligned
u
&
\,=\,
F(x)
\,=\,
\tfrac{1}{2}\,x_1^2
+
\tfrac{1}{6}\,x_1^3\,f_{1,1,1}
+
{\rm O}_{x'}(1)\,{\rm O}_x(2)
+
{\rm O}_x(4),
\\
v
&
\,=\,
G(y)
\,=\,
\tfrac{1}{2}\,y_1^2
+
\tfrac{1}{6}\,y_1^3\,g_{1,1,1}
+
{\rm O}_{y'}(1)\,{\rm O}_y(2)
+
{\rm O}_y(4),
\endaligned
\]

\begin{Assertion}
One can normalize $g_{1,1,1} := 0$.
\end{Assertion}

\proof
With free $b_1 \in \R$, use the map belonging to $G_\stab^2$:
\[
y_1
:=
x_1+b_1\,u,
\ \ \ \ \ 
y_2
:=
x_2,
\ \ \ \ \ 
\dots\dots,
\ \ \ \ \ 
y_n
:=
x_n,
\ \ \ \ \ 
v
:=
u.
\]
So $y' =x'$, hence ${\rm O}_{y'}(1)\, {\rm O}_y(2) = 
{\rm O}_{x'}(1)\, {\rm O}_x(2)$, and the fundamental
equation~({\ref{eqFG}}) reads:
\[
\aligned
0
&
\,\equiv\,
-\,\tfrac{1}{2}\,x_1^2
-
\tfrac{1}{6}\,f_{1,1,1}\,x_1^3
-
{\rm O}_{x'}(1)\,{\rm O}_x(2)
-
{\rm O}_x(4)
\\
&
\ \ \ \ \
+
\tfrac{1}{2}\,
\big(
x_1+b_1\tfrac{1}{2}\,x_1^2+{\rm O}_x(3)
\big)^2
+
\tfrac{1}{6}\,
g_{1,1,1}\,
\big(
x_1+{\rm O}_x(2)
\big)^3
+
{\rm O}_{x'}(1)\,{\rm O}_x(2)
+
{\rm O}_x(4)
\\
&
\,\equiv\,
\tfrac{1}{6}\,x_1^3\,
\big[
-f_{1,1,1}+3\,b_1+g_{1,1,1}
\big]
+
{\rm O}_{x'}(1)\,{\rm O}_x(2)
+
{\rm O}_x(4).
\endaligned
\]

No monomial $x_1^3$ can appear in remainders. Hence the coefficient
of $x_1^3$ must vanish. This means that $g_{1,1,1} := f_{1,1,1}
- 3\, b_1$ necessarily. But since $b_1$ is a free
parameter in the affine transformation, 
we can choose $b_1 := \frac{1}{3}\,
f_{1,1,1}$ to normalize $g_{1,1,1} := 0$.
\endproof

To normalize further, we can restart from this $v = G(y)$ 
having $g_{1,1,1} = 0$, call it $u = F(x)$ with $f_{1,1,1} = 0$,
and again normalize the new target $v = G(y)$ with
$g_{1,1,1} = 0$. In other words,
both hypersurfaces are normalized similarly (as always):
\[
\aligned
u
&
\,=\,
F(x)
\,=\,
\tfrac{1}{2}\,x_1^2
+
0
+
{\rm O}_{x'}(1)\,{\rm O}_x(2)
+
{\rm O}_x(4),
\\
v
&
\,=\,
G(y)
\,=\,
\tfrac{1}{2}\,y_1^2
+
0
+
{\rm O}_{y'}(1)\,{\rm O}_y(2)
+
{\rm O}_y(4).
\endaligned
\]

Furthermore, before taking account
of the normalizations $f_{1,1,1} := 0$
and $g_{1,1,1} := 0$, remind that 
the current stability group $G_\stab^2$ is:
\[
\left[
\begin{array}{ccccc}
a_{1,1} & \red{\bf 0} & \cdots & \red{\bf 0} & b_1
\\
a_{2,1} & a_{2,2} & \cdots & a_{2,n} & b_2
\\
\vdots & \ddots & \vdots & \vdots
\\
a_{n,1} & a_{n,2} & \cdots & a_{n,n} & b_n
\\
\red{\bf 0} & \red{\bf 0} & \cdots & \red{\bf 0} & 
\red{a_{1,1}^2}
\end{array}
\right]^{\green{\bf 2}}.
\]
This being a subgroup of $\GL(n+1,\R)$, its 
block-trigonal determinant
must be nonzero, whence:
\[
a_{1,1}
\,\neq\,
0
\,\neq\,
\left\vert
\begin{array}{ccc}
a_{2,2} & \cdots & a_{2,n}
\\
\vdots & \ddots & \vdots 
\\
a_{n,2} & \cdots & a_{n,n}
\end{array}
\right\vert.
\]

Next, let cubic terms of the form $x_1^2\, {\rm O}_{x'}(1)$ appear:
\[
u
\,=\,
F(x)
\,=\,
\tfrac{1}{2}\,x_1^2
+
0
+
\tfrac{1}{2}\,
x_1^2\,
\big(
\varphi_2\,x_2+\cdots+\varphi_n\,x_n
\big)
+
x_1\,{\rm O}_{x'}(2)
+
{\rm O}_{x'}(3)
+
{\rm O}_x(4).
\]

\begin{Assertion}
The remaining cubic terms $x_1\,{\rm O}_{x'}(2) + {\rm O}_{x'}(3)
\equiv 0$ are zero.
\end{Assertion}

\proof
It suffices to show:
\reqnomode\usetagform{EngelLie}
\begin{align}
0
&
\overset{\text{\bf ?}}{\,=\,}
F_{x_1x_ix_j}(0)
\tag{(\forall\,\,2\,\leqslant\,i,\,j\,\leqslant\,n),}
\\
0
&
\overset{\text{\bf ?}}{\,=\,}
F_{x_ix_jx_k}(0)
\tag{(\forall\,\,2\,\leqslant\,i,\,j,\,k\,\leqslant\,n).}
\end{align}

By previous normalizations, we have $0 = F_{x_1 x_i}(0)$
for all $i = 2, \dots, n$.
Differentiating~({\ref{F-Hrank-1}}) with respect to $x_1$ and
to $x_k$ yields these vanishings:
\begin{footnotesize}
\begin{align*}
F_{x_1x_ix_j}(0)
&
\,=\,
\frac{F_{x_1x_1x_i}(0)\,\zero{F_{x_1x_j}(0)}
+\zero{F_{x_1x_i}(0)}\,F_{x_1x_1x_j}(0)}{
F_{x_1x_1}(0)}
-
\frac{\zero{F_{x_1x_i}(0)}\,\zero{F_{x_1x_j}(0)}\,
F_{x_1x_1x_1}(0)}{F_{x_1x_1}(0)\,F_{x_1x_1}(0)}
\,=\,
0,
\\
F_{x_ix_jx_k}(0)
&
\,=\,
\frac{F_{x_1x_ix_k}(0)\,\zero{F_{x_1x_j}(0)}
+\zero{F_{x_1x_i}(0)}\,F_{x_1x_jx_k}(0)}{
F_{x_1x_1}(0)}
-
\frac{\zero{F_{x_1x_i}(0)}\,\zero{F_{x_1x_j}(0)}\,
F_{x_1x_1x_k}(0)}{F_{x_1x_1}(0)\,F_{x_1x_1}(0)}
\,=\,
0.
\qedhere
\end{align*}
\end{footnotesize}
\endproof

Thus, the two hypersurfaces are:
\[
\aligned
u
&
\,=\,
F(x)
\,=\,
\tfrac{1}{2}\,x_1^2
+
0
+
\tfrac{1}{2}\,
x_1^2\,
\big(
\varphi_2\,x_2+\cdots+\varphi_n\,x_n
\big)
+
{\rm O}_x(4).
\\
v
&
\,=\,
G(y)
\,=\,
\tfrac{1}{2}\,
y_1^2
+
0
+
\tfrac{1}{2}\,
y_1^2\,
\big(
\psi_2\,y_2+\cdots+\psi_n\,y_n
\big)
+
{\rm O}_y(4).
\endaligned
\]

\begin{Assertion}
The property $0 = \varphi_2 = \cdots = \varphi_n$
is equivalent to $0 = \psi_2 = \cdots = \psi_n$.
\end{Assertion}

\proof
The general map of $G_\stab^2$ 
which stabilizes the normalization up to order $2$
writes:
\[
\aligned
y_1
&
\,=\,
a_{1,1}\,x_1
\ \ \ \ \ \ \ \ \ \ \ \ \ \ \ \ \ \ \ \ \ \ \ \ \ \ \ \ \ \ \ \ \ \ \
\ \ \ \ \ \,
+b_1\,u,
\\
y_2
&
\,=\,
a_{2,1}\,x_1
+
a_{2,2}\,x_2
+\cdots+
a_{2,n}\,x_n
+
b_2\,u,
\\
\cdots
&
\cdots\cdots\cdots\cdots\cdots\cdots\cdots\cdots\cdots\cdots\cdots
\cdots\cdots
\\
y_n
&
\,=\,
a_{n,1}\,x_1
+
a_{n,2}\,x_2
+\cdots+
a_{n,n}\,x_n
+
b_n\,u,
\\
v
&
\,=\,
\ \ \ \ \ \ \ \ \ \ \ \ \ \ \ \ \ \ \ \ \ \ \ \ \ \ \ \ \ \ \ \ \ \ \
\ \ \ \ \ \ \ \ \ \ \ \ \ \ \ \ \ \ \ \
a_{1,1}^2\,u.
\endaligned
\]
Therefore, the fundamental equation~({\ref{eqFG}}) writes:
\[
\aligned
0
&
\,\equiv\,
-\,a_{1,1}^2\,
\Big[
\tfrac{1}{2}\,x_1^2
+
\tfrac{1}{2}\,x_1^2\,
\big(
\varphi_2\,x_2
+\cdots+
\varphi_n\,x_n
\big)
+
{\rm O}_x(4)
\Big]
\\
&
\ \ \ \ \ 
+
\tfrac{1}{2}\,
\big(
a_{1,1}\,x_1
+
b_1\,\tfrac{1}{2}\,x_1^2
+
{\rm O}_x(3)
\big)^2
\\
&
\ \ \ \ \ 
+
\tfrac{1}{2}\,
\big(
a_{1,1}\,x_1
+
{\rm O}_x(2)
\big)^2\,
\Big[
\smallsum{2\leqslant j\leqslant n}\,
\psi_j\,
\big(
a_{j,1}\,x_1
+
a_{j,2}\,x_2
+\cdots+
a_{j,n}\,x_n
+
{\rm O}_x(2)
\big)
\Big]
+
{\rm O}_x(4).
\endaligned
\]

Of course, the coefficient of $x_1^2$ is zero.
Next, picking the coefficients, of $x_1^3$
and of $x_1^2 x_2$, \dots, $x_1^2 x_n$, we get:
\leqnomode\usetagform{default}
\begin{align}
\label{b1-psi}
b_1
&
\,=\,
-\,a_{1,1}\,
\smallsum{2\leqslant j\leqslant n}\,
\psi_j\,a_{j,1},
\\
\varphi_2
&
\,=\,
\smallsum{2\leqslant j\leqslant n}\,
\psi_j\,a_{j,2},\,\,
\dots\dots\dots,\,\,
\varphi_n
\,=\,
\smallsum{2\leqslant j\leqslant n}\,
\psi_j\,a_{j,n}.
\notag
\end{align}
With the invertible $(n-1) \times (n-1)$ matrix 
$(a_{j,k})$, we thus have $\varphi = {}^{\TT} a \cdot \psi$,
hence $\varphi = 0$ iff $\psi = 0$.
\endproof

Before we discuss the two distinct affinely invariant cases
$\varphi = 0$ and $\varphi \neq 0$, we must examine 
infinitesimal affine automorphisms.

\SectionHead{Tangency at Order $2$}
{tangency-order-2}

Take a hypersurface normalized at order $2$:
\[
u
\,=\,
\tfrac{1}{2}\,x_1^2
+
{\rm O}_x(3)
\,\,=\,
F(x).
\]
A general affine vector field:
\[
\aligned
L
&
\,=\,\,\,\,\,\,
\Big(
T_1
+
A_{1,1}\,x_1
+\cdots+
A_{1,n}\,x_n
+
B_1\,u
\Big)\,
\tfrac{\partial}{\partial x_1}
\,+
\\
&
\ \ \ \ \
+
\Big(
T_2
+
A_{2,1}\,x_1
+\cdots+
A_{2,n}\,x_n
+
B_2\,u
\Big)\,
\tfrac{\partial}{\partial x_2}
\,+
\\
&
\ \ \ \ \
+
\cdots\cdots\cdots\cdots\cdots\cdots\cdots\cdots\cdots\cdots\cdots
\cdots\cdots
\,+
\\
&
\ \ \ \ \
+
\Big(
T_n
+
A_{n,1}\,x_1
+\cdots+
A_{n,n}\,x_n
+
B_n\,u
\Big)\,
\tfrac{\partial}{\partial x_n}
\,+
\\
&
\ \ \ \ \
+
\Big(
T_0
+
C_1\,x_1
+\cdots+
C_n\,x_n
+
D\,u
\Big)\,
\tfrac{\partial}{\partial u},
\endaligned
\]
is tangent to the hypersurface $\{u = F(x)\}$
if and only if:
\leqnomode\usetagform{default}
\begin{align}
\label{L-tangent-u-Fx}
L\big(
-\,u+F(x)
\big)
\Big\vert_{u:=F(x)}
\,\equiv\,
0,
\end{align}
identically in $\C\{x_1, \dots, x_n\}$.

Here, neglecting terms of order $\geqslant 2$, this equation gives:
\[
\aligned
0
&
\,\equiv\,
-\,T_0
-
C_1\,x_1
-
C_2\,x_2
-\cdots-
C_n\,x_n
-
D\,{\rm O}_x(2)
\\
&
\ \ \ \ \
+
\big(
T_1
+
{\rm O}_x(1)
\big)\,
\big(
x_1
+
{\rm O}_x(2)
\big)
\\
&
\ \ \ \ \
+
\big(
T_2
+
{\rm O}_x(1)
\big)\,
{\rm O}_x(2)
\\
&
\ \ \ \ \
+
\cdots\cdots\cdots\cdots\cdots\cdots
\\
&
\ \ \ \ \
+
\big(
T_n
+
{\rm O}_x(1)
\big)\,
{\rm O}_x(2),
\endaligned
\]
whence:
\[
0
\,=\,
T_0,
\ \ \ \ \
C_1
\,=\,
T_1,
\ \ \ \ \
C_2
\,=\,
0,
\ \ \ \ \
\dots\dots,
\ \ \ \ \
C_n
\,=\,
0.
\]
Thus:
\leqnomode\usetagform{default}
\begin{align}
\label{L-order-2}
L
&
\,=\,\,\,\,\,\,
\Big(
T_1
+
A_{1,1}\,x_1
+\cdots+
A_{1,n}\,x_n
+
B_1\,u
\Big)\,
\tfrac{\partial}{\partial x_1}
\,+
\notag
\\
&
\ \ \ \ \
+
\Big(
T_2
+
A_{2,1}\,x_1
+\cdots+
A_{2,n}\,x_n
+
B_2\,u
\Big)\,
\tfrac{\partial}{\partial x_2}
\,+
\notag
\\
&
\ \ \ \ \
+
\cdots\cdots\cdots\cdots\cdots\cdots\cdots\cdots\cdots\cdots\cdots
\cdots\cdots
\,+
\\
&
\ \ \ \ \
+
\Big(
T_n
+
A_{n,1}\,x_1
+\cdots+
A_{n,n}\,x_n
+
B_n\,u
\Big)\,
\tfrac{\partial}{\partial x_n}
\,+
\notag
\\
&
\ \ \ \ \
+
\Big(
\ \ \ \ \ \ \ \ \ \
T_1\,x_1
\ \ \ \ \ \ \ \ \ \ \ \ \ \ \ \ \ \ \ \ \ \ \ \ \ \ \ \ \ \ 
+
D\,u
\Big)\,
\tfrac{\partial}{\partial u}.
\notag
\end{align}

\SectionHead{The Product Case $H^n \cong H^1 \times \R^{n-1}$}
{product-case-Hn-H1-Rn-1}

First, examine the {\em affinely invariant} case $\varphi_2 = 
\cdots = \varphi_n = 0$:
\[
u
\,=\,
\tfrac{1}{2}\,x_1^2
+
0
+
{\rm O}_x(4).
\]
Notice then that ({\ref{b1-psi}}) becomes $b_1 = 0$.

\begin{Proposition}
If such a hypersurface is affinely homogeneous, then
$F(x) = F(x_1)$ is independent of $x_2, \dots, x_n$.
\end{Proposition}

\proof
Now, examine the tangency equation~({\ref{L-tangent-u-Fx}})
using $L$ from~({\ref{L-order-2}}) modulo ${\rm O}_x(3)$:
\[
\aligned
0
&
\,\equiv\,
-\,
T_1\,x_1
-
D\,\big(
\tfrac{1}{2}\,x_1^2
+
{\rm O}_x(4)
\big)
\\
&
\ \ \ \ \ 
+
\Big(
T_1
+
A_{1,1}\,x_1
+
A_{1,2}\,x_2
+\cdots+
A_{1,n}\,x_n
+
{\rm O}_x(2)
\Big)\,
\big(
x_1
+
{\rm O}_x(3)
\big)
\\
&
\ \ \ \ \ 
+
\Big(
T_2
+
{\rm O}_x(1)
\Big)\,
{\rm O}_x(3)
\\
&
\ \ \ \ \ 
+
\cdots\cdots\cdots\cdots\cdots\cdots
\\
&
\ \ \ \ \ 
+
\Big(
T_n
+
{\rm O}_x(1)
\Big)\,
{\rm O}_x(3).
\endaligned
\]
The coefficients of $x_1^2$, of $x_1x_2$, \dots, of $x_1x_n$
must vanish, which gives:
\[
D
\,=\,
2\,A_{1,1},
\ \ \ \ \ 
A_{1,2}
\,=\,
0,
\ \ \ \ \ 
\dots\dots,
\ \ \ \ \ 
A_{1,n}
\,=\,
0.
\]

Thus:
\leqnomode\usetagform{default}
\begin{align}
\label{L-order-3}
L
&
\,=\,\,\,\,\,\,
\Big(
T_1
+
A_{1,1}\,x_1
\ \ \ \ \ \ \ \ \ \ \ \ \ \ \ \ \ \ \ \ \ \ \ \ \ \ \ \ \ \ \ \ \ \ \
\ \ \ \ \ \ \ 
+
B_1\,u
\Big)\,
\tfrac{\partial}{\partial x_1}
\,+
\notag
\\
&
\ \ \ \ \
+
\Big(
T_2
+
A_{2,1}\,x_1
+
A_{2,2}\,x_2
+\cdots+
A_{2,n}\,x_n
+
B_2\,u
\Big)\,
\tfrac{\partial}{\partial x_2}
\,+
\notag
\\
&
\ \ \ \ \
+
\cdots\cdots\cdots\cdots\cdots\cdots\cdots\cdots\cdots\cdots\cdots
\cdots\cdots\cdots\cdots\cdots
\,+
\\
&
\ \ \ \ \
+
\Big(
T_n
+
A_{n,1}\,x_1
+
A_{n,2}\,x_2
+\cdots+
A_{n,n}\,x_n
+
B_n\,u
\Big)\,
\tfrac{\partial}{\partial x_n}
\,+
\notag
\\
&
\ \ \ \ \
+
\Big(
\ \ \ \ \ \ \ \ \ \
T_1\,x_1
\ \ \ \ \ \ \ \ \ \ \ \ \ \ \ \ \ \ \ \ \ \ \ \ \ \ \ \ \ \ \ \ \ \ \
\ \ \ \ \ \
+
2\,A_{1,1}\,u
\Big)\,
\tfrac{\partial}{\partial u}.
\notag
\end{align}

Next, let order $4$ terms appear:
\[
u
\,=\,
\tfrac{1}{2}\,x_1^2
+
0
+
F^4(x)
+
{\rm O}_x(5)
\ \ \ \ \ \ \ \ \ \ \ \ \ \ \ \ \ \ \ \
\text{where}
\ \ \ \ \ \ \ \ \ \ \ \ \ \ \ \ \ \ \ \
F^4(x)
\,:=\,
\sum_{
\sigma_1,\dots,\sigma_n\geqslant 0
\atop
\sigma_1+\cdots+\sigma_n=4}\,
\frac{x_1^{\sigma_1}\cdots x_n^{\sigma_n}}{
\sigma_1!\,\cdots\,\sigma_n!}\,
f_{\sigma_1,\dots,\sigma_n}.
\]

\begin{Assertion}
This $F^4(x) = F^4(x_1)$ is independent of
$(x_2, \dots, x_n)$.
\end{Assertion}

\proof
Write\big/examine the tangency equation up to order $3$:
\[
\aligned
0
&
\,\equiv\,
-\,
T_1\,x_1
-
2\,A_{1,1}\,
\big(
\tfrac{1}{2}\,x_1^2
+
{\rm O}_x(4)
\big)
\\ 
&
\ \ \ \ \
+
\Big(
T_1
+
A_{1,1}\,x_1
\ \ \ \ \ \ \ \ \ \ \ \ \ \ \ \ \ \ \ \ \ \ \ \ \ \
+
B_1\,
\big(
\tfrac{1}{2}\,x_1^2
+
{\rm O}_x(4)
\big)
\Big)\,
\big(
x_1
+
F_{x_1}^4
+
{\rm O}_x(4)
\big)
\\ 
&
\ \ \ \ \
+
\big(
T_2
+
{\rm O}_x(1)
\big)\,
\big(
F_{x_2}^4
+
{\rm O}_x(4)
\big)
\\ 
&
\ \ \ \ \
+
\cdots\cdots\cdots\cdots\cdots\cdots\cdots\cdots\cdots
\\ 
&
\ \ \ \ \
+
\big(
T_n
+
{\rm O}_x(1)
\big)\,
\big(
F_{x_n}^4
+
{\rm O}_x(4)
\big).
\endaligned
\]

This equation, viewed modulo ${\rm O}_x(4)$,
is a polynomial in $(x_1, x_2, \dots, x_n)$
of degree $\leqslant 3$ which must vanish identically:
\leqnomode\usetagform{default}
\begin{align}
\label{T1-Fx1-4-B1}
0
&
\,\equiv\,
\ \ \
T_1\,F_{x_1}^4
+
B_1\,\tfrac{1}{2}\,x_1^3
\notag
\\
&
\ \ \ \ \
+
T_2\,F_{x_2}^4
\\
&
\ \ \ \ \
+
\cdots\cdots
\notag
\\
&
\ \ \ \ \
+
T_n\,F_{x_n}^4.
\notag
\end{align}

\begin{Notation}
{\bf [Coefficient picking]}
For $(\tau_1, \tau_2, \dots, \tau_n) \in \N^n$,
given a converging power series: 
\[
E(x)
\,=\, 
E(x_1,x_2,\dots,x_n)
\,=\,
\sum_{\sigma_1,\sigma_2,\dots,\sigma_n\geqslant 0}\,
\frac{x_1^{\sigma_1}x_2^{\sigma_2}\cdots x_n^{\sigma_n}}{
\sigma_1!\,\sigma_2!\,\cdots\,\sigma_n!}\,
E_{x_1^{\sigma_1}x_2^{\sigma_2}\cdots x_n^{\sigma_n}}(0),
\]
denote:
\[
\big[
x_1^{\tau_1}x_2^{\tau_2}\cdots x_n^{\tau_n}
\big]
\big(
E(x)
\big)
\,:=\,
\frac{1}{\tau_1!\,\tau_2!\,\cdots\,\tau_n!}\,
E_{x_1^{\tau_1}x_2^{\tau_2}\cdots x_n^{\tau_n}}(0).
\]
\end{Notation}

Therefore, in~({\ref{T1-Fx1-4-B1}}),
the coefficients of {\em all}
monomials $x_1^{\tau_1} x_2^{\tau_2} \cdots x_n^{\tau_n}$
with $\tau_1 + \tau_2 + \cdots + \tau_n = 3$
must vanish.

Disregarding the monomial $x_1^3$ in order not
to let $B_1$ intervene, {\em i.e.}
taking all
$\tau_1 + \tau_2 + \cdots + \tau_n = 3$ 
with $\tau_1 \neq 3$, we get:
\[
\aligned
0
&
\,=\,
\ \ \ 
T_1\,
\Big(
\big[
x_1^{\tau_1}x_2^{\tau_2}\cdots x_n^{\tau_n}
\big]
\big(
F_{x_1}^4
\big)
\Big)
\\
&
\ \ \ \ \
+
T_2\,
\Big(
\big[
x_1^{\tau_1}x_2^{\tau_2}\cdots x_n^{\tau_n}
\big]
\big(
F_{x_2}^4
\big)
\Big)
\\
&
\ \ \ \ \
+
\cdots\cdots\cdots\cdots\cdots\cdots\cdots\cdots\cdot
\\
&
\ \ \ \ \
+
T_n\,
\Big(
\big[
x_1^{\tau_1}x_2^{\tau_2}\cdots x_n^{\tau_n}
\big]
\big(
F_{x_n}^4
\big)
\Big).
\endaligned
\]

The following goes almost without saying.

\begin{Observation}
{\bf [Transitivity]}
Since $T_0 = 0$, at the origin $0 = (0, 0, \dots, 0) \in 
\R_{x,u}^{n+1}$, the value of $L$ is:
\[
L
\big\vert_0
\,=\,
T_1\,
\frac{\partial}{\partial x_1}
\Big\vert_0
+
T_2\,
\frac{\partial}{\partial x_2}
\Big\vert_0
+\cdots+
T_n\,
\frac{\partial}{\partial x_n}
\Big\vert_0,
\]
and since also: 
\[
T_0H^n
\,=\,
\Span\,
\Big(
\frac{\partial}{\partial x_1}
\Big\vert_0,\,
\frac{\partial}{\partial x_2}
\Big\vert_0,\,
\dots,\,
\frac{\partial}{\partial x_n}
\Big\vert_0
\Big),
\]
no linear relation between $T_1, T_2, \dots, T_n$ only
can hold for $H^n \subset \R^{n+1}$ to be affinely homogeneous.\qed
\end{Observation}

Consequently, the coefficients of $T_1$, of $T_2$, \dots,
of $T_n$ above must all vanish. 
Restricting attention to $\tau_1 = 0$, 
and taking all $(\tau_2, \dots, \tau_n) \in \N^{n-1}$
with $\tau_2 + \cdots + \tau_n = 3$,
we obtain:
\[
\aligned
0
&
\,=\,
\ \ \
T_1\,
\Big(
\big[
x_2^{\tau_2}\cdots x_n^{\tau_n}
\big]
\big(
F_{x_1}^4
\big)
\Big)
\\
&
\ \ \ \ \
+
T_2\,
\Big(
\big[
x_2^{\tau_2}\cdots x_n^{\tau_n}
\big]
\big(
F_{x_2}^4
\big)
\Big)
\\
&
\ \ \ \ \
+
\cdots\cdots\cdots\cdots\cdots\cdots\cdots
\\
&
\ \ \ \ \
+
T_n\,
\Big(
\big[
x_2^{\tau_2}\cdots x_n^{\tau_n}
\big]
\big(
F_{x_n}^4
\big)
\Big),
\endaligned
\]
whence:
\[
0
\,=\,
\big[
x_2^{\tau_2}\cdots x_n^{\tau_n}
\big]
\big(
F_{x_2}^4
\big),
\ \ \ \ \ 
\dots\dots,
\ \ \ \ \ 
\big[
x_2^{\tau_2}\cdots x_n^{\tau_n}
\big]
\big(
F_{x_n}^4
\big),
\]
which is equivalent to:
\[
0
\,\equiv\,
F_{x_2}^4,
\ \ \ \ \ 
\dots\dots,
\ \ \ \ \ 
0
\,\equiv\,
F_{x_n}^4.
\]
Thus $F^4(x)$ is independent of $x_2, \dots, x_n$,
and can be denoted $E^4(x_1)$.
\endproof

We can therefore let homogeneous terms of order $5$ appear:
\[
u
\,=\,
\tfrac{1}{2}\,
x_1^2
+
0
+
E^4(x_1)
+
F^5(x)
+
{\rm O}_x(6).
\]
Directly, let us reason by induction. With $m \geqslant 5$,
assume:
\[
u
\,=\,
\tfrac{1}{2}\,
x_1^2
+
0
+
E^4(x_1)
+\cdots+
E^{m-1}(x_1)
+
F^m(x)
+
{\rm O}_x(m+1).
\]
To terminate the proof of the proposition, we must show that
$F^m(x) \equiv E^m(x_1)$ is independent of $x_2, \dots, x_n$.

Using $L$ from~({\ref{L-order-3}}), write\big/examine
tangency equations up to order $m-1$:
\[
\footnotesize
\aligned
0
&
\,\equiv\,
-\,
T_1\,x_1
-
2\,A_{1,1}\,
\Big(
\tfrac{1}{2}\,x_1^2
+
E^4
+\cdots+
E^{m-1}
+
{\rm O}_x(m)
\Big)
\\
&
\ \ \ \ \
+
\Big[
T_1
+
A_{1,1}\,x_1
+
B_1\,
\big(
\tfrac{1}{2}\,x_1^2
+
E^4
+\cdots+
E^{m-2}
+
{\rm O}_x(m-1)
\big)
\Big]\,
\Big[
x_1
+
E_{x_1}^4
+\cdots+
E_{x_1}^{m-1}
+
F_{x_1}^m
+
{\rm O}_x(m)
\Big]
\\
&
\ \ \ \ \
+
\big[
T_2
+
{\rm O}_x(1)
\big]\,
\big[
F_{x_2}^m
+
{\rm O}_x(m)
\big]
\\
&
\ \ \ \ \
+
\cdots\cdots\cdots\cdots\cdots\cdots\cdots\cdots\cdots
\cdot
\\
&
\ \ \ \ \
+
\big[
T_n
+
{\rm O}_x(1)
\big]\,
\big[
F_{x_n}^m
+
{\rm O}_x(m)
\big].
\endaligned
\]
In the first two lines, all appearing monomials $x_1^{\tau_1}
x_2^{\tau_2} \cdots x_n^{\tau_n}$ with $\tau_1 + \tau_2 + \cdots +
\tau_n \leqslant m-1$ are such that $\tau_1 \geqslant 1$, except the
ones in $T_1\, F_{x_1}^m$.

Therefore, when applying the coefficients-taking operators
$[x_2^{\tau_2} \cdots x_n^{\tau_n}]
(\centersmallbullet)$ for all $(\tau_2, \dots, \tau_n) \in \N^{n-1}$
with $\tau_2 + \cdots + \tau_n = m-1$, it remains only:
\[
\aligned
0
&
\,=\,
\ \ \
T_1\,
\Big(
\big[
x_2^{\tau_2}\cdots x_n^{\tau_n}
\big]
\big(
F_{x_1}^m
\big)
\Big)
\\
&
\ \ \ \ \
+
T_2\,
\Big(
\big[
x_2^{\tau_2}\cdots x_n^{\tau_n}
\big]
\big(
F_{x_2}^m
\big)
\Big)
\\
&
\ \ \ \ \
+
\cdots\cdots\cdots\cdots\cdots\cdots\cdots
\\
&
\ \ \ \ \
+
T_n\,
\Big(
\big[
x_2^{\tau_2}\cdots x_n^{\tau_n}
\big]
\big(
F_{x_n}^m
\big)
\Big),
\endaligned
\]
whence:
\[
0
\,=\,
\big[
x_2^{\tau_2}\cdots x_n^{\tau_n}
\big]
\big(
F_{x_2}^m
\big),
\ \ \ \ \ 
\dots\dots,
\ \ \ \ \ 
\big[
x_2^{\tau_2}\cdots x_n^{\tau_n}
\big]
\big(
F_{x_n}^m
\big),
\]
which is equivalent to:
\[
0
\,\equiv\,
F_{x_2}^m,
\ \ \ \ \ 
\dots\dots,
\ \ \ \ \ 
0
\,\equiv\,
F_{x_n}^m.
\]
Thus $F^m(x)$ is independent of $x_2, \dots, x_n$,
can be denoted $E^m(x_1)$, and a plain induction
on $m \longrightarrow \infty$ concludes
the proof of the proposition:
\[
u
\,=\,
\tfrac{1}{2}\,x_1^2
+
\sum_{m=4}^\infty\,
E^m(x_1)
\,=:\,
F(x_1).
\qedhere
\]
\endproof

\SectionHead{Interlude: Expansion of $F$ in Homogeneous 
(In)Dependent Monomials}
{interlude-dependent-monomials}

Expand:
\[
F(x)
\,=\,
\sum_{\sigma_1,\dots,\sigma_n\geqslant 0}\,
x_1^{\sigma_1}\cdots x_n^{\sigma_n}\,
F_{\sigma_1,\dots,\sigma_n},
\]
with $F_\sigma := \frac{1}{\sigma!}\, F_{x^\sigma} (0)$.
In accordance with 
Terminology~{\ref{Terminology-border-dependent-jets}},
let us introduce

\begin{Terminology}
A monomial 
$x_1^{\sigma_1} x_2^{\sigma_2} \cdots x_n^{\sigma_n}\,
F_{\sigma_1, \sigma_2, \dots, \sigma_n}$ will be said:

\smallskip\noindent$\bullet$\,
{\sl independent} if $\sigma_2 + \cdots + \sigma_n \leqslant 1$;

\smallskip\noindent$\bullet$\,
{\sl border-dependent} if 
$\sigma_2 + \cdots + \sigma_n = 2$;

\smallskip\noindent$\bullet$\,
{\sl body-dependent} if 
$\sigma_2 + \cdots + \sigma_n \geqslant 3$.

\end{Terminology}

Then $F$ decomposes into $3$ parts:
\[
F
\,=\,
\sum_{\sigma_2+\cdots+\sigma_n\leqslant 1}\,
x^\sigma\,F_\sigma
+
\sum_{\sigma_2+\cdots+\sigma_n=2}\,
x^\sigma\,F_\sigma
+
\sum_{\sigma_2+\cdots+\sigma_n\geqslant 3}\,
x^\sigma\,F_\sigma.
\]

If one wants to emphasize independent monomials only,
the dependent monomials can be gathered as a plain remainder:
\[
F
\,=\,
\sum_{\sigma_1,\sigma_2,\dots,\sigma_n\geqslant 0
\atop
\sigma_2+\cdots+\sigma_n\leqslant 1}\,
x_1^{\sigma_1}x_2^{\sigma_2}\cdots x_n^{\sigma_n}\,
F_{\sigma_1,\sigma_2,\dots,\sigma_n}
+
{\rm O}_{x_2,\dots,x_n}(2).
\]

If one wants to also show border-dependent monomials,
the body-dependent monomials can be gathered as a plain remainder:
\[
F
\,=\,
\sum_{\sigma_1,\sigma_2,\dots,\sigma_n\geqslant 0
\atop
\sigma_2+\cdots+\sigma_n\leqslant 1}\,
x_1^{\sigma_1}x_2^{\sigma_2}\cdots x_n^{\sigma_n}\,
F_{\sigma_1,\sigma_2,\dots,\sigma_n}
+
\sum_{\sigma_1,\sigma_2,\dots,\sigma_n\geqslant 0
\atop
\sigma_2+\cdots+\sigma_n=2}\,
x_1^{\sigma_1}x_2^{\sigma_2}\cdots x_n^{\sigma_n}\,
F_{\sigma_1,\sigma_2,\dots,\sigma_n}
+
{\rm O}_{x_2,\dots,x_n}(3).
\]

The decomposition of $F(x)$ as a sum of homogeneous terms
will be constantly used:
\[
F(x)
\,=\,
\sum_{m=2}^\infty\,
\sum_{\sigma_1+\sigma_2+\cdots+\sigma_n=m}\,
x_1^{\sigma_1}x_2^{\sigma_2}\cdots x_n^{\sigma_n}\,
F_{\sigma_1,\sigma_2,\dots,\sigma_n}
\,\,=:\,
\sum_{m=2}^\infty\,
F^m(x).
\]
Then each $F^m(x)$ can be subjected to similar decompositions:
\[
\aligned
F^m(x)
&
\,=\,
\sum_{\sigma_1+\sigma_2+\cdots+\sigma_n=m
\atop
\ \ \ 
\sigma_2+\cdots+\sigma_n\leqslant 1}\,
x^\sigma\,F_\sigma
+
\sum_{\sigma_1+\sigma_2+\cdots+\sigma_n=m
\atop
\ \ \ 
\sigma_2+\cdots+\sigma_n=2}\,
x^\sigma\,F_\sigma
+
\sum_{\sigma_1+\sigma_2+\cdots+\sigma_n=m
\atop
\ \ \ 
\sigma_2+\cdots+\sigma_n\geqslant 3}\,
x^\sigma\,F_\sigma
\\
&
\,=\,
\sum_{\sigma_1+\sigma_2+\cdots+\sigma_n=m
\atop
\ \ \ 
\sigma_2+\cdots+\sigma_n\leqslant 1}\,
x^\sigma\,F_\sigma
+
\sum_{\sigma_1+\sigma_2+\cdots+\sigma_n=m
\atop
\ \ \ 
\sigma_2+\cdots+\sigma_n=2}\,
x^\sigma\,F_\sigma
+
{\rm O}_{x_2,\dots,x_n}(3)
\\
&
\,=\,
\sum_{\sigma_1+\sigma_2+\cdots+\sigma_n=m
\atop
\ \ \ 
\sigma_2+\cdots+\sigma_n\leqslant 1}\,
x^\sigma\,F_\sigma
+
{\rm O}_{x_2,\dots,x_n}(2).
\endaligned
\]
Lastly, recall the abbreviation:
\[
x'
\,:=\,
(x_2,\dots,x_n).
\]

\SectionHead{The Nondegenerate Order $3$ Case and a Starting Induction}
{nondegenerate-case-starting-induction}

Disregarding the product case of
Section~{\ref{product-case-Hn-H1-Rn-1}}, 
assume that 
$(\varphi_2, \dots, \varphi_n) \neq (0, \dots, 0)$.
A natural affine coordinate change:
\[
u
\,=\,
\tfrac{1}{2}\,x_1^2
+
\tfrac{1}{2}\,
x_1^2\,
\big(
\underbrace{
\varphi_2\,x_2+\cdots+\varphi_n\,x_n}_{
=:\,\,{\sf new}\,\,x_2}
\big)
+
{\rm O}_x(4),
\]
enables to normalize:
\[
u
\,=\,
\tfrac{1}{2}\,x_1^2
+
\tfrac{1}{2}\,x_1^2\,x_2
+
{\rm O}_x(4).
\]

Now, set up a general inductive reasoning. For some integer $\nu$
with $3 \leqslant \nu \leqslant n$, assume that, modulo dependent
monomials which, by 
Section~{\ref{interlude-dependent-monomials}}, 
can be all gathered as a remainder ${\rm O}_{x'}(2)$, 
the hypersurface equation
is of the form:
\[
u
\,=\,
\frac{x_1^2}{2}
+
\frac{x_1^2x_2}{2}
+\cdots+
\frac{x_1^{\nu-1}x_{\nu-1}}{(\nu-1)!}
+
{\rm O}_{x'}(2)
+
{\rm O}_x(\nu+1).
\]
Then let appear all independent monomials of order $\nu+1$:
\[
\aligned
u
&
\,=\,
\frac{x_1^2}{2}
+
\frac{x_1^2x_2}{2}
+\cdots+
\frac{x_1^{\nu-1}x_{\nu-1}}{(\nu-1)!}
\\
&
\ \ \ \ \
+
\varphi_1\,
\frac{x_1^\nu x_1}{\nu!}
+\cdots+
\varphi_{\nu-1}\,
\frac{x_1^\nu x_{\nu-1}}{\nu!}
+
\varphi_\nu\,
\frac{x_1^\nu x_\nu}{\nu!}
+\cdots+
\varphi_n\,
\frac{x_1^\nu x_n}{\nu!}
+
{\rm O}_{x'}(2)
+
{\rm O}_x(\nu+2).
\endaligned
\]
Again, if $(\varphi_\nu, \dots, \varphi_n) \neq (0, \dots, 0)$,
a natural affine coordinate change:
\[
\aligned
u
&
\,=\,
\frac{x_1^2}{2}
+
\frac{x_1^2x_2}{2}
+\cdots+
\frac{x_1^{\nu-1}x_{\nu-1}}{(\nu-1)!}
\\
&
\ \ \ \ \
+
\frac{x_1^\nu}{\nu!}\,
\Big(
\underbrace{
\varphi_1\,x_1
+\cdots+
\varphi_{\nu-1}\,x_{\nu-1}
+
\varphi_\nu\,x_\nu
+\cdots+
\varphi_n\,x_n
}_{=:\,\,{\sf new}\,\,x_\nu}
\Big)
+
{\rm O}_{x'}(2)
+
{\rm O}_x(\nu+2),
\endaligned
\]
leads to:
\[
u
\,=\,
\frac{x_1^2}{2}
+
\frac{x_1^2x_2}{2}
+\cdots+
\frac{x_1^{\nu-1}x_{\nu-1}}{(\nu-1)!}
+
\frac{x_1^\nu x_\nu}{\nu!}
+
{\rm O}_{x'}(2)
+
{\rm O}_x(\nu+2),
\]
and then inductively to:
\[
u
\,=\,
\frac{x_1^2}{2}
+
\frac{x_1^2x_2}{2}
+\cdots+
\frac{x_1^\nu x_\nu}{\nu!}
+\cdots+
\frac{x_1^n x_n}{n!}
+
{\rm O}_{x'}(2)
+
{\rm O}_x(n+2),
\]
These are the hypersurfaces we mainly want to study:
they involve all variables $x_1, x_2, \dots, x_n$.

But before going further, we must examine\big/study
the circumstance where $\varphi_\nu = \cdots = \varphi_n = 0$:
\[
\aligned
u
&
\,=\,
\frac{x_1^2}{2}
+
\frac{x_1^2x_2}{2}
+\cdots+
\frac{x_1^{\nu-1}x_{\nu-1}}{(\nu-1)!}
\\
&
\ \ \ \ \
+
\varphi_1\,
\frac{x_1^\nu x_1}{\nu!}
+\cdots+
\varphi_{\nu-1}\,
\frac{x_1^\nu x_{\nu-1}}{\nu!}
+
0
+\cdots+
0
+
{\rm O}_{x'}(2)
+
{\rm O}_x(\nu+2).
\endaligned
\]
This degenerate branch will again lead to a product situation.

Since the arguments in the next 
Section~{\ref{product-case-Hn-Hnu-1-R-n-nu}}
will again involve applying an affine vector field $L$
to the equation $0 = - u + F(x)$, and since $L$ is a
first-order derivation,
we need to know the border-dependent monomials as well.
Recall that, by Section~{\ref{interlude-dependent-monomials}},
body-dependent monomials that are 
not border-dependent can be gathered as ${\rm O}_{x'}(3)$.

To organize properly the thought, we 
consider simultaneously the
two cases ${\bf 0} \cdot \frac{x_1^\nu x_\nu}{\nu!}$ and
${\bf 1} \cdot \frac{x_1^\nu x_\nu}{\nu!}$, 
by setting up an

\begin{InductionHypothesis}
\label{Induct-Hyp}
{\em For some integer $\nu$ with $3 \leqslant \nu \leqslant n$,
the hypersurface equation writes:}
\[
\aligned
u
\,=\,
\frac{x_1^2}{2}
+
\frac{x_1^2x_2}{2}
&
+
\sum_{m=3}^{\nu-1}\,
\Big(
\frac{x_1^mx_m}{m!}
+
x_1^{m-1}\,
\sum_{i,j\geqslant 2
\atop
i+j=m+1}\,
\tfrac{1}{2}\,
\frac{x_ix_j}{(i-1)!\,(j-1)!}
+
R^{m+1}(x_1,x_2,\dots,x_{m-1})
\Big)
\\
&
\ \ \ \ \
+
\left\{
\begin{smallmatrix}
{\bf 0}
\\
{\bf 1}
\end{smallmatrix}
\right\}
\frac{x_1^\nu x_\nu}{\nu!}
+
x_1^{\nu-1}\,
\sum_{i,j\geqslant 2}\,
\tfrac{1}{2}\,
x_ix_j\,
\Lambda_{i,j}^\nu
+
R^{\nu+1}(x_1,x_2,\dots,x_\nu,\dots,x_n)
+
{\rm O}_x(\nu+2),
\endaligned
\]
{\em 
where $R^{m+1}$ is homogeneous of order $m+1$
in $(x_1, x_2, \dots, x_{m-1})$, 
and is of order $\geqslant 3$ in $(x_2, \dots, x_{m-1})$,
where the $\Lambda_{i,j}^\nu = 
\Lambda_{j,i}^\nu$ are unknown constants,
and where $R^{\nu+1}$ is homogeneous of order $\nu+1$ 
in $(x_1, x_2, \dots, x_\nu, \dots, x_n)$,
and is of order $\geqslant 3$ in $(x_2, \dots, x_\nu, \dots,
x_n)$.
}
\end{InductionHypothesis}

We therefore assume that up to $m = \nu-1$, the values
of the border-dependent jets have been computed, 
as they appear within the large parentheses. 
For $\nu=3$, the formula holds true with empty
$\sum_{m=3}^{\nu-1}$.
To close the induction on $\nu$, 
we must show that $R^{\nu+1} = R^{\nu+1}(x_1,\dots,x_{\nu-1})$
is independent of $x_\nu, x_{\nu+1}, \dots, x_n$, 
and we must determine the values
of the $\Lambda_{i,j}^\nu$. In the nondegenerate
branch ${\bf 1} \cdot \frac{x_1^\nu x_\nu}{\nu!}$
that we will treat later,
we will have to show that $\Lambda_{i,j}^\nu = 0$
whenever $i+j \neq \nu + 1$, and
that $\Lambda_{i,j}^\nu = \frac{1}{(i-1)! (j-1)!}$
when $i + j = \nu + 1$.

Abbreviate the homogeneous terms of $F(x)$ normalized 
up to order $\leqslant \nu$, namely up to
$m = \nu-1$ in the sum, as:
\[
u
\,=\,
N^2(x_1)
+
N^3(x_1,x_2)
+\cdots+
N^\nu(x_1,x_2,\dots,x_{\nu-1})
+
{\rm O}_x(\nu+1),
\]
where, for $3 \leqslant m \leqslant \nu - 1$:
\[
N^{m+1}
\,:=\,
\frac{x_1^mx_m}{m!}
+
x_1^{m-1}\,
\sum_{i,j\geqslant 2
\atop
i+j=m+1}\,
\tfrac{1}{2}\,
\frac{x_ix_j}{(i-1)!\,(j-1)!}
+
R^{m+1}(x_1,x_2,\dots,x_{m-1}),
\]
and abbreviate also the full dependent 
remainder homogeneous of order $\nu+1$
after $\big\{ \begin{smallmatrix} {\bf 0} \\ {\bf 1}
\end{smallmatrix} \big\}\, \frac{x_1^\nu x_\nu}{\nu !}$ as:
\[
S^{\nu+1}
\big(
x_1,x_2,\dots,x_\nu,\dots,x_n
\big)
\,:=\,
x_1^{\nu-1}
\sum_{i,j\geqslant 2}\,
\tfrac{1}{2}\,
x_ix_j\,
\Lambda_{i,j}^\nu
+
R^{\nu+1}
\big(
x_1,x_2,\dots,x_\nu,\dots,x_n
\big),
\]
which is of order $\geqslant 2$ in $(x_2, \dots, x_\nu, \dots,
x_n)$, so that:
\[
u
\,=\,
F(x)
\,=\,
N^2
+\cdots+
N^\nu
+
\left\{
\begin{smallmatrix}
{\bf 0}
\\
{\bf 1}
\end{smallmatrix}
\right\}
\frac{x_1^\nu x_\nu}{\nu!}
+
S^{\nu+1}
+
{\rm O}_x(\nu+2).
\]

\begin{Assertion}
\label{Ass-S-nu-1}
This function $S^{\nu+1}$ is independent of $x_\nu,
x_{\nu+1}, \dots, x_n$.
\end{Assertion}

\proof
For any two indices $\nu \leqslant k, \ell \leqslant n$,
we must have by our constant Hessian rank $1$ hypothesis:
\[
\footnotesize
\aligned
0
&
\,\equiv\,
F_{x_1x_1}(x)\cdot F_{x_kx_\ell}(x)
-
F_{x_1x_k}(x)\cdot F_{x_1x_\ell}(x)
\\
&
\,\equiv\,
\big(1+{\rm O}_x(1)\big)
\cdot
\Big(
S_{x_kx_\ell}^{\nu+1}
+
{\rm O}_x(\nu)
\Big)
-
\Big(
\big\{
\begin{smallmatrix}
{\bf 0}
\\
{\bf 1}
\end{smallmatrix}
\big\}\,
\tfrac{x_1^{\nu-1}}{(\nu-1)!}\,
\tfrac{\partial x_\nu}{\partial x_k}
+
S_{x_1x_k}^{\nu+1}
+
{\rm O}_x(\nu)
\Big)
\cdot
\Big(
\big\{
\begin{smallmatrix}
{\bf 0}
\\
{\bf 1}
\end{smallmatrix}
\big\}\,
\tfrac{x_1^{\nu-1}}{(\nu-1)!}\,
\tfrac{\partial x_\nu}{\partial x_\ell}
+
S_{x_1x_\ell}^{\nu+1}
+
{\rm O}_x(\nu)
\Big)
\\
&
\,\equiv\,
S_{x_kx_\ell}^{\nu+1}
+
{\rm O}_x(\nu)
-
{\rm O}_x(2\nu-2),
\endaligned
\]
since $2\nu-2 \geqslant \nu$ as $\nu \geqslant 3$,
and this yields:
\[
0
\,\equiv\,
S_{x_kx_\ell}^{\nu+1},
\]
Since $S^{\nu+1}$ is of order $\geqslant 2$ in
$(x_2, \dots, x_n)$, this concludes.
\endproof

\SectionHead{The Product Case $H^n \cong H^{\nu-1} \times 
\R^{n-\nu}$}
{product-case-Hn-Hnu-1-R-n-nu}

Here, we treat the degenerate branch
${\bf 0} \cdot \frac{x_1^\nu x_\nu}{\nu!}$. 
Thus:
\[
u
\,=\,
F(x)
\,=\,
N^2(x_1)
+\cdots+
N^\nu(x_1,\dots,x_{\nu-1})
+
{\bf 0}
+
S^{\nu+1}(x_1,\dots,x_{\nu-1})
+
{\rm O}_x(\nu+2).
\]
Notice that the term $S^{\nu+1}$
(homogeneous of order $\nu+1$)
depends only on $x_1, \dots, x_{\nu-1}$,
as does the preceding one $N^\nu$.

\begin{Proposition}
If such a hypersurface is affinely homogeneous, then
$F = F(x_1,\dots,x_{\nu-1})$ is independent of $x_\nu, \dots, x_n$.
\end{Proposition}

\proof
Let homogeneous terms of order $\nu+2$ appear:
\[
u
\,=\,
F
\,=\,
N^2
+\cdots+
N^\nu
+
S^{\nu+1}
+
F^{\nu+2}(x_1,\dots,x_{\nu-1},x_\nu,\dots,x_n)
+
{\rm O}(\nu+3).
\]
We claim that $F^{\nu+2}$ is independent of $x_\nu, \dots, x_n$.

Indeed, using $L$ from~({\ref{L-order-3}}), write\big/examine
tangency equations up to order $\nu + 1$:
\[
\!\!\!\!\!\!\!\!\!\!\!\!\!\!\!\!\!\!\!\!\!\!\!\!\!
\footnotesize
\aligned
0
&
\,\equiv\,
-\,
T_1\,x_1
-
2\,A_{1,1}\,
\Big(
N^2+\cdots+N^\nu
+
S^{\nu+1}
+
{\rm O}_x(\nu+2)
\Big)
\\
&
+
\Big[
T_1
+
A_{1,1}x_1
\ \ \ \ \ \ \ \ \ \ \ \ \ \ \ \ \ \ \ \ \ \ \ \ \ \ \ \ \ \ \ \ \ \ \
\ \ \ \ \ \ \ \ \ \ \ \ \ \ \ 
+
B_1
\big(
N^2+\cdots+N^\nu
+
{\rm O}_x(\nu+1)
\big)
\Big]
\cdot
\Big[
N_{x_1}^2+\cdots+N_{x_1}^\nu
+
S_{x_1}^{\nu+1}
+
F_{x_1}^{\nu+2}
+
{\rm O}_x(\nu+2)
\Big]
\\
&
+
\Big[
T_2
+
A_{2,1}x_1+\cdots+A_{2,\nu}x_\nu+\cdots+A_{2,n}x_n
+
B_2
\big(
N^2+\cdots+N^{\nu-1}+{\rm O}_x(\nu)
\big)
\Big]
\cdot
\Big[
N_{x_2}^3+\cdots+N_{x_2}^\nu
+
S_{x_2}^{\nu+1}
+
F_{x_2}^{\nu+2}
+
{\rm O}_x(\nu+2)
\Big]
\\
&
+
\cdots\cdots\cdots\cdots\cdots\cdots\cdots\cdots\cdots
\cdots\cdots\cdots\cdots\cdots\cdots\cdots\cdots\cdots
\cdots\cdots\cdots\cdots\cdots\cdots\cdots\cdots\cdots
\cdots\cdots\cdots\cdots\cdots\cdots\cdots\cdots\cdots
\cdots\cdots\cdots
\cdot
\\
&
+
\Big[
T_{\nu-1}
+
A_{\nu-1,1}x_1+\cdots+A_{\nu-1,\nu}x_\nu+\cdots+A_{\nu-1,n}x_n
+
B_{\nu-1}
\big(
N^2+{\rm O}_x(3)
\big)
\Big]
\cdot
\Big[
N_{x_{\nu-1}}^\nu
+
S_{x_{\nu-1}}^{\nu+1}
+
F_{x_{\nu-1}}^{\nu+2}
+
{\rm O}_x(\nu+2)
\Big]
\\
&
+
\Big[
T_\nu
+
{\rm O}_x(1)
\Big]
\cdot
\Big[
F_{x_\nu}^{\nu+2}
+
{\rm O}_x(\nu+2)
\Big]
\\
&
+
\cdots\cdots\cdots\cdots\cdots\cdots\cdots\cdots\cdots
\cdots\cdots
\cdot
\\
&
+
\Big[
T_n
+
{\rm O}_x(1)
\Big]
\cdot
\Big[
F_{x_n}^{\nu+2}
+
{\rm O}_x(\nu+2)
\Big].
\endaligned
\]

For all $(\tau_\nu, \dots, \tau_n) \in 
\N^{n-\nu+1}$ with $\tau_\nu + \cdots + \tau_n = \nu + 1$, 
to this equation, 
apply the coefficients-picking operators
$[x_\nu^{\tau_\nu} \cdots x_n^{\tau_n}] ( \centersmallbullet )$.
Since $N^2$, $N^3$, \dots, $N^\nu$,
$S^{\nu+1}$ depend only on $(x_1, \dots, x_{\nu-1})$,
we obtain:
\[
\aligned
0
&
\,\equiv\,
\ \ \ 
T_1\,
\Big(
[x_\nu^{\tau_\nu}\cdots x_n^{\tau_n}]
\big(F_{x_1}^{\nu+2}\big)
\Big)
\\
&
\ \ \ \ \ 
+
\cdots\cdots\cdots\cdots\cdots\cdots\cdots\cdots
\\
&
\ \ \ \ \ 
+
T_{\nu-1}\,
\Big(
[x_\nu^{\tau_\nu}\cdots x_n^{\tau_n}]
\big(F_{x_{\nu-1}}^{\nu+2}\big)
\Big)
\\
&
\ \ \ \ \ 
+
T_\nu\,
\Big(
[x_\nu^{\tau_\nu}\cdots x_n^{\tau_n}]
\big(F_{x_\nu}^{\nu+2}\big)
\Big)
\\
&
\ \ \ \ \ 
+
\cdots\cdots\cdots\cdots\cdots\cdots\cdots
\cdot
\\
&
\ \ \ \ \ 
+
T_n\,
\Big(
[x_\nu^{\tau_\nu}\cdots x_n^{\tau_n}]
\big(F_{x_n}^{\nu+2}\big)
\Big),
\endaligned
\]
whence:
\[
0
\,=\,
\big[
x_\nu^{\tau_\nu}\cdots x_n^{\tau_n}
\big]
\big(
F_{x_\nu}^{\nu+2}
\big),
\ \ \ \ \ 
\dots\dots,
\ \ \ \ \ 
\big[
x_\nu^{\tau_\nu}\cdots x_n^{\tau_n}
\big]
\big(
F_{x_n}^{\nu+2}
\big),
\]
which is equivalent to:
\[
0
\,\equiv\,
F_{x_\nu}^{\nu+2},
\ \ \ \ \ 
\dots\dots,
\ \ \ \ \ 
0
\,\equiv\,
F_{x_n}^{\nu+2}.
\]
Thus $F^{\nu+2}$ is independent of $x_\nu, \dots, x_n$,
and can be denoted $E^{\nu+2}(x_1,\dots,x_{\nu-1})$.

Next, let homogeneous terms of order $\nu + 3$ appear:
\[
u
\,=\,
F
\,=\,
N^2
+\cdots+
N^\nu
+
S^{\nu+1}
+
E^{\nu+2}
+
F^{\nu+3}(x_1,\dots,x_{\nu-1},x_\nu,\dots,x_n)
+
{\rm O}(\nu+4).
\]
By writing the tangency equation up to order $\nu + 2$,
we realize similarly that 
$F^{\nu+3}$ is independent of $x_\nu, \dots, x_n$.
Proceeding by induction, we conclude that:
\[
u
\,=\,
F
\,=\,
N^2
+\cdots+
N^\nu
+
S^{\nu+1}
+
\sum_{m=\nu+2}^\infty\,
E^m(x_1,\dots,x_{\nu-1}).
\qedhere
\]
\endproof

We now disregard this degenerate case. 

\SectionHead{The Nondegenerate Case ${\bf 1} \cdot 
\frac{x_1^\nu x_\nu}{\nu!}$}
{nondegenerate-case-x1nu-xnu}

At last, we can start to treat the most interesting branch
${\bf 1} \cdot \frac{x_1^\nu x_\nu}{\nu!}$. Thus,
as already stated by the Induction 
Hypothesis~{\ref{Induct-Hyp}}, 
we start from:
\[
\aligned
u
\,=\,
\frac{x_1^2}{2}
+
\frac{x_1^2x_2}{2}
&
+
\sum_{m=3}^{\nu-1}\,
\Big(
\frac{x_1^mx_m}{m!}
+
x_1^{m-1}\,
\sum_{i,j\geqslant 2
\atop
i+j=m+1}\,
\tfrac{1}{2}\,
\frac{x_ix_j}{(i-1)!\,(j-1)!}
+
{\rm O}_{x_2,\dots,x_{m-1}}(3)
\Big)
\\
&
\ \ \ \ \
+
\frac{x_1^\nu x_\nu}{\nu!}
+
x_1^{\nu-1}\,
\sum_{2\leqslant i,j\leqslant\nu-1}\,
\tfrac{1}{2}\,
x_ix_j\,
\Lambda_{i,j}^\nu
+
{\rm O}_{x_2,\dots,x_{\nu-1}}(3)
+
{\rm O}_x(\nu+2),
\endaligned
\]
where we already know from Assertion~{\ref{Ass-S-nu-1}} that
$\Lambda_{i,j}^\nu = 0$ when $i \geqslant \nu$ or $j \geqslant \nu$,
and that the body-dependent remainder $R^{\nu+1}$ is independent of
$x_\nu, x_{\nu+1}, \dots, x_n$, hence is of order $ \geqslant 3$ in
$x_2, \dots, x_{\nu-1}$.

Here is the statement we promised.

\begin{Proposition}
\label{Prp-Lambda-i-i-nu-zero}
One has $\Lambda_{i,j}^\nu = 0$ whenever $i + j \neq \nu + 1$,
and $\Lambda_{i,j}^\nu = \frac{1}{(i-1)!\,(j-1)!}$
when $i + j = \nu+1$.
\end{Proposition}

\proof
Firstly:
\[
F_{x_1x_1}
\,=\,
1
+
{\rm O}_{x_2,\dots,x_{\nu-1}}(1)
+
{\rm O}_x(\nu).
\]
Secondly, for any two indices $2 \leqslant i,j \leqslant \nu-1$:
\[
\aligned
F_{x_1x_i}
&
\,=\,
\frac{x_1^{i-1}}{(i-1)!}
+
{\rm O}_{x_2,\dots,x_{\nu-1}}(1)
+
{\rm O}_x(\nu),
\\
F_{x_1x_j}
&
\,=\,
\frac{x_1^{j-1}}{(j-1)!}
+
{\rm O}_{x_2,\dots,x_{\nu-1}}(1)
+
{\rm O}_x(\nu).
\endaligned
\]
Here, we can also write ${\rm O}_{x'}(1)$ instead
of the more precise (but not useful) 
${\rm O}_{x_2, \dots, x_{\nu-1}} (1)$.

Thirdly, to compute $F_{x_i x_j}$, consider two subcases:

\smallskip\noindent$\bullet$\,
when $4 \leqslant i+j \leqslant \nu$:
\[
F_{x_ix_j}
\,=\,
x_1^{i+j-2}\,
\tfrac{1}{(i-1)!\,(j-1)!}
+
x_1^{\nu-1}\,
\Lambda_{i,j}^\nu
+
{\rm O}_{x'}(1)
+
{\rm O}_x(\nu);
\]

\smallskip\noindent$\bullet$\,
when $\nu+1 \leqslant i+j \leqslant 2\nu - 2$:
\[
F_{x_ix_j}
\,=\,
0
+
x_1^{\nu-1}\,
\Lambda_{i,j}^\nu
+
{\rm O}_{x'}(1)
+
{\rm O}_x(\nu).
\]

The Hessian vanishing writes:
\[
\aligned
0
&
\,\equiv\,
F_{x_1x_1}\cdot F_{x_ix_j}
-
F_{x_1x_i}\cdot F_{x_1x_j}
\\
&
\,\equiv\,
\Big(
1
+
{\rm O}_{x'}(1)
+
{\rm O}_x(\nu)
\Big)
\Big(
\big\{
\begin{smallmatrix}
{\bf 0} \\ {\bf 1}
\end{smallmatrix}
\big\}\,
\tfrac{1}{(i-1)!(j-1)!}\,
x_1^{i+j-2}
+
x_1^{\nu-1}\,\Lambda_{i,j}^\nu
+
{\rm O}_{x'}(1)
+
{\rm O}_x(\nu)
\Big)
\\
&
\ \ \ \ \
-
\Big(
\tfrac{x_1^{i-1}}{(i-1)!}
+
{\rm O}_{x'}(1)
+
{\rm O}_x(\nu)
\Big)
\Big(
\tfrac{x_1^{j-1}}{(j-1)!}
+
{\rm O}_{x'}(1)
+
{\rm O}_x(\nu)
\Big),
\endaligned
\]
and it can be expanded in the various cases.

\smallskip\noindent$\bullet$\,
When $4 \leqslant i + j \leqslant \nu$, it gives
$\Lambda_{i,j}^\nu = 0$.

\smallskip\noindent$\bullet$\,
For $i + j = \nu+1$, it gives 
$\Lambda_{i,j}^\nu = \frac{1}{(i-1)!(j-1)!}$.

\smallskip\noindent$\bullet$\,
When $\nu+2 \leqslant i+j \leqslant 2\nu-2$, it gives
$\Lambda_{i,j}^\nu = 0$.
\endproof

This closes induction from $\nu-1$ to $\nu$,
while $3 \leqslant \nu \leqslant n$,
{\em cf.} Induction Hypothesis~{\ref{Induct-Hyp}}.

\SectionHead{Summary and Beyond}
{summary-and-beyond}

We started from any hypersurface $H^n \subset \R^{n+1}$
whose Hessian has constant rank $1$ and we showed,
generally, that one can let appear
monomials $\frac{x_1^2}{2}$, 
$\frac{x_1^2x_2}{2!}$, $\frac{x_1^3x_3}{3!}$,
$\frac{x_1^4x_4}{4!}$, $\frac{x_1^5x_5}{5!}$,
\dots, until the process stops, and
we proved the following

\begin{Theorem}
\label{Thm-product-nH}
Let $H^n \subset \R^{n+1}$ be a local affinely homogeneous
hypersurface
having constant Hessian rank $1$. Then there exists
an integer $1 \leqslant n_H \leqslant n$ and
affine coordinates $(x_1, \dots, x_n)$ in which:
\[
H^n
\,=\,
H^{n_H}
\times
\R_{x_{n_H+1},\dots,x_n}^{n-n_H-1}
\]
is a product of an affinely homogeneous
hypersurface $H^{n_H} \subset \R^{n_H + 1}$
times a `dumb' $\R^{n-n_H-1}$,
and is graphed as:
\[
\aligned
u
&
\,=\,
\frac{x_1^2}{2}
+
\frac{x_1^2x_2}{2}
+
\sum_{m=3}^{n_H}\,
\Big(
\frac{x_1^mx_m}{m!}
+
x_1^{m-1}
\sum_{i,j\geqslant 2
\atop
i+j=m+1}\,
\tfrac{1}{2}\,
\frac{x_ix_j}{(i-1)!(j-1)!}
+
{\rm O}_{x_2,\dots,x_{m-1}}(3)
\Big)
\\
& 
\ \ \ \ \ \ \ \ \ \ \ \ \ \ \ \ \ \ \ \ \ \ \ \ \ 
+
\sum_{m=n_H+2}^\infty\,
E^m(x_1,\dots,x_{n_H}),
\endaligned
\]
with graphing function $F = F(x_1, \dots, x_{n_H})$ 
independent of $x_{n_H+1}, \dots, x_n$.\qed
\end{Theorem}

Disregarding the product cases (branches) 
$n_H = 1, \dots, n_H = n-1$
which amount to similar
considerations in lower dimensions, we will from now on
study the class of {\sl nondegenerate}
hypersurfaces, those involving {\em all} variables
$(x_1, \dots, x_n)$:
\[
u
\,=\,
\frac{x_1^2}{2}
+
\frac{x_1^2x_2}{2}
+
\sum_{m=3}^n\,
\Big(
\frac{x_1^mx_m}{m!}
+
x_1^{m-1}
\sum_{i,j\geqslant 2
\atop
i+j=m+1}\,
\tfrac{1}{2}\,
\frac{x_ix_j}{(i-1)!(j-1)!}
+
{\rm O}_{x_2,\dots,x_{m-1}}(3)
\Big)
+
{\rm O}_x(n+2).
\]
The next (substantial) task is to examine further the remainder
${\rm O}_x(n+2)$.

We gather all remainders ${\rm O}_{x_2, \dots, x_{m-1}}(3)$
plainly as ${\rm O}_{x'} (3)$, where $x' := (x_2, \dots, x_n)$,
and we let appear the independent homogeneous
monomials of order $n+2$, namely:
\[
F_{n+2,0,\dots,0}\,
\frac{x_1^{n+1}x_1}{(n+2)!}
+
F_{n+1,1,\dots,0}\,
\frac{x_1^{n+1}x_2}{(n+1)!}
+\cdots+
F_{n+1,0,\dots,1}\,
\frac{x_1^{n+1}x_n}{(n+1)!}.
\]
A reasoning similar to that of the proof
of Proposition~{\ref{Prp-Lambda-i-i-nu-zero}}
shows that border-dependent monomials of homogeneous
order $n+2$
have the same form, hence the equation of the
hypersurface writes:
\[
\aligned
u
&
\,=\,
\frac{x_1^2}{2}
+
\frac{x_1^2x_2}{2}
+
\sum_{m=3}^n\,
\Big(
\frac{x_1^mx_m}{m!}
+
x_1^{m-1}
\sum_{i,j\geqslant 2
\atop
i+j=m+1}\,
\tfrac{1}{2}\,
\frac{x_ix_j}{(i-1)!(j-1)!}
\Big)
\\
&
\ \ \ \ \
+
F_{n+2,0,\dots,0}\,
\frac{x_1^{n+1}x_1}{(n+2)!}
+
F_{n+1,1,\dots,0}\,
\frac{x_1^{n+1}x_2}{(n+1)!}
+\cdots+
F_{n+1,0,\dots,1}\,
\frac{x_1^{n+1}x_n}{(n+1)!}
\\
&
\ \ \ \ \
+
x_1^n
\sum_{i,j\geqslant 2
\atop
i+j=n+2}\,
\tfrac{1}{2}\,
\frac{x_ix_j}{(i-1)!(j-1)!}
+
{\rm O}_{x'}(3)
+
{\rm O}_x(n+3).
\endaligned
\]

\begin{Question}
{\sl How to normalize these order $n+2$
independent coefficients
$F_{n+2,0,\dots,0}$, $F_{n+1,1,\dots,0}$,
\dots, $F_{n+1,0,\dots,1}$?}
\end{Question}

\SectionHead{Normalization of Order $n+2$ Terms}
{normalization-n-2}

Consider therefore {\em another} similar (nondegenerate) hypersurface,
in coordinates $(y_1, \dots, y_n, v)$:
\[
\aligned
v
&
\,=\,
\frac{y_1^2}{2}
+
\frac{y_1^2y_2}{2}
+
\sum_{m=3}^n\,
\Big(
\frac{y_1^my_m}{m!}
+
y_1^{m-1}
\sum_{i,j\geqslant 2
\atop
i+j=m+1}\,
\tfrac{1}{2}\,
\frac{y_iy_j}{(i-1)!(j-1)!}
\Big)
\\
&
\ \ \ \ \
+
G_{n+2,0,\dots,0}\,
\frac{y_1^{n+1}y_1}{(n+2)!}
+
G_{n+1,1,\dots,0}\,
\frac{y_1^{n+1}y_2}{(n+1)!}
+\cdots+
G_{n+1,0,\dots,1}\,
\frac{y_1^{n+1}y_n}{(n+1)!}
\\
&
\ \ \ \ \
+
y_1^n
\sum_{i,j\geqslant 2
\atop
i+j=n+2}\,
\tfrac{1}{2}\,
\frac{y_iy_j}{(i-1)!(j-1)!}
+
{\rm O}_{y'}(3)
+
{\rm O}_y(n+3).
\endaligned
\]
From Section~{\ref{normalization-order-2}},
we know what is the stability subgroup at order $2$.

\smallskip\noindent{\bf Step~1.}
Determine subgroup reduction stabilizing
terms of order $\leqslant n+1$:
\[
\overset{G_\stab^2}{
\left[
\begin{array}{ccccc}
a_{1,1} & \red{\bf 0} & \cdots & \red{\bf 0} & b_1
\\
a_{2,1} & a_{2,2} & \cdots & a_{2,n} & b_2
\\
\vdots & \vdots & \ddots & \vdots & \vdots
\\
a_{n,1} & a_{n,2} & \cdots & a_{n,n} & b_n
\\
\red{\bf 0} & \red{\bf 0} & \cdots & \red{\bf 0} & 
\red{a_{1,1}^2}
\end{array}
\right]^{\green{\bf 2}}}
\,\,\,\leadsto\,\,\,
\overset{G_\stab^{n+1}}{
\left[
\begin{array}{ccccc}
? & \red{\bf 0} & \cdots & \red{\bf 0} & ?
\\
? & ? & \cdots & ? & ?
\\
\vdots & \vdots & \ddots & \vdots & \vdots
\\
? & ? & \cdots & ? & ?
\\
\red{\bf 0} & \red{\bf 0} & \cdots & \red{\bf 0} & 
\red{a_{1,1}^{n+1}}
\end{array}
\right]^{\green{\bf 2}}}.
\]

\smallskip\noindent{\bf Step~2.}
Determine how this subgroup acts on order $n+2$ terms,
and normalize those coefficients among
$G_{n+2,1,\dots,0}$, $G_{n+1,1,\dots,0}$, \dots,
$G_{n+1,0,\dots,1}$ which can be normalized.

\smallskip

This computational task being nontrivial,
let us start concretely, in low dimensions. 
The calculations of the next Section will not be detailed,
and the remainders ${\rm O}_x( \ast)$ 
will not be written.

\SectionHead{Stabilizing Order $n+1$ Terms 
in Dimensions $n = 2, 3, 4, 5, 6$}
{stabilizing-G-n-2-3-4-5-6}

\noindent$\bullet$\,
In dimension $n = 2$:
\[
\aligned
u
&
\,=\,
\frac{x_1^2}{2}
\\
&
\ \ \ \ \
+
\frac{x_1^2x_2}{2}
\\
&
\ \ \ \ \
+
F_{4,0}\,\frac{x_1^3x_1}{24}
+
F_{3,1}\,\frac{x_1^3x_2}{6},
\endaligned
\]
the stability group is:
\[
G_{\stab}^3
\colon
\ \ \ \ \
\left[
\begin{array}{ccc}
a_{1,1} & \red{\bf 0} & \red{-a_{1,1}a_{2,1}}
\\
a_{2,1} & \red{1} & b_2
\\
\red{\bf 0} & \red{\bf 0} & \red{a_{1,1}^2}
\end{array}
\right]^{\green{\bf 3}},
\]
and its action gives:
\[
\aligned
0
&
\overset{\green{40}}{\,\,=\,\,}
-\,\frac{1}{24}\,a_{1,1}^2\,F_{4,0}
+
\frac{1}{24}\,a_{1,1}^4\,G_{4,0}
+
\frac{1}{6}\,a_{1,1}^3\,a_{2,1}\,
G_{3,1}
+
\frac{1}{8}\,a_{1,1}^2\,a_{2,1}^2
+
\frac{1}{4}\,a_{1,1}^2\,\boxed{b_2},
\\
0
&
\overset{\green{31}}{\,\,=\,\,}
-\,\frac{1}{6}\,a_{1,1}^2\,F_{3,1}
+
\frac{1}{6}\,a_{1,1}^3\,G_{3,1}.
\endaligned
\]
The free group parameter $b_2$ can be used
to normalize $G_{4,0} := 0$.

\medskip\noindent$\bullet$\,
In dimension $n = 3$:
\[
\aligned
u
&
\,=\,
\frac{x_1^2}{2}
\\
&
\ \ \ \ \
+
\frac{x_1^2x_2}{2}
\\
&
\ \ \ \ \
+
\frac{x_1^3x_3}{6}
+
\frac{x_1^2x_2^2}{2}
\\
&
\ \ \ \ \
+
F_{5,0,0}\,\frac{x_1^4x_1}{120}
+
F_{4,1,0}\,\frac{x_1^4x_2}{24}
+
F_{4,0,1}\,\frac{x_1^4x_3}{24}
+
\frac{x_1^3x_2x_3}{2},
\endaligned
\]
the stability group is:
\[
G_{\stab}^4
\colon
\ \ \ \ \
\left[
\def\arraystretch{1.25}
\begin{array}{cccc}
a_{1,1} & \red{\bf 0} & \red{\bf 0} & \red{-a_{1,1}a_{2,1}}
\\
a_{2,1} & \red{1} & \red{\bf 0} & 
\red{-\frac{1}{2}a_{2,1}^2-\frac{2}{3}a_{1,1}a_{3,1}}
\\
a_{3,1} & \red{\bf 0} & \red{\frac{1}{a_{1,1}}} & b_3
\\
\red{\bf 0} & \red{\bf 0} & \red{\bf 0} & \red{a_{1,1}^2}
\end{array}
\right]^{\green{\bf 4}},
\]
and its action gives:
\[
\aligned
0
&
\overset{\green{500}}{\,\,=\,\,}
-\,\frac{1}{120}\,a_{1,1}^2\,F_{5,0,0}
+
\frac{1}{120}\,a_{1,1}^5\,G_{5,0,0}
+
\frac{1}{24}\,a_{1,1}^4\,a_{2,1}\,G_{4,1,0}
+
\frac{1}{12}\,a_{1,1}^3\,a_{2,1}\,a_{3,1}
\\
&
\ \ \ \ \ \ \ \ \ \ \ \ \ \ \ \ \ \ \ \ \ \ \ \ \ \ \ \ \ \ \ \ \ \ \
\ \ \ \ \ \ \ \ \ \ \ \ \ \ \ \ \ \ \ \ \ \ \ \ \ \ \ \ \ \ \ \ \ \ \
\ \ \ \ \ \ \ \ \ \ \ \ \ \ \ \ \ \ \ \ 
+
\frac{1}{24}\,a_{1,1}^4\,a_{3,1}\,G_{4,0,1}
+
\frac{1}{12}\,a_{1,1}^3\,\boxed{b_3}
\\
0
&
\overset{\green{410}}{\,\,=\,\,}
-\,
\frac{1}{24}\,a_{1,1}^2\,F_{4,1,0}
+
\frac{1}{24}\,a_{1,1}^4\,G_{4,1,0},
\\
0
&
\overset{\green{401}}{\,\,=\,\,}
-\,\frac{1}{24}\,a_{1,1}^2\,F_{4,0,1}
+
\frac{1}{24}\,a_{1,1}^3\,G_{4,0,1}
+
\frac{1}{12}\,a_{1,1}^2\,\boxed{a_{2,1}}.
\endaligned
\]
The free group parameter $a_{2,1}$
can be used to normalize $G_{4,0,1} := 0$,
and the free group parameter
$b_3$ can be used to normalize $G_{5,0,0} := 0$.

\medskip\noindent$\bullet$\,
In dimension $n = 4$:
\[
\aligned
u
&
\,=\,
\frac{x_1^2}{2}
\\
&
\ \ \ \ \
+
\frac{x_1^2x_2}{2}
\\
&
\ \ \ \ \
+
\frac{x_1^3x_3}{6}
+
\frac{x_1^2x_2^2}{2}
\\
&
\ \ \ \ \
+
\frac{x_1^4x_4}{24}
+
\frac{x_1^2x_2^3}{2}
+
\frac{x_1^3x_2x_3}{2}
\\
&
\ \ \ \ \
+
F_{6,0,0,0}\,\frac{x_1^5x_1}{720}
+
F_{5,1,0,0}\,\frac{x_1^5x_2}{120}
+
F_{5,0,1,0}\,\frac{x_1^5x_3}{120}
+
F_{5,0,0,1}\,\frac{x_1^5x_4}{120}
+
\frac{x_1^2x_2^4}{2}
+
x_1^3x_2^2x_3
+
\frac{x_1^4x_2x_4}{6}
+
\frac{x_1^4x_3^2}{8},
\endaligned
\]
the stability group is:
\[
G_{\stab}^5
\colon
\ \ \ \ \
\left[
\def\arraystretch{1.25}
\begin{array}{ccccc}
a_{1,1} & \red{\bf 0} & \red{\bf 0} & \red{\bf 0} & 
\red{-a_{1,1}a_{2,1}}
\\
a_{2,1} & \red{1} & \red{\bf 0} & \red{\bf 0} &
\red{-\frac{1}{2}a_{2,1}^2-\frac{2}{3}a_{1,1}a_{3,1}}
\\
a_{3,1} & \red{\bf 0} & \red{\frac{1}{a_{1,1}}} & \red{\bf 0} &
\red{-\frac{1}{2}a_{1,1}a_{4,1}-a_{2,1}a_{3,1}} 
\\
a_{4,1} & \red{\bf 0} & \red{\frac{-2a_{2,1}}{a_{1,1}^2}} &
\red{\frac{1}{a_{1,1}^2}} & b_4 
\\
\red{\bf 0} & \red{\bf 0} & \red{\bf 0} & \red{\bf 0} & \red{a_{1,1}^2}
\end{array}
\right]^{\green{\bf 5}},
\]
and its action gives:
\[
\aligned
0
&
\overset{\green{6000}}{\,\,=\,\,}
-\,\frac{1}{720}\,a_{1,1}^2\,F_{6,0,0,0}
+
\frac{1}{720}\,a_{1,1}^6\,G_{6,0,0,0}
+
\frac{1}{72}\,a_{1,1}^4\,a_{3,1}^2
+
\frac{1}{48}\,a_{1,1}^4\,a_{2,1}\,a_{4,1}
\\
&
\ \ \ \ \ \ \ \ \ \ \ \ \ \ \ \ \ \ \ \ \ \
+
\frac{1}{120}\,G_{5,1,0,0}\,a_{1,1}^5\,a_{2,1}
+
\frac{1}{120}\,G_{5,0,1,0}\,a_{1,1}^5\,a_{3,1}
+
\frac{1}{120}\,G_{5,0,0,1}\,a_{1,1}^5\,a_{4,1}
+
\frac{1}{48}\,a_{1,1}^4\,\boxed{b_4}
\\
0
&
\overset{\green{5100}}{\,\,=\,\,}
-\,
\frac{1}{120}\,a_{1,1}^2\,F_{5,1,0,0}
+
\frac{1}{120}\,a_{1,1}^5\,G_{5,1,0,0},
\\
0
&
\overset{\green{5010}}{\,\,=\,\,}
-\,
\frac{1}{120}\,a_{1,1}^2\,F_{5,0,1,0}
+
\frac{1}{120}\,a_{1,1}^4\,G_{5,0,1,0}
-
\frac{1}{24}\,a_{1,1}^2\,a_{2,1}^2
-
\frac{1}{60}\,a_{1,1}^3\,a_{2,1}
+
\frac{1}{36}\,a_{1,1}^3\,\boxed{a_{3,1}},
\\
0
&
\overset{\green{5001}}{\,\,=\,\,}
-\,\frac{1}{120}\,a_{1,1}^2\,F_{5,0,0,1}
+
\frac{1}{120}\,a_{1,1}^3\,G_{5,0,0,1}
+
\frac{1}{24}\,a_{1,1}^2\,\boxed{a_{2,1}}.
\endaligned
\]
The free group parameters $a_{2,1}$, $a_{3,1}$, $b_4$
can be used to normalize $G_{5,0,0,1} := 0$, 
$G_{5,0,1,0} := 0$, $G_{6,0,0,0} := 0$.

\medskip\noindent$\bullet$\,
In dimension $n = 5$:
\[
\aligned
u
&
\,=\,
\frac{x_1^2}{2}
\\
&
\ \ \ \ \
+
\frac{x_1^2x_2}{2}
\\
&
\ \ \ \ \
+
\frac{x_1^3x_3}{6}
+
\frac{x_1^2x_2^2}{2}
\\
&
\ \ \ \ \
+
\frac{x_1^4x_4}{24}
+
\frac{x_1^2x_2^3}{2}
+
\frac{x_1^3x_2x_3}{2}
\\
&
\ \ \ \ \
+
\frac{x_1^5x_5}{120}
+
\frac{x_1^2x_2^4}{2}
+
x_1^3x_2^2x_3
+
\frac{x_1^4x_2x_4}{6}
+
\frac{x_1^4x_3^2}{8}
\\
&
\ \ \ \ \
+
F_{7,0,0,0,0}\,\frac{x_1^6x_1}{5040}
+
F_{6,1,0,0,0}\,\frac{x_1^6x_2}{720}
+
F_{6,0,1,0,0}\,\frac{x_1^6x_3}{720}
+
F_{6,0,0,1,0}\,\frac{x_1^6x_4}{720}
+
F_{6,0,0,0,1}\,\frac{x_1^6x_5}{720}
\\
&
\ \ \ \ \
+
\tfrac{1}{2}\,x_1^2x_2^5
+
\tfrac{5}{3}\,x_1^3x_2^3x_3
+
\tfrac{5}{8}\,x_1^4x_2x_3^2
+
\tfrac{5}{12}\,x_1^4x_2^2x_4
+
\tfrac{1}{12}\,x_1^5x_3x_4
+
\tfrac{1}{24}\,x_1^5x_2x_5,
\endaligned
\]
the stability group is:
\[
G_{\stab}^6
\colon
\ \ \ \ \
\left[
\def\arraystretch{1.25}
\begin{array}{cccccc}
a_{1,1} & \red{\bf 0} & \red{\bf 0} & \red{\bf 0} & \red{\bf 0} & 
\red{-a_{1,1}a_{2,1}}
\\
a_{2,1} & \red{1} & \red{\bf 0} & \red{\bf 0} & \red{\bf 0} &
\red{-\frac{1}{2}a_{2,1}^2-\frac{2}{3}a_{1,1}a_{3,1}}
\\
a_{3,1} & \red{\bf 0} & \red{\frac{1}{a_{1,1}}} & \red{\bf 0} & 
\red{\bf 0} &
\red{-a_{2,1}a_{3,1}-\frac{1}{2}a_{1,1}a_{4,1}} 
\\
a_{4,1} & \red{\bf 0} & \red{-\frac{2a_{2,1}}{a_{1,1}^2}} & 
\red{\frac{1}{a_{1,1}^2}} & \red{\bf 0} & 
\red{-a_{2,1}a_{4,1}-\frac{2}{3}a_{3,1}^2-\frac{2}{5}a_{1,1}a_{5,1}}
\\
a_{5,1} & \red{\bf 0} & 
\red{\frac{5a_{2,1}^2}{a_{1,1}^3}
-\frac{10}{3}\frac{a_{3,1}}{a_{1,1}^2}} & 
\red{-\frac{5a_{2,1}}{a_{1,1}^3}} &
\red{\frac{1}{a_{1,1}^3}} & 
b_5
\\
\red{\bf 0} & \red{\bf 0} & \red{\bf 0} & \red{\bf 0} & \red{\bf 0} & 
\red{a_{1,1}^2}
\end{array}
\right]^{\green{\bf 6}},
\]
and its action gives:
\[
\footnotesize
\aligned
0
&
\overset{\green{70000}}{\,\,=\,\,}
-\,\frac{1}{5040}\,a_{1,1}^2\,F_{7,0,0,0,0}
+
\frac{1}{5040}\,a_{1,1}^7\,G_{7,0,0,0,0}
+
\frac{1}{144}\,a_{1,1}^5\,a_{3,1}\,a_{4,1}
+
\frac{1}{240}\,a_{1,1}^5\,a_{2,1}\,a_{5,1}
\\
&
\ \ \ \ \
+
\frac{1}{720}\,G_{6,1,0,0,0}\,a_{1,1}^6\,a_{2,1}
+
\frac{1}{720}\,G_{6,0,1,0,0}\,a_{1,1}^6\,a_{3,1}
+
\frac{1}{720}\,G_{6,0,0,1,0}\,a_{1,1}^6\,a_{4,1}
+
\frac{1}{720}\,G_{6,0,0,0,1}\,a_{1,1}^6\,a_{5,1}
+
\frac{1}{240}\,a_{1,1}^5\,\boxed{b_5}
\\
0
&
\overset{\green{61000}}{\,\,=\,\,}
-\,
\frac{1}{720}\,a_{1,1}^2\,F_{6,1,0,0,0}
+
\frac{1}{720}\,a_{1,1}^6\,G_{6,1,0,0,0},
\\
0
&
\overset{\green{60100}}{\,\,=\,\,}
-\,
\frac{1}{720}\,a_{1,1}^2\,F_{6,0,1,0,0}
+
\frac{1}{720}\,a_{1,1}^5\,G_{6,0,1,0,0}
-
\frac{1}{36}\,a_{1,1}^3\,a_{2,1}\,a_{3,1}
+
\frac{1}{48}\,a_{1,1}^2\,a_{2,1}^3
+
\frac{1}{144}\,a_{1,1}^3\,a_{2,1}^2\,G_{6,0,0,0,1}
\\
&
\ \ \ \ \ 
-
\frac{1}{216}\,a_{1,1}^4\,a_{3,1}\,G_{6,0,0,0,1}
-
\frac{1}{360}\,a_{1,1}^4\,a_{2,1}\,G_{6,0,0,1,0}
+
\frac{1}{144}\,a_{1,1}^4\,\boxed{a_{4,1}},
\\
0
&
\overset{\green{60010}}{\,\,=\,\,}
-\,
\frac{1}{720}\,a_{1,1}^2\,F_{6,0,0,1,0}
+
\frac{1}{720}\,a_{1,1}^4\,G_{6,0,0,1,0}
-
\frac{1}{32}\,a_{1,1}^2\,a_{2,1}^2
-
\frac{1}{144}\,a_{1,1}^3\,a_{2,1}\,G_{6,0,0,0,1}
+
\frac{1}{72}\,a_{1,1}^3\,\boxed{a_{3,1}},
\\
0
&
\overset{\green{60001}}{\,\,=\,\,}
-\,\frac{1}{720}\,a_{1,1}^2\,F_{6,0,0,0,1}
+
\frac{1}{720}\,a_{1,1}^3\,G_{6,0,0,0,1}
+
\frac{1}{80}\,a_{1,1}^2\,\boxed{a_{2,1}}.
\endaligned
\]
The free group parameters $a_{2,1}$, $a_{3,1}$, $a_{4,1}$, $b_4$
can be used to normalize $G_{6,0,0,0,1} := 0$, 
$G_{6,0,0,1,0} := 0$, $G_{6,0,1,0,0} := 0$, $G_{7,0,0,0,0} := 0$.

\medskip\noindent$\bullet$\,
In dimension $n = 6$:
\[
\scriptsize
\aligned
u
&
\,=\,
\frac{x_1^2}{2}
\\
&
\ \ \ \ \
+
\frac{x_1^2x_2}{2}
\\
&
\ \ \ \ \
+
\frac{x_1^3x_3}{6}
+
\frac{x_1^2x_2^2}{2}
\\
&
\ \ \ \ \
+
\frac{x_1^4x_4}{24}
+
\frac{x_1^2x_2^3}{2}
+
\frac{x_1^3x_2x_3}{2}
\\
&
\ \ \ \ \
+
\frac{x_1^5x_5}{120}
+
\frac{x_1^2x_2^4}{2}
+
x_1^3x_2^2x_3
+
\frac{x_1^4x_2x_4}{6}
+
\frac{x_1^4x_3^2}{8}
\\
&
\ \ \ \ \
+
\frac{x_1^6x_6}{720}
+
\frac{x_1^2x_2^5}{2}
+
\tfrac{5}{3}\,x_1^3x_2^3x_3
+
\tfrac{5}{8}\,x_1^4x_2x_3^2
+
\tfrac{5}{12}\,x_1^4x_2^2x_4
+
\tfrac{1}{12}\,x_1^5x_3x_4
+
\tfrac{1}{24}\,x_1^5x_2x_5
\\
&
\ \ \ \ \
+
F_{8,0,0,0,0,0}\,\frac{x_1^7x_1}{40320}
+
F_{7,1,0,0,0,0}\,\frac{x_1^7x_2}{5040}
+
F_{7,0,1,0,0,0}\,\frac{x_1^7x_3}{5040}
+
F_{7,0,0,1,0,0}\,\frac{x_1^7x_4}{5040}
+
F_{7,0,0,0,1,0}\,\frac{x_1^7x_5}{5040}
+
F_{7,0,0,0,0,1}\,\frac{x_1^7x_6}{5040}
\\
&
\ \ \ \ \
+
\tfrac{1}{8}\,x_1^5x_3^3
+
\tfrac{1}{72}\,x_1^6x_4^2
+
\tfrac{5}{2}\,x_1^3x_2^4x_3
+
\tfrac{15}{8}\,x_1^4x_2^2x_3^2
+
\tfrac{5}{6}\,x_1^4x_2^3x_4
+
\tfrac{1}{8}\,x_1^5x_2^2x_5
+
\tfrac{1}{48}\,x_1^6x_3x_5
+
\tfrac{1}{2}\,x_1^5x_2x_3x_4
+
\tfrac{1}{2}\,x_1^2x_2^6
+
\tfrac{1}{120}\,x_1^6x_2x_6,
\endaligned
\]
the stability group is:
\[
\footnotesize
\aligned
\left[
\def\arraystretch{1.25}
\begin{array}{ccccccc}
a_{1,1} & \red{\bf 0} & \red{\bf 0} & \red{\bf 0} & \red{\bf 0} & 
\red{\bf 0} & \red{-a_{1,1}a_{2,1}}
\\
a_{2,1} & \red{1} & \red{\bf 0} & \red{\bf 0} & \red{\bf 0} & 
\red{\bf 0} &
\red{-\frac{1}{2}a_{2,1}^2-\frac{2}{3}a_{1,1}a_{3,1}}
\\
a_{3,1} & \red{\bf 0} & \red{\frac{1}{a_{1,1}}} & \red{\bf 0} & 
\red{\bf 0} & \red{\bf 0} & 
\red{-a_{2,1}a_{3,1}-\frac{1}{2}a_{1,1}a_{4,1}} 
\\
a_{4,1} & \red{\bf 0} & \red{-\frac{2a_{2,1}}{a_{1,1}^2}} & 
\red{\frac{1}{a_{1,1}^2}} & \red{\bf 0} & \red{\bf 0} & 
\red{-a_{2,1}a_{4,1}-\frac{2}{3}a_{3,1}^2-\frac{2}{5}a_{1,1}a_{5,1}}
\\
a_{5,1} & \red{\bf 0} & 
\red{\frac{5a_{2,1}^2}{a_{1,1}^3}
-\frac{10}{3}\frac{a_{3,1}}{a_{1,1}^2}} & 
\red{-\frac{5a_{2,1}}{a_{1,1}^3}} &
\red{\frac{1}{a_{1,1}^3}} & 
\red{\bf 0} &
\red{-a_{2,1}a_{5,1}-\frac{1}{3}a_{1,1}a_{6,1}
-\frac{5}{3}a_{3,1}a_{4,1}}
\\
a_{6,1} & \red{\bf 0} & 
\red{20\frac{a_{2,1}a_{3,1}}{a_{1,1}^3}
-15\frac{a_{2,1}^3}{a_{1,1}^4}-5\frac{a_{4,1}}{a_{2,1}^2}} & 
\red{\frac{45}{2}\frac{a_{2,1}^2}{a_{1,1}^4}
-10\frac{a_{3,1}}{a_{1,1}^3}} &
\red{-9\frac{a_{2,1}}{a_{1,1}^4}} &
\red{\frac{1}{a_{1,1}^4}} & 
b_6
\\
\red{\bf 0} & \red{\bf 0} & \red{\bf 0} & \red{\bf 0} & 
\red{\bf 0} & \red{\bf 0} & 
\red{a_{1,1}^2}
\end{array}
\right]^{\green{\bf 7}},
\endaligned
\]
and its action gives:
\[
\!\!\!\!\!\!\!\!\!\!\!\!\!\!\!\!\!\!\!\!\!\!\!\!\!
\scriptsize
\aligned
0
&
\overset{\green{800000}}{\,\,=\,\,}
-\,\frac{1}{40320}\,a_{1,1}^2\,F_{8,0,0,0,0,0}
+
\frac{1}{40320}\,a_{1,1}^7\,G_{8,0,0,0,0}
+
\\
&
\ \ \ \ \
+
\frac{1}{5040}\,G_{7,1,0,0,0,0}\,a_{1,1}^7\,a_{2,1}
+
\frac{1}{5040}\,G_{7,0,1,0,0,0}\,a_{1,1}^7\,a_{3,1}
+
\frac{1}{5040}\,G_{7,0,0,1,0,0}\,a_{1,1}^7\,a_{4,1}
+
\frac{1}{5040}\,G_{7,0,0,0,1,0}\,a_{1,1}^7\,a_{5,1}
+
\frac{1}{5040}\,G_{7,0,0,0,0,1}\,a_{1,1}^7\,a_{6,1}
\\
&
\ \ \ \ \
+
\frac{1}{1152}\,a_{1,1}^6\,a_{4,1}^2
+
\frac{1}{1440}\,a_{1,1}^6\,a_{2,1}\,a_{6,1}
+
\frac{1}{720}\,a_{1,1}^6\,a_{3,1}\,a_{5,1}
+
\frac{1}{1440}\,a_{1,1}^6\,\boxed{b_6}
\\
0
&
\overset{\green{710000}}{\,\,=\,\,}
-\,
\frac{1}{5040}\,a_{1,1}^2\,F_{7,1,0,0,0,0}
+
\frac{1}{5040}\,a_{1,1}^7\,G_{7,1,0,0,0,0},
\\
0
&
\overset{\green{701000}}{\,\,=\,\,}
-\,
\frac{1}{5040}\,a_{1,1}^2\,F_{7,0,1,0,0,0}
+
\frac{1}{5040}\,a_{1,1}^6\,G_{7,0,1,0,0,0}
-
\frac{1}{144}\,a_{1,1}^4\,a_{2,1}\,a_{4,1}
+
\frac{1}{48}\,a_{1,1}^3\,a_{2,1}^2\,a_{3,1}
-
\frac{1}{96}\,a_{1,1}^2\,a_{2,1}^4
-
\frac{1}{216}\,a_{1,1}^4\,a_{3,1}^2
\\
&
\ \ \ \ \
-
\frac{1}{2520}\,a_{1,1}^5\,a_{2,1}\,G_{7,0,0,1,0,0}
+
\frac{1}{252}\,a_{1,1}^4\,a_{2,1}\,a_{3,1}\,G_{7,0,0,0,0,1}
-
\frac{1}{336}\,a_{1,1}^3\,a_{2,1}^3\,G_{7,0,0,0,0,1}
-
\frac{1}{1008}\,a_{1,1}^4\,a_{4,1}\,G_{7,0,0,0,0,1}
-
\frac{1}{1008}\,a_{1,1}^4\,a_{2,1}^2\,G_{7,0,0,0,1,0}
\\
&
\ \ \ \ \
-
\frac{1}{1512}\,a_{1,1}^5\,a_{3,1}\,G_{7,0,0,0,1,0}
+
\frac{1}{720}\,a_{1,1}^5\,\boxed{a_{5,1}},
\\
0
&
\overset{\green{700100}}{\,\,=\,\,}
-\,
\frac{1}{5040}\,a_{1,1}^2\,F_{7,0,0,1,0,0}
+
\frac{1}{5040}\,a_{1,1}^5\,G_{7,0,0,1,0,0}
+
\frac{1}{224}\,a_{1,1}^3\,a_{2,1}^2\,G_{7,0,0,0,0,1}
-
\frac{1}{504}\,a_{1,1}^4\,a_{3,1}\,G_{7,0,0,0,0,1}
+
\frac{1}{48}\,a_{1,1}^2\,a_{2,1}^3
\\
&
\ \ \ \ \
-
\frac{1}{48}\,a_{1,1}^3\,a_{2,1}\,a_{3,1}
-
\frac{1}{1008}\,a_{1,1}^4\,a_{2,1}\,G_{7,0,0,0,1,0}
+
\frac{1}{288}\,a_{1,1}^4\,\boxed{a_{4,1}},
\\
0
&
\overset{\green{700010}}{\,\,=\,\,}
-\,\frac{1}{5040}\,a_{1,1}^2\,F_{7,0,0,0,1,0}
+
\frac{1}{5040}\,a_{1,1}^4\,G_{7,0,0,0,1,0}
-
\frac{1}{560}\,a_{1,1}^3\,a_{2,1}\,G_{7,0,0,0,0,1}
-
\frac{1}{80}\,a_{1,1}^2\,a_{2,1}^2
+
\frac{1}{240}\,a_{1,1}^3\,\boxed{a_{3,1}},
\\
0
&
\overset{\green{700001}}{\,\,=\,\,}
-\,\frac{1}{5040}\,a_{1,1}^2\,F_{7,0,0,0,0,1}
+
\frac{1}{5040}\,a_{1,1}^3\,G_{7,0,0,0,0,1}
+
\frac{1}{360}\,a_{1,1}^2\,\boxed{a_{2,1}}.
\endaligned
\]
The free group parameters $a_{2,1}$, $a_{3,1}$, $a_{4,1}$, 
$a_{5,1}$, $b_4$
can be used to normalize $G_{7,0,0,0,0,1} := 0$, 
$G_{7,0,0,0,1,0} := 0$, $G_{7,0,0,1,0,0} := 0$, 
$G_{7,0,1,0,0,0} := 0$, 
$G_{8,0,0,0,0,0} := 0$.

\medskip

Instead of attempting to {\em dominate}
the combinatorics of such formulas in any dimension
$n \geqslant 2$, we will {\em infinitesimalize} 
the determination of the stability group at order
$n+1$, and also, we will {\em infinitesimalize} 
its action on coefficients of order $n+2$.

\SectionHead{Tangency at Order $2$ in Dimension $n$}
{tangency-order-2-dimension-n}

Before starting, 
in any dimension $n \geqslant 2$ and in continuation
with Section~{\ref{tangency-order-2}},
let us examine tangency of $L$
up to order $2$ to the
hypersurface $u = \frac{1}{2}\, x_1^2 + {\rm O}_x(3)$.
Thus, in~({\ref{L-order-2}}), 
put $T_1 := 0$, \dots, $T_n := 0$:
\[
\aligned
L
&
\,=\,\,\,\,\,\,
\Big(
A_{1,1}\,x_1
+\cdots+
A_{1,n}\,x_n
+
B_1\,u
\Big)\,
\tfrac{\partial}{\partial x_1}
\,+
\notag
\\
&
\ \ \ \ \
+
\Big(
A_{2,1}\,x_1
+\cdots+
A_{2,n}\,x_n
+
B_2\,u
\Big)\,
\tfrac{\partial}{\partial x_2}
\,+
\notag
\\
&
\ \ \ \ \
+
\cdots\cdots\cdots\cdots\cdots\cdots\cdots\cdots\cdots\cdots\cdots
\cdot\cdot
\,+
\\
&
\ \ \ \ \
+
\Big(
A_{n,1}\,x_1
+\cdots+
A_{n,n}\,x_n
+
B_n\,u
\Big)\,
\tfrac{\partial}{\partial x_n}
\,+
\notag
\\
&
\ \ \ \ \
+
\Big(
\ \ \ \ \ \ \ \ \ \ \ \ \ \ \ \ \ \ \ \ \ \ \ \ \ \ \ \ \ \ \ \ \ \ \ 
\ \ \ \ 
+
D\,u
\Big)\,
\tfrac{\partial}{\partial u}.
\endaligned
\]

\begin{Lemma}
\label{Lm-tangency-L-order-2-dim-n}
Tangency $\pi^2 \big( L (-u+F) \big\vert_{u=F} \big)$ 
up to order $2$ holds iff the coefficients matrix
of $L$ reads:
\[
\left[
\begin{array}{ccccc}
A_{1,1} & \red{\bf 0} & \cdots & \red{\bf 0} & B_1
\\
A_{2,1} & A_{2,2} & \cdots & A_{2,n} & B_2
\\
\vdots & \vdots & \ddots & \vdots & \vdots
\\
A_{n,1} & A_{n,2} & \cdots & A_{n,n} & B_n
\\
\red{\bf 0} & \red{\bf 0} & \cdots & \red{\bf 0} & \red{2A_{1,1}}
\end{array}
\right]^{\green{\bf 2}}.
\]
\end{Lemma}

\proof
Write the graph as $u = \frac{1}{2}\, x_1^2 + {\rm O}_x(3)$,
and compute modulo ${\rm O}_x(3)$:
\[
\aligned
0
&
\,\equiv\,
-\,D\,
\big(
\tfrac{1}{2}\,x_1^2
+
{\rm O}_x(3)
\big)
\\
&
\ \ \ \ \
+
\Big(
A_{1,1}\,x_1
+
A_{1,2}\,x_2
+\cdots+
A_{1,n}\,x_n
+
B_1\,{\rm O}_x(2)
\Big)\,
\big(
x_1
+
{\rm O}_x(2)
\big)
\\
&
\ \ \ \ \
+
\Big(
A_{2,1}\,x_1
+
A_{2,2}\,x_2
+\cdots+
A_{2,n}\,x_n
+
B_2\,{\rm O}_x(2)
\Big)\,
\big(
0
+
{\rm O}_x(2)
\big)
\\
&
\ \ \ \ \
+
\cdots\cdots\cdots\cdots\cdots\cdots\cdots\cdots\cdots\cdots\cdots
\cdots\cdots\cdots\cdots\cdots\cdots\cdots\cdots
\\
&
\ \ \ \ \
+
\Big(
A_{n,1}\,x_1
+
A_{n,2}\,x_2
+\cdots+
A_{n,n}\,x_n
+
B_n\,{\rm O}_x(2)
\Big)\,
\big(
0
+
{\rm O}_x(2)
\big).
\endaligned
\]
The coefficients of $x_1^2$, of $x_1x_2$, \dots, of
$x_1 x_n$ must vanish, which concludes:
\[
D
\,=\,
\red{2\,A_{1,1}},
\ \ \ \ \
A_{1,2}
\,=\,
\red{\bf 0},
\ \ \ \ \
\dots\dots,
\ \ \ \ \ 
A_{1,n}
\,=\,
\red{\bf 0}.
\qedhere
\]
\endproof

\SectionHead{Tangency in Dimensions $n = 2, 3, 4, 5, 6$}
{tangency-low-dimensions}

Before endeavoring to treat the crucial order $n+1$ 
in general dimension 
$n \geqslant 2$,
let us show what formulas exist in low dimensions.

\medskip\noindent$\bullet$\,
In dimension $n = 2$, with:
\[
\aligned
u
&
\,=\,
\tfrac{x_1^2}{2}
\\
&
\ \ \ \ \
+
\tfrac{x_1^2x_2}{2},
\endaligned
\ \ \ \ \ \ \ \ \ \ \ \ \ \ \ \ \ \ \ \
\aligned
L
&
\,=\,
\big(
A_{1,1}\,x_1+A_{1,2}\,x_2+B_1\,u
\big)\,\partial_{x_1}
\\
&
\ \ \ \ \
+
\big(
A_{2,1}\,x_1+A_{2,2}\,x_2+B_2\,u
\big)\,\partial_{x_2}
\\
&
\ \ \ \ \
+
\big(
C_1\,x_1+C_2\,x_2+D\,u
\big)\,\partial_u,
\endaligned
\]
the tangency equation in orders $\leqslant 3$:
\[
\aligned
0
&
\,\equiv\,
-\,C_1\,x_1
-
C_2\,x_2
-
D\,\big(
\tfrac{x_1^2}{2}
+
\tfrac{x_1^2x_2}{2}
\big)
\\
&
\ \ \ \ \
+
\big(
A_{1,1}\,x_1+A_{1,2}\,x_2+B_1\,\tfrac{x_1^2}{2}
\big)\,
\big(x_1+x_1x_2
\big)
\\
&
\ \ \ \ \
+
\big(
A_{2,1}\,x_1+A_{2,2}\,x_2
\big)\,
\big(
\tfrac{x_1^2}{2}
\big),
\endaligned
\]
gives at order $1$ as we know:
\[
C_1
\,:=\,
0,
\ \ \ \ \ \ \ \ \ \ \ \ \ \ \ \ \ \ \ \
C_2
\,:=\,
0,
\]
then at order $2$ as we know:
\[
D
\,:=\,
2\,A_{1,1},
\ \ \ \ \ \ \ \ \ \ \ \ \ \ \ \ \ \ \ \
A_{1,2}
\,:=\,
0,
\]
and it remains:
\[
\aligned
0
&
\,\equiv\,
-
2\,A_{1,1}\,\big(
\tfrac{x_1^2}{2}
+
\tfrac{x_1^2x_2}{2}
\big)
\\
&
\ \ \ \ \
+
\big(
A_{1,1}\,x_1+B_1\,\tfrac{x_1^2}{2}
\big)\,
\big(x_1+x_1x_2
\big)
\\
&
\ \ \ \ \
+
\big(
A_{2,1}\,x_1+A_{2,2}\,x_2
\big)\,
\big(
\tfrac{x_1^2}{2}
\big)
\\
&
\,\equiv\,
x_1^3\,
\big[
\tfrac{1}{2}\,B_1
+
\tfrac{1}{2}\,A_{2,1}
\big]
+
x_1^2x_2\,
\big[
-\,A_{1,1}
+
A_{1,1}
+
\tfrac{1}{2}\,A_{2,2}
\big],
\endaligned
\]
which gives at order $3$:
\[
B_1
\,:=\,
-\,A_{2,1},
\ \ \ \ \ \ \ \ \ \ \ \ \ \ \ \ \ \ \ \
A_{2,2}
\,:=\,
0.
\]

The reductions of the coefficients matrix of $L$ read:
\[
\left[
\begin{array}{ccc}
A_{1,1} & A_{1,2} & B_1
\\
A_{2,1} & A_{2,2} & B_2
\\
C_1 & C_2 & D
\end{array}
\right]^{\green{0}}
\leadsto
\left[
\begin{array}{ccc}
A_{1,1} & A_{1,2} & B_1
\\
A_{2,1} & A_{2,2} & B_2
\\
\red{\bf 0} & \red{\bf 0} & D
\end{array}
\right]^{\green{1}}
\leadsto
\left[
\begin{array}{ccc}
A_{1,1} & \red{\bf 0} & B_1
\\
A_{2,1} & A_{2,2} & B_2
\\
\red{\bf 0} & \red{\bf 0} & \red{2A_{1,1}}
\end{array}
\right]^{\green{2}}
\leadsto
\left[
\begin{array}{ccc}
A_{1,1} & \red{\bf 0} & \red{-A_{2,1}}
\\
A_{2,1} & \red{\bf 0} & B_2
\\
\red{\bf 0} & \red{\bf 0} & D
\end{array}
\right]^{\green{3}}.
\]
 
\medskip\noindent$\bullet$\,
In dimension $n = 3$, with:
\[
\aligned
u
&
\,=\,
\tfrac{x_1^2}{2}
\\
&
\ \ \ \ \
+
\tfrac{x_1^2x_2}{2}
\\
&
\ \ \ \ \
+
\tfrac{x_1^3x_3}{6}
+
\tfrac{x_1^2x_2^2}{2},
\endaligned
\ \ \ \ \ \ \ \ \ \ \ \ \ \ \ \ \ \ \ \
\aligned
L
&
\,=\,
\big(
A_{1,1}\,x_1+A_{1,2}\,x_2+A_{1,3}\,x_3+B_1\,u
\big)\,\partial_{x_1}
\\
&
\ \ \ \ \
+
\big(
A_{2,1}\,x_1+A_{2,2}\,x_2+A_{2,3}\,x_3+B_2\,u
\big)\,\partial_{x_2}
\\
&
\ \ \ \ \
+
\big(
A_{3,1}\,x_1+A_{3,2}\,x_2+A_{3,3}\,x_3+B_3\,u
\big)\,\partial_{x_3}
\\
&
\ \ \ \ \
+
\big(
C_1\,x_1+C_2\,x_2+C_3\,x_3+D\,u
\big)\,\partial_u,
\endaligned
\]
the tangency equation up to order $4$, 
whose expansion can be done,
is:
\[
\aligned
0
&
\,\equiv\,
-\,C_1\,x_1-C_2\,x_2-C_3\,x_3
-
D\,\big(
\tfrac{x_1^2}{2}
+
\tfrac{x_1^2x_2}{2}
+
\tfrac{x_1^3x_3}{6}
\big)
\\
&
\ \ \ \ \
+
\big(
A_{1,1}x_1+A_{1,2}x_2+A_{1,3}x_3
+
B_1
\big[
\tfrac{x_1^2}{2}+\tfrac{x_1^2x_2}{2}
\big]
\big)\,
\big(
x_1+x_1x_2+\tfrac{x_1^2x_3}{2}
\big)
\\
&
\ \ \ \ \
+
\big(
A_{2,1}x_1+A_{2,2}x_2+A_{2,3}x_3
+
B_2
\big[
\tfrac{x_1^2}{2}
\big]
\big)\,
\big(
\tfrac{x_1^2}{2}+x_1^2x_2
\big)
\\
&
\ \ \ \ \
+
\big(
A_{3,1}x_1+A_{3,2}x_2+A_{3,3}x_3
\big)\,
\big(
\tfrac{x_1^3}{6}
\big).
\endaligned
\]

Starting from order $2$ thanks to 
Lemma~{\ref{Lm-tangency-L-order-2-dim-n}},
the matrix reductions read:
\[
\left[
\begin{array}{cccc}
A_{1,1} & \red{\bf 0} & \red{\bf 0} & B_1
\\
A_{2,1} & A_{2,2} & A_{2,3} & B_2
\\
A_{3,1} & A_{3,2} & A_{3,3} & B_3
\\
\red{\bf 0} & \red{\bf 0} & \red{\bf 0} & \red{2A_{1,1}}
\end{array}
\right]^{\green{2}}
\,\,\,\leadsto\,\,\,
\left[
\begin{array}{cccc}
A_{1,1} & \red{\bf 0} & \red{\bf 0} & \red{-A_{2,1}}
\\
A_{2,1} & \red{\bf 0} & \red{\bf 0} & B_2
\\
A_{3,1} & A_{3,2} & A_{3,3} & B_3
\\
\red{\bf 0} & \red{\bf 0} & \red{\bf 0} & \red{2A_{1,1}}
\end{array}
\right]^{\green{3}}
\,\,\,\leadsto\,\,\,
\left[
\begin{array}{cccc}
A_{1,1} & \red{\bf 0} & \red{\bf 0} & \red{-A_{2,1}}
\\
A_{2,1} & \red{\bf 0} & \red{\bf 0} & \red{-\frac{2}{3}A_{3,1}}
\\
A_{3,1} & \red{\bf 0} & \red{-A_{1,1}} & B_3
\\
\red{\bf 0} & \red{\bf 0} & \red{\bf 0} & \red{2A_{1,1}}
\end{array}
\right]^{\green{4}}.
\]

\medskip\noindent$\bullet$\,
In dimension $n = 4$, with:
\[
\aligned
u
&
\,=\,
\tfrac{x_1^2}{2}
\\
&
\ \ \ \ \
+
\tfrac{x_1^2x_2}{2}
\\
&
\ \ \ \ \
+
\tfrac{x_1^3x_3}{6}
+
\tfrac{x_1^2x_2^2}{2}
\\
&
\ \ \ \ \
+
\tfrac{x_1^4x_4}{24}
+
\tfrac{x_1^2x_2^3}{2}
+
\tfrac{x_1^3x_2x_3}{2},
\endaligned
\ \ \ \ \ \ \ \ \ \ \ \ \ \ \ \ \ \ \ \
\aligned
L
&
\,=\,
\big(
A_{1,1}x_1+A_{1,2}x_2+A_{1,3}x_3+A_{1,4}x_4+B_1u
\big)\partial_{x_1}
\\
&
\ \ \ \ \
+
\big(
A_{2,1}x_1+A_{2,2}x_2+A_{2,3}x_3+A_{2,4}x_4+B_2u
\big)\partial_{x_2}
\\
&
\ \ \ \ \
+
\big(
A_{3,1}x_1+A_{3,2}x_2+A_{3,3}x_3+A_{3,4}x_4+B_3u
\big)\partial_{x_3}
\\
&
\ \ \ \ \
+
\big(
A_{4,1}x_1+A_{4,2}x_2+A_{4,3}x_3+A_{4,4}x_4+B_4u
\big)\partial_{x_4}
\\
&
\ \ \ \ \
+
\big(
C_1x_1+C_2x_2+C_3x_3+C_4x_4+Du
\big)\partial_u,
\endaligned
\]
one can convince oneself that up to order $4$ included,
the same equations appear as in dimension $n = 3$,
hence one can replace the coefficients 
of $L$ from what has been obtained just above, so that
the tangency equation up to order $5$ reads:
\[
\aligned
0
&
\,\equiv\,
-\,2\,A_{1,1}\,
\big(
\tfrac{x_1^2}{2}
+
\tfrac{x_1^2x_2}{2}
+
\tfrac{x_1^3x_3}{6}
+
\tfrac{x_1^4x_4}{24}
\big)
\\
&
\ \ \ \ \
+
\big(
A_{1,1}\,x_1
-
A_{2,1}\,
\big[
\tfrac{x_1^2}{2}
+
\tfrac{x_1^2x_2}{2}
+
\tfrac{x_1^3x_3}{6}
\big]
\big)\,
\big(
x_1+x_1x_2+\tfrac{x_1^2x_3}{2}+\tfrac{x_1^3x_4}{6}
\big)
\\
&
\ \ \ \ \
+
\big(
A_{2,1}\,x_1
-
\tfrac{2}{3}\,A_{3,1}\,
\big[
\tfrac{x_1^2}{2}
+
\tfrac{x_1^2x_2}{2}
\big]
\big)\,
\big(
\tfrac{x_1^2}{2}
+
x_1^2x_2
+
\tfrac{x_1^3x_3}{2}
\big)
\\
&
\ \ \ \ \
+
\big(
A_{3,1}\,x_1
-
A_{1,1}\,x_3
+
B_3
\big[
\tfrac{x_1^2}{2}
\big]
\big)\,
\big(
\tfrac{x_1^3}{6}
+
\tfrac{x_1^3x_2}{2}
\big)
\\
&
\ \ \ \ \
+
\big(
A_{4,1}\,x_1+A_{4,2}\,x_2+A_{4,3}\,x_3+A_{4,4}\,x_4
\big)\,
\big(
\tfrac{x_1^4}{24}
\big).
\endaligned
\]
By expanding this equation, one can see 
that the matrix reduction at order
$5$ reads:
\[
\left[
\def\arraystretch{1.25}
\begin{array}{ccccc}
A_{1,1} & \red{\bf 0} & \red{\bf 0} & \red{\bf 0} & 
\red{-\frac{2}{2}A_{2,1}}
\\
A_{2,1} & \red{\bf 0} & \red{\bf 0} & \red{\bf 0} & 
\red{-\frac{2}{3}A_{3,1}}
\\
A_{3,1} & \red{\bf 0} & \red{-A_{1,1}} & \red{\bf 0} & B_3
\\
A_{4,1} & A_{4,2} & A_{4,3} & A_{4,4} & B_4
\\
\red{\bf 0} & \red{\bf 0} & \red{\bf 0} & \red{\bf 0} & \red{2A_{1,1}}
\end{array}
\right]^{\green{4}}
\,\,\,\leadsto\,\,\,
\left[
\def\arraystretch{1.25}
\begin{array}{ccccc}
A_{1,1} & \red{\bf 0} & \red{\bf 0} & \red{\bf 0} & 
\red{-\frac{2}{2}A_{2,1}}
\\
A_{2,1} & \red{\bf 0} & \red{\bf 0} & \red{\bf 0} & 
\red{-\frac{2}{3}A_{3,1}}
\\
A_{3,1} & \red{\bf 0} & \red{-A_{1,1}} & \red{\bf 0} & 
\red{-\frac{2}{4}A_{4,1}}
\\
A_{4,1} & \red{\bf 0} & \red{-2A_{2,1}} & \red{-2A_{1,1}} & B_4
\\
\red{\bf 0} & \red{\bf 0} & \red{\bf 0} & \red{\bf 0} & \red{2A_{1,1}}
\end{array}
\right]^{\green{5}}.
\]

\medskip\noindent$\bullet$\,
In dimension $n = 5$, with:
\[
\aligned
u
&
\,=\,
\tfrac{x_1^2}{2}
\\
&
\ \ \ \ \
+
\tfrac{x_1^2x_2}{2}
\\
&
\ \ \ \ \
+
\tfrac{x_1^3x_3}{6}
+
\tfrac{x_1^2x_2^2}{2}
\\
&
\ \ \ \ \
+
\tfrac{x_1^4x_4}{24}
+
\tfrac{x_1^2x_2^3}{2}
+
\tfrac{x_1^3x_2x_3}{2}
\\
&
\ \ \ \ \
+
\tfrac{x_1^5x_5}{120}
+
\tfrac{x_1^2x_2^4}{2}
+
x_1^3x_2^2x_3
+
\tfrac{x_1^4x_2x_4}{6}
+
\tfrac{x_1^4x_3^2}{8},
\endaligned
\]
the tangency equation up to order $5$ is:
\[
\aligned
0
&
\,\equiv\,
-\,2\,A_{1,1}\,
\big(
\tfrac{x_1^2}{2}
+
\tfrac{x_1^2x_2}{2}
+
\tfrac{x_1^3x_3}{6}
+
\tfrac{x_1^4x_4}{24}
+
\tfrac{x_1^5x_5}{120}
\big)
\\
&
\ \ \ \ \
+
\big(
A_{1,1}\,x_1
-
A_{2,1}\,
\big[
\tfrac{x_1^2}{2}
+
\tfrac{x_1^2x_2}{2}
+
\tfrac{x_1^3x_3}{6}
+
\tfrac{x_1^4x_4}{24}
\big]
\big)\,
\big(
x_1+x_1x_2+\tfrac{x_1^2x_3}{2}+\tfrac{x_1^3x_4}{6}
+
\tfrac{x_1^4x_5}{24}
\big)
\\
&
\ \ \ \ \
+
\big(
A_{2,1}\,x_1
-
\tfrac{2}{3}\,A_{3,1}\,
\big[
\tfrac{x_1^2}{2}
+
\tfrac{x_1^2x_2}{2}
+
\tfrac{x_1^3x_3}{6}
\big]
\big)\,
\big(
\tfrac{x_1^2}{2}
+
x_1^2x_2
+
\tfrac{x_1^3x_3}{2}
+
\tfrac{x_1^4x_4}{6}
\big)
\\
&
\ \ \ \ \
+
\big(
A_{3,1}\,x_1
-
A_{1,1}\,x_3
-
\tfrac{2}{4}\,A_{4,1}\,
\big[
\tfrac{x_1^2}{2}
+
\tfrac{x_1^2x_2}{2}
\big]
\big)\,
\big(
\tfrac{x_1^3}{6}
+
\tfrac{x_1^3x_2}{2}
+
\tfrac{x_1^4x_3}{4}
\big)
\\
&
\ \ \ \ \
+
\big(
A_{4,1}\,x_1
-
2\,A_{2,1}\,x_3
-
2\,A_{1,1}\,x_4
+
B_4\,
\big[
\tfrac{x_1^2}{2}
\big]
\big)\,
\big(
\tfrac{x_1^4}{24}
+
\tfrac{x_1^4x_2}{6}
\big)
\\
&
\ \ \ \ \
+
\big(
A_{5,1}\,x_1+A_{5,2}\,x_2+A_{5,3}\,x_3+A_{5,4}\,x_4+A_{5,5}\,x_5
\big)\,
\big(
\tfrac{x_1^5}{120}
\big),
\endaligned
\]
and the matrix reduction is:
\[
\!\!\!\!\!\!\!\!\!\!\!\!\!\!\!
\left[
\def\arraystretch{1.25}
\begin{array}{cccccc}
A_{1,1} & \red{\bf 0} & \red{\bf 0} & \red{\bf 0} & \red{\bf 0} & 
\red{-\frac{2}{2}A_{2,1}}
\\
A_{2,1} & \red{\bf 0} & \red{\bf 0} & \red{\bf 0} & \red{\bf 0} & 
\red{-\frac{2}{3}A_{3,1}}
\\
A_{3,1} & \red{\bf 0} & \red{-A_{1,1}} & \red{\bf 0} & \red{\bf 0} & 
\red{-\frac{2}{4}A_{4,1}}
\\
A_{4,1} & \red{\bf 0} & \red{-2A_{2,1}} & \red{-2A_{1,1}} & 
\red{\bf 0} & B_4
\\
A_{5,1} & A_{5,2} & A_{5,3} & A_{5,4} & A_{5,5} & B_5
\\
\red{\bf 0} & \red{\bf 0} & \red{\bf 0} & \red{\bf 0} & \red{\bf 0} & 
\red{2A_{1,1}}
\end{array}
\right]^{\green{5}}
\,\,\,\leadsto\,\,\,
\left[
\def\arraystretch{1.25}
\begin{array}{cccccc}
A_{1,1} & \red{\bf 0} & \red{\bf 0} & \red{\bf 0} & \red{\bf 0} & 
\red{-\frac{2}{2}A_{2,1}}
\\
A_{2,1} & \red{\bf 0} & \red{\bf 0} & \red{\bf 0} & \red{\bf 0} & 
\red{-\frac{2}{3}A_{3,1}}
\\
A_{3,1} & \red{\bf 0} & \red{-A_{1,1}} & \red{\bf 0} & \red{\bf 0} & 
\red{-\frac{2}{4}A_{4,1}}
\\
A_{4,1} & \red{\bf 0} & \red{-2A_{2,1}} & \red{-2A_{1,1}} & 
\red{\bf 0} & \red{-\frac{2}{5}A_{5,1}}
\\
A_{5,1} & \red{\bf 0} & \red{-\frac{10}{3}A_{3,1}} & 
\red{-5A_{2,1}} & \red{-3A_{1,1}} & B_5
\\
\red{\bf 0} & \red{\bf 0} & \red{\bf 0} & \red{\bf 0} & \red{\bf 0} & 
\red{2A_{1,1}}
\end{array}
\right]^{\green{6}}.
\]

\medskip\noindent$\bullet$\,
In dimension $n = 6$, with:
\[
\aligned
u
&
\,=\,
\tfrac{x_1^2}{2}
\\
&
\ \ \ \ \
+
\tfrac{x_1^2x_2}{2}
\\
&
\ \ \ \ \
+
\tfrac{x_1^3x_3}{6}
+
\tfrac{x_1^2x_2^2}{2}
\\
&
\ \ \ \ \
+
\tfrac{x_1^4x_4}{24}
+
\tfrac{x_1^2x_2^3}{2}
+
\tfrac{x_1^3x_2x_3}{2}
\\
&
\ \ \ \ \
+
\tfrac{x_1^5x_5}{120}
+
\tfrac{x_1^2x_2^4}{2}
+
x_1^3x_2^2x_3
+
\tfrac{x_1^4x_2x_4}{6}
+
\tfrac{x_1^4x_3^2}{8}
\\
&
\ \ \ \ \
+
\tfrac{x_1^6x_6}{720}
+
\tfrac{x_1^2x_2^5}{2}
+
\tfrac{5}{3}x_1^3x_2^3x_3
+
\tfrac{5}{8}x_1^4x_2x_3^2
+
\tfrac{5}{12}x_1^4x_2^2x_4
+
\tfrac{1}{12}x_1^5x_3x_4
+
\tfrac{1}{24}x_1^5x_2x_5,
\endaligned
\]
the tangency equation up to order $6$ is:
\[
\footnotesize
\aligned
0
&
\,\equiv\,
-\,2\,A_{1,1}\,
\big(
\tfrac{x_1^2}{2}
+
\tfrac{x_1^2x_2}{2}
+
\tfrac{x_1^3x_3}{6}
+
\tfrac{x_1^4x_4}{24}
+
\tfrac{x_1^5x_5}{120}
+
\tfrac{x_1^6x_6}{720}
\big)
\\
&
\ \ \ \ \
+
\big(
A_{1,1}\,x_1
-
A_{2,1}\,
\big[
\tfrac{x_1^2}{2}
+
\tfrac{x_1^2x_2}{2}
+
\tfrac{x_1^3x_3}{6}
+
\tfrac{x_1^4x_4}{24}
+
\tfrac{x_1^5x_5}{120}
\big]
\big)\,
\big(
x_1+x_1x_2+\tfrac{x_1^2x_3}{2}+\tfrac{x_1^3x_4}{6}
+
\tfrac{x_1^4x_5}{24}
+
\tfrac{x_1^5x_6}{120}
\big)
\\
&
\ \ \ \ \
+
\big(
A_{2,1}\,x_1
-
\tfrac{2}{3}\,A_{3,1}\,
\big[
\tfrac{x_1^2}{2}
+
\tfrac{x_1^2x_2}{2}
+
\tfrac{x_1^3x_3}{6}
+
\tfrac{x_1^4x_4}{24}
\big]
\big)\,
\big(
\tfrac{x_1^2}{2}
+
x_1^2x_2
+
\tfrac{x_1^3x_3}{2}
+
\tfrac{x_1^4x_4}{6}
+
\tfrac{x_1^5x_5}{24}
\big)
\\
&
\ \ \ \ \
+
\big(
A_{3,1}x_1
-
A_{1,1}x_3
-
\tfrac{2}{4}\,A_{4,1}\,
\big[
\tfrac{x_1^2}{2}
+
\tfrac{x_1^2x_2}{2}
+
\tfrac{x_1^3x_3}{6}
\big]
\big)\,
\big(
\tfrac{x_1^3}{6}
+
\tfrac{x_1^3x_2}{2}
+
\tfrac{x_1^4x_3}{4}
+
\tfrac{x_1^5x_4}{12}
\big)
\\
&
\ \ \ \ \
+
\big(
A_{4,1}\,x_1
-
2\,A_{2,1}\,x_3
-
2\,A_{1,1}\,x_4
-
\tfrac{2}{5}\,A_{5,1}\,
\big[
\tfrac{x_1^2}{2}
+
\tfrac{x_1^2x_2}{2}
\big]
\big)\,
\big(
\tfrac{x_1^4}{24}
+
\tfrac{x_1^4x_2}{6}
+
\tfrac{x_1^5x_3}{12}
\big)
\\
&
\ \ \ \ \
+
\big(
A_{5,1}\,x_1
-
\tfrac{10}{3}\,A_{3,1}\,x_3
-
5\,A_{2,1}\,x_4
-
3\,A_{1,1}\,x_5
+
B_5\,
\big[
\tfrac{x_1^2}{2}
\big]
\big)\,
\big(
\tfrac{x_1^5}{120}
+
\tfrac{x_1^5x_2}{24}
\big)
\\
&
\ \ \ \ \
+
\big(
A_{6,1}\,x_1+A_{6,2}\,x_2+A_{6,3}\,x_3+A_{6,4}\,x_4+A_{6,5}\,x_5
+A_{6,6}\,x_6
\big)\,
\big(
\tfrac{x_1^6}{720}
\big),
\endaligned
\]
and the matrix reduction is:
\[
\!\!\!\!\!\!\!\!\!\!\!\!\!\!\!\!\!\!\!\!\!\!\!\!\!
\footnotesize
\aligned
\left[
\def\arraystretch{1.25}
\begin{array}{ccccccc}
A_{1,1} & \red{\bf 0} & \red{\bf 0} & \red{\bf 0} & \red{\bf 0} & 
 \red{\bf 0} & \red{-\frac{2}{2}A_{2,1}}
\\
A_{2,1} & \red{\bf 0} & \red{\bf 0} & \red{\bf 0} & \red{\bf 0} & 
 \red{\bf 0} & \red{-\frac{2}{3}A_{3,1}}
\\
A_{3,1} & \red{\bf 0} & \red{-A_{1,1}} & \red{\bf 0} & \red{\bf 0} & 
\red{\bf 0} & \red{-\frac{2}{4}A_{4,1}}
\\
A_{4,1} & \red{\bf 0} & \red{-2A_{2,1}} & \red{-2A_{1,1}} & 
\red{\bf 0} & \red{\bf 0} & \red{-\frac{2}{5}A_{5,1}}
\\
A_{5,1} & \red{\bf 0} & \red{-\frac{10}{3}A_{3,1}} & 
\red{-5A_{2,1}} & \red{-3A_{1,1}} & \red{\bf 0} & B_5
\\
A_{6,1} & A_{6,2} & A_{6,3} & A_{6,4} & A_{6,5} &
A_{6,6} & B_6
\\
\red{\bf 0} & \red{\bf 0} & \red{\bf 0} & \red{\bf 0} & \red{\bf 0} & 
 \red{\bf 0} & \red{2A_{1,1}}
\end{array}
\right]^{\green{6}}
\leadsto
\left[
\def\arraystretch{1.25}
\begin{array}{ccccccc}
A_{1,1} & \red{\bf 0} & \red{\bf 0} & \red{\bf 0} & \red{\bf 0} & 
 \red{\bf 0} & \red{-\frac{2}{2}A_{2,1}}
\\
A_{2,1} & \red{\bf 0} & \red{\bf 0} & \red{\bf 0} & \red{\bf 0} & 
 \red{\bf 0} & \red{-\frac{2}{3}A_{3,1}}
\\
A_{3,1} & \red{\bf 0} & \red{-A_{1,1}} & \red{\bf 0} & \red{\bf 0} & 
\red{\bf 0} & \red{-\frac{2}{4}A_{4,1}}
\\
A_{4,1} & \red{\bf 0} & \red{-2A_{2,1}} & \red{-2A_{1,1}} & 
\red{\bf 0} & \red{\bf 0} & \red{-\frac{2}{5}A_{5,1}}
\\
A_{5,1} & \red{\bf 0} & \red{-\frac{10}{3}A_{3,1}} & 
\red{-5A_{2,1}} & \red{-3A_{1,1}} & \red{\bf 0} & 
\red{-\frac{2}{6}A_{6,1}}
\\
A_{6,1} & \red{\bf 0} & \red{-5A_{4,1}} & \red{-10A_{3,1}} & 
\red{-9A_{2,1}} & \red{-4A_{1,1}} & B_6
\\
\red{\bf 0} & \red{\bf 0} & \red{\bf 0} & \red{\bf 0} & \red{\bf 0} & 
 \red{\bf 0} & \red{2A_{1,1}}
\end{array}
\right]^{\green{7}}.
\endaligned
\]

\SectionHead{Projections $\pi^m (\centersmallbullet)$ and
$\pi_{\ind}^m (\centersmallbullet)$}
{projection-pi-m-ind}

For any integer $m \geqslant 1$ and any power series vanishing
at the origin:
\[
E(x)
\,=\,
E(x_1,\dots,x_n)
\,=\,
\sum_{\sigma_1,\dots,\sigma_n\geqslant 0
\atop
\sigma_1+\cdots+\sigma_n\geqslant 1}\,
x_1^{\sigma_1}\cdots x_n^{\sigma_n}\,
E_{\sigma_1,\dots,\sigma_n},
\]
define:
\[
\pi^m
\big(E(x)\big)
\,:=\,
\sum_{\sigma_1+\cdots+\sigma_n\leqslant m}\,
x_1^{\sigma_1}\cdots x_n^{\sigma_n}\,
E_{\sigma_1,\dots,\sigma_n}.
\]
Several times later, the following elementary fact will be useful. 

\begin{Observation}
\label{Obs-pi-product-m}
Given two integers $k_1, k_2 \in \N$ and 
two power series $H_1(x) \in {\rm O}_x(k_1)$ and
$H_2(x) \in {\rm O}_x(k_2)$, namely: 
\[
H_1(x)
\,=\,
\sum_{\sigma_1+\cdots+\sigma_n\geqslant k_1}\, 
H_{1,\sigma_1,\dots,\sigma_n}\,
x_1^{\sigma_1}\cdots x_n^{\sigma_n},
\ \ \ \ \ \ \ \ \ \ \ \ \ \ \ \ \ \ \ \
H_2(x)
\,=\,
\sum_{\sigma_1+\cdots+\sigma_n\geqslant k_2}\, 
H_{2,\sigma_1,\dots,\sigma_n}\,
x_1^{\sigma_1}\cdots x_n^{\sigma_n},
\]
for any (homogeneous) order $m \geqslant 0$,
one has:
\[
\pi^m
\big(
H_1
\cdot
H_2
\big)
\,=\,
\pi^m
\Big(
\pi^{m-k_2}
(H_1)
\cdot
\pi^{m-k_1}
(H_2)
\Big),
\]
the right-hand side being understood as $0$ when 
$m < \min\, (k_1, k_2)$.\qed
\end{Observation}

As our hypersurfaces $H^n \subset \R^{n+1}$ have constant
Hessian rank $1$, independent monomials are of interest.
Accordingly, define:
\[
\pi_{\ind}\big(E(x)\big)
\,:=\,
\sum_{i=0}^\infty\,
\Big(
E_{i+1,0,\dots,0}\,x_1^{i+1}
+
E_{i,1,\dots,0}\,x_1^ix_2
+\cdots+
E_{i,0,\dots,1}\,x_1^ix_n.
\Big),
\]
and also:
\[
\pi_{\ind}^m
\,:=\,
\pi^m\circ\pi_{\ind}
\,=\,
\pi_{\ind}\circ\pi^m,
\]
that is:
\[
\pi_{\ind}^m
\big(E(x)\big)
\,:=\,
\sum_{0\leqslant i\leqslant m-1}\,
\Big(
E_{i+1,0,\dots,0}\,x_1^{i+1}
+
E_{i,1,\dots,0}\,x_1^ix_2
+\cdots+
E_{i,0,\dots,1}\,x_1^ix_n.
\Big).
\]

\begin{Observation}
\label{Obs-pi-product-m-independent}
Given two integers $k_1, k_2 \in \N$ and 
two power series $H_1(x) \in {\rm O}_x(k_1)$ and
$H_2(x) \in {\rm O}_x(k_2)$, 
for any (homogeneous) order $m \geqslant 0$,
one has:
\[
\pi_{\ind}^m
\big(
H_1
\cdot
H_2
\big)
\,=\,
\pi_{\ind}^m
\Big(
\pi_{\ind}^{m-k_2}
(H_1)
\cdot
\pi_{\ind}^{m-k_1}
(H_2)
\Big),
\]
the right-hand side being understood as $0$ when 
$m < \min\, (k_1, k_2)$.\qed
\end{Observation}

Soon, we will need a notation to select {\em homogeneous}
monomials
of order exactly equal to some fixed integeer $m \geqslant 1$ 
(notice that the index $m$
is lower-case now):
\[
\pi_m
\big(E(x)\big)
\,:=\,
\sum_{\sigma_1+\cdots+\sigma_n=m}\,
x_1^{\sigma_1}\cdots x_n^{\sigma_n}\,
E_{\sigma_1,\dots,\sigma_n},
\]
and also to select the {\em independent} homogeneous ones:
\[
\pi_m^{\ind}
\big(E(x)\big)
\,:=\,
E_{m,0,\dots,0}\,x_1^m
+
E_{m-1,1,\dots,0}\,x_1^{m-1}x_2
+\cdots+
E_{m-1,0,\dots,1}\,x_1^{m-1}x_n.
\]

\SectionHead{Tangency at Order $n+1$ in Dimension $n$}
{tangency-order-n-1}

Now, launch induction on dimension. 
In any dimension $n-1$, in coordinates
$(x_1, \dots, x_{n-1}, u)$, 
we know that the hypersurface equation writes
up to order $n-1+1$ as:
\[
\aligned
u
&
\,=\,
\frac{x_1^2}{2!}
+
\frac{x_1^2x_2}{2!}
+
\sum_{m=3}^{n-2}\,
\Big(
\frac{x_1^mx_m}{m!}
+
x_1^{m-1}
\sum_{i,j\geqslant 2
\atop
i+j=m+1}\,
\tfrac{1}{2}\,
\frac{x_ix_j}{(i-1)!(j-1)!}
\Big)
\\
&
\ \ \ \ \
+
\frac{x_1^{n-1}x_{n-1}}{(n-1)!}
+
x_1^{n-2}
\sum_{i,j\geqslant 2
\atop
i+j=n}\,
\tfrac{1}{2}\,
\frac{x_ix_j}{(i-1)!(j-1)!}
\\
&
\ \ \ \ \
+
{\rm O}_{x_2,\dots,x_{n-1}}(2)
+
{\rm O}_{x_1,x_2,\dots,x_{n-1}}(n+1),
\endaligned
\]
with explicit independent and border-dependent monomials.

A linear (vanishing at the origin) affine vector field writes:
\[
L
\,=\,
X_1\,\frac{\partial}{\partial x_1}
+\cdots+
X_{n-1}\,\frac{\partial}{\partial x_{n-1}}
+
U\,\frac{\partial}{\partial u},
\]
with:
\[
\aligned
X_1
&
\,=\,
A_{1,1}\,x_1
+\cdots+
A_{1,n-1}\,x_{n-1}
+
B_1\,u,
\\
\cdots
&
\cdots\cdots\cdots\cdots\cdots\cdots\cdots\cdots\cdots\cdots\cdots
\cdots\cdots
\\
X_{n-1}
&
\,=\,
A_{n-1,1}\,x_1
+\cdots+
A_{n-1,n-1}\,x_{n-1}
+
B_{n-1}\,u,
\\
U
&
\,=\,
C_1\,x_1
+\cdots+
C_{n-1}\,x_{n-1}
+
D\,u.
\endaligned
\]

We consider tangency of $L$ to such a hypersurface
$u = F(x_1, \dots, x_{n-1})$ up to order
$n-1$, namely the condition:
\[
0
\,=\,
\pi^{n-1}
\Big(
L\big(-u+F\big)
\big\vert_{u=F}
\Big),
\]
which is equivalent to:
\[
0
\,=\,
\pi_{\ind}^{n-1}
\Big(
L\big(-u+F\big)
\big\vert_{u=F}
\Big).
\]
We start considerations from 
Lemma~{\ref{Lm-tangency-L-order-2-dim-n}} at order $2$,
viewed in dimension $n-1$.

\begin{InductionHypothesis}
\label{Induction-Hypothesis-L-n-1}
The vector space of fields $L$ which are tangent to
order $\leqslant n-1$ is of dimension $n-1 + 1$,
parametrized by $A_{1,1}$, \dots, $A_{n-1,1}$, $B_{n-1}$,
with the other constants defined by:
\[
\underset{(1\leqslant i\leqslant n-1)}{
A_{i,2}\,=\,0,}
\ \ \ \ \ \ \ 
\underset{(3\leqslant j\leqslant n-1)}{
A_{1,j}\,=\,A_{2,j}\,=\,0,}
\ \ \ \ \ \ \ 
\underset{(3\leqslant j\leqslant i\leqslant n-1)}{
A_{i,j}\,=\,{\textstyle{\frac{-(j-2)}{i-j+1}\,
\binom{i}{j}}}\,
A_{i-j+1,1},}
\ \ \ \ \ \ \ 
\underset{(3\leqslant j\leqslant n-2)}{
B_j\,=\,-\tfrac{2}{j+1}\,A_{j+1,1},}
\]
or equivalently, the matrix of the coefficients of $L$ is:
\[
\left[
\def\arraystretch{1.25}
\begin{array}{ccccccc}
A_{1,1} & \red{\bf 0} & \red{\bf 0} & \red{\bf 0} &
\red{\cdots} & \red{\bf 0} & \red{-\tfrac{2}{2}A_{2,1}}
\\
A_{2,1} & \red{\bf 0} & \red{\bf 0} & \red{\bf 0} &
\red{\cdots} & \red{\bf 0} & \red{-\tfrac{2}{3}A_{3,1}}
\\
A_{3,1} & \red{\bf 0} & \red{-A_{1,1}} & \red{\bf 0} &
\red{\cdots} & \red{\bf 0} & \red{-\tfrac{2}{4}A_{4,1}}
\\
A_{4,1} & \red{\bf 0} & \red{-2A_{2,1}} & \red{-A_{1,1}} &
\red{\cdots} & \red{\bf 0} & \red{-\tfrac{2}{5}A_{4,1}}
\\
\vdots & \red{\vdots} & \red{\vdots} & \red{\vdots} & 
\red{\ddots} & \red{\vdots} & \red{\vdots}
\\
A_{n-2,1} & \red{\bf 0} &
\red{\frac{-1}{n-4}\binom{n-2}{3}A_{n-4,1}} &
\red{\frac{-2}{n-5}\binom{n-2}{4}A_{n-5,1}} & \red{\cdots} &
\red{\bf 0} & \red{-\frac{2}{n-1}A_{n-1,1}}
\\
A_{n-1,1} & \red{\bf 0} &
\red{\frac{-1}{n-3}\binom{n-1}{3}A_{n-3,1}} &
\red{\frac{-2}{n-4}\binom{n-1}{4}A_{n-4,1}} & \red{\cdots} &
\red{-(n-3)A_{1,1}} & B_{n-1}
\\
\red{\bf 0} & \red{\bf 0} & \red{\bf 0} & \red{\bf 0} &
\red{\cdots} & \red{\bf 0} & \red{2\,A_{1,1}}
\end{array}
\right].
\]
\end{InductionHypothesis}

The goal is to show that similar expressions hold
in dimension $n$. Thus, in coordinates
$(x_1, \dots, x_{n-1}, x_n, u)$, consider a hypersurface
whose equation writes up to order $n+1$ as:
\leqnomode\usetagform{default}
\begin{align}
\label{u-F-n-plus-1}
u
&
\,=\,
\frac{x_1^2}{2!}
+
\frac{x_1^2x_2}{2!}
+
\sum_{m=3}^{n-2}\,
\Big(
\frac{x_1^mx_m}{m!}
+
x_1^{m-1}
\sum_{i,j\geqslant 2
\atop
i+j=m+1}\,
\tfrac{1}{2}\,
\frac{x_ix_j}{(i-1)!(j-1)!}
\Big)
\notag
\\
&
\ \ \ \ \
+
\frac{x_1^{n-1}x_{n-1}}{(n-1)!}
+
x_1^{n-2}
\sum_{i,j\geqslant 2
\atop
i+j=n}\,
\tfrac{1}{2}\,
\frac{x_ix_j}{(i-1)!(j-1)!}
\\
&
\ \ \ \ \
+
\frac{x_1^nx_n}{n!}
+
x_1^{n-1}
\sum_{i,j\geqslant 2
\atop
i+j=n+1}\,
\tfrac{1}{2}\,
\frac{x_ix_j}{(i-1)!(j-1)!}
\notag
\\
&
\ \ \ \ \
+
{\rm O}_{x_2,\dots,x_{n-1},x_n}(2)
+
{\rm O}_{x_1,x_2,\dots,x_{n-1},x_n}(n+2),
\notag
\end{align}
and consider a linear vector field:
\[
L
\,=\,
X_1\,\frac{\partial}{\partial x_1}
+\cdots+
X_{n-1}\,\frac{\partial}{\partial x_{n-1}}
+
X_n\,\frac{\partial}{\partial x_n}
+
U\,\frac{\partial}{\partial u},
\]
with:
\[
\aligned
X_1
&
\,=\,
A_{1,1}\,x_1
+\cdots+
A_{1,n-1}\,x_{n-1}
+
B_1\,u,
\\
\cdots
\cdot
&
\cdots\cdots\cdots\cdots\cdots\cdots\cdots\cdots\cdots\cdots\cdots
\cdots\cdot\cdot
\\
X_{n-1}
&
\,=\,
A_{n-1,1}\,x_1
+\cdots+
A_{n-1,n}\,x_n
+
B_{n-1}\,u,
\\
X_n
&
\,=\,
A_{n,1}\,x_1
+\cdots+
A_{n,n}\,x_n
+
B_n\,u,
\\
U
&
\,=\,
C_1\,x_1
+\cdots+
C_n\,x_n
+
D\,u.
\endaligned
\]

Similarly, we consider tangency of $L$ to such a hypersurface $u =
F(x_1, \dots, x_{n-1}, x_n)$ up to order $n$, namely the condition:
\[
0
\,=\,
\pi^n
\Big(
L\big(-u+F\big)
\big\vert_{u=F}
\Big),
\]
which is equivalent to:
\leqnomode\usetagform{default}
\begin{align}
\label{pi-ind-n-tangency-L}
0
\,=\,
\pi_{\ind}^n
\Big(
L\big(-u+F\big)
\big\vert_{u=F}
\Big).
\end{align}

Putting $x_n := 0$, taking only $\pi_{\sf ind}^{n-1} 
(\centersmallbullet)$,
applying the Induction 
Hypothesis~{\ref{Induction-Hypothesis-L-n-1}},
we see that the matrix of coefficients of
$L$ incorporates known elements in dimension $n-1$:
\[
\left[
\def\arraystretch{1.25}
\begin{array}{cccccccc}
A_{1,1} & \red{\bf 0} & \red{\bf 0} & \red{\bf 0} &
\red{\cdots} & \red{\bf 0} & A_{1,n} & \red{-\tfrac{2}{2}A_{2,1}}
\\
A_{2,1} & \red{\bf 0} & \red{\bf 0} & \red{\bf 0} &
\red{\cdots} & \red{\bf 0}& A_{2,n} & \red{-\tfrac{2}{3}A_{3,1}}
\\
A_{3,1} & \red{\bf 0} & \red{-A_{1,1}} & \red{\bf 0} &
\red{\cdots} & \red{\bf 0} & A_{3,n} & \red{-\tfrac{2}{4}A_{4,1}}
\\
A_{4,1} & \red{\bf 0} & \red{-2A_{2,1}} & \red{-2A_{1,1}} &
\red{\cdots} & \red{\bf 0} & A_{4,n} & \red{-\tfrac{2}{5}A_{4,1}}
\\
\vdots & \red{\vdots} & \red{\vdots} & \red{\vdots} & 
\red{\ddots} & \red{\vdots} & \vdots & \red{\vdots}
\\
A_{n-2,1} & \red{\bf 0} &
\red{\frac{-1}{n-4}\binom{n-2}{3}A_{n-4,1}} &
\red{\frac{-2}{n-5}\binom{n-2}{4}A_{n-5,1}} & \red{\cdots} &
\red{\bf 0} & A_{n-2,n} & \red{-\frac{2}{n-1}A_{n-1,1}}
\\
A_{n-1,1} & \red{\bf 0} &
\red{\frac{-1}{n-3}\binom{n-1}{3}A_{n-3,1}} &
\red{\frac{-2}{n-4}\binom{n-1}{4}A_{n-4,1}} & \red{\cdots} &
\red{-(n-3)A_{1,1}} & A_{n-1,n} & B_{n-1}
\\
A_{n,1} & A_{n,2} & A_{n,3} & A_{n,4} & \cdots & A_{n-1,n} & A_{n,n} & 
B_n
\\
\red{\bf 0} & \red{\bf 0} & \red{\bf 0} & \red{\bf 0} &
\red{\cdots} & \red{\bf 0} & C_n & \red{2\,A_{1,1}}
\end{array}
\right].
\]

Without putting $x_n := 0$, 
by examining $\pi_{\sf ind}^{n-1} 
(\centersmallbullet)$, one can realize
that:
\[
A_{1,n}
\,=\,
\cdots
\,=\,
A_{n-1,n}
\,=\,
\red{\bf 0},
\ \ \ \ \ \ \ \ \ \ \ \ \ \ \ \ \ \ \ \
C_n
\,=\,
\red{\bf 0},
\]
and the detailed verification of these vanishings 
is implicitly contained
in the proof of the next Lemma~{\ref{Lm-tangency-order-n}}. 
So the matrix is:
\[
\left[
\def\arraystretch{1.25}
\begin{array}{cccccccc}
A_{1,1} & \red{\bf 0} & \red{\bf 0} & \red{\bf 0} &
\red{\cdots} & \red{\bf 0} & \red{\bf 0} & \red{-\tfrac{2}{2}A_{2,1}}
\\
A_{2,1} & \red{\bf 0} & \red{\bf 0} & \red{\bf 0} &
\red{\cdots} & \red{\bf 0}& \red{\bf 0} & \red{-\tfrac{2}{3}A_{3,1}}
\\
A_{3,1} & \red{\bf 0} & \red{-A_{1,1}} & \red{\bf 0} &
\red{\cdots} & \red{\bf 0} & \red{\bf 0} & \red{-\tfrac{2}{4}A_{4,1}}
\\
A_{4,1} & \red{\bf 0} & \red{-2A_{2,1}} & \red{-2A_{1,1}} &
\red{\cdots} & \red{\bf 0} & \red{\bf 0} & \red{-\tfrac{2}{5}A_{4,1}}
\\
\vdots & \red{\vdots} & \red{\vdots} & \red{\vdots} & 
\red{\ddots} & \red{\vdots} & \vdots & \red{\vdots}
\\
A_{n-2,1} & \red{\bf 0} &
\red{\frac{-1}{n-4}\binom{n-2}{3}A_{n-4,1}} &
\red{\frac{-2}{n-5}\binom{n-2}{4}A_{n-5,1}} & \red{\cdots} &
\red{\bf 0} & \red{\bf 0} & \red{-\frac{2}{n-1}A_{n-1,1}}
\\
A_{n-1,1} & \red{\bf 0} &
\red{\frac{-1}{n-3}\binom{n-1}{3}A_{n-3,1}} &
\red{\frac{-2}{n-4}\binom{n-1}{4}A_{n-4,1}} & \red{\cdots} &
\red{-(n-3)A_{1,1}} & \red{\bf 0} & B_{n-1}
\\
A_{n,1} & A_{n,2} & A_{n,3} & A_{n,4} & \cdots & A_{n-1,n} & A_{n,n} & 
B_n
\\
\red{\bf 0} & \red{\bf 0} & \red{\bf 0} & \red{\bf 0} &
\red{\cdots} & \red{\bf 0} & \red{\bf 0} & \red{2\,A_{1,1}}
\end{array}
\right],
\]
and in order to complete the induction, 
we must determine the values of
$A_{n,1}$, $A_{n,2}$, $A_{n,3}$, $A_{n,4}$,
\dots, $A_{n-1,n}$, $A_{n,n}$, $B_{n-1}$,
as follows.

\begin{Lemma}
\label{Lm-tangency-order-n}
The tangency condition~{\eqref{pi-ind-n-tangency-L}} 
at order $n+1$ forces the announced values:
\[
B_{n-1}
\,=\,
-\,\frac{2}{n}\,
A_{n,1},
\ \ \ \ \ \ \ \ \ \ \ \ \ \ \ \ \ \ \ \
A_{n,2}
\,=\,
\red{\bf 0},
\ \ \ \ \ \ \ \ \ \ \ \ \ \ \ \ \ \ \ \
\underset{(3\leqslant k\leqslant n)}{
A_{n,k}
\,=\,
-\,
\frac{k-2}{n-k+1}\,
\binom{n}{k}\,
A_{n-k+1,1}}.
\]
\end{Lemma}

\proof
So, let us abbreviate the hypersurface equation~{\eqref{u-F-n-plus-1}}
including
independent and border-dependent monomials 
up to order $n+1$ as:
\[
u
\,=\,
\overline{F}^{n+1}(x)
+
{\rm O}_{x'}(2)
+
{\rm O}_x(n+2),
\]
and start to examine the tangency equation:
\leqnomode\usetagform{default}
\begin{footnotesize}
\begin{align}
\label{tangency-equation-U-Xi-Fxi-n-1}
0
&
\,\equiv\,
\pi_{\ind}^{n+1}
\Big(
L
\big(
-\,u+F
\big)
\big\vert_{u=F}
\Big)
\notag
\\
&
\,\equiv\,
\pi_{\ind}^{n+1}
\Big(
-\,U(x,u)
+
\sum_{i=1}^n\,
X_i(x,u)
\cdot
F_{x_i}(x)
\Big\vert_{u=F(x)}
\Big)
\notag
\\
&
\,\equiv\,
\pi_{\ind}^{n+1}
\bigg(
-\,U
\Big(
x,
\overline{F}^{n+2}(x)+{\rm O}_{x'}(2)+{\rm O}_x(n+2)
\Big)
\notag
\\
&
\ \ \ \ \ \ \ \ \ \ \ \ \ \ \ \ \ \ \ \
+
\sum_{i=1}^n\,
X_i
\Big(
x,
\overline{F}^{n+2}(x)+{\rm O}_{x'}(2)+{\rm O}_x(n+2)
\Big)
\cdot
\Big(
\overline{F}_{x_i}^{n+1}+{\rm O}_{x'}(1)+{\rm O}_x(n+1)
\Big)
\bigg)
\notag
\\
&
\,\equiv\,
\pi_{\ind}^{n+1}
\bigg(
-\,U
\big(x,
\overline{F}^{n+2}(x)
\big)
+
\sum_{i=1}^n\,
X_i
\big(
x,\overline{F}^{n+2}(x)
\big)
\cdot
\overline{F}_{x_i}^{n+1}(x)
\bigg),
\end{align}
\end{footnotesize}where 
we use $0 = \pi_{\sf ind}^{n+1} \big( {\rm O}_{x'}(1) \big) =
\pi^{n+1} \big( {\rm O}_x(n+2)\big)$, and we use the fact that the
$X_i\big( x, \overline{F}^{n+1}(x) \big)$ are all ${\rm O}_x(1)$,
hence multiplied by the last remainder ${\rm O}_x(n+1)$ become ${\rm
O}_x(n+2)$.

In two steps, 
we want to apply Observation~{\ref{Obs-pi-product-m-independent}} 
to the products
$X_i \cdot \overline{F}_{x_i}^{n+2}$ above for $i = 1, \dots, n$.
Firstly, each $X_i$ is an ${\rm O}_x(1)$, so:
\[
\pi_{\ind}^{n+1}
\big(
X_i
\cdot
\overline{F}_{x_i}^{n+1}
\big)
\,=\,
\pi_{\ind}^{n+1}
\Big(
\pi_{\ind}^{n+1}(X_i)
\cdot
\pi_{\ind}^n
\big(
\overline{F}_{x_i}^{n+1}
\big)
\Big).
\]

Secondly, the vanishing orders at $0$
of the $\overline{F}_{x_i}^{n+2}$
increase, as one may verify, by 
differentiating~{\ref{u-F-n-plus-1}},
that:
\[
\aligned
\pi_{\ind}^n
\big(
\overline{F}_{x_1}^{n+1}
\big)
&
\,=\,
\frac{x_1^1}{1!}
\bigg\vert
+
\frac{x_1^1x_2}{2!}
+\cdots+
\frac{x_1^{n-1}x_n}{(n-1)!},
\\
\pi_{\ind}^n
\big(
\overline{F}_{x_2}^{n+1}
\big)
&
\,=\,
\frac{x_1^2}{2!}
\bigg\vert
+
\frac{x_1^2x_2}{1!\,1!}
+\cdots+
\frac{x_1^{n-1}x_{n-1}}{1!\,(n-2)!},
\\
\pi_{\ind}^n
\big(
\overline{F}_{x_3}^{n+1}
\big)
&
\,=\,
\frac{x_1^3}{3!}
\bigg\vert
+
\frac{x_1^3x_2}{2!\,1!}
+\cdots+
\frac{x_1^{n-1}x_{n-2}}{2!\,(n-3)!},
\\
\cdots\cdots\cdots\cdot
&
\cdots\cdots\cdots\cdots\cdots\cdots\cdots\cdots\cdots\cdots\cdots
\cdots\cdots
\\
\pi_{\ind}^n
\big(
\overline{F}_{x_{n-3}}^{n+1}
\big)
&
\,=\,
\frac{x_1^{n-3}}{(n-3)!}
\bigg\vert
+
\frac{x_1^{n-3}x_2}{(n-4)!\,1!}
+
\frac{x_1^{n-2}x_3}{(n-4)!\,2!}
+
\frac{x_1^{n-1}x_4}{(n-4)!\,3!},
\\
\pi_{\ind}^n
\big(
\overline{F}_{x_{n-2}}^{n+1}
\big)
&
\,=\,
\frac{x_1^{n-2}}{(n-2)!}
\bigg\vert
+
\frac{x_1^{n-2}x_2}{(n-3)!\,1!}
+
\frac{x_1^{n-1}x_3}{(n-3)!\,2!}
\\
\pi_{\ind}^n
\big(
\overline{F}_{x_{n-1}}^{n+1}
\big)
&
\,=\,
\frac{x_1^{n-1}}{(n-1)!}
\bigg\vert
+
\frac{x_1^{n-1}x_2}{(n-2)!\,1!}
\\
\pi_{\ind}^n
\big(
\overline{F}_{x_n}^{n+1}
\big)
&
\,=\,
\frac{x_1^n}{n!};
\endaligned
\]
here for later use, the vertical bar separates the pure monomials
$x_1^\ast$
from the monomials $x_1^\ast x_2, \dots, x_1^\ast x_n$.
Therefore, applying again 
Observation~{\ref{Obs-pi-product-m-independent}}, 
it suffice to know the independent parts of
$X_1$, $X_2$, $X_3$, \dots, $X_{n-3}$,
$X_{n-2}$, $X_{n-1}$, $X_n$ up to decreasing orders:
\[
\aligned
\pi_{\ind}^{n+1}
\big(
X_1
\cdot
\overline{F}_{x_1}^{n+1}
\big)
&
\,=\,
\pi_{\ind}^{n+1}
\Big(
\pi_{\ind}^n(X_1)
\cdot
\pi_{\ind}^n
\big(\overline{F}_{x_1}^{n+1}\big)
\Big),
\\
\pi_{\ind}^{n+1}
\big(
X_2
\cdot
\overline{F}_{x_2}^{n+1}
\big)
&
\,=\,
\pi_{\ind}^{n+1}
\Big(
\pi_{\ind}^{n-1}(X_2)
\cdot
\pi_{\ind}^n
\big(\overline{F}_{x_2}^{n+1}\big)
\Big),
\\
\pi_{\ind}^{n+1}
\big(
X_3
\cdot
\overline{F}_{x_3}^{n+1}
\big)
&
\,=\,
\pi_{\ind}^{n+1}
\Big(
\pi_{\ind}^{n-2}(X_3)
\cdot
\pi_{\ind}^n
\big(\overline{F}_{x_3}^{n+1}\big)
\Big),
\\
\cdots\cdots\cdots\cdots\cdots\cdot\cdot
&
\cdots\cdots\cdots\cdots\cdots\cdots\cdots\cdots\cdots\cdots
\cdot\cdot
\\
\pi_{\ind}^{n+1}
\big(
X_{n-3}
\cdot
\overline{F}_{x_{n-3}}^{n+1}
\big)
&
\,=\,
\pi_{\ind}^{n+1}
\Big(
\pi_{\ind}^4(X_{n-3})
\cdot
\pi_{\ind}^n
\big(\overline{F}_{x_{n-3}}^{n+1}\big)
\Big),
\\
\pi_{\ind}^{n+1}
\big(
X_{n-2}
\cdot
\overline{F}_{x_{n-2}}^{n+1}
\big)
&
\,=\,
\pi_{\ind}^{n+1}
\Big(
\pi_{\ind}^3(X_{n-2})
\cdot
\pi_{\ind}^n
\big(\overline{F}_{x_{n-2}}^{n+1}\big)
\Big),
\\
\pi_{\ind}^{n+1}
\big(
X_{n-1}
\cdot
\overline{F}_{x_{n-1}}^{n+1}
\big)
&
\,=\,
\pi_{\ind}^{n+1}
\Big(
\pi_{\ind}^2(X_{n-1})
\cdot
\pi_{\ind}^n
\big(\overline{F}_{x_{n-1}}^{n+1}\big)
\Big),
\\
\pi_{\ind}^{n+1}
\big(
X_n
\cdot
\overline{F}_{x_n}^{n+1}
\big)
&
\,=\,
\pi_{\ind}^{n+1}
\Big(
\pi_{\ind}^1(X_n)
\cdot
\pi^n
\big(\overline{F}_{x_n}^{n+1}\big)
\Big),
\endaligned
\]

Remind that in all the $X_i$, before taking
$\pi_{\ind}^{n-i+1} (X_i)$ as above, 
we have to replace $u$ by
$\overline{F}^{n+1}(x)$ 
which, for $3 \leqslant i \leqslant n-2$, gives:
\[
\aligned
X_i
&
\,=\,
X_i
\big(
x,
\overline{F}^{n+1}(x)
\big)
\\
&
\,=\,
A_{i,1}x_1
-
\sum_{3\leqslant j\leqslant i}\,
\tfrac{j-2}{i-j+1}\,
{\textstyle{\binom{i}{j}}}\,
A_{i-j+1}\,
x_j
-
\tfrac{2}{i+1}\,
A_{i+1,1}\,
\overline{F}^{n+2}(x),
\endaligned
\]
with similar formulas for $i = 1, 2, n-1$.
Serendipitously, for all $3 \leqslant j \leqslant i$,
the values of the coefficients of the $x_j$
will be unimportant, hence we will abbreviate:
\[
\pi_{\ind}^{n-i+1}
(X_i)
\,=\,
A_{i,1}\,x_1
-
\ast\,x_3
-\cdots-
\ast\,x_i
-
\tfrac{2}{i+1}\,A_{i+1,1}
\Big[
\tfrac{x_1^2}{2!}
\Big\vert
+
\tfrac{x_1^2x_2}{2!}
+\cdots+
\tfrac{x_1^{n-i}x_{n-i}}{(n-i)!}
\Big].
\]
Serendipitously also, in $X_n$, 
no contribution of $B_n\, u \equiv B_n\, 
\overline{F}^{n+1}$ occurs,
since 
$\overline{F}^{n+1} = {\rm O}_x(2)$:
\[
\aligned
\pi_{\ind}^1(X_n)
&
\,=\,
\pi_{\ind}^1
\Big(
A_{n,1}\,x_1
+\cdots+
A_{n,n}\,x_n
+
B_n\,{\rm O}_x(2)
\Big)
\\
&
\,=\,
A_{n,1}\,x_1
+\cdots+
A_{n,n}\,x_n.
\endaligned
\]

At last, we can write in length all the terms
of the tangency equation~{\eqref{tangency-equation-U-Xi-Fxi-n-1}}, 
as follows without mentionining
$\pi_{\ind}^{n+1} (\centersmallbullet)$:
\[
\!\!\!\!\!\!\!\!\!\!\!\!\!\!\!\!\!\!\!\!\!\!\!\!\!
\footnotesize
\aligned
0
&
\,\equiv\,
-\,2\,A_{1,1}\,
\Big(
\tfrac{x_1^2}{2!}
+
\tfrac{x_1^2x_2}{2!}
+\cdots+
\tfrac{x_1^nx_n}{n!}
\Big)
\\
&
\ \ \ \ \
+
\Big(
A_{1,1}x_1
-
\tfrac{2}{2}A_{2,1}
\Big[
\tfrac{x_1^2}{2!}
\Big\vert
+
\tfrac{x_1^2x_2}{2!}
+\cdots+
\tfrac{x_1^{n-1}x_{n-1}}{(n-1)!}
\Big]
\Big)\,
\Big(
\tfrac{x_1^1}{1!}
\Big\vert
+
\tfrac{x_1^2x_2}{1!}
+\cdots+
\tfrac{x_1^{n-2}x_{n-1}}{(n-2)!}
+
\tfrac{x_1^{n-1}x_n}{(n-1)!}
\Big)
\\
&
\ \ \ \ \ 
+
\Big(
A_{2,1}x_1
-
\tfrac{2}{3}A_{3,1}
\Big[
\tfrac{x_1^2}{2!}
\Big\vert
+
\tfrac{x_1^2x_2}{2!}
+\cdots+
\tfrac{x_1^{n-2}x_{n-2}}{(n-2)!}
\Big]
\Big)
\Big(
\tfrac{x_1^2}{2!}
\Big\vert
+
\tfrac{x_1^2x_2}{1!\,1!}
+\cdots+
\tfrac{x_1^{n-2}x_{n-2}}{1!\,(n-3)!}
+
\tfrac{x_1^{n-1}x_{n-1}}{1!\,(n-2)!}
\Big)
\\
&
\ \ \ \ \ 
+
\Big(
A_{3,1}x_1
-
\ast x_3
-
\tfrac{2}{4}A_{4,1}
\Big[
\tfrac{x_1^2}{2!}
\Big\vert
+
\tfrac{x_1^2x_2}{2!}
+\cdots+
\tfrac{x_1^{n-3}x_{n-3}}{(n-3)!}
\Big]
\Big)
\Big(
\tfrac{x_1^3}{3!}
\Big\vert
+
\tfrac{x_1^3x_2}{2!\,1!}
+\cdots+
\tfrac{x_1^{n-2}x_{n-3}}{2!\,(n-4)!}
+
\tfrac{x_1^{n-1}x_{n-2}}{2!\,(n-3)!}
\Big)
\\
&
\ \ \ \ \ 
+
\cdots\cdots\cdots\cdots\cdots\cdots\cdots\cdots\cdots\cdots\cdots
\cdots\cdots\cdots\cdots\cdots\cdots\cdots\cdots\cdots\cdots\cdots
\cdots\cdots\cdots\cdots\cdots\cdots\cdots\cdots\cdots
\\
&
\ \ \ \ \ 
+
\Big(
A_{i,1}x_1
-\ast x_3-\cdots-\ast x_i
-
\tfrac{2}{i+1}A_{i+1,1}
\Big[
\tfrac{x_1^2}{2!}
\Big\vert
+
\tfrac{x_1^2x_2}{2!}
+\cdots+
\tfrac{x_1^{n-i}x_{n-i}}{(n-i)!}
\Big]
\Big)
\Big(
\tfrac{x_1^i}{i!}
\Big\vert
+
\tfrac{x_1^ix_2}{(i-1)!1!}
+\cdots+
\tfrac{x_1^{n-2}x_{n-i}}{(i-1)!(n-i-1)!}
+
\tfrac{x_1^{n-1}x_{n-i+1}}{(i-1)!(n-i)!}
\Big)
\\
&
\ \ \ \ \ 
+
\cdots\cdots\cdots\cdots\cdots\cdots\cdots\cdots\cdots\cdots\cdots
\cdots\cdots\cdots\cdots\cdots\cdots\cdots\cdots\cdots\cdots\cdots
\cdots\cdots\cdots\cdots\cdots\cdots\cdots\cdots\cdots
\\
&
\ \ \ \ \ 
+
\Big(
A_{n-3,1}x_1
-
\ast x_3-\cdots-\ast x_{n-3}
-
\tfrac{2}{n-2}A_{n-2,1}
\Big[
\tfrac{x_1^2}{2!}
\Big\vert
+
\tfrac{x_1^2x_2}{2!}
+
\tfrac{x_1^3x_3}{3!}
\Big]
\Big)
\Big(
\tfrac{x_1^{n-3}}{(n-3)!}
\Big\vert
+
\tfrac{x_1^{n-3}x_2}{(n-4)!1!}
+
\tfrac{x_1^{n-2}x_3}{(n-4)!2!}
+
\tfrac{x_1^{n-1}x_4}{(n-4)!3!}
\Big)
\\
&
\ \ \ \ \ 
+
\Big(
A_{n-2,1}x_1
-
\ast x_3-\cdots-\ast x_{n-2}
-
\tfrac{2}{n-1}A_{n-1,1}
\Big[
\tfrac{x_1^2}{2!}
\Big\vert
+
\tfrac{x_1^2x_2}{2!}
\Big]
\Big)
\Big(
\tfrac{x_1^{n-2}}{(n-2)!}
\Big\vert
+
\tfrac{x_1^{n-2}x_2}{(n-3)!1!}
+
\tfrac{x_1^{n-1}x_3}{(n-3)!2!}
\Big)
\\
&
\ \ \ \ \ 
+
\Big(
A_{n-1,1}x_1
-
\ast x_3-\cdots-\ast x_{n-1}
+
B_{n-1}
\Big[
\tfrac{x_1^2}{2!}
\Big]
\Big)
\Big(
\tfrac{x_1^{n-1}}{(n-1)!}
\Big\vert
+
\tfrac{x_1^{n-1}x_2}{(n-2)!1!}
\Big)
\\
&
\ \ \ \ \ 
+
\Big(
A_{n,1}x_1
+\cdots+
A_{n,k}x_k
+\cdots+
A_{n,n}x_n
\Big)
\Big(
\tfrac{x_1^n}{n!}
\Big).
\endaligned
\]

We know by the induction assumption that
$\pi_{\ind}^n (\centersmallbullet)$ applied to this gives zero.
Thus, when applying $\pi_{\ind}^{n+1} (\centersmallbullet)$,
it suffices to collect all independent monomials of order
$n+1$, namely $x_1^n x_1$, $x_1^n x_2$, \dots, 
$x_1^n x_n$, and to determine the coefficients
of these monomials, which should all vanish.

Let us explain how to determine the coefficient
of a general `intermediate' monomial $x_1^k x_k$ with
$2 \leqslant k \leqslant n-1$:

\smallskip\noindent$\bullet$\,
the line $i = n-k$ contributes {\em two} terms;

\smallskip\noindent$\bullet$\,
the line $i = n-k+1$ contributes {\em one} term;

\smallskip\noindent$\bullet$\,
the line $i = k$ contributes {\em one} term; 

\smallskip

\noindent
and that is all. Therefore, we obtain:
\[
x_1^kx_k\,
\Big(
-\tfrac{2}{n-k+1}
A_{n-k+1,1}\,
\Big[\,
\tfrac{1}{2}\,
\tfrac{1}{(n-k-1)!(k-1)!}
+
\tfrac{1}{k!\,(n-k)!}
\Big]
+
A_{n-k+1,1}\,
\tfrac{1}{(n-k)!(k-1)!}
+
\tfrac{1}{n!}\,
A_{n,k}
\Big),
\]

The coefficients of $x_1^n x_1$, of
$x_1^n x_2$, can
be seen\big/determined directly from what is written, 
and the coefficients of $x_1^n x_{n-1}$, of
$x_1^n x_n$ as well, so that we obtain:
\[
\aligned
0
&
\,\equiv\,
x_1^{n+1}\,
\Big(
\tfrac{1}{2!}\,
\tfrac{1}{(n-1)!}\,
B_{n-1}
+
\tfrac{1}{n!}\,
A_{n,1}
\Big)
\\
&
\ \ \ \ \ 
+
x_1^nx_2\,
\Big(
\Big[
\zero{
-
\tfrac{2}{n-1}\,
\tfrac{1}{2!}\,\tfrac{1}{(n-3)!1!}
-
\tfrac{2}{n-1}\,
\tfrac{1}{2!(n-2)!}
+
\tfrac{1}{(n-2)!1!}}
\Big]
A_{n-1,1}
+
\tfrac{1}{n!}
A_{n,2}
\Big)
\\
&
\ \ \ \ \ 
+
\cdots\cdots\cdots\cdots\cdots\cdots\cdots\cdots\cdots\cdots\cdots
\cdots\cdots\cdots\cdots\cdots\cdots\cdots\cdots\cdots\cdots\cdots
\\
&
\ \ \ \ \ 
+
x_1^kx_k\,
\Big(
\Big[
-\tfrac{2}{n-k+1}\,
\tfrac{1}{2!}\,
\tfrac{1}{(n-k+1)!(k-1)!}
-\tfrac{2}{n-k+1}\,
\tfrac{1}{k!(n-k)!}
+
\tfrac{1}{(n-k)!(k-1)!}
\Big]
A_{n-k+1,1}
+
\tfrac{1}{n!}\,
A_{n,k}
\Big)
\\
&
\ \ \ \ \ 
+
\cdots\cdots\cdots\cdots\cdots\cdots\cdots\cdots\cdots\cdots\cdots
\cdots\cdots\cdots\cdots\cdots\cdots\cdots\cdots\cdots\cdots\cdots
\\
&
\ \ \ \ \ 
+
x_1^{n-1}x_{n-1}\,
\Big(
\Big[
-\tfrac{2}{2}\,
\tfrac{1}{2!}\,
\tfrac{1}{0!(n-2)!}
-
\tfrac{2}{2}\,
\tfrac{1}{(n-1)!1!}
+
\tfrac{1}{1!(n-2)!}
\Big]\,
A_{2,1}
+
\tfrac{1}{n!}\,
A_{n,n-1}
\Big)
\\
&
\ \ \ \ \ 
+
x_1^nx_n\,
\Big(
\Big[
-2\,\tfrac{1}{n!}
+
\tfrac{1}{(n-1)!}
\Big]
+
\tfrac{1}{n!}\,
A_{n,n}
\Big)
\endaligned
\]
So, we obtain the announced values for
$B_{n-1}$, $A_{n,2}$, \dots, $A_{n,k}$, \dots, $A_{n,n-1}$, 
$A_{n,n}$.
\endproof

In conclusion, the induction on the dimension $n$ is complete.\qed

\SectionHead{Infinitesimal Action at Order $n+2$}
{infinitesimal-action-order-n-2}

Beyond dimension $n \leqslant 6$
(Section~{\ref{stabilizing-G-n-2-3-4-5-6}}), 
we have no explicit
formula for the subgroup of $\GL(n+1,\R)$ which
stabilizes the normalization~{\eqref{u-F-n-plus-1}}
up to order $n+1$ and which would be
valid for any $n \geqslant 2$.
We will therefore proceed in an infinitesimal manner.
This will be less expensive, computationally speaking.

Thanks to Section~{\ref{tangency-order-n-1}}, we may 
take a vector field tangent 
to~{\eqref{u-F-n-plus-1}}
up to order $\leqslant n+1$:
\leqnomode\usetagform{default}
\begin{align}
\label{Lstab-order-n-1}
L
&
\,=\,\ \ \
\big(
A_{1,1}\,x_1
\ \ \ \ \ \ \ \ \ \ \ \ \ \ \ \ \ \ \ \ \ \ \ \ \ \ \ \ \ \ \ \ \ \ \
\ \ \ \ \ \ \ \ \ \ \ \ \ \ \ \ \ \ \ \ \ \ \ \ \ \ \ \ \ \ \ \ \ \ \ 
\ \ 
-
\tfrac{2}{1+1}\,A_{2,1}\,u
\big)\,
\partial_{x_1}
\notag
\\
&
\ \ \ \ \
+
\big(
A_{2,1}\,x_1
\ \ \ \ \ \ \ \ \ \ \ \ \ \ \ \ \ \ \ \ \ \ \ \ \ \ \ \ \ \ \ \ \ \ \
\ \ \ \ \ \ \ \ \ \ \ \ \ \ \ \ \ \ \ \ \ \ \ \ \ \ \ \ \ \ \ \ \ \ \
\ \
-
\tfrac{2}{2+1}\,A_{3,1}\,u
\big)\,
\partial_{x_2}
\notag
\\
&
\ \ \ \ \
+
\big(
A_{3,1}\,x_1
+
0
+
A_{3,1}\,x_3
\ \ \ \ \ \ \ \ \ \ \ \ \ \ \ \ \ \ \ \ \ \ \ \ \ \ \ \ \ \ \ \ \ \ \
\ \ \ \ \ \ \ \ \ \ \ \ \ \
-
\tfrac{2}{3+1}\,A_{4,1}\,u
\big)\,
\partial_{x_3}
\notag
\\
&
\ \ \ \ \
+
\cdots\cdots\cdots\cdots\cdots\cdots\cdots\cdots\cdots\cdots\cdots
\cdots\cdots\cdots\cdots\cdots\cdots\cdots\cdots\cdots\cdots\cdot
\\
&
\ \ \ \ \
+
\big(
A_{n-1,1}x_1
+
0
+
A_{n-1,1}\,x_3
+\cdots+
A_{n-1,n-1}x_{n-1}
\ \ \ \ \
-
\tfrac{2}{n}\,A_{n,1}\,u
\big)\,
\partial_{x_{n-1}}
\notag
\\
&
\ \ \ \ \
+
\big(
A_{n,1}x_1
+
0
\ \ \ \ \
+
A_{n,3}\,x_3
+\cdots+
A_{n,n-1}x_{n-1}
+
A_{n,n}x_n
+
B_n\,u
\big)\,
\partial_{x_n}
\notag
\\
&
\ \ \ \ \
+
\big(
\ \ \ \ \ \ \ \ \ \ \ \ \ \ \ \ \ \ \ \ \ \ \ \ \ \ \ \ \ \ \ \ \ \ \
\ \ \ \ \ \ \ \ \ \ \ \ \ \ \ \ \ \ \ \ \ \ \ \ \ \ \ \ \ \ \ \ \ \ \
\ \ \ \ \ \ \ \ \ \ \ \ \ \ \ \ \ \ \ 
2\,A_{1,1}\,u
\big)\,
\partial_u,
\notag
\end{align}
where, for all $3 \leqslant j \leqslant i \leqslant n$:
\leqnomode\usetagform{default}
\begin{align}
\label{values-A-i-j}
A_{i,j}
\,=\,
-\,\frac{(j-2)}{i-j+1}\,
\binom{i}{j}.
\end{align}

With small $\varepsilon \approx 0$, this $L$ has the approximate flow:
\begin{footnotesize}
\leqnomode\usetagform{default}
\begin{align}
\label{y-x-epsilon}
y_1
&
\,=\,
x_1
+
\varepsilon\,
\Big(
A_{1,1}\,x_1
\ \ \ \ \ \ \ \ \ \ \ \ \ \ \ \ \ \ \ \ \ \ \ \ \ \ \ \ \ \ \ \ \ \ \
\ \ \ \ \ \ \ \ \ \ \ \ \ \ \ \ \ \ \ \ \ \ \ \ \ \ \ \ \ \ \ \ \ \ \ 
\ \ 
-
\tfrac{2}{1+1}\,A_{2,1}\,u
\Big)
+
{\rm O}(\varepsilon^2),
\notag
\\
y_2
&
\,=\,
x_2
+
\varepsilon\,
\Big(
A_{2,1}\,x_1
\ \ \ \ \ \ \ \ \ \ \ \ \ \ \ \ \ \ \ \ \ \ \ \ \ \ \ \ \ \ \ \ \ \ \
\ \ \ \ \ \ \ \ \ \ \ \ \ \ \ \ \ \ \ \ \ \ \ \ \ \ \ \ \ \ \ \ \ \ \
\ \
-
\tfrac{2}{2+1}\,A_{3,1}\,u
\Big)
+
{\rm O}(\varepsilon^2),
\notag
\\
y_3
&
\,=\,
x_3
+
\varepsilon\,
\Big(
A_{3,1}\,x_1
+
0
+
A_{3,1}\,x_3
\ \ \ \ \ \ \ \ \ \ \ \ \ \ \ \ \ \ \ \ \ \ \ \ \ \ \ \ \ \ \ \ \ \ \
\ \ \ \ \ \ \ \ \ \ \ \ \ \
-
\tfrac{2}{3+1}\,A_{4,1}\,u
\Big)
+
{\rm O}(\varepsilon^2),
\notag
\\
\cdot
&
\cdots\cdots\cdots\cdots\cdots\cdots\cdots\cdots\cdots\cdots\cdots
\cdots\cdots\cdots\cdots\cdots\cdots\cdots\cdots\cdots\cdots\cdots
\cdots\cdots\cdots\cdots\cdots
\cdot
\\
y_{n-1}
&
\,=\,
x_{n-1}
+
\varepsilon\,
\Big(
A_{n-1,1}x_1
+
0
+
A_{n-1,1}\,x_3
+\cdots+
A_{n-1,n-1}x_{n-1}
\ \
-
\tfrac{2}{n}\,A_{n,1}\,u
\Big)
+
{\rm O}(\varepsilon^2),
\notag
\\
y_n
&
\,=\,
x_n
+
\varepsilon\,
\Big(
A_{n,1}x_1
+
0
\ \ \ \ \ \
+
A_{n,3}\,x_3
+\cdots+
A_{n,n-1}x_{n-1}
+
A_{n,n}x_n
+
B_n\,u
\Big)
+
{\rm O}(\varepsilon^2),
\notag
\\
v
&
\,=\,
u\ \
+
\varepsilon\,
\Big(
\ \ \ \ \ \ \ \ \ \ \ \ \ \ \ \ \ \ \ \ \ \ \ \ \ \ \ \ \ \ \ \ \ \ \
\ \ \ \ \ \ \ \ \ \ \ \ \ \ \ \ \ \ \ \ \ \ \ \ \ \ \ \ \ \ \ \ \ \ \
\ \ \ \ \ \ \ \ \ \ \ \ \ \ \ \ \ \ \ \ \ \ \ \ 
2\,A_{1,1}\,u
\Big)
+
{\rm O}(\varepsilon^2).
\notag
\end{align}
\end{footnotesize}

For $1 \leqslant i \leqslant n-1$, we may write the intermediate
lines as:
\leqnomode\usetagform{default}
\begin{align}
\label{y-i-x-i-epsilon}
\ \ \ \ \ 
y_i
&
\,=\,
x_i
+
\varepsilon\,
\Big(
A_{i,i}\,x_1
+
0
+
\sum_{3\leqslant j\leqslant i}\,
A_{i,j}\,x_j
+
0+\cdots+0
-
\tfrac{2}{i+1}\,A_{i+1,1}\,u
\Big)
+
{\rm O}(\varepsilon^2).
\end{align}

With its independent monomials of order $n+2$,
the hypersurface equation on the left is:
\leqnomode\usetagform{default}
\begin{align}
\label{u-F-order-n-2}
u
&
\,=\,
\frac{x_1^2}{2}
+
\frac{x_1^2x_2}{2}
+
\sum_{m=3}^n\,
\Big(
\frac{x_1^mx_m}{m!}
+
x_1^{m-1}
\sum_{i,j\geqslant 2
\atop
i+j=m+1}\,
\tfrac{1}{2}\,
\frac{x_ix_j}{(i-1)!(j-1)!}
\Big)
\notag
\\
&
\ \ \ \ \
+
F_{n+2,0,\dots,0}\,
\frac{x_1^{n+1}x_1}{(n+2)!}
+
F_{n+1,1,\dots,0}\,
\frac{x_1^{n+1}x_2}{(n+1)!}
+\cdots+
F_{n+1,0,\dots,1}\,
\frac{x_1^{n+1}x_n}{(n+1)!}
\\
&
\ \ \ \ \
+
x_1^n
\sum_{i,j\geqslant 2
\atop
i+j=n+2}\,
\tfrac{1}{2}\,
\frac{x_ix_j}{(i-1)!(j-1)!}
+
{\rm O}_{x'}(3)
+
{\rm O}_x(n+3).
\notag
\end{align}
Similarly,
the hypersurface equation on the right is:
\leqnomode\usetagform{default}
\begin{align}
\label{y-G-order-n-2}
v
&
\,=\,
\frac{y_1^2}{2}
+
\frac{y_1^2y_2}{2}
+
\sum_{m=3}^n\,
\Big(
\frac{y_1^my_m}{m!}
+
y_1^{m-1}
\sum_{i,j\geqslant 2
\atop
i+j=m+1}\,
\tfrac{1}{2}\,
\frac{y_iy_j}{(i-1)!(j-1)!}
\Big)
\notag
\\
&
\ \ \ \ \
+
G_{n+2,0,\dots,0}\,
\frac{y_1^{n+1}y_1}{(n+2)!}
+
G_{n+1,1,\dots,0}\,
\frac{y_1^{n+1}y_2}{(n+1)!}
+\cdots+
G_{n+1,0,\dots,1}\,
\frac{y_1^{n+1}y_n}{(n+1)!}
\\
&
\ \ \ \ \
+
y_1^n
\sum_{i,j\geqslant 2
\atop
i+j=n+2}\,
\tfrac{1}{2}\,
\frac{y_iy_j}{(i-1)!(j-1)!}
+
{\rm O}_{y'}(3)
+
{\rm O}_y(n+3).
\notag
\end{align}

Now, assume that $(y_1, \dots, y_n, v)$ are replaced
from~{\eqref{y-x-epsilon}}
in $0 = - v + G(y)$, that 
$u$ is replaced by $F(x_1, \dots, x_n)$
from~{\eqref{u-F-order-n-2}}, and give a name
to the concerned fundamental equation (to be studied in a while):
\[
0
\,\equiv\,
E(x_1,\dots,x_n,\varepsilon)
\,=\,
E(x,\varepsilon).
\]
Consider terms of order $\leqslant n+1$:
\[
\pi^{n+1}
\bigg(
\sum_{\sigma_1,\dots,\sigma_n}\,
E_{\sigma_1,\dots,\sigma_n}(\varepsilon)\,
x_1^{\sigma_1}\cdots x_n^{\sigma_n}
\bigg)
\,=\,
\sum_{\sigma_1+\dots+\sigma_n\leqslant n+1}\,
E_{\sigma_1,\dots,\sigma_n}(\varepsilon)\,
x_1^{\sigma_1}\cdots x_n^{\sigma_n}.
\]

\begin{Assertion}
\label{Ass-tangency-L-order-n-1}
The fact that $L$ is tangent up to order $n+1$ implies:
\[
0
\,\equiv\,
\pi^{n+1}
\big(
-v+G(y)
\big)
\,=\,
\pi^{n+1}\,
\big(
E(x,\varepsilon)
\big).
\eqno\qed
\]
\end{Assertion}

Hence, look at order $n+2$ terms, namely at:
\[
\aligned
0
&
\,\equiv\,
\pi^{n+2}\big(E(x,\varepsilon)\big)
\,=\,
\pi_{n+2}\big(E(x,\varepsilon)\big)
\\
&
\,=:\,
E_{n+2,0,\dots,0}(\varepsilon)\,
x_1^{n+1}x_1
+
E_{n+1,1,\dots,0}(\varepsilon)\,
x_1^{n+1}x_2
+\cdots+
E_{n+1,0,\dots,1}(\varepsilon)\,
x_1^{n+1}x_n.
\endaligned
\]
For $\varepsilon = 0$, the map is the identity, hence without
computation we know that:
\[
0
\,\equiv\,
\sum_{\nu_1+\cdots+\nu_n=n+2}\,
\bigg(
-
\frac{F_{\nu_1,\dots,\nu_n}}{\nu_1!\cdots\nu_n!}
+
\frac{G_{\nu_1,\dots,\nu_n}}{\nu_1!\cdots\nu_n!}
+
\varepsilon\,T_{\nu_1,\dots,\nu_n}
+
{\rm O}(\varepsilon^2)
\bigg)
+
{\rm O}_{x_1,\dots,x_n}(n+3),
\]
and our key goal is to determine these terms $T_{\nu_1, \dots,
\nu_n}$ that are of order $1$ in $\varepsilon$\,\,---\,\,more 
precisely,
to compute the independent ones:
\leqnomode\usetagform{default}
\begin{align}
\label{names-T}
T_{n+2,0,\dots,0},
\ \ \ \ \
T_{n+1,1,\dots,0},
\ \ \ \ \
\dots,
\ \ \ \ \
T_{n+1,0,\dots,1}.
\end{align}

Equivalently, we may write the fundamental equation
by specifying the dependent monomials as the
remainder ${\rm O}_{x_2, \dots, x_n} (2)$:
\leqnomode\usetagform{default}
\begin{align}
\label{F-G-epsilon}
0
&
\,\equiv\,
E(x,\varepsilon)
\notag
\\
&
\,\equiv\,
x_1^{n+1}x_1
\bigg(
-
\frac{F_{n+2,0,\dots,0}}{(n+2)!}
+
\frac{G_{n+2,0,\dots,0}}{(n+2)!}
+
\varepsilon\,T_{n+2,0,\dots,0}
+
{\rm O}(\varepsilon^2)
\bigg)
\notag
\\
&
\ \ \ \ \
+
x_1^{n+1}x_2
\bigg(
-
\frac{F_{n+1,1,\dots,0}}{(n+1)!}
+
\frac{G_{n+1,1,\dots,0}}{(n+1)!}
+
\varepsilon\,T_{n+1,1,\dots,0}
+
{\rm O}(\varepsilon^2)
\bigg)
\\
&
\ \ \ \ \
+
\cdots\cdots\cdots\cdots\cdots\cdots\cdots\cdots\cdots\cdots\cdots
\cdots\cdots\cdots\cdots\cdots\cdots\cdots
\notag
\\
&
\ \ \ \ \
+
x_1^{n+1}x_n
\bigg(
-
\frac{F_{n+1,0,\dots,1}}{(n+1)!}
+
\frac{G_{n+1,0,\dots,1}}{(n+1)!}
+
\varepsilon\,T_{n+1,0,\dots,1}
+
{\rm O}(\varepsilon^2)
\bigg)
\notag
\\
&
\ \ \ \ \
+
{\rm O}_{x_2,\dots,x_n}(2)
+
{\rm O}_{x_1,x_2,\dots,x_n}(n+3).
\notag
\end{align}

Introduce the operator:
\[
\Pi_{n+2}^{\ind}
(\centersmallbullet)
\,:=\,
\frac{d}{d\varepsilon}\Big\vert_{\varepsilon=0}
\Big(
\pi_{n+2}^{\ind}
(\centersmallbullet)
\Big)
\,=\,
\pi_{n+2}^{\ind}
\Big(
\frac{d}{d\varepsilon}\Big\vert_{\varepsilon=0}
(\centersmallbullet)
\Big),
\]
which selects what we want to compute:
\[
\Pi_{n+2}^{\ind}
\big(E(x,\varepsilon)\big)
\,=\,
\frac{x_1^{n+1}x_1}{(n+2)!}\,
T_{n+2,0,\dots,0}
+
\frac{x_1^{n+1}x_2}{(n+1)!}\,
T_{n+1,1,\dots,0}
+\cdots+
\frac{x_1^{n+1}x_n}{(n+1)!}\,
T_{n+1,0,\dots,1}.
\]
Also, for some use about remainders, we will need
to consider all independent monomials of order $\leqslant n+2$
(notice the placement\big/switch of indices):
\[
\Pi_{\ind}^{n+2}
(\centersmallbullet)
\,:=\,
\frac{d}{d\varepsilon}\Big\vert_{\varepsilon=0}
\Big(
\pi_{\ind}^{n+2}
(\centersmallbullet)
\Big)
\,=\,
\pi_{\ind}^{n+2}
\Big(
\frac{d}{d\varepsilon}\Big\vert_{\varepsilon=0}
(\centersmallbullet)
\Big).
\]

Thus in summary, we must apply $\Pi_{n+2} (\centersmallbullet)$
to all (numerous) terms of~{\eqref{y-G-order-n-2}}.
We start by the remainders.

\begin{Lemma}
One has:
\[
\frac{d}{d\varepsilon}\Big\vert_{\varepsilon=0}
\big(
{\rm O}_{y_2,\dots,y_n}(3)
\big)
\,=\,
{\rm O}_{x_2,\dots,x_n}(2).
\]
\end{Lemma}

\proof
Take any monomial $y_2^{\nu_2} \cdots y_n^{\nu_n}$
with $\nu_2 + \cdots + \nu_n \geqslant 3$,
and abbreviate~{\eqref{y-x-epsilon}} as:
\[
y_i
\,=\,
x_i
+
\varepsilon\,R_i(x)
+
{\rm O}(\varepsilon^2)
\eqno
{\scriptstyle{(2\,\leqslant\,i\,\leqslant\,n)}},
\]
whence:
\[
\footnotesize
\aligned
\big(
x_2+\varepsilon\,R_2+{\rm O}(\varepsilon^2)
\big)^{\nu_2}
&
\cdots
\big(
x_n+\varepsilon\,R_n+{\rm O}(\varepsilon^2)
\big)^{\nu_n}
\,=\,
\\
&
\,=\,
\Big(
x_2^{\nu_2}
+
\nu_2\,x_2^{\nu_2-1}\,\varepsilon\,R_2
+
{\rm O}(\varepsilon^2)
\Big)
\cdots
\Big(
x_n^{\nu_n}
+
\nu_n\,x_n^{\nu_n-1}\,\varepsilon\,R_n
+
{\rm O}(\varepsilon^2)
\Big)
\\
&
\,=\,
x_2^{\nu_2}\cdots x_n^{\nu_n}
+
\varepsilon\,
\Big[
\nu_2\,x_2^{\nu_2-1}R_2x_3^{\nu_3}\cdots x_n^{\nu_n}
+\cdots+
x_2^{\nu_2}\cdots x_{n-1}^{\nu_{n-1}}
\nu_n x_n^{\nu_n-1}R_n
\Big]
+
{\rm O}(\varepsilon^2),
\endaligned
\]
and here $\nu_2 - 1 + \nu_3 + \cdots + \nu_n \geqslant 2$,
\dots, $\nu_2 + \cdots + \nu_{n-1} + \nu_n - 1 \geqslant 2$
as well.
\endproof

Thus, the remainder ${\rm O}_{y_2, \dots, y_n} (3)$
in~{\eqref{y-G-order-n-2}} has contribution equal to $0$
in $\Pi_{\ind}^{n+2}(E)$.

\smallskip

Next, still in~{\eqref{y-G-order-n-2}},
consider border-dependent monomials from the sum
$\sum_{m=3}^n$, but only for $m \leqslant n-1$ at first.

\begin{Lemma}
For every $3 \leqslant m \leqslant n-1$
and for all $i, j \geqslant 2$ with $i + j = m+1$,
one has:
\[
0
\,=\,
\Pi_{n+2}^{\ind}
\big(
y_1^{m-1}y_iy_j
\big).
\]
\end{Lemma}

\proof
To simplify~{\eqref{y-i-x-i-epsilon}}
after replacement of $u$ by $F$, 
abbreviate:
\[
\alpha_i
\,:=\,
-\,\tfrac{2}{i+1}\,
\tfrac{A_{i+1}}{(n+1-m)!},
\]
and use $\cdots$ to denote
monomials of order $\leqslant n-m+1$:
\[
\footnotesize
\aligned
y_i
&
\,=\,
x_i
+
\varepsilon\,
\bigg\{
A_{i,1}
+
\sum_{3\leqslant j\leqslant i}\,
A_{i,j}\,x_j
-
\tfrac{2}{i+1}\,
A_{i+1,1}\,
\Big[
\cdots
+
\frac{x_1^{n+1-m}x_{n+1-m}}{(n+1-m)!}
+
{\rm O}_x(n-m+3)
\Big]
\bigg\}
+
{\rm O}(\varepsilon^2)
\\
&
\,=\,
x_i
+
\varepsilon\,
\Big\{
\cdots
+
\alpha_i\,
x_1^{n+1-m}\,x_{n+1-m}
+
{\rm O}_x(n-m+3)
\Big\}
+
{\rm O}(\varepsilon^2),
\endaligned
\]
write the product:
\[
\aligned
y_1^{m-1}y_iy_j
&
\,=\,
\Big(
x_1
+
\varepsilon\,
\Big\{
\cdots
+
\alpha_1\,x_1^{n+1-m}\,x_{n+1-m}
+
{\rm O}_x(n-m+3)
\Big\}
+
{\rm O}(\varepsilon^2)
\Big)^{m-1}
\\
&
\ \ \ \ \
\cdot
\bigg(
x_i
+
\varepsilon\,
\Big\{
\cdots
+
\alpha_i\,x_1^{n+1-m}\,x_{n+1-m}
+
{\rm O}_x(n-m+3)
\Big\}
+
{\rm O}(\varepsilon^2)
\bigg)^1
\\
&
\ \ \ \ \
\cdot
\bigg(
x_j
+
\varepsilon\,
\Big\{
\cdots
+
\alpha_j\,x_1^{n+1-m}\,x_{n+1-m}
+
{\rm O}_x(n-m+3)
\Big\}
+
{\rm O}(\varepsilon^2)
\bigg)^1,
\endaligned
\]
and expand:
\[
\aligned
\pi_{n+2}
\big(
y_1^{m-1}y_iy_j
\big)
\,=\,
\zero{
\pi_{n+2}
\big(
x_1^{m-1}x_ix_j
\big)}
+
\varepsilon\,
\Big\{
&
x_1^{m-2}\,\tbinom{m-1}{1}\,\alpha_1\,x_1^{n+1-m}\,x_{n+1-m}\,x_i\,x_j
\\
&
+
x_1^{m-1}\,\alpha_i\,x_1^{n+1-m}\,x_{n+1-m}\,x_j
\\
&
+
x_1^{m-1}\,x_i\,\alpha_j\,x_1^{n+1-m}\,x_{n+1-m}
\Big\}
+
{\rm O}(\varepsilon^2).
\endaligned
\]
To conclude, 
observe that $x_i x_j$ are dependent monomials,
and since $n-m+1 \geqslant 2$, observe that the
two monomials $x_{n+1-m}\,x_j$ and 
$x_{n+1-m}\,x_i$ are also dependent.
\endproof

It therefore remains to compute:
\[
\aligned
\Pi_{n+2}^{\ind}
\bigg(
-\,v
&
+
\frac{y_1^2}{2!}
+
\frac{y_1^2y_2}{2!}
+
\frac{y_1^3y_3}{3!}
+\cdots+
\frac{y_1^{n-1}y_{n-1}}{(n-1)!}
\\
&
+
\frac{y_1^ny_n}{n!}
+
y_1^{n-1}
\sum_{i,j\geqslant 2
\atop
i+j=n+1}\,
\tfrac{1}{2}\,
\frac{y_iy_j}{(i-1)!(j-1)!}
\\
&
+
G_{n+2,0,\dots,0}\,
\frac{y_1^{n+1}y_1}{(n+2)!}
+
G_{n+1,1,\dots,0}\,
\frac{y_1^{n+1}y_2}{(n+1)!}
+\cdots+
G_{n+1,0,\dots,1}\,
\frac{y_1^{n+1}y_n}{(n+1)!}
\\
&
\ \ \ \ \ \ \ \ \ \ \ \ \ \ \ \ \ \ \ \ \ \ \ \ \ \ \ \ \ \ \ \ \ \ \
\ \ \ \ \ \ \ \ \ \ \ \ \ \ \ \ \ \ \ \ \ \ \ \ \ \ \ \ \ \ \ \ \ \ \
+
y_1^n
\sum_{i,j\geqslant 2
\atop
i+j=n+2}\,
\tfrac{1}{2}\,
\frac{y_iy_j}{(i-1)!(j-1)!}
\bigg).
\endaligned
\]
Since $\Pi_{n+2}^{\ind} (\centersmallbullet)$ is linear,
we can proceed termwise.

Compute the first term:
\[
\small
\aligned
\Pi_{n+2}^{\ind}
(-v)
&
\,=\,
\pi_{n+2}^{\ind}
\bigg(
\frac{d}{d\varepsilon}\Big\vert_{\varepsilon=0}
\Big(
-u
-
\varepsilon\,2\,A_{1,1}\,u
+
{\rm O}(\varepsilon^2)
\Big)
\bigg)
\\
&
\,=\,
\pi_{n+2}^{\ind}
\bigg(
-2\,A_{1,1}\,
\Big[
\frac{x_1^2}{2!}
+
\frac{x_1^2x_2}{2!}
+\cdots+
\frac{x_1^nx_n}{n!}
\\
&
\ \ \ \ \
+
F_{n+2,0,\dots,0}\,
\frac{x_1^{n+1}x_1}{(n+2)!}
+
F_{n+1,1,\dots,0}\,
\frac{x_1^{n+1}x_2}{(n+1)!}
+\cdots+
F_{n+1,0,\dots,1}\,
\frac{x_1^{n+1}x_n}{(n+1)!}
+
{\rm O}_{x_1,\dots,x_n}(n+3)
\Big]
\bigg),
\endaligned
\]
and obtain:
\leqnomode\usetagform{default}
\begin{align}
\label{1-term-n-plus-2}
\Pi_{n+2}^{\ind}(-v)
&
\,=\,
-\,\frac{2}{(n+2)!}\,
F_{n+2,0,\dots,0}\,
A_{1,1}\,x_1^{n+2}
-
\frac{2}{(n+1)!}\,
F_{n+1,1,\dots,0}\,
A_{1,1}\,x_1^{n+1}x_2
-\cdots-
\notag
\\
&
\ \ \ \ \
-\cdots-
\frac{2}{(n+1)!}\,
F_{n+1,0,\dots,1}\,
A_{1,1}\,x_1^{n+1}x_n.
\end{align}

Next, for $\kappa \in \N$, abbreviating:
\[
\mathcal{R}^\kappa
\,:=\,
{\rm O}_{x_2,\dots,x_n}(2)
+
{\rm O}_{x_1,x_2,\dots,x_n}(\kappa),
\]
compute the second term:
\[
\!\!\!\!\!\!\!\!\!\!\!\!\!\!\!
\footnotesize
\aligned
\Pi_{n+2}^{\ind}
\big(
\tfrac{y_1^2}{2}
\big)
&
\,=\,
\Pi_{n+2}^{\ind}
\bigg(
\tfrac{1}{2}
\Big(
x_1
+
\varepsilon
\Big\{
A_{1,1}x_1
-
\tfrac{2}{2}A_{2,1}
\Big[
\tfrac{x_1^2}{2}
+\cdots+
\tfrac{x_1^{n-1}x_{n-1}}{(n-1)!}
+
\tfrac{x_1^nx_n}{n!}
+
{\rm O}_{x'}(2)
+
{\rm O}_x(n+2)
\Big]
\Big\}
+
{\rm O}(\varepsilon^2)
\Big)^2
\bigg)
\\
&
\,=\,
\Pi_{n+2}^{\ind}
\bigg(
\tfrac{1}{2}
\Big(
x_1
+
\varepsilon
\Big\{
\cdots
-
\tfrac{2}{2}A_{2,1}\,
\tfrac{x_1^nx_n}{n!}
+
\mathcal{R}^{n+2}
\Big\}
+
{\rm O}(\varepsilon^2)
\Big)^2
\bigg),
\endaligned
\]
expand the square, select the coefficient of $\varepsilon^1$,
select the (single) independent monomial, and obtain:
\leqnomode\usetagform{default}
\begin{align}
\label{2-term-n-plus-2}
\Pi_{n+2}^{\ind}
\big(
\tfrac{y_1^2}{2}
\big)
\,=\,
\tfrac{1}{2}\,2\,x_1\,
\big(
-\tfrac{2}{2}
A_{2,1}
\big)\,
\tfrac{x_1^nx_n}{n!}
\,=\,
-\,\tfrac{1}{n!}\,
A_{2,1}\,x_1^{n+1}x_n.
\end{align}

Treat the third term:
\[
\aligned
\Pi_{n+2}^{\ind}
\big(
\tfrac{y_1^2y_2}{2}
\big)
\,=\,
\Pi_{n+2}^{\ind}
\bigg(
\tfrac{1}{2!}
&
\Big(
x_1
+
\varepsilon\,
\Big[
\cdots
-
\tfrac{2}{2}\,A_{2,1}\,
\tfrac{x_1^{n-1}x_{n-1}}{(n-1)!}
+
\mathcal{R}^{n+1}
\Big]
+
{\rm O}(\varepsilon^2)
\Big)^2
\\
\cdot
&
\Big(
x_2
+
\varepsilon\,
\Big[
\cdots
-
\tfrac{2}{3}\,A_{3,1}\,
\tfrac{x_1^{n-1}x_{n-1}}{(n-1)!}
+
\mathcal{R}^{n+1}
\Big]
+
{\rm O}(\varepsilon^2)
\Big)^1
\bigg)
\endaligned
\]
that is:
\leqnomode\usetagform{default}
\begin{align}
\label{3-term-n-plus-2}
\ \ \ \ \ \ \ \ \ \
\Pi_{n+2}^{\ind}
\big(
\tfrac{y_1^2y_2}{2}
\big)
\,=\,
\tfrac{1}{2}\,
x_1^2\,
\big(
-\tfrac{2}{3}\,A_{3,1}
\big)\,
\tfrac{x_1^{n-1}x_{n-1}}{(n-1)!}
\,=\,
-\,\tfrac{2}{3}\,
\tfrac{1}{2!\,(n-1)!}\,
A_{3,1}\,
x_1^{n+1}x_{n+1}.
\end{align}

Similarly:
\[
\aligned
\Pi_{n+2}^{\ind}
\big(
\tfrac{y_1^3y_3}{3!}
\big)
\,=\,
\Pi_{n+2}^{\ind}
\bigg(
\tfrac{1}{3!}
&
\Big(
x_1
+
\varepsilon\,
\Big[
\cdots
-
\tfrac{2}{2}\,A_{2,1}\,
\tfrac{x_1^{n-2}x_{n-2}}{(n-2)!}
+
\mathcal{R}^{n}
\Big]
+
{\rm O}(\varepsilon^2)
\Big)^3
\\
\cdot
&
\Big(
x_3
+
\varepsilon\,
\Big[
\cdots
-
\tfrac{2}{4}\,A_{4,1}\,
\tfrac{x_1^{n-2}x_{n-2}}{(n-2)!}
+
\mathcal{R}^{n}
\Big]
+
{\rm O}(\varepsilon^2)
\Big)^1
\bigg)
\endaligned
\]
that is:
\[
\Pi_{n+2}^{\ind}
\big(
\tfrac{y_1^3y_3}{3!}
\big)
\,=\,
\tfrac{1}{3!}\,
x_1^3\,
\big(
-\tfrac{2}{4}\,A_{4,1}
\big)\,
\tfrac{x_1^{n-2}x_{n-2}}{(n-2)!}
\,=\,
-\,\tfrac{2}{4}\,
\tfrac{1}{3!\,(n-2)!}\,
A_{4,1}\,
x_1^{n+1}x_{n-2}.
\]

Now, consider general $m$ with $3 \leqslant m \leqslant n-1$:
\[
\aligned
\Pi_{n+2}^{\ind}
\big(
\tfrac{y_1^my_m}{m!}
\big)
\,=\,
\Pi_{n+2}^{\ind}
\bigg(
\tfrac{1}{m!}
&
\Big(
x_1
+
\varepsilon\,
\Big[
\cdots
-
\tfrac{2}{2}\,A_{2,1}\,
\tfrac{x_1^{n-m+1}x_{n-m+1}}{(n-m+1)!}
+
\mathcal{R}^{n-m+3}
\Big]
+
{\rm O}(\varepsilon^2)
\Big)^m
\\
\cdot
&
\Big(
x_m
+
\varepsilon\,
\Big[
\cdots
-
\tfrac{2}{m+1}\,A_{m+1,1}\,
\tfrac{x_1^{n-m+1}x_{n-m+1}}{(n-m+1)!}
+
\mathcal{R}^{n-m+3}
\Big]
+
{\rm O}(\varepsilon^2)
\Big)^1
\bigg),
\endaligned
\]
that is:
\leqnomode\usetagform{default}
\begin{align}
\label{4-term-n-plus-2}
\Pi_{n+2}^{\ind}
\big(
\tfrac{y_1^my_m}{m!}
\big)
&
\,=\,
\tfrac{1}{m!}\,
x_1^m\,
\big(
-\tfrac{2}{m+1}\,A_{m+1,1}
\big)\,
\tfrac{x_1^{n-m+1}x_{n-m+1}}{(n-m+1)!}
\notag
\\
&
\,=\,
-\
\tfrac{2}{m+1}\,
\tfrac{1}{m!(n-m+1)!}\,
A_{m+1,1}\,
x_1^{n+1}x_{n-m+1}.
\end{align}

For $m = n$, the result is different, two monomials will be obtained:
\[
\aligned
\Pi_{n+2}^{\ind}
\big(
\tfrac{1}{n!}\,y_1^ny_n
\big)
&
\,=\,
\Pi_{n+2}^{\ind}
\bigg(
\tfrac{1}{n!}\,
\Big(
x_1
+
\varepsilon\,
\Big[
\cdots
-
\tfrac{2}{2}\,A_{2,1}\,
\tfrac{x_1^2}{2!}
+
\mathcal{R}^3
\Big]
+
{\rm O}(\varepsilon^2)
\Big)^n
\\
&
\ \ \ \ \ \ \ \ \ \ \ \ \ \ 
\cdot
\Big(
x_n
+
\varepsilon\,
\Big[
\cdots
+
B_n\,\tfrac{x_1^2}{2}
+
\mathcal{R}^3
\Big]
+
{\rm O}(\varepsilon^2)
\Big)^1
\bigg)^n
\\
&
\,=\,
\tfrac{1}{n!}\,x_1^n\,B_n\,\tfrac{x_1^2}{2!}
+
\tfrac{1}{n!}\,
{\textstyle{\binom{n}{1}}}\,
x_1^{n-1}\,
\big(
-\tfrac{2}{2}\,A_{2,1}
\big)\,
\tfrac{x_1^2}{2!}\,x_n,
\endaligned
\]
that is:
\leqnomode\usetagform{default}
\begin{align}
\label{5-term-n-plus-2}
\Pi_{n+2}^{\ind}
\big(
\tfrac{1}{n!}\,y_1^ny_n
\big)
\,=\,
\tfrac{1}{n!\,2!}\,B_n\,x_1^{n+2}
-
\tfrac{1}{(n-1)!\,2!}\,
A_{2,1}\,x_1^{n+1}\,x_n.
\end{align}

Next, with any $i, j \geqslant 2$ satisfying $i+j = n+1$,
whence $i, j \leqslant n-1$:
\[
\aligned
\Pi_{n+2}^{\ind}
\big(
y_1^{n-1}
\tfrac{1}{2}\,
\tfrac{y_iy_j}{(i-1)!(j-1)!}
\big)
&
\,=\,
\tfrac{1}{2}\,
\tfrac{1}{(i-1)!(j-1)!}\,
\Pi_{n+2}^{\ind}
\bigg(
\Big(
x_1
+
\varepsilon\,
\Big[
\cdots
-
\tfrac{2}{2}\,A_{2,1}\,\tfrac{x_1^2}{2}
+
\mathcal{R}^3
\Big]
+
{\rm O}(\varepsilon^2)
\Big)^{n-1}
\\
&
\ \ \ \ \ \ \ \ \ \ \ \ \ \ \ \ \ \ \ \ \ \ \ \ \ \ \ \ \ \ \ \ \ \ \ \
\cdot
\Big(
x_i
+
\varepsilon\,
\Big[
\cdots
-
\tfrac{2}{i+1}\,A_{i+1,1}\,\tfrac{x_1^2}{2}
+
\mathcal{R}^3
\Big]
+
{\rm O}(\varepsilon^2)
\Big)^1
\\
&
\ \ \ \ \ \ \ \ \ \ \ \ \ \ \ \ \ \ \ \ \ \ \ \ \ \ \ \ \ \ \ \ \ \ \ \
\cdot
\Big(
x_j
+
\varepsilon\,
\Big[
\cdots
-
\tfrac{2}{j+1}\,A_{j+1,1}\,\tfrac{x_1^2}{2}
+
\mathcal{R}^3
\Big]
+
{\rm O}(\varepsilon^2)
\Big)^1
\bigg)
\\
&
\,=\,
0
\\
&
\ \ \ \ \
+
\tfrac{1}{2}\,
\tfrac{1}{(i-1)!(j-1)!}\,
x_1^{n-1}\,
\big(
-\tfrac{2}{i+1}\,A_{i+1,1}
\big)\,
\tfrac{x_1^2}{2}\,x_j
\\
&
\ \ \ \ \
+
\tfrac{1}{2}\,
\tfrac{1}{(i-1)!(j-1)!}\,
x_1^{n-1}\,
\big(
-\tfrac{2}{j+1}\,A_{j+1,1}
\big)\,
\tfrac{x_1^2}{2}\,x_i,
\endaligned
\]
that is:
\[
\Pi_{n+2}^{\ind}
\big(
y_1^{n-1}
\tfrac{1}{2}\,
\tfrac{y_iy_j}{(i-1)!(j-1)!}
\big)
\,=\,
-\,
\tfrac{1}{2(i+1)}\,
\tfrac{1}{(i-1)!(j-1)!}\,
A_{i+1,1}\,x_1^{n+1}x_j
-
\tfrac{1}{2(j+1)}\,
\tfrac{1}{(i-1)!(j-1)!}\,
A_{j+1,1}\,x_1^{n+1}x_i.
\]
In fact, we must perform a sum:
\[
\aligned
\Pi_{n+2}^{\ind}
\bigg(
y_1^{n-1}
\sum_{i,j\geqslant 2
\atop
i+j=n+1}\,
\tfrac{1}{2}\,
\frac{y_iy_j}{(i-1)!(j-1)!}
\bigg)
&
\,=\,
-\,\tfrac{1}{2}
\sum_{j=2}^{n-1}\,
\frac{1}{n-j+2}\,
\frac{1}{(n-j)!(j-1)!}\,
A_{n-j+2,1}\,
x_1^{n+1}x_j
\\
&
\ \ \ \ \
-\,\tfrac{1}{2}
\sum_{i=2}^{n-1}\,
\frac{1}{n-i+2}\,
\frac{1}{(i-1)!(n-i)!}\,
A_{n-i+2,1}\,
x_1^{n+1}x_i
\endaligned
\]
and by exchanging $j \longleftrightarrow i$ in the 
first sum, the two sums are {\em equal}, whence:
\leqnomode\usetagform{default}
\begin{align}
\label{6-term-n-plus-2}
\ \ \ \ \ \ \ \ \ \
\Pi_{n+2}^{\ind}
\Big(
y_1^{n-1}
\sum_{i,j\geqslant 2
\atop
i+j=n+1}\,
\tfrac{1}{2}\,
\tfrac{y_iy_j}{(i-1)!(j-1)!}
\Big)
\,=\,
-\,
\sum_{i=2}^{n-1}\,
\tfrac{1}{n-i+2}\,
\tfrac{1}{(i-1)!(n-i)!}\,
A_{n-i+2,1}\,
x_1^{n+1}x_i.
\end{align}

The next term:
\[
\aligned
\Pi_{n+2}^{\ind}
\Big(
G_{n+2,0,\dots,0}\,
\tfrac{y_1^{n+2}}{(n+2)!}
\Big)
&
\,=\,
\tfrac{G_{n+2,0,\dots,0}}{(n+2)!}\,
\Pi_{n+2}^{\ind}
\bigg(
\Big(
x_1
+
\varepsilon\,
\Big[
A_{1,1}\,x_1
+
\mathcal{R}^2
\Big]
+
{\rm O}(\varepsilon^2)
\Big)^{n+2}
\bigg)
\\
&
\,=\,
\tfrac{G_{n+2,0,\dots,0}}{(n+2)!}\,
{\textstyle{\binom{n+2}{1}}}\,
x_1^{n+1}\,A_{1,1}\,x_1
\endaligned
\]
has value:
\leqnomode\usetagform{default}
\begin{align}
\label{7-term-n-plus-2}
\Pi_{n+2}^{\ind}
\Big(
G_{n+2,0,\dots,0}\,
\tfrac{y_1^{n+2}}{(n+2)!}
\Big)
\,=\,
\tfrac{G_{n+2,0,\dots,0}}{(n+1)!}\,
A_{1,1}\,
x_1^{n+2}.
\end{align}

Then:
\[
\aligned
\Pi_{n+2}^{\ind}
\Big(
G_{n+1,1,\dots,0}\,
\tfrac{y_1^{n+1}y_2}{(n+1)!}
\Big)
&
\,=\,
\tfrac{G_{n+1,1,\dots,0}}{(n+1)!}\,
\Pi_{n+2}^{\ind}
\bigg(
\Big(
x_1
+
\varepsilon\,
\Big[
A_{1,1}\,x_1
+
\mathcal{R}^2
\Big]
+
{\rm O}(\varepsilon^2)
\Big)^{n+1}
\bigg)
\\
&
\ \ \ \ \ \ \ \ \ \ \ \ \ \ \ \ \ \ \ \ \ \ \ \ \ \ \ \ \ \ \ \ 
\cdot
\Big(
x_2
+
\varepsilon\,
\Big[
A_{2,1}\,x_1
+
\mathcal{R}^2
\Big]
+
{\rm O}(\varepsilon^2)
\Big)^1
\bigg)
\\
&
\,=\,
\tfrac{G_{n+1,1,\dots,0}}{(n+1)!}\,
x_1^{n+1}\,A_{2,1}\,x_1
+
\tfrac{G_{n+1,1,\dots,0}}{(n+1)!}\,
x_1^n\,
{\textstyle{\binom{n+1}{1}}}\,
A_{1,1}\,x_1x_2,
\endaligned
\]
that is:
\leqnomode\usetagform{default}
\begin{align}
\label{8-term-n-plus-2}
\ \ \ \ \ \ \ \ \ \
\Pi_{n+2}^{\ind}
\Big(
G_{n+1,1,\dots,0}\,
\tfrac{y_1^{n+1}y_2}{(n+1)!}
\Big)
\,=\,
\tfrac{G_{n+1,1,\dots,0}}{(n+1)!}\,
A_{2,1}\,x_1^{n+2}
+
\tfrac{G_{n+1,1,\dots,0}}{n!}\,
A_{1,1}\,x_1^{n+1}x_2.
\end{align}

Next:
\[
\aligned
\Pi_{n+2}^{\ind}
\Big(
G_{n+1,0,1,\dots,0}\,
\tfrac{y_1^{n+1}y_3}{(n+1)!}
\Big)
&
\,=\,
\tfrac{G_{n+1,0,1,\dots,0}}{(n+1)!}\,
\Pi_{n+2}^{\ind}
\bigg(
\Big(
x_1
+
\varepsilon\,
\Big[
A_{1,1}\,x_1
+
\mathcal{R}^2
\Big]
+
{\rm O}(\varepsilon^2)
\Big)^{n+1}
\bigg)
\\
&
\ \ \ \ \ \ \ \ \ \ \ \ \ \ \ \ \ \ \ \ \ \ \ \ \ \ \ \ \ \ \ \ 
\cdot
\Big(
x_3
+
\varepsilon\,
\Big[
A_{3,1}\,x_1
-
A_{1,1}\,x_3
+
\mathcal{R}^2
\Big]
+
{\rm O}(\varepsilon^2)
\Big)^1
\bigg)
\\
&
\,=\,
\tfrac{G_{n+1,0,1,\dots,0}}{(n+1)!}\,
\Big\{
x_1^{n+1}\,A_{3,1}\,x_1
-
x_1^{n+1}\,A_{1,1}\,x_3
+
x_1^n\,
{\textstyle{\binom{n+1}{1}}}\,
A_{1,1}\,x_1x_3
\Big\},
\endaligned
\]
that is:
\leqnomode\usetagform{default}
\begin{align}
\label{9-term-n-plus-2}
\ \ \ \ \ \ \ \ \ \
\Pi_{n+2}^{\ind}
\Big(
G_{n+1,0,1,\dots,0}\,
\tfrac{y_1^{n+1}y_3}{(n+1)!}
\Big)
\,=\,
\tfrac{G_{n+1,0,1,\dots,0}}{(n+1)!}\,
\Big\{
A_{3,1}\,x_1^{n+2}
+
\tfrac{n}{(n+1)!}\,
A_{1,1}\,x_1^{n+1}\,x_3
\Big\}.
\end{align}

Then for general $k$ with $3 \leqslant k \leqslant n$, 
compute\,\,---\,\,in $0\cdots1\cdots0$, the $1$ is at position 
$k-1$\,\,---\,\,:
\[
\footnotesize
\aligned
\Pi_{n+2}^{\ind}
\Big(
G_{n+1,0\cdots1\cdots1}\,
\tfrac{y_1^{n+1}y_k}{(n+1)!}
\Big)
&
\,=\,
\tfrac{G_{n+1,0\cdots1\cdots0}}{(n+1)!}\,
\Pi_{n+2}^{\ind}
\bigg(
\Big(
x_1
+
\varepsilon\,
\Big[
A_{1,1}\,x_1
+
\mathcal{R}^2
\Big]
+
{\rm O}(\varepsilon^2)
\Big)^{n+1}
\\
&
\ \ \ \ \ \ \ \ \ \ \ \ \ \ \ \ \ \ \ \ \ \ \ \ \ \ \ \ \ \ \ \ \ \ \
\cdot
\Big(
x_k
+
\varepsilon\,
\Big[
A_{k,1}\,x_1
+
\smallsum{3\leqslant j\leqslant k}\,
A_{k,j}\,x_j
+
\mathcal{R}^2
\Big]
+
{\rm O}(\varepsilon^2)
\Big)^1
\bigg)
\\
&
\,=\,
\tfrac{G_{n+1,0\cdots1\cdots0}}{(n+1)!}\,
\bigg\{
x_1^{n+1}\,
\Big(
A_{k,1}\,x_1
+
\smallsum{3\leqslant j\leqslant k}\,
A_{k,j}\,x_j
\Big)
+
x_1^n\,
{\textstyle{\binom{n+1}{1}}}\,
A_{1,1}\,x_1x_k
\bigg\}
\\
&
\,=\,
\tfrac{G_{n+1,0\cdots1\cdots0}}{(n+1)!}\,
\bigg\{
A_{k,1}\,x_1^{n+2}
+
\smallsum{3\leqslant j\leqslant k-1}\,
A_{k,j}\,x_1^{n+1}x_j
+
\big(
A_{k,k}
+
(n+1)\,A_{1,1}
\big)\,
x_1^{n+1}x_k
\bigg\}.
\endaligned
\]
Now, replace the $A_{k,j}$ by their values
from~{\eqref{values-A-i-j}}, especially:
\[
A_{k,k}
\,=\,
-\,(k-2)\,A_{1,1},
\]
and obtain:
\begin{footnotesize}
\leqnomode\usetagform{default}
\begin{align}
\label{10-term-n-plus-2}
\Pi_{n+2}^{\ind}
\Big(
G_{n+1,0\cdots1\cdots0}\,
\tfrac{y_1^{n+1}y_k}{(n+1)!}
\Big)
\,=\,
\tfrac{1}{(n+1)!}\,
G_{n+1,0\cdots1\cdots0}\,
A_{k,1}\,x_1^{n+2}
&
-
\smallsum{3\leqslant j\leqslant k-1}\,
\tfrac{1}{(n+1)!}\,
\tfrac{j-2}{k-j+1}\,
\tbinom{k}{j}\,
G_{n+1,0\cdots1\cdots0}\,
A_{k-j+1,1}\,x_1^{n+1}x_j
\notag
\\
&
+
\tfrac{n-k+3}{(n+1)!}\,
G_{n+1,0\cdots1\cdots0}\,
A_{1,1}\,x_1^{n+1}x_k.
\end{align}
\end{footnotesize}

Finally, for any $i, j \geqslant 2$ with $i+j = n+2$:
\[
\aligned
\Pi_{n+2}^{\ind}
\Big(
y_1^n\,\tfrac{1}{2}\,
\tfrac{y_iy_j}{(i-1)!(j-1)!}
\Big)
&
\,=\,
\tfrac{1}{2}\,
\tfrac{1}{(i-1)!(j-1)!}\,
\Pi_{n+2}^{\ind}
\bigg(
\Big(
x_1
+
\varepsilon\,
\Big[
A_{1,1}\,x_1
+
\mathcal{R}^2
\Big]
+
{\rm O}(\varepsilon^2)
\Big)^n
\\
&
\ \ \ \ \ \ \ \ \ \ \ \ \ \ \ \ \ \ \ \ \ \ \ \ \ \ \ \ \ \ \ \ \ \ \
\cdot
\Big(
x_i
+
\varepsilon\,
\Big[
A_{i,1}\,x_1
+
\smallsum{3\leqslant k\leqslant i}\,
A_{i,k}\,x_k
+
\mathcal{R}^2
\Big]
+
{\rm O}(\varepsilon^2)
\Big)^1
\\
&
\ \ \ \ \ \ \ \ \ \ \ \ \ \ \ \ \ \ \ \ \ \ \ \ \ \ \ \ \ \ \ \ \ \ \
\cdot
\Big(
x_j
+
\varepsilon\,
\Big[
A_{j,1}\,x_1
+
\smallsum{3\leqslant l\leqslant i}\,
A_{j,l}\,x_l
+
\mathcal{R}^2
\Big]
+
{\rm O}(\varepsilon^2)
\Big)^1
\bigg),
\endaligned
\]
that is:
\[
\Pi_{n+2}^{\ind}
\Big(
y_1^n\,\tfrac{1}{2}\,
\tfrac{y_iy_j}{(i-1)!(j-1)!}
\Big)
\,=\,
\tfrac{1}{2}\,
\tfrac{1}{(i-1)!(j-1)!}\,
A_{i,1}\,x_1^{n+1}x_j
+
\tfrac{1}{2}\,
\tfrac{1}{(i-1)!(j-1)!}\,
A_{j,1}\,x_1^{n+1}x_i.
\]

In fact, we must perform a sum:
\[
\aligned
\Pi_{n+2}^{\ind}
\bigg(
y_1^n
\sum_{i,j\geqslant 2
\atop
i+j=n+2}\,
\tfrac{1}{2}\,
\frac{y_iy_j}{(i-1)!(j-1)!}
\bigg)
&
\,=\,
-\,\tfrac{1}{2}
\sum_{j=2}^n\,
\frac{1}{(n-j+1)!(j-1)!}\,
A_{n-j+2,1}\,
x_1^{n+1}x_j
\\
&
\ \ \ \ \
-\,\tfrac{1}{2}
\sum_{i=2}^n\,
\frac{1}{(i-1)!(n-i+1)!}\,
A_{n-i+2,1}\,
x_1^{n+1}x_i,
\endaligned
\]
and by exchanging $j \longleftrightarrow i$ in the 
first sum, the two sums are {\em equal}, whence:
\leqnomode\usetagform{default}
\begin{align}
\label{11-term-n-plus-2}
\ \ \ \ \ \ \ \ \ \
\Pi_{n+2}^{\ind}
\Big(
y_1^n
\sum_{i,j\geqslant 2
\atop
i+j=n+2}\,
\tfrac{1}{2}\,
\tfrac{y_iy_j}{(i-1)!(j-1)!}
\Big)
\,=\,
-\,
\sum_{i=2}^n\,
\tfrac{1}{(i-1)!(n-i+1)!}\,
A_{n-i+2,1}\,
x_1^{n+1}x_i.
\end{align}

Now that we have computed all terms
{\eqref{1-term-n-plus-2}},
{\eqref{2-term-n-plus-2}},
{\eqref{3-term-n-plus-2}},
{\eqref{4-term-n-plus-2}},
{\eqref{5-term-n-plus-2}},
{\eqref{6-term-n-plus-2}},
{\eqref{7-term-n-plus-2}},
{\eqref{8-term-n-plus-2}},
{\eqref{9-term-n-plus-2}},
{\eqref{10-term-n-plus-2}},
{\eqref{11-term-n-plus-2}},
we can sum them, namely:
\[
\aligned
\!\!\!\!\!\!\!\!\!\!\!\!\!\!\!
\Pi_{n+2}
&
\big(-v+G(y)\big)
\,=\,
\\
&
\,=\,
-\,\tfrac{2}{(n+2)!}\,
F_{n+2,0,\dots,0}\,
A_{1,1}\,x_1^{n+2}
-
\tfrac{2}{(n+1)!}\,
F_{n+1,1,\dots,0}\,A_{1,1}\,x_1^{n+1}x_2
-\cdots-
\tfrac{2}{(n+1)!}\,
F_{n+1,0,\dots,1}\,A_{1,1}\,x_1^{n+1}x_n
\\
&
\ \ \ \ \
-\,
\frac{1}{n!}\,A_{2,1}\,x_1^{n+1}x_n
\\
&
\ \ \ \ \ 
-\,\tfrac{1}{3}\,
\tfrac{1}{(n-1)!}\,
A_{3,1}\,x_1^{n+1}x_{n-1}
-\cdots-
\tfrac{2}{m+1}\,
\tfrac{1}{m!(n-m+1)!}\,
A_{m+1,1}\,x_1^{n+1}x_{n-m+1}
-\cdots-
\tfrac{2}{n}\,
\tfrac{1}{(n-1)!2!}\,
A_{n,1}\,
x_1^{n+1}x_2
\\
&
\ \ \ \ \ 
+
\tfrac{1}{n!2}\,B_n\,x_1^{n+2}
-
\tfrac{1}{(n-1)!2}\,
A_{2,1}\,x_1^{n+1}x_n
-
\sum_{i=2}^{n-1}\,
\tfrac{1}{n-i+2}\,
\tfrac{1}{(i-1)!(n-i)!}\,
A_{n+2-i,1}\,
x_1^{n+1}x_i
\\
&
\ \ \ \ \ 
+
\tfrac{1}{(n+1)!}\,
G_{n+2,0,\dots,0}\,A_{1,1}\,x_1^{n+2}
\\
&
\ \ \ \ \ 
+
\tfrac{1}{(n+1)!}\,
G_{n+1,1,\dots,0}\,
A_{2,1}\,x_1^{n+2}
+
\tfrac{1}{n!}\,
G_{n+1,1,0,\dots,0}\,
A_{1,1}\,x_1^nx_2
\\
&
\ \ \ \ \ 
+
\sum_{3\leqslant k\leqslant n}\,
\bigg\{
\tfrac{1}{(n+1)!}\,
G_{n+1,0\cdots1\cdots0}\,
A_{k,1}\,
x_1^{n+2}
-
\smallsum{3\leqslant j\leqslant k-1}\,
\tfrac{1}{(n+1)!}\,
\tfrac{j-2}{k-j+1}\,
\tbinom{k}{j}\,
G_{n+1,0\cdots1\cdots0}\,
A_{k-j+1,1}\,
x_1^{n+1}x_j
\\
&
\ \ \ \ \ \ \ \ \ \ \ \ \ \ \ \ \ \ \ \ \ \ \ \ \ \ \ \ \ \ \ \ \ \ \
\ \ \ \ \ \ \ \ \ \ \ \ \ \ \ \ \ \ \ \ \ \ \ \ \ \ \ \ \ \ \ \ \ \ \
\ \ \ \ \ \ \ \ \ \ \ \ \ \ \ \ \ \ \ \ \ \ \ \ \ \ \ \ \ \ \ \ \ \ 
+
\tfrac{n-k+3}{(n+1)!}\,
G_{n+1,0\cdots1\cdots0}\,
A_{1,1}\,
x_1^{n+1}k_x
\bigg\}
\\
&
\ \ \ \ \ 
+
\sum_{i=2}^n\,
\tfrac{1}{(i-1)!(n-i+1)!}\,
A_{n-i+2,1}\,
x_1^{n+1}x_i,
\endaligned
\]
and collect the coefficients of the independent monomials:
\[
\footnotesize
\aligned
0
&
\,\equiv\,
x_1^{n+2}
\Big\{
-\tfrac{2}{(n+2)!}\,F_{n+2,0,\dots,0}\,A_{1,1}
+
\tfrac{1}{(n+1)!}\,G_{n+2,0,\dots,0}\,A_{1,1}
+
\smallsum{2\leqslant k\leqslant n}\,
\tfrac{1}{(n+1)!}\,
G_{n+1,0\cdots1\cdots}\,
A_{k,1}
+
\frac{1}{n!2}\,
B_n
\Big\}
\\
&
\ \ \ \ \
+
x_1^{n+1}x_2
\Big\{
-\tfrac{2}{(n+1)!}\,F_{n+1,1,\dots,0}\,A_{1,1}
+
\tfrac{1}{n!}\,G_{n+1,1,\dots,0}\,A_{1,1}
\underbrace{
-
\tfrac{2}{n}\,
\tfrac{1}{(n-1)!2!}\,A_{n,1}
-
\tfrac{1}{n}\,\tfrac{1}{1!(n-2)!}\,
A_{n,1}
+
\tfrac{1}{1!}\,\tfrac{1}{(n-1)!}\,A_{n,1}}_{
=\,\,A_{n,1}\,\tfrac{1}{n!2!}\,\big[-2-(n-1)2+n\cdot 2\big]
\,\,=\,\,A_{n,1}\,\tfrac{1}{n!2!}\,\big[\zero{-2+2}\big]
\,\,=\,\,0}
\Big\}
\endaligned
\]
\[
\aligned
+
x_1^{n+1}x_3
\Big\{
&
-\tfrac{2}{(n+1)!}\,F_{n+1,01\cdots0}\,A_{1,1}
+
\tfrac{(n-3+3)}{(n+1)!}\,
G_{n+1,01\cdots0}\,A_{1,1}
+
\ast\,A_{2,1}+\cdots+\ast\,A_{n-2,1}
\\
&
\ \ \ \ \ \ \ \ \ \ \ \ \ \ \ \ \ \ \ \ \ \ \ \ \ \ \ \ \ \ \ \ \ \ \
+
\underbrace{
A_{n-1,1}
\big[
-\tfrac{2}{n-1}\,\tfrac{1}{(n-2)!3!}
-
\tfrac{1}{n-1}\,\tfrac{1}{2!(n-3)!}
+
\tfrac{1}{2!(n-2)!}
\big]}_{
=\,\,A_{n-1,1}\,\tfrac{1}{(n-1)!3!}\,
\big[-2-(n-2)3+(n-1)3\big]
\,\,=\,\,
A_{n-1,1}\,\tfrac{1}{(n-1)!3!}\,\big[-2+3\big]
}
\Big\}
\endaligned
\]
\[
\aligned
+
x_1^{n+1}x_4
\Big\{
&
-\tfrac{2}{(n+1)!}\,F_{n+1,001\cdots0}\,A_{1,1}
+
\tfrac{n-4+3}{(n+1)!}\,G_{n+1,001\cdots0}\,A_{1,1}
+
\ast\,A_{2,1}+\cdots+\ast\,A_{n-3,1}
\\
&
\ \ \ \ \ \ \ \ \ \ \ \ \ \ \ \ \ \ \ \ \ \ \ \ \ \ \ \ \ \ \ \ \ \ \
+
\underbrace{
A_{n-2,1}
\big[
-\tfrac{2}{n-2}\,\tfrac{1}{(n-3)!4!}
-
\tfrac{1}{n-2}\,\tfrac{1}{3!(n-4)!}
+
\tfrac{1}{3!(n-3)!}
\big]}_{
=\,\,
A_{n-2,1}\,\tfrac{1}{(n-2)!4!}\,
\big[
-2-(n-3)4+(n-2)4
\big]
\,\,=\,\,
A_{n-2,1}\,\tfrac{1}{(n-2)!4!}\,
\big[-2+4\big]
}
\Big\}
\endaligned
\]
\[
\aligned
+
x_1^{n+1}x_5
\Big\{
&
-\tfrac{2}{(n+1)!}\,F_{n+1,0001\cdots0}\,A_{1,1}
+
\tfrac{n-5+3}{(n+1)!}\,G_{n+1,0001\cdots0}\,A_{1,1}
+
\ast\,A_{2,1}+\cdots+\ast\,A_{n-4,1}
\\
&
\ \ \ \ \ \ \ \ \ \ \ \ \ \ \ \ \ \ \ \ \ \ \ \ \ \ \ \ \ \ \ \ \ \ \
+
\underbrace{
A_{n-3,1}
\big[
-\tfrac{2}{n-3}\,\tfrac{1}{(n-4)!5!}
-
\tfrac{1}{n-3}\,\tfrac{1}{4!(n-5)!}
+
\tfrac{1}{4!(n-4)!}
\big]}_{
=\,\,
A_{n-3,1}\,\tfrac{1}{(n-3)!5!}\,
\big[
-2-(n-4)5+(n-3)5
\big]
\,\,=\,\,
A_{n-3,1}\,\tfrac{1}{(n-3)!5!}\,
\big[-2+5\big]
}
\Big\}
\endaligned
\]
\[
+
\cdots\cdots\cdots\cdots\cdots\cdots\cdots\cdots\cdots\cdots\cdots
\cdots\cdots\cdots\cdots\cdots\cdots\cdots\cdots\cdots\cdots\cdots
\cdots\cdots\cdots\cdots\cdots
\]
\[
\aligned
+
x_1^{n+1}x_k
\Big\{
-\tfrac{2}{(n+1)!}\,F_{n+1,0\cdots1\cdots0}\,A_{1,1}
+
\tfrac{n-k+3}{(n+1)!}\,
G_{n+1,0\cdots1\cdots0}\,A_{1,1}
&
+
\ast\,A_{2,1}+\cdots+\ast\,A_{n-k+1,1}
\\
&
+
A_{n-k+2,1}\,\tfrac{1}{(n-k+2)!k!}\,
\big[
-2+k
\big]
\Big\}
\endaligned
\]
\[
+
\cdots\cdots\cdots\cdots\cdots\cdots\cdots\cdots\cdots\cdots\cdots
\cdots\cdots\cdots\cdots\cdots\cdots\cdots\cdots\cdots\cdots\cdots
\cdots\cdots\cdots\cdots\cdots
\]
\[
\aligned
+
x_1^{n+1}x_{n-1}
\Big\{
-\tfrac{2}{(n+1)!}\,
F_{n+1,0\cdots10}\,A_{1,1}
+
&
\tfrac{4}{(n+1)!}\,
G_{n+1,0\cdots10}\,
A_{1,1}
+
\ast\,A_{2,1}
\\
&
\ \ \ \ \ \ \ \ \ \ \ \ \ \ \ \ \ \ \ \ \ \ \ \ \ \ \ \ \ \ \ \ \ \ \
+
A_{3,1}\,
\tfrac{1}{3!(n-1)!}\,
\big[-2+n-1\big]
\Big\}
\endaligned
\]
\[
\aligned
+
x_1^{n+1}x_n
\Big\{
&
-\tfrac{2}{(n+1)!}\,F_{n+1,0\cdots01}\,A_{1,1}
+
\tfrac{3}{(n+1)!}\,G_{n+1,0\cdots01}\,A_{1,1}
\\
&
\ \ \ \ \ \ \ \ \ \ \ \ \ \ \ \ \ \ \ \ \ \ \ \ \ \ \ \ \ \ \ \ \ \ \
\ \ \ \ \ \ \ \ \ \ \ \ \ \ \ \ \ \ \ \ \ \ \ \ \ \ \ \ \ \ \ \ \ \ \
+
\underbrace{
A_{2,1}\,
\big[
-\tfrac{1}{n!}
-
\tfrac{1}{(n-1)!2!}
-
\emptyset
+
\tfrac{1}{(n-1)!1!}
\big]}_{
=\,\,A_{2,1}\,\tfrac{1}{2!n!}\,
\big[-2-n+2n\big]
\,\,=\,\,
A_{2,1}\,\tfrac{1}{2!n!}\,\big[n-2\big]
}
\Big\}.
\endaligned
\]
With this expression, we have therefore
computed the quantities $T_{n+2, 0,\dots,0}$,
$T_{n+1,1,\dots,0}$, $T_{n+1,0,1,\dots,0}$, \dots, $T_{n+1,0,\dots,1}$
introduced in~{\eqref{names-T}}.

\SectionHead{Normalizations of Order $(n+2)$ Monomials}
{normalizations-order-n-plus-2}

Now, coming back to~{\eqref{F-G-epsilon}},
thanks to these expressions of $T_{n+2, 0,\dots,0}$,
$T_{n+1,1,\dots,0}$, $T_{n+1,0,1,\dots,0}$, \dots, 
$T_{n+1,0,\dots,1}$, 
the coefficients of 
$x_1^{n+2}$, of $x_1^{n+1}x_2$, of $x_1^{n+1}x_3$, \dots,
of $x_1^{n+1}x_n$ must vanish, hence we obtain
for the first three:
\[
\footnotesize
\aligned
0
&
\,\equiv\,
-\,\tfrac{1}{(n+2)!}\,F_{n+2,0,\dots,0}
+
\tfrac{1}{(n+2)!}\,G_{n+2,0,\dots,0}
\\
&
\ \ \ \ \
+
\varepsilon\,
\Big\{
-\tfrac{2}{(n+2)!}\,
F_{n+2,0,\dots,0}\,A_{1,1}
+
\tfrac{1}{(n+1)!}\,
G_{n+2,0,\dots,0}\,A_{1,1}
+
\smallsum{2\leqslant k\leqslant n}\,
\tfrac{1}{(n+1)!}\,
G_{n+1,0\cdots1\cdots0}\,A_{k,1}
+
\tfrac{1}{n!\,2}\,
B_n
\Big\}
+
{\rm O}(\varepsilon^2),
\\
0
&
\,\equiv\,
-\,\tfrac{1}{(n+1)!}\,
F_{n+1,1,\dots,0}
+
\tfrac{1}{(n+1)!}\,
G_{n+1,1,\dots,0},
\endaligned
\]
\[
\footnotesize
\aligned
0
&
\,\equiv\,
-\,
\tfrac{1}{(n+1)!}\,
F_{n+1,0,1,\dots,0}
+
\tfrac{1}{(n+1)!}\,
G_{n+1,0,1,\dots,0}
\\
&
\ \ \ \ \
+
\varepsilon\,
\Big\{
-\tfrac{2}{(n+1)!}\,F_{n+1,0,1,\dots,0}\,A_{1,1}
+
\tfrac{n-3+3}{(n+1)!}\,
G_{n+1,0,1,\dots,0}\,A_{1,1}
+
\ast\,A_{2,1}+\cdots+\ast\,A_{n-2,1}
+
\tfrac{-2+3}{(n-1)!3!}\,A_{n-1,1}
\Big\}
+
{\rm O}(\varepsilon^2),
\endaligned
\]
for general $k$ with $3 \leqslant k \leqslant n-1$:
\[
\footnotesize
\aligned
0
&
\,\equiv\,
-\,
\tfrac{1}{(n+1)!}\,
F_{n+1,0\cdots1\cdots0}
+
\tfrac{1}{(n+1)!}\,
G_{n+1,0\cdots1\cdots0}
\\
&
\ \ \ \ \ 
+
\varepsilon\,
\Big\{
-\tfrac{2}{(n+1)!}\,F_{n+1,0\cdots1\cdots0}\,A_{1,1}
+
\tfrac{n-k+3}{(n+1)!}\,
G_{n+1,0\cdots1\cdots0}\,A_{1,1}
+\ast\,A_{2,1}+\cdots+\ast\,A_{n-k+1,1}
\\
&
\ \ \ \ \ \ \ \ \ \ \ \ \ \ \ \ \ \ \ \ \ \ \ \ \ \ \ \ \ \ \ \ \ \ \
\ \ \ \ \ \ \ \ \ \ \ \ \ \ \ \ \ \ \ \ \ \ \ \ \ \ \ \ \ \ \ \ \ \ \
\ \ \ \ \ \ \ \ \ \ \ \ \ \ \ \ \ \ \ \ \ \ \ \ \ \
+
\tfrac{-2+k}{(n-k+2)!k!}\,A_{n-k+2,1}
\Big\}
+
{\rm O}(\varepsilon^2),
\endaligned
\]
and for the last two:
\[
\footnotesize
\aligned
0
&
\,\equiv\,
-\,
\tfrac{1}{(n+1)!}\,
F_{n+1,0,\dots,1,0}
+
\tfrac{1}{(n+1)!}\,
G_{n+1,0,\dots,1,0}
\\
&
\ \ \ \ \ 
+
\varepsilon\,
\Big\{
-\tfrac{2}{(n+1)!}\,
F_{n+1,0,\dots,1,0}\,A_{1,1}
+
\tfrac{4}{(n+1)!}\,
G_{n+1,0,\dots,1,0}\,A_{1,1}
+
\ast\,A_{2,1}
+
\tfrac{-2+n-1}{3!\,(n-1)!}
\Big\}
+
{\rm O}(\varepsilon^2),
\\
0
&
\,\equiv\,
-\,
\tfrac{1}{(n+1)!}\,
F_{n+1,0,\dots,0,1}
+
\tfrac{1}{(n+1)!}\,
G_{n+1,0,\dots,0,1}
\\
&
\ \ \ \ \ 
+
\varepsilon\,
\Big\{
-\tfrac{2}{(n+1)!}\,
F_{n+1,0,\dots,0,1}\,A_{1,1}
+
\tfrac{3}{(n+1)!}\,
G_{n+1,0,\dots,0,1}\,A_{1,1}
+
\tfrac{-2+n}{2!\,n!}\,A_{2,1}
\Big\}
+
{\rm O}(\varepsilon^2).
\endaligned
\]

Now, we want to solve the $G_\smallbullet$ in terms of the
$F_\smallbullet$. 
For some $\nu = (\nu_1, \dots, \nu_n) \in \N^n$,
assume given an equation with constant $\alpha$, $\beta$,
$\Lambda_\nu$:
\[
0
\,\equiv\,
-\,F_\nu
+
G_\nu
+
\varepsilon\,
\big\{
\alpha\,F_\nu
+
\beta\,G_\nu
+
\Lambda_\nu
\big\}
+
{\rm O}(\varepsilon^2).
\]
To determine the term $S_\nu$ with:
\[
G_\nu
\overset{\text{\bf ?}}{\,=\,}
F_\nu
+
\varepsilon\,S_\nu
+
{\rm O}(\varepsilon^2),
\]
replace and identify:
\[
\aligned
0
&
\,\equiv\,
\zero{-\,F_\nu+F_\nu}
+
\varepsilon\,S_\nu
+
{\rm O}(\varepsilon^2)
+
\varepsilon\,
\Big\{
\alpha\,F_\nu
+
\beta\,
\big(
G_\nu+\varepsilon\,S_\nu+{\rm O}(\varepsilon^2)
\big)
+
\Lambda_\nu
\Big\}
+
{\rm O}(\varepsilon^2)
\\
&
\,\equiv\,
\varepsilon\,
\big\{
S_\nu
+
(\alpha+\beta)\,G_\nu
+
\Lambda_\nu
\big\}
+
{\rm O}(\varepsilon^2),
\endaligned
\]
whence:
\[
S_\nu
\,:=\,
-\,(\alpha+\beta)\,G_\nu
-
\Lambda_\nu.
\]
Same kinds of formulas exist for a linear system
involving several $F_\nu$, $G_\nu$, as above.

We therefore obtain:
\[
\footnotesize
\aligned
G_{n+2,0,\dots,0}
&
\,=\,
F_{n+2,0,\dots,0}
-
\varepsilon\,
\Big\{
n\,F_{n+2,0,\dots,0}\,A_{1,1}
+
\smallsum{2\leqslant k\leqslant n}\,
(n+2)\,F_{n+1,0\cdots1\cdots0}\,A_{k,1}
+
\tbinom{n+2}{2}\,B_n
\Big\}
+
{\rm O}(\varepsilon^2),
\\
G_{n+1,1,\dots,0}
&
\,=\,
F_{n+1,1,\dots,0}
-
\varepsilon\,
\Big\{
(n-1)\,F_{n+1,1,\dots,0}\,A_{1,1}
\Big\}
+
{\rm O}(\varepsilon^2)
\\
G_{n+1,0,1,,\dots,0}
&
\,=\,
F_{n+1,0,1,\dots,0}
-
\varepsilon\,
\Big\{
(n-2)\,F_{n+1,0,1,\dots,0}\,A_{1,1}
+
\ast\,A_{2,1}+\cdots+\ast\,A_{n-2,1}
+
\tfrac{-2+3}{(n-1)!3!}\,
(n+1)!\,
A_{n-1,1}
\Big\}
+
{\rm O}(\varepsilon^2),
\endaligned
\]
for general $k$ with $3 \leqslant k \leqslant n-1$:
\[
\aligned
G_{n+1,0\cdots1\cdots0}
\,=\,
F_{n+1,0\cdots1\cdots0}
-
\varepsilon\,
\Big\{
(n-k+1)\,F_{n+1,0\cdots1\cdots0}\,A_{1,1}
&
+
\ast\,A_{2,1}+\cdots+\ast\,A_{n-k+1,1}
\\
&
+
\tfrac{-2+k}{(n-k+2)!k!}\,(n+1)!\,
A_{n-k+2,1}
\Big\}
+
{\rm O}(\varepsilon^2),
\endaligned
\]
and:
\[
\aligned
G_{n+1,0,\dots,1,0}
&
\,=\,
F_{n+1,0,\dots,1,0}
-
\varepsilon\,
\Big\{
2\,F_{n+1,0,\dots,1,0}\,A_{1,1}
+
\ast\,A_{2,1}
+
\tfrac{-2+n-1}{3!(n-1)!}\,
(n+1)!\,
A_{3,1}
\Big\}
+
{\rm O}(\varepsilon^2),
\\
G_{n+1,0,\dots,0,1}
&
\,=\,
F_{n+1,0,\dots,0,1}
-
\varepsilon\,
\Big\{
1\,F_{n+1,0,\dots,0,1}\,A_{1,1}
+
\tfrac{-2+n}{2!n!}\,
(n+1)!\,
A_{2,1}
\Big\}
+
{\rm O}(\varepsilon^2).
\endaligned
\]

Since it is more natural, 
we start from the last equation.

\begin{Proposition}
\label{Prp-normalizations-n-plus-2}
One can normalize on the right, and then on the left also:
\[
\aligned
G_{n+1,0,\dots,0,1}
&
\,:=\,0,
&
\ \ \ \ \ \ \ \ \ \ \ \ \ \ \ \ \ \ \ \ \ \ \ \ \ \
F_{n+1,0,\dots,0,1}
&
\,:=\,0,
\\
G_{n+1,0,\dots,1,0}
&
\,:=\,0,
&
\ \ \ \ \ \ \ \ \ \ \ \ \ \ \ \ \ \ \ \ \ \ \ \ \ \
F_{n+1,0,\dots,0,1}
&
\,:=\,0,
\\
\cdots\cdots\cdots
&
\cdots\cdots
&
\ \ \ \ \ \ \ \ \ \ \ \ \ \ \ \ \ \ \ \ \ \ \ \ \ \
\cdots\cdots\cdots
&
\cdots\cdots
\\
G_{n+1,0,1,\dots,0}
&
\,:=\,0,
&
\ \ \ \ \ \ \ \ \ \ \ \ \ \ \ \ \ \ \ \ \ \ \ \ \ \
F_{n+1,0,1,\dots,0}
&
\,:=\,0,
\\
G_{n+2,0,0,\dots,0}
&
\,:=\,0,
&
\ \ \ \ \ \ \ \ \ \ \ \ \ \ \ \ \ \ \ \ \ \ \ \ \ \
F_{n+2,0,0,\dots,0}
&
\,:=\,0,
\endaligned
\]
while the coefficient of $\frac{1}{(n+1)!}\,
x_1^{n+1}x_2$ is a relative invariant:
\[
G_{n+1,10\cdots0}
\,\,\propto\,\,
F_{n+1,10\cdots0}.
\]
\end{Proposition}

\proof
Indeed, the variable $A_{2,1}$ is free to make the first
normalization, then the variable $A_{3,1}$ is free to
make the second one, and so on
until $B_n$ is free make the last normalization.
Notice that the system is triangular.
Two other views of these normalizations 
will be given in the next two 
Section~{\ref{normalizations-n-2-via-jet-prolongations}}
and~{\ref{alternative-normalizations-order-n-plus-2}}.
\endproof

\SectionHead{Alternative, More Direct Normalizations at Order $(n+2)$}
{alternative-normalizations-order-n-plus-2}

As we know from 
Section~{\ref{tangency-order-n-1}}, 
the matrix of coefficients of a vector field
$L$ tangent up to order $n+1$ is:
\[
\!\!\!\!\!\!\!\!\!\!\!\!\!\!\!
\scriptsize
\aligned
\left[
\def\arraystretch{1.25}
\begin{array}{ccccccccc}
A_{1,1} & \red{\bf 0} & \red{\bf 0} & \red{\bf 0} &
\red{\cdots} & \red{\bf 0} & \red{\bf 0} & \red{\bf 0} & 
\red{-\tfrac{2}{2}A_{2,1}}
\\
A_{2,1} & \red{\bf 0} & \red{\bf 0} & \red{\bf 0} &
\red{\cdots} & \red{\bf 0} & \red{\bf 0} & \red{\bf 0} & 
\red{-\tfrac{2}{3}A_{3,1}}
\\
A_{3,1} & \red{\bf 0} & \red{-A_{1,1}} & \red{\bf 0} &
\red{\cdots} & \red{\bf 0} & \red{\bf 0} & \red{\bf 0} & 
\red{-\tfrac{2}{4}A_{4,1}}
\\
A_{4,1} & \red{\bf 0} & \red{-2A_{2,1}} & \red{-2A_{1,1}} &
\red{\cdots} & \red{\bf 0} & \red{\bf 0} & \red{\bf 0} & 
\red{-\tfrac{2}{5}A_{4,1}}
\\
\vdots & \red{\vdots} & \red{\vdots} & \red{\vdots} & 
\red{\ddots} & \red{\vdots} & \vdots & \red{\vdots} & 
\red{\vdots}
\\
A_{n-2,1} & \red{\bf 0} &
\red{\frac{-1}{n-4}\binom{n-2}{3}A_{n-4,1}} &
\red{\frac{-2}{n-5}\binom{n-2}{4}A_{n-5,1}} & \red{\cdots} &
\red{\frac{-(n-4)}{1}\binom{n-2}{n-2}A_{1,1}} & 
\red{\bf 0} & \red{\bf 0} & 
\red{-\frac{2}{n-1}A_{n-1,1}}
\\
A_{n-1,1} & \red{\bf 0} &
\red{\frac{-1}{n-3}\binom{n-1}{3}A_{n-3,1}} &
\red{\frac{-2}{n-4}\binom{n-1}{4}A_{n-4,1}} & \red{\cdots} &
\red{\frac{-(n-4)}{2}\binom{n-1}{n-2}A_{2,1}} & 
\red{\frac{-(n-3)}{1}\binom{n-1}{n-1}A_{1,1}} & 
\red{\bf 0} & 
\red{-\frac{2}{n}A_{n,1}}
\\
A_{n,1} & \red{\bf 0} & 
\red{\frac{-1}{n-2}\binom{n}{3}A_{n-2,1}} &
\red{\frac{-2}{n-3}\binom{n}{4}A_{n-3,1}} & \cdots & 
\red{\frac{-(n-4)}{3}\binom{n}{n-2}A_{3,1}} & 
\red{\frac{-(n-3)}{2}\binom{n}{n-1}A_{2,1}} & 
\red{\frac{-(n-2)}{1}\binom{n}{n}A_{1,1}} & 
B_n
\\
\red{\bf 0} & \red{\bf 0} & \red{\bf 0} & \red{\bf 0} &
\red{\cdots} & \red{\bf 0} & \red{\bf 0} & \red{\bf 0} & 
\red{2\,A_{1,1}}
\end{array}
\right],
\endaligned
\]
So with these coefficients, 
$0 \equiv \pi_{\ind}^{n+1} 
\big( L \big( - u + F(x) \big) \big\vert_{u=F(x)} \big)$.

Then we apply the derivation $L (\centersmallbullet)$ to
the hypersurface equation~{\ref{u-F-order-n-2}} 
written up to order $n+2$, 
and we want to compute $\pi_{n+2}^{\ind} \big(
L \big( - u + F(x) \big) \big\vert_{u=F(x)} \big)$,
namely we want to compute
$\pi_{n+2}^{\ind} (\centersmallbullet)$ of:
\[
\!\!\!\!\!\!\!\!\!\!\!\!\!\!\!\!\!\!\!\!\!\!\!\!\!
\footnotesize
\aligned
0
&
\,\equiv\,
-\,2\,A_{1,1}\,
\Big(
\tfrac{x_1^2}{2!}
+
\tfrac{x_1^2x_2}{2!}
+\cdots+
\tfrac{x_1^nx_n}{n!}
+
F_{n+2,0,\dots,0}\tfrac{x_1^{n+1}x_1}{(n+2)!}
+
F_{n+1,1,\dots,0}\tfrac{x_1^{n+1}x_2}{(n+1)!}
+\cdots+
F_{n+1,0,\dots,1}\tfrac{x_1^{n+1}x_n}{(n+1)!}
\Big)
\\
&
\ \ \ \ \
+
\Big(
A_{1,1}x_1
-
\tfrac{2}{2}A_{2,1}
\Big[
\tfrac{x_1^2}{2!}
\Big\vert
+
\tfrac{x_1^2x_2}{2!}
+\cdots+
\tfrac{x_1^nx_n}{n!}
\Big]
\Big)\,
\Big(
\tfrac{x_1^1}{1!}
\Big\vert
+
\tfrac{x_1^2x_2}{1!}
+\cdots+
\tfrac{x_1^{n-1}x_n}{(n-1)!}
\\
&
\ \ \ \ \ \ \ \ \ \ \ \ \ \ \ \ \ \ \ \ \ \ \ \ \ \ \ \ \ \ \ \ \ \ \
\ \ \ \ \ \ \ \ \ \ \ \ \ \ \ \ \ \ \ \ \ \ \ \ \ \ \ \ \ \ \ \ \ \ \
\ \ \ \ \ \ \ \ \ \ \ \ \ \ \ \ \ \ \ \ \ \ \ \ \ 
+
F_{n+2,0,\dots,0}\tfrac{x_1^{n+1}}{(n+1)!}
+
F_{n+1,1,\dots,0}\tfrac{x_1^nx_2}{n!}
+\cdots+
F_{n+1,0,\dots,1}\tfrac{x_1^nx_n}{n!}
\Big)
\\
&
\ \ \ \ \ 
+
\Big(
A_{2,1}x_1
-
\tfrac{2}{3}A_{3,1}
\Big[
\tfrac{x_1^2}{2!}
\Big\vert
+
\tfrac{x_1^2x_2}{2!}
+\cdots+
\tfrac{x_1^{n-1}x_{n-1}}{(n-1)!}
\Big]
\Big)
\Big(
\tfrac{x_1^2}{2!}
\Big\vert
+
\tfrac{x_1^2x_2}{1!\,1!}
+\cdots+
\tfrac{x_1^{n-1}x_{n-1}}{1!\,(n-2)!}
+
F_{n+1,1,\dots,0}\tfrac{x_1^{n+1}}{(n+1)!}
+
\tfrac{x_1^nx_n}{1!(n-1)!}
\Big)
\\
&
\ \ \ \ \ 
+
\Big(
A_{3,1}x_1
-
\ast A_{1,1}x_3
-
\tfrac{2}{4}A_{4,1}
\Big[
\tfrac{x_1^2}{2!}
\Big\vert
+
\tfrac{x_1^2x_2}{2!}
+\cdots+
\tfrac{x_1^{n-2}x_{n-2}}{(n-2)!}
\Big]
\Big)
\Big(
\tfrac{x_1^3}{3!}
\Big\vert
+
\tfrac{x_1^3x_2}{2!\,1!}
+\cdots+
\tfrac{x_1^{n-1}x_{n-2}}{2!\,(n-3)!}
+
F_{n+1,0,1,\dots,0}\tfrac{x_1^{n+1}}{(n+1)!}
+
\tfrac{x_1^nx_{n-1}}{2!(n-2)!}
\Big)
\\
&
\ \ \ \ \ 
+
\cdots\cdots\cdots\cdots\cdots\cdots\cdots\cdots\cdots\cdots\cdots
\cdots\cdots\cdots\cdots\cdots\cdots\cdots\cdots\cdots\cdots\cdots
\cdots\cdots\cdots\cdots\cdots\cdots\cdots\cdots\cdots
\\
&
\ \ \ \ \ 
+
\Big(
A_{i,1}x_1
-\ast A_{i-2,1}x_3-\cdots-\ast A_{2,1}x_{i-1}-(i-2)A_{1,1}x_i
-
\tfrac{2}{i+1}A_{i+1,1}
\Big[
\tfrac{x_1^2}{2!}
\Big\vert
+
\tfrac{x_1^2x_2}{2!}
+\cdots+
\tfrac{x_1^{n-i+1}x_{n-i+1}}{(n-i+1)!}
\Big]
\Big)
\\
&
\ \ \ \ \ \ \ \ \ \ \ \ \ \ \ \ \ \ \ \ \ \ \ \ \ \ \ \ \ \ \ \ \ \ \
\ \ \ \ \ \ \ \ \ \ \ \ \ \ \ \ \ \ \ \ \ \ \ \ \ \ \ \ \ \ \ \ \ \ \
\ \ \ \ \ \ \ \ \ \ \ \ \ \ \ \ \
\Big(
\tfrac{x_1^i}{i!}
\Big\vert
+
\tfrac{x_1^ix_2}{(i-1)!1!}
+\cdots+
\tfrac{x_1^{n-1}x_{n-i+1}}{(i-1)!(n-i)!}
+
F_{n+1,0\cdots1\cdots0}\tfrac{x_1^{n+1}}{(n+1)!}
+
\tfrac{x_1^nx_{n-i+2}}{(i-1)!(n-i+1)!}
\Big)
\\
&
\ \ \ \ \ 
+
\cdots\cdots\cdots\cdots\cdots\cdots\cdots\cdots\cdots\cdots\cdots
\cdots\cdots\cdots\cdots\cdots\cdots\cdots\cdots\cdots\cdots\cdots
\cdots\cdots\cdots\cdots\cdots\cdots\cdots\cdots\cdots
\\
&
\ \ \ \ \ 
+
\Big(
A_{n-3,1}x_1
-
\ast A_{n-5,1}x_3-\cdots-\ast A_{2,1}x_{n-4}-(n-5)A_{1,1}x_{n-3}
-
\tfrac{2}{n-2}A_{n-2,1}
\Big[
\tfrac{x_1^2}{2!}
\Big\vert
+
\tfrac{x_1^2x_2}{2!}
+
\tfrac{x_1^3x_3}{3!}
+
\tfrac{x_1^4x_4}{4!}
\Big]
\Big)
\\
&
\ \ \ \ \ \ \ \ \ \ \ \ \ \ \ \ \ \ \ \ \ \ \ \ \ \ \ \ \ \ \ \ \ \ \
\ \ \ \ \ \ \ \ \ \ \ \ \ \ \ \ \ \ \ \ \ \ \ \ \ \ \ \ \ \ \ \ \ \ \
\ \ \ \ \ \ \ \ \ \ \ \ \ \ \ \ \
\Big(
\tfrac{x_1^{n-3}}{(n-3)!}
\Big\vert
+
\tfrac{x_1^{n-3}x_2}{(n-4)!1!}
+
\tfrac{x_1^{n-2}x_3}{(n-4)!2!}
+
\tfrac{x_1^{n-1}x_4}{(n-4)!3!}
+
F_{n+1,0,\dots,1,0,0,0}\tfrac{x_1^{n+1}}{(n+1)!}
+
\tfrac{x_1^nx_5}{(n-4)!4!}
\Big)
\\
&
\ \ \ \ \ 
+
\Big(
A_{n-2,1}x_1
-
\ast A_{n-4,1}x_3-\cdots-\ast A_{2,1}x_{n-3}-(n-4)A_{1,1}x_{n-2}
-
\tfrac{2}{n-1}A_{n-1,1}
\Big[
\tfrac{x_1^2}{2!}
\Big\vert
+
\tfrac{x_1^2x_2}{2!}
+
\tfrac{x_1^3x_3}{3!}
\Big]
\Big)
\\
&
\ \ \ \ \ \ \ \ \ \ \ \ \ \ \ \ \ \ \ \ \ \ \ \ \ \ \ \ \ \ \ \ \ \ \
\ \ \ \ \ \ \ \ \ \ \ \ \ \ \ \ \ \ \ \ \ \ \ \ \ \ \ \ \ \ \ \ \ \ \
\ \ \ \ \ \ \ \ \ \ \ \ \ \ \ \ \
\Big(
\tfrac{x_1^{n-2}}{(n-2)!}
\Big\vert
+
\tfrac{x_1^{n-2}x_2}{(n-3)!1!}
+
\tfrac{x_1^{n-1}x_3}{(n-3)!2!}
+
F_{n+1,0,\dots,1,0,0}\tfrac{x_1^{n+1}}{(n+1)!}
+
\tfrac{x_1^nx_4}{(n-3)!3!}
\Big)
\\
&
\ \ \ \ \ 
+
\Big(
A_{n-1,1}x_1
-
\ast A_{n-3,1}x_3-\cdots-\ast A_{2,1}x_{n-2}-(n-3)A_{1,1}x_{n-1}
-
\tfrac{2}{n}A_{n,1}
\Big[
\tfrac{x_1^2}{2!}
\Big\vert
+
\tfrac{x_1^2x_2}{2!}
\Big]
\Big)
\\
&
\ \ \ \ \ \ \ \ \ \ \ \ \ \ \ \ \ \ \ \ \ \ \ \ \ \ \ \ \ \ \ \ \ \ \
\ \ \ \ \ \ \ \ \ \ \ \ \ \ \ \ \ \ \ \ \ \ \ \ \ \ \ \ \ \ \ \ \ \ \
\ \ \ \ \ \ \ \ \ \ \ \ \ \ \ \ \
\Big(
\tfrac{x_1^{n-1}}{(n-1)!}
\Big\vert
+
\tfrac{x_1^{n-1}x_2}{(n-2)!1!}
+
F_{n+1,0,\dots,1,0}\tfrac{x_1^{n+1}}{(n+1)!}
+
\tfrac{x_1^nx_3}{(n-2)!2!}
\Big)
\\
&
\ \ \ \ \ 
+
\Big(
A_{n,1}x_1
-
\ast A_{n-2,1}x_3-\cdots-\ast A_{2,1}x_{n-1}-(n-2)A_{1,1}x_n
+
B_n
\Big[
\tfrac{x_1^2}{2!}
\Big]
\Big)
\Big(
\tfrac{x_1^n}{n!}
\Big\vert
+
\tfrac{x_1^nx_2}{(n-1)!1!}
+
F_{n+1,0,\dots,0,1}\tfrac{x_1^{n+1}}{(n+1)!}
+
\tfrac{x_1^nx_2}{(n-1)!2!}
\Big).
\endaligned
\]

By careful inspection of the products and after relevant
simplifications, we find exactly the same 
factors (up to sign) of $\varepsilon^1$ 
as in the equations preceding
Proposition~{\ref{Prp-normalizations-n-plus-2}}:
\[
\aligned
0
&
\,\equiv\,
\tfrac{x_1^{n+2}}{(n+2)!}
\Big\{
n\,F_{n+2,0,\dots,0}\,A_{1,1}
+
\smallsum{2\leqslant k\leqslant n}\,
(n+2)\,F_{n+1,0\cdots1\cdots0}\,A_{k,1}
+
\tbinom{n+2}{2}\,B_n
\Big\}
\\
&
\ \ \ \ \ 
+
\tfrac{x_1^{n+1}x_2}{(n+1)!}
\Big\{
(n-1)\,F_{n+1,1,\dots,0}\,A_{1,1}
\Big\}
\\
&
\ \ \ \ \ 
+
\tfrac{x_1^{n+1}x_3}{(n+1)!}
\Big\{
(n-2)\,F_{n+1,0,1,\dots,0}\,A_{1,1}
+
\ast\,A_{2,1}+\cdots+\ast\,A_{n-2,1}
+
\tfrac{-2+3}{(n-1)!3!}\,
(n+1)!\,
A_{n-1,1}
\Big\}
\\
&
\ \ \ \ \
+
\cdots\cdots\cdots\cdots\cdots\cdots\cdots\cdots\cdots
\cdots\cdots\cdots\cdots\cdots\cdots\cdots\cdots\cdots
\cdots\cdots\cdots\cdots\cdots\cdots\cdots\cdots\cdots
\cdot
\\
&
\ \ \ \ \
+
\tfrac{x_1^{n+1}x_k}{(n+1)!}
\Big\{
(n-k+1)\,F_{n+1,0\cdots1\cdots0}\,A_{1,1}
+
\ast\,A_{2,1}+\cdots+\ast\,A_{n-k+1,1}
+
\tfrac{-2+k}{(n-k+2)!k!}\,(n+1)!\,
A_{n-k+2,1}
\Big\}
\\
&
\ \ \ \ \
+
\cdots\cdots\cdots\cdots\cdots\cdots\cdots\cdots\cdots
\cdots\cdots\cdots\cdots\cdots\cdots\cdots\cdots\cdots
\cdots\cdots\cdots\cdots\cdots\cdots\cdots
\\
&
\ \ \ \ \ 
+
\tfrac{x_1^{n+1}x_{n-1}}{(n+1)!}
\Big\{
2\,F_{n+1,0,\dots,1,0}\,A_{1,1}
+
\ast\,A_{2,1}
+
\tfrac{-2+n-1}{3!(n-1)!}\,
(n+1)!\,
A_{3,1}
\Big\}
\\
&
\ \ \ \ \ 
+
\tfrac{x_1^{n+1}x_{n-1}}{(n+1)!}
\Big\{
1\,F_{n+1,0,\dots,0,1}\,A_{1,1}
+
\tfrac{-2+n}{2!n!}\,
(n+1)!\,
A_{2,1}
\Big\}.
\endaligned
\]
In conclusion, this computation is a {\sl shortcut}
to the longer computation done previously with infinitesimal
$\varepsilon$.

\begin{Examples}
In dimension $n = 2$:
\[
\aligned
\tfrac{x_1^4}{4!}\,
\Big\{
&
2\,F_{4,0}\,A_{1,1}
+
4\,F_{3,1}\,A_{2,1}
+
6\,B_2
\Big\},
\\
+
\tfrac{x_1^3x_2}{3!}\,
\Big\{
&
F_{3,1}\,A_{1,1}
\Big\}.
\endaligned
\]

\noindent
In dimension $n = 3$:
\[
\aligned
\tfrac{x_1^5}{5!}\,
\Big\{
&
3\,F_{5,0,0}\,A_{1,1}
+
5\,F_{4,1,0}\,A_{2,1}
+
5\,F_{4,0,1}\,A_{3,1}
+
10\,B_3
\Big\},
\\
+
\tfrac{x_1^4x_2}{4!}\,
\Big\{
&
2\,F_{4,1,0}\,A_{1,1}
\Big\},
\\
+
\tfrac{x_1^4x_3}{4!}\,
\Big\{
&
F_{4,0,1}\,A_{1,1}
+
2\,A_{2,1}
\Big\}.
\endaligned
\]

\noindent
In dimension $n = 4$:
\[
\aligned
\tfrac{x_1^6}{6!}\,
\Big\{
&
4\,F_{6,0,0,0}\,A_{1,1}
+
6\,F_{5,1,0,0}\,A_{2,1}
+
6\,F_{5,0,1,0}\,A_{3,1}
+
6\,F_{5,0,0,1}\,A_{4,1}
+
15\,B_4
\Big\},
\\
+
\tfrac{x_1^5x_2}{5!}\,
\Big\{
&
3\,F_{5,1,0,0}\,A_{1,1}
\Big\},
\\
+
\tfrac{x_1^5x_3}{5!}\,
\Big\{
&
2\,F_{5,0,1,0}\,A_{1,1}
-
2\,F_{5,0,0,1}\,A_{2,1}
+
\tfrac{10}{3}\,A_{3,1}
\Big\},
\\
+
\tfrac{x_1^5x_4}{5!}\,
\Big\{
&
1\,F_{5,0,0,1}\,A_{1,1}
+
5\,A_{2,1}
\Big\}.
\endaligned
\]

\noindent
In dimension $n = 5$:
\[
\aligned
\tfrac{x_1^7}{7!}\,
\Big\{
&
5\,F_{7,0,0,0,0}\,A_{1,1}
+
7\,F_{6,1,0,0,0}\,A_{2,1}
+
7\,F_{6,0,1,0,0}\,A_{3,1}
+
7\,F_{6,0,0,1,0}\,A_{4,1}
+
7\,F_{6,0,0,0,1}\,A_{5,1}
+
21\,B_5
\Big\},
\\
+
\tfrac{x_1^6x_2}{6!}\,
\Big\{
&
4\,F_{6,1,0,0,0}\,A_{1,1}
\Big\},
\\
+
\tfrac{x_1^6x_3}{6!}\,
\Big\{
&
3\,F_{6,0,1,0,0}\,A_{1,1}
-
2\,F_{6,0,0,1,0}\,A_{2,1}
-
\tfrac{10}{3}\,F_{6,0,0,0,1}\,A_{3,1}
+
5\,A_{4,1}
\Big\},
\\
+
\tfrac{x_1^6x_4}{6!}\,
\Big\{
&
2\,F_{6,0,0,1,0}\,A_{1,1}
-
5\,F_{6,0,0,0,1}\,A_{2,1}
+
10\,A_{3,1}
\Big\},
\\
+
\tfrac{x_1^6x_4}{6!}\,
\Big\{
&
1\,F_{6,0,0,0,1}\,A_{1,1}
+
9\,A_{2,1}
\Big\}.
\endaligned
\]

\noindent
In dimension $n = 6$:
\[
\aligned
\tfrac{x_1^8}{8!}\,
\Big\{
&
6\,F_{8,0,0,0,0,0}\,A_{1,1}
+
8\,F_{7,1,0,0,0,0}\,A_{2,1}
+
8\,F_{7,0,1,0,0,0}\,A_{3,1}
+
8\,F_{7,0,0,1,0,0}\,A_{4,1}
\\
&
\ \ \ \ \ \ \ \ \ \ \ \ \ \ \ \ \ \ \ \ \ \ \ \ \ \ \ \ \ \ \ \ \ \ \
\ \ \ \ \ \ \ \ \ \ \ \ \ \ \ \ \ \ \ \ \ \ \ \ \ \ \ \ \ \ \ \ \ \ \
+
8\,F_{7,0,0,0,1,0}\,A_{5,1}
+
8\,F_{7,0,0,0,0,1}\,A_{6,1}
+
28\,B_6
\Big\},
\\
+
\tfrac{x_1^7x_2}{7!}\,
\Big\{
&
5\,F_{7,1,0,0,0,0}\,A_{1,1}
\Big\},
\\
+
\tfrac{x_1^7x_3}{7!}\,
\Big\{
&
4\,F_{7,0,1,0,0,0}\,A_{1,1}
-
2\,F_{7,0,0,1,0,0}\,A_{2,1}
-
\tfrac{10}{3}\,F_{7,0,0,0,1,0}\,A_{3,1}
-
5\,F_{7,0,0,0,0,1}\,A_{4,1}
+
7\,A_{5,1}
\Big\},
\\
+
\tfrac{x_1^7x_4}{7!}\,
\Big\{
&
3\,F_{7,0,0,1,0,0}\,A_{1,1}
-
5\,F_{7,0,0,0,1,0}\,A_{2,1}
-
10\,F_{7,0,0,0,0,1}\,A_{3,1}
+
\tfrac{35}{2}\,A_{4,1}
\Big\},
\\
+
\tfrac{x_1^7x_5}{7!}\,
\Big\{
&
2\,F_{7,0,0,0,1,0}\,A_{1,1}
-
9\,F_{7,0,0,0,0,1}\,A_{2,1}
+
21\,A_{3,1}
\Big\},
\\
+
\tfrac{x_1^7x_6}{7!}\,
\Big\{
&
1\,F_{7,0,0,0,0,1}\,A_{1,1}
+
14\,A_{2,1}
\Big\}.
\endaligned
\]
\end{Examples}

\SectionHead{Normalizations at Order $(n+2)$ Via Jet Prolongations}
{normalizations-n-2-via-jet-prolongations}

As before, consider $u = u(x_1, \dots, x_n)$ as
a function of $(x_1, \dots, x_n)$. The letter $F$ will not be used.
Then to each partial
derivative $u_{x_1^{\nu_1} \cdots x_n^{\nu_n}} (x_1, \dots, x_n)$,
one can associate an independent coordinate (variable),
denoted similarly $u_{x_1^{\nu_1} \cdots x_n^{\nu_n}}$,
or sometimes more simply $u_{\nu_1,\dots,\nu_n}$.
Background appears {\em e.g.} 
in~{\cite{Chen-Merker-2019}}. 

For each integer $\kappa \geqslant 0$,
introduce the {\sl jet space of order $\kappa$},
namely $\R^{n + \binom{n+\kappa}{\kappa}}$ equipped with the
coordinates:
\[
\Big(
x_1,\dots,x_n,u,\,
\big(
u_{x_1^{\nu_1}\cdots x_n^{\nu_n}}
\big)_{1\leqslant\nu_1+\cdots+\nu_n\leqslant\kappa}
\Big).
\]
An abbreviation will be used:
\[
u^{(\kappa)}
\,:=\,
\big(
u_{x_1^{\nu_1}\cdots x_n^{\nu_n}}
\big)_{1\leqslant\nu_1+\cdots+\nu_n\leqslant\kappa}.
\]

Introduce also, for $i = 1, \dots, n$, 
the {\sl total differentiation operators:}
\[
{\sf D}_{x_i}
\,:=\,
\frac{\partial}{\partial x_i}
+
u_{x_i}\,\frac{\partial}{\partial u}
+
\sum_{m=1}^\infty\,
\sum_{\nu_1+\cdots+\nu_n=m}\,
u_{x_ix_1^{\nu_1}\cdots x_n^{\nu_n}}\,
\frac{\partial}{
\partial 
u_{x_1^{\nu_1}\cdots x_n^{\nu_n}}},
\]
which commute one with another.

Now, given a general vector field in the $(x,u)$-space:
\[
L
\,=\,
\sum_{i=1}^n\,
X_i(x,u)\,\frac{\partial}{\partial x_i}
+
U(x,u)\,\frac{\partial}{\partial u},
\]
its prolongation~{\cite{Chen-Merker-2019}} to the infinite jet space:
\[
L^{(\infty)}
\,=\,
L
+
\sum_{\nu_1+\cdots+\nu_n\geqslant 1}\,
U_{\nu_1,\dots,\nu_n}\,
\frac{\partial}{\partial u_{x_1^{\nu_1}\cdots x_n^{\nu_n}}},
\]
expresses how the (differentiated)
flow of $L$ acts on higher order jets,
and its coefficients $U_{\nu_1, \dots, \nu_n}$
are uniquely determined by the formulas:
\[
U_{\nu_1,\dots,\nu_n}
\,:=\,
{\sf D}_{x_1^{\nu_1}}\cdots{\sf D}_{x_n^{\nu_n}}
\Big(
U
-
\smallsum{1\leqslant i\leqslant n}\,
X_i\cdot u_{x_i}
\Big)
+
\smallsum{1\leqslant i\leqslant n}\,
X^i\cdot u_{x_ix_1^{\nu_1}\cdots x_n^{\nu_n}}.
\]
It is known that they depend on jet coordinates of order
$\leqslant \nu_1 + \cdots + \nu_n$:
\[
U_{\nu_1,\dots,\nu_n}
\,=\,
U_{\nu_1,\dots,\nu_n}
\Big(
x,u,\,u^{(\nu_1+\cdots+\nu_n)}
\Big).
\]

Now, according to Section~{\ref{tangency-order-n-1}}
({\em see} also the matrix at the
beginning of 
Section~{\ref{alternative-normalizations-order-n-plus-2}}),
the general vector field which stabilizes
the normalizations up order $n+1$,
{\em i.e.} which is tangent up to order $n+1$
to~{\eqref{u-F-n-plus-1}}, writes:
\[
\aligned
L
&
\,=\,
\Big(
A_{1,1}x_1
\ \ \ \ \ \ \ \ \ \ \ \ \ \ \ \ \ \ \ \ \ \ \ \ \ \ \ \ \ \ \ \ \ \ \
\ \ \ \ \ \ \ \ \ \ \ \ \ \ \ \ \ \ \ \ \ \ \ \ \ \ \ \ \ \ \ \ \ \ \
\ \ \ \ \ \ \ \ \ \ \ \ \ \ \ \ \ \ \ \ \ \ \ \ \ \ \ \ \ \
-
\tfrac{2}{2}A_{2,1}u
\Big)\frac{\partial}{\partial x_1}
\\
&
\ \ \ \ \
+
\Big(
A_{2,1}x_1
\ \ \ \ \ \ \ \ \ \ \ \ \ \ \ \ \ \ \ \ \ \ \ \ \ \ \ \ \ \ \ \ \ \ \
\ \ \ \ \ \ \ \ \ \ \ \ \ \ \ \ \ \ \ \ \ \ \ \ \ \ \ \ \ \ \ \ \ \ \
\ \ \ \ \ \ \ \ \ \ \ \ \ \ \ \ \ \ \ \ \ \ \ \ \ \ \ 
-
\tfrac{2}{3}A_{3,1}u
\Big)\frac{\partial}{\partial x_2}
\\
&
\ \ \ \ \
+
\Big(
A_{3,1}x_1
-
A_{1,1}x_3
\ \ \ \ \ \ \ \ \ \ \ \ \ \ \ \ \ \ \ \ \ \ \ \ \ \ \ \ \ \ \ \ \ \ \
\ \ \ \ \ \ \ \ \ \ \ \ \ \ \ \ \ \ \ \ \ \ \ \ \ \ \ \ \ \ \ \ \ \ \
\ \ \ \ \ \ \ \ \ \ \ \
-
\tfrac{2}{4}A_{4,1}u
\Big)\frac{\partial}{\partial x_3}
\\
&
\ \ \ \ \
+
\cdots\cdots\cdots\cdots\cdots\cdots\cdots\cdots\cdots
\cdots\cdots\cdots\cdots\cdots\cdots\cdots\cdots\cdots
\cdots\cdots\cdots\cdots\cdots\cdots\cdots
\cdot\cdot
\\
&
\ \ \ \ \
+
\Big(
A_{n-1,1}x_1
-
\tfrac{1}{n-3}\tbinom{n-1}{3}A_{n-3,1}x_3
-\cdots-
\tfrac{n-3}{1}\tbinom{n-1}{n-1}A_{1,1}x_{n-1}
\ \ \ \ \ \ \ \ \ \ \ \ \ \ \
-
\tfrac{2}{n}A_{n,1}u
\Big)\frac{\partial}{\partial x_{n-1}}
\\
&
\ \ \ \ \
+
\Big(
A_{n,1}x_1
-
\tfrac{1}{n-2}\tbinom{n}{3}A_{n-2,1}x_3
-\cdots-
\tfrac{n-3}{2}\tbinom{n}{n-1}A_{2,1}x_{n-1}
-
\tfrac{n-2}{1}\tbinom{n}{n}A_{1,1}x_n
+
B_nu
\Big)\frac{\partial}{\partial x_n}
\\
&
\ \ \ \ \
+
\Big(
\ \ \ \ \ \ \ \ \ \ \ \ \ \ \ \ \ \ \ \ \ \ \ \ \ \ \ \ \ \ \ \ \ \ \
\ \ \ \ \ \ \ \ \ \ \ \ \ \ \ \ \ \ \ \ \ \ \ \ \ \ \ \ \ \ \ \ \ \ \
\ \ \ \ \ \ \ \ \ \ \ \ \ \ \ \ \ \ \ \ \ \ \ \ \ \ \ \ \ \ \ \ \ \ \
\ \ \ \ \
+
2A_{1,1}u
\Big)\frac{\partial}{\partial u}.
\endaligned
\]
Then prolong it to the jet space of order $n+2$:
\[
\small
\aligned
L^{(n+2)}
\,=\,
L
+
\sum_{\nu_1+\cdots+\nu_n=1}\,
U_{\nu_1,\dots,\nu_n}
\big(
x,u,\,u^{(1)}
\big)
+\cdots
&
+
\sum_{\nu_1+\cdots+\nu_n=n+1}\,
U_{\nu_1,\dots,\nu_n}
\big(
x,u,\,u^{(n+1)}
\big)
\frac{\partial}{\partial u_{\nu_1,\dots,\nu_n}}
\\
&
+
\sum_{\nu_1+\cdots+\nu_n=n+2}\,
U_{\nu_1,\dots,\nu_n}
\big(
x,u,\,u^{(n+2)}
\big)
\frac{\partial}{\partial u_{\nu_1,\dots,\nu_n}}.
\endaligned
\]
The following (admitted) statement can be established
in a general theoretical context.

\begin{Proposition}
If $L$ is tangent to~{\eqref{u-F-order-n-2}} up to order
$n+1$, then at the origin $(x,u) = (0,0)$, it holds:
\[
0
\,=\,
U_{\nu_1,\dots,\nu_n}
\big(
0,0,\,
u^{(\nu_1+\cdots+\nu_n)}
\big)
\eqno
{\scriptstyle{(\forall\,\nu_1+\cdots+\nu_n\,\leqslant\,n+1)}}.
\qed
\]
\end{Proposition}

Furthermore, at order equal to $n+2$, considering only independent
jets, still at the origin, it can be shown
(in a general theoretical context) that one recovers 
(up to sign and a change of notation) 
the expressions appearing in
Section~{\ref{alternative-normalizations-order-n-plus-2}}.

\begin{Proposition}
If $L$ is tangent to~{\eqref{u-F-order-n-2}} up to order
$n+1$, then at the origin $(x,u) = (0,0)$, it holds:
\[
\footnotesize
\aligned
U_{n+2,0,\dots,0}
\big(
0,0,u^{(n+2)}
\big)
&
\,=\,
-\,n\,u_{n+2,0,\dots,0}\,A_{1,1}
-
\smallsum{2\leqslant k\leqslant n}\,
(n+2)\,u_{n+1,0\cdots1\cdots0}\,A_{k,1}
-
\tbinom{n+2}{2}\,B_n,
\\
U_{n+1,1,\dots,0}
\big(
0,0,u^{(n+2)}
\big)
&
\,=\,
-\,(n-1)\,u_{n+1,1,\dots,0}\,A_{1,1},
\notag
\\
U_{n+1,0,1,\dots,0}
\big(
0,0,u^{(n+2)}
\big)
&
\,=\,
-\,
(n-2)\,u_{n+1,0,1,\dots,0}\,A_{1,1}
-
\ast\,A_{2,1}-\cdots-\ast\,A_{n-2,1}
-
\tfrac{-2+3}{(n-1)!3!}\,
(n+1)!\,
A_{n-1,1}
\\
\cdots\cdots\cdots\cdots\cdots\cdots\cdots
\cdot
&
\cdots\cdots\cdots\cdots\cdots\cdots\cdots\cdots\cdots
\cdots\cdots\cdots\cdots\cdots\cdots\cdots\cdots\cdots
\cdots\cdots\cdots\cdots\cdots\cdots\cdots\cdots\cdots
\\
U_{n+1,0\cdots1\cdots0}
\big(
0,0,u^{(n+2)}
\big)
&
\,=\,
-\,
(n-k+1)\,u_{n+1,0\cdots1\cdots0}\,A_{1,1}
-
\ast\,A_{2,1}-\cdots-\ast\,A_{n-k+1,1}
-
\tfrac{-2+k}{(n-k+2)!k!}\,(n+1)!\,
A_{n-k+2,1}
\\
\cdots\cdots\cdots\cdots\cdots\cdots\cdots
\cdot
&
\cdots\cdots\cdots\cdots\cdots\cdots\cdots\cdots\cdots
\cdots\cdots\cdots\cdots\cdots\cdots\cdots\cdots\cdots
\cdots\cdots\cdots\cdots\cdots\cdots\cdots
\\
U_{n+1,0,\dots,1,0}
\big(
0,0,u^{(n+2)}
\big)
&
\,=\,
-\,
2\,F_{n+1,0,\dots,1,0}\,A_{1,1}
-
\ast\,A_{2,1}
-
\tfrac{-2+n-1}{3!(n-1)!}\,
(n+1)!\,
A_{3,1}
\\
U_{n+1,0,\dots,0,1}
\big(
0,0,u^{(n+2)}
\big)
&
\,=\,
-\,
1\,F_{n+1,0,\dots,0,1}\,A_{1,1}
-
\tfrac{-2+n}{2!n!}\,
(n+1)!\,
A_{2,1}.
\qed
\endaligned
\]
\end{Proposition}

\begin{Example}
For instance, when $n = 3$, the vector field from
Section~{\ref{tangency-low-dimensions}}:
\[
\aligned
L
&
\,=\,
\big(
A_{1,1}\,x_1
\ \ \ \ \ \ \ \ \ \ \ \ \ \ \ \ \ \ \ \
-
\tfrac{2}{2}
A_{2,1}\,u
\big)\,\tfrac{\partial}{\partial x_1}
\\
&
\ \ \ \ \
+
\big(
A_{2,1}\,x_1
\ \ \ \ \ \ \ \ \ \ \ \ \ \ \ \ \ 
-
\tfrac{2}{3}
A_{3,1}\,u
\big)\,\tfrac{\partial}{\partial x_2}
\\
&
\ \ \ \ \
+
\big(
A_{3,1}\,x_1
-
A_{1,1}\,x_3
-
\tfrac{2}{4}
A_{4,1}\,u
\big)\,\tfrac{\partial}{\partial x_3}
\\
&
\ \ \ \ \
+
\big(
\ \ \ \ \ \ \ \ \ \ \ \ \ \ \ \ \ \ \ \ \ \ \ \ \ \ \ \ \ \ \ \ 
2\,A_{1,1}\,u
\big)\,\tfrac{\partial}{\partial u},
\endaligned
\]
has prolongation of order $5$ above the origin $(x,u) = (0,0)$
given by:
\[
\aligned
L^{(5)}
&
\,=\,
-\,
\big(
3\,u_{5,0,0}\,A_{1,1}
+
5\,u_{4,1,0}\,A_{2,1}
+
5\,u_{4,0,1}\,A_{3,1}
+
10\,B_3
\big)\,
\tfrac{\partial}{\partial u_{5,0,0}}
\\
&
\ \ \ \ \
-\,
\big(
2\,u_{4,1,0}\,A_{1,1}
\big)\,
\tfrac{\partial}{\partial u_{4,1,0}}
\\
&
\ \ \ \ \
-
\big(
u_{4,0,1}\,A_{1,1}
+
2\,A_{2,1}
\big)\,
\tfrac{\partial}{\partial u_{4,0,1}}.
\endaligned
\]
So with $n = 3$, 
in the space $\R^3 \in (u_{5,0,0}, u_{4,1,0}, u_{4,0,1})$
of pure $(n+2)$-jets, we have a linear space
of vector fields
parametrized by the $4$ free coefficients:
\[
A_{1,1},\,A_{2,1},\,A_{3,1},\,B_3.
\]
This is in fact the Lie algebra of the action of
the stability group in orders $\leqslant n+1$
acting on pure $n+2$ order monomials. 

Without computing the flows of these vector fields
(which would amount to recover
the formulas of Section~{\ref{stabilizing-G-n-2-3-4-5-6}}),
we can realize (once again) that the two normalizations:
\[
F_{4,0,1}
\,:=\,
0,
\ \ \ \ \ \ \ \ \ \ \ \ \ \ \ \ \ \ \ \ \ \ \ \ \ \
F_{5,0,0}
\,:=\,
0,
\]
are possible, as follows.

\smallskip\noindent$\bullet$\,
By taking $A_{2,1} := 1$ and the others zero, we have:
\[
L^{(5)}
\,=\,
-\,\big(5\,u_{4,1,0}\,\big)\,
\tfrac{\partial}{\partial u_{5,0,0}}
-
0\,
\tfrac{\partial}{\partial u_{4,1,0}}
-
2\,
\tfrac{\partial}{\partial u_{4,0,1}},
\]
and the flow of the
constant vector field along the $u_{4,0,1}$-axis
clearly crosses the axis $\{ u_{4,0,1} = 0 \}$.
 
\smallskip\noindent$\bullet$\,
Assuming $u_{4,0,1} = 0$, by taking $B_3 := 1$ and the others zero,
we have:
\[
L^{(5)}
\,=\,
-\,10\,
\tfrac{\partial}{\partial u_{5,0,0}}
-
0\,
\tfrac{\partial}{\partial u_{4,1,0}}
-
0\,
\tfrac{\partial}{\partial u_{4,0,1}},
\]
a vector field that stabilizes $\{ u_{4,0,1} = 0\}$,
and whose flow crosses the axis $\{ u_{5,0,0} = 0 \}$.\qed
\end{Example}

\SectionHead{Normalizations of Order $(n+3)$ Monomials}
{normalizations-order-n-plus-3}

Thanks to Proposition~{\ref{Prp-normalizations-n-plus-2}},
several order $n+2$ independent monomials can be normalized to zero:
\leqnomode\usetagform{default}
\begin{align}
\label{normalize-F-order-n-2}
\ \ \ \ \ \ \ \ \ \ \ \ \ \
F_{n+2,0,\dots,0}
\,:=\,
0,
\ \ \ \ \ \ \
F_{n+1,0,1,\dots,0}
\,:=\,
0,
\ \ \ \ \ \ \
\dots\dots,
\ \ \ \ \ \ \
F_{n+1,0,\dots,1}
\,:=\,0,
\end{align}
hence if we come back
to~{\eqref{u-F-order-n-2}}
and let appear the next order $n+3$ monomials\,\,---\,\,disregarding 
as before the body-dependent ones\,\,---, we have
the equation:
\[
\aligned
u
&
\,=\,
\frac{x_1^2}{2}
+
\frac{x_1^2x_2}{2}
+
\sum_{m=3}^n\,
\Big(
\frac{x_1^mx_m}{m!}
+
x_1^{m-1}
\sum_{i,j\geqslant 2
\atop
i+j=m+1}\,
\tfrac{1}{2}\,
\frac{x_ix_j}{(i-1)!(j-1)!}
\Big)
\\
&
\ \ \ \ \
+
0
+
F_{n+1,1,\dots,0}\,
\frac{x_1^{n+1}x_2}{(n+1)!}
+
0
+\cdots+
0
+
x_1^n
\sum_{i,j\geqslant 2
\atop
i+j=n+2}\,
\tfrac{1}{2}\,
\frac{x_ix_j}{(i-1)!(j-1)!}
\\
&
\ \ \ \ \
+
F_{n+3,0,\dots,0}\,
\frac{x_1^{n+2}x_1}{(n+3)!}
+
F_{n+2,1,\dots,0}\,
\frac{x_1^{n+2}x_2}{(n+2)!}
+
F_{n+2,0,1,\dots,0}\,
\frac{x_1^{n+2}x_3}{(n+2)!}
+\cdots+
F_{n+2,0,\dots,1}\,
\frac{x_1^{n+2}x_n}{(n+2)!}
\\
&
\ \ \ \ \ \ \ \ \ \ \ \ \ \ \ \ \ \ \ \ \ \ \ \ \ \ \ \ \ \
+
F_{n+1,1,\dots,0}\,
\frac{x_1^{n+1}x_2x_2}{n!}
+
x_1^{n+1}
\sum_{i,j\geqslant 2
\atop
i+j=n+3}\,
\tfrac{1}{2}\,
\frac{x_ix_j}{(i-1)!(j-1)!}
+
{\rm O}_{x'}(3)
+
{\rm O}_x(n+4),
\endaligned
\]
where the border-dependent monomials
written in the last line can be
checked to include the supplementary monomial
$F_{n+1,1,\dots,0}\,
\frac{x_1^nx_2x_2}{n!}$
by reasoning as in 
Proposition~{\ref{Prp-Lambda-i-i-nu-zero}}.

Next, we come back to the expression of a general affine
vector field $L$ tangent up to order $n+1$.
Taking account of 
the current normalization~{\eqref{normalize-F-order-n-2}},
by looking at the equation at the end of 
Section~{\ref{alternative-normalizations-order-n-plus-2}},
we see that the tangency of $L$ up to order $n+2$
requires:
\[
A_{2,1}
\,:=\,0,
\ \ \ \ \ \ \
A_{3,1}
\,:=\,0,
\ \ \ \ \ \ \
\dots\dots,
\ \ \ \ \ \ \
A_{n-1,1}
\,:=\,0,
\ \ \ \ \ \ \
B_n
\,:=\,
0,
\]
hence coming back to~{\eqref{Lstab-order-n-1}}
we get the following vector field:
\[
\aligned
L
&
\,=\,\ \ \
\big(
A_{1,1}\,x_1
\ \ \ \ \ \ \ \ \ \ \ \ \ \ \ \ \ \ \ \ \ \ \ \ \ \ \ \ \ \ \ \ \ \ \
\ \ \ \ \ \ \ \ \ \ \ \ \ \ \ \ \ \ \ \ \ \ \ \ \ \ \ \ \ \ \ \ \ \ \ 
\ \ \ \ \ \ \ \ \ \ \ \ \ \ \ \ \ \ \ \ \ \ \ \ \ 
\big)\,
\partial_{x_1}
\\
&
\ \ \ \ \
+
\big(
0
\ \ \ \ \ \ \ \ \ \ \ \ \ \ \ \ \ \ \ \ \ \ \ \ \ \ \ \ \ \ \ \ \ \ \
\ \ \ \ \ \ \ \ \ \ \ \ \ \ \ \ \ \ \ \ \ \ \ \ \ \ \ \ \ \ \ \ \ \ \
\ \ \ \ \ \ \ \ \ \ \ \ \ \ \ \ \ \ \ \ \ \ \ \ \ \ \ \ \ \ \ \
\ \
\big)\,
\partial_{x_2}
\\
&
\ \ \ \ \
+
\big(
\ \ \ \ \ \ \ \ \ \ \ \ \ 
-
A_{1,1}\,x_3
\ \ \ \ \ \ \ \ \ \ \ \ \ \ \ \ \ \ \ \ \ \ \ \ \ \ \ \ \ \ \ \ \ \ \
\ \ \ \ \ \ \ \ \ \ \ \ \ \ \ \ \ \ \ \ \ \ \ \ \ \ \ \ \ \ \ \ \ \ \
\ \ \ \ \ \ \ 
\big)\,
\partial_{x_3}
\\
&
\ \ \ \ \
+
\cdots\cdots\cdots\cdots\cdots\cdots\cdots\cdots\cdots\cdots\cdots
\cdots\cdots\cdots\cdots\cdots\cdots\cdots\cdots\cdots\cdots\cdot
\\
&
\ \ \ \ \
+
\big(
\ \ \ \ \ \ \ \ \ \ \ \ \ \ \ \ \ \ \ \ \ \ \ \ \ \ \ \ \ \ \ \ \ \ \
\ \ \ \ \ \ \ \ \ \ \ \ \ \ \ \
-(n-3)A_{1,1}x_{n-1}
\ \ \ \ \
-
\tfrac{2}{n}\,A_{n,1}\,u
\big)\,
\partial_{x_{n-1}}
\\
&
\ \ \ \ \
+
\big(
A_{n,1}x_1
\ \ \ \ \ \ \ \ \ \ \ \ \ \ \ \ \ \ \ \ \ \ \ \ \ \ \ \ \ \ \ \ \ \ \
\ \ \ \ \ \ \ \ \ \ \ \ \ \ \ \ \ \ \ \
-(n-2)A_{1,1}x_n
\ \ \ \ \ \ \ \ \ \ \ \ \
\big)\,
\partial_{x_n}
\\
&
\ \ \ \ \
+
\big(
\ \ \ \ \ \ \ \ \ \ \ \ \ \ \ \ \ \ \ \ \ \ \ \ \ \ \ \ \ \ \ \ \ \ \
\ \ \ \ \ \ \ \ \ \ \ \ \ \ \ \ \ \ \ \ \ \ \ \ \ \ \ \ \ \ \ \ \ \ \
\ \ \ \ \ \ \ \ \ \ \ \ \ \ \ \ \ \ \ \ \ \ \ \ \
2\,A_{1,1}\,u
\big)\,
\partial_u,
\endaligned
\]
whose flow, according to Assertion~{\ref{Ass-tangency-L-order-n-1}}
(at the next order), stabilizes monomials of 
order $\leqslant n+2$.

Now, to see the action on order $n+3$ monomials, we
have to apply $\pi_{n+3}^{\ind} (\centersmallbullet)$
to the tangency equation, which gives (without writing
terms which will not contribute):
\[
\footnotesize
\aligned
0
&
\,\equiv\,
-\,2\,A_{1,1}
\Big(
\cdots
+
F_{n+3,0,\dots,0}
\tfrac{x_1^{n+2}x_1}{(n+3)!}
+
F_{n+2,1,\dots,0}
\tfrac{x_1^{n+2}x_2}{(n+2)!}
+
F_{n+2,0,1,\dots,0}
\tfrac{x_1^{n+2}x_3}{(n+2)!}
+\cdots+
F_{n+2,0,\dots,1}
\tfrac{x_1^{n+2}x_n}{(n+2)!}
\Big)
\\
&
\ \ \ \ \
+
\big(
A_{1,1}x_1
\big)
\Big(
\cdots
+
F_{n+3,0,\dots,0}
\tfrac{x_1^{n+2}}{(n+2)!}
+
F_{n+2,1,\dots,0}
\tfrac{x_1^{n+1}x_2}{(n+1)!}
+
F_{n+2,0,1,\dots,0}
\tfrac{x_1^{n+1}x_3}{(n+1)!}
+\cdots+
F_{n+2,0,\dots,1}
\tfrac{x_1^{n+1}x_n}{(n+2)!}
\Big)
\\
&
\ \ \ \ \
+
0
\\
&
\ \ \ \ \
+
\big(-A_{1,1}x_3\big)
\Big(
\cdots
+
F_{n+2,0,1,\dots,0}
\tfrac{x_1^{n+2}}{(n+2)!}
+
x_1^{n+1}\,
\tfrac{x_n}{2!(n-1)!}
\Big)
\\
&
\ \ \ \ \
+
\big(-2A_{1,1}x_4\big)
\Big(
\cdots
+
F_{n+2,0,0,1,\dots,0}
\tfrac{x_1^{n+2}}{(n+2)!}
+
x_1^{n+1}\,
\tfrac{x_{n-1}}{3!(n-2)!}
\Big)
\\
&
\ \ \ \ \
+
\cdots\cdots\cdots\cdots\cdots\cdots\cdots\cdots\cdots
\cdots\cdots\cdots
\\
&
\ \ \ \ \
+
\big(
-(n-4)A_{1,1}x_{n-2}
\big)
\Big(
\cdots
+
F_{n+2,0,\dots,1,0,0}
\tfrac{x_1^{n+2}}{(n+2)!}
+
x_1^{n+1}
\tfrac{x_5}{(n-3)!4!}
\Big)
\\
&
\ \ \ \ \
+
\Big(
-(n-3)A_{1,1}x_{n-1}-\tfrac{2}{n}A_{n,1}
\big[
\tfrac{x_1^2}{2!}
\big\vert
+\tfrac{x_1^2x_2}{2!}+\tfrac{x_1^3x_3}{3!}
\big]
\Big)
\Big(
\tfrac{x_1^{n-1}}{(n-1)!}
\Big\vert
+
\tfrac{x_1^{n-1}x_2}{1!(n-2)!}
+
\tfrac{x_1^nx_3}{2!(n-2)!}
+
F_{n+2,0,\dots,1,0}
\tfrac{x_1^{n+2}}{(n+2)!}
+
\tfrac{x_1^{n+1}x_4}{3!(n-2)!}
\Big)
\\
&
\ \ \ \ \
+
\big(
A_{n,1}x_1
-
(n-2)A_{1,1}x_n
\big)
\Big(
\tfrac{x_1^n}{n!}
\Big\vert
+
\tfrac{x_1^nx_2}{1!(n-1)!}
+
F_{n+2,0,\dots,1}
\tfrac{x_1^{n+2}}{(n+2)!}
+
\tfrac{x_1^{n+1}x_3}{2!(n-1)!}
\Big).
\endaligned
\]
We then collect all the independent monomials of order $n+3$:
\[
\!\!\!\!\!\!\!\!\!\!\!\!\!\!\!
\scriptsize
\aligned
0
&
\equiv
x_1^{n+3}
\Big\{
-\tfrac{2}{(n+3)!}
F_{n+3,0,\dots,0}
A_{1,1}
+
\tfrac{1}{(n+2)!}
F_{n+3,0,\dots,0}
A_{1,1}
+
\tfrac{1}{(n+2)!}
F_{n+2,0,\dots,1}
A_{n,1}
\Big\}
\\
&
\ \ \ \ \
x_1^{n+2}x_2
\Big\{
\tfrac{-2}{(n+2)!}
F_{n+2,1,\dots,0}
A_{1,1}
+
\tfrac{1}{(n+1)!}
F_{n+2,1,\dots,0}
A_{1,1}
\Big\}
\\
&
\ \ \ \ \
x_1^{n+2}x_3
\Big\{
\tfrac{-2}{(n+2)!}
F_{n+2,0,1,\dots,0}
A_{1,1}
+
\tfrac{1}{(n+1)!}
F_{n+2,0,1,\dots,0}
A_{1,1}
-
\tfrac{1}{(n+2)!}
F_{n+2,0,1,\dots,0}
A_{1,1}
-
\tfrac{2}{n}
\tfrac{1}{3!(n-1)!}
A_{n,1}
-
\tfrac{2}{n}
\tfrac{1}{2!}
\tfrac{1}{(n-2)!2!}
A_{n,1}
+
\tfrac{1}{2!(n-1)!}
A_{n,1}
\Big\}
\\
&
\ \ \ \ \
x_1^{n+2}x_4
\Big\{
\tfrac{-2}{(n+2)!}
F_{n+2,0,0,1,\dots,0}
A_{1,1}
+
\tfrac{1}{(n+1)!}
F_{n+2,0,0,1,\dots,0}
A_{1,1}
-
\tfrac{2}{(n+2)!}
F_{n+2,0,0,1,\dots,0}
A_{1,1}
\Big\}
\\
&
\ \ \ \ \
\cdots\cdots\cdots\cdots\cdots\cdots\cdots\cdots\cdots
\cdots\cdots\cdots\cdots\cdots\cdots\cdots\cdots\cdots
\\
&
\ \ \ \ \
x_1^{n+2}x_{n-2}
\Big\{
\tfrac{-2}{(n+2)!}
F_{n+2,0,\dots,1,0,0}
A_{1,1}
+
\tfrac{1}{(n+1)!}
F_{n+2,0,\dots,1,0,0}
A_{1,1}
-
\tfrac{n-4}{(n+2)!}
F_{n+2,0,\dots,1,0,0}
A_{1,1}
\Big\}
\\
&
\ \ \ \ \
x_1^{n+2}x_{n-1}
\Big\{
\tfrac{-2}{(n+2)!}
F_{n+2,0,\dots,1,0}
A_{1,1}
+
\tfrac{1}{(n+1)!}
F_{n+2,0,\dots,1,0}
A_{1,1}
-
\tfrac{n-3}{(n+2)!}
F_{n+2,0,\dots,1,0}
A_{1,1}
\Big\}
\\
&
\ \ \ \ \
x_1^{n+2}x_n
\Big\{
\tfrac{-2}{(n+2)!}
F_{n+2,0,\dots,1}
A_{1,1}
+
\tfrac{1}{(n+1)!}
F_{n+2,0,\dots,1}
A_{1,1}
-
\tfrac{n-2}{(n+2)!}
F_{n+2,0,\dots,1}
A_{1,1}
\Big\}.
\endaligned
\]
After simplifications, this contracts as:
\[
\aligned
0
&
\equiv
x_1^{n+3}
\Big\{
\tfrac{n+1}{(n+3)!}
F_{n+3,0,\dots,0}
A_{1,1}
+
\tfrac{1}{(n+2)!}
F_{n+2,0,\dots,1}
A_{n,1}
\Big\}
\\
&
\ \ \ \ \
x_1^{n+2}x_2
\Big\{
\tfrac{n}{(n+2)!}
F_{n+2,1,\dots,0}
A_{1,1}
\Big\}
\\
&
\ \ \ \ \
x_1^{n+2}x_3
\Big\{
\tfrac{n-1}{(n+2)!}
F_{n+2,0,1,\dots,0}
A_{1,1}
+
\tfrac{1}{3!n!}
A_{n,1}
\Big\}
\\
&
\ \ \ \ \
x_1^{n+2}x_4
\Big\{
\tfrac{n-2}{(n+2)!}
F_{n+2,0,0,1,\dots,0}
A_{1,1}
\Big\}
\\
&
\ \ \ \ \
\cdots\cdots\cdots\cdots\cdots\cdots\cdots\cdots\cdots
\cdots\cdots\cdots\cdots\cdots\cdots\cdots\cdots\cdots
\\
&
\ \ \ \ \
x_1^{n+2}x_{n-2}
\Big\{
\tfrac{4}{(n+2)!}
F_{n+2,0,\dots,1,0,0}
A_{1,1}
\Big\}
\\
&
\ \ \ \ \
x_1^{n+2}x_{n-1}
\Big\{
\tfrac{3}{(n+2)!}
F_{n+2,0,\dots,1,0}
A_{1,1}
\Big\}
\\
&
\ \ \ \ \
x_1^{n+2}x_n
\Big\{
\tfrac{2}{(n+2)!}
F_{n+2,0,\dots,1}
A_{1,1}
\Big\}.
\endaligned
\]

Thus, we see that we can normalize:
\[
F_{n+2,0,1,\dots,0}
\,:=\,
0,
\]
while all the other independent coefficients of order $n+3$
are relative invariants.

Finally, to stabilize the obtained order $(n+3)$ normalizations:
\[
\aligned
u
&
\,=\,
\frac{x_1^2}{2}
+
\frac{x_1^2x_2}{2}
+
\sum_{m=3}^n\,
\Big(
\frac{x_1^mx_m}{m!}
+
x_1^{m-1}
\sum_{i,j\geqslant 2
\atop
i+j=m+1}\,
\tfrac{1}{2}\,
\frac{x_ix_j}{(i-1)!(j-1)!}
\Big)
\\
&
\ \ \ \ \
+
0
+
F_{n+1,1,\dots,0}\,
\frac{x_1^{n+1}x_2}{(n+1)!}
+
0
+\cdots+
0
+
x_1^n
\sum_{i,j\geqslant 2
\atop
i+j=n+2}\,
\tfrac{1}{2}\,
\frac{x_ix_j}{(i-1)!(j-1)!}
\\
&
\ \ \ \ \
+
F_{n+3,0,\dots,0}\,
\frac{x_1^{n+2}x_1}{(n+3)!}
+
F_{n+2,1,\dots,0}\,
\frac{x_1^{n+2}x_2}{(n+2)!}
+
0
+
F_{n+2,0,0,1,\dots,0}\,
\frac{x_1^{n+2}x_4}{(n+2)!}
+\cdots+
F_{n+2,0,\dots,1}\,
\frac{x_1^{n+2}x_n}{(n+2)!}
\\
&
\ \ \ \ \ \ \ \ \ \ \ \ \ \ \ \ \ \ \ \ \ \ \ \ \ \ \ \ \ \
+
F_{n+1,1,\dots,0}\,
\frac{x_1^{n+1}x_2x_2}{n!}
+
x_1^{n+1}
\sum_{i,j\geqslant 2
\atop
i+j=n+3}\,
\tfrac{1}{2}\,
\frac{x_ix_j}{(i-1)!(j-1)!}
+
{\rm O}_{x'}(3)
+
{\rm O}_x(n+4),
\endaligned
\]
we deduce that:
\[
A_{n,1}
\,:=\,
0.
\]
The isotropy Lie algebra is $1$-dimensional, represented
by the matrix:
\[
\left[
\def\arraystretch{1.25}
\begin{array}{cccccccc}
A_{1,1} & \red{\bf 0} & \red{\bf 0} & \red{\bf 0} &
\red{\cdots} & \red{\bf 0} & \red{\bf 0} & \red{\bf 0}
\\
\red{\bf 0} & \red{\bf 0} & \red{\bf 0} & \red{\bf 0} &
\red{\cdots} & \red{\bf 0}& \red{\bf 0} & \red{\bf 0}
\\
\red{\bf 0} & \red{\bf 0} & \red{-A_{1,1}} & \red{\bf 0} &
\red{\cdots} & \red{\bf 0} & \red{\bf 0} & \red{\bf 0}
\\
\red{\bf 0} & \red{\bf 0} & \red{\bf 0} & \red{-2A_{1,1}} &
\red{\cdots} & \red{\bf 0} & \red{\bf 0} & \red{\bf 0}
\\
\vdots & \red{\vdots} & \red{\vdots} & \red{\vdots} & 
\red{\ddots} & \red{\vdots} & \vdots & \red{\vdots}
\\
\red{\bf 0} & \red{\bf 0} & \red{\bf 0} &
\red{\bf 0} & \red{\cdots} &
\red{-(n-3)A_{1,1}} & \red{\bf 0} & \red{\bf 0}
\\
\red{\bf 0} & \red{\bf 0} & \red{\bf 0} & \red{\bf 0} & \cdots & 
\red{\bf 0} & \red{-(n-2)A_{1,1}} & 
\red{\bf 0}
\\
\red{\bf 0} & \red{\bf 0} & \red{\bf 0} & \red{\bf 0} &
\red{\cdots} & \red{\bf 0} & \red{\bf 0} & \red{2\,A_{1,1}}
\end{array}
\right].
\]

\SectionHead{Orders $(n+4)$ and $(n+3)$}
{orders-n-3-n-4}

The corresponding matrix Lie group:
\[
\left[\!
\begin{array}{cccccccc}
a_{1,1} & 0 & 0 & 0 & \cdots & 0 & 0 & 0
\\ 
0 & 1 & 0 & 0 & \cdots & 0 & 0 & 0
\\
0 & 0 & \frac{1}{a_{1,1}} & 0 & \cdots & 0 & 0 & 0
\\
0 & 0 & 0 & \frac{1}{a_{1,1}^2} & \cdots & 0 & 0 & 0
\\
\vdots & \vdots & \vdots & \vdots & \ddots & \vdots & \vdots & \vdots
\\
0 & 0 & 0 & 0 & \cdots & \frac{1}{a_{1,1}^{n-3}} & 0 & 0
\\
0 & 0 & 0 & 0 & \cdots & 0 & \frac{1}{a_{1,1}^{n-2}} & 0
\\
0 & 0 & 0 & 0 & \cdots & 0 & 0 & a_{1,1}^2
\end{array}
\!\right],
\]
consists of the plain dilations:
\[
y_1
\,=\,
a_{1,1}\,x_1,
\ \ \ \ \ \
y_2
\,=\,
0,
\ \ \ \ \ \
y_3
\,=\,
\tfrac{1}{a_{1,1}}\,
x_3,
\ \ \ \ \ \
\dots,
\ \ \ \ \ \
y_n
\,=\,
\tfrac{1}{a_{1,1}^{n-2}}\,
x_n,
\ \ \ \ \ \
v
\,=\,
a_{1,1}^2\,u.
\]

\begin{Observation}
All power series coefficients $F_{\sigma_1,\dots,\sigma_n}$ 
are relative invariants.
\end{Observation}

\proof
So the possible remaining maps
fixing the origin which 
send a hypersurface normalized as above:
\[
u
\,=\,
\sum_{\sigma_1,\dots,\sigma_n}\,
x_1^{\sigma_1}
\cdots
x_n^{\sigma_n}\,
F_{\sigma_1,\dots,\sigma_n},
\] 
to a similarly normalized hypersurface:
\[
v
\,=\,
\sum_{\sigma_1,\dots,\sigma_n}\,
y_1^{\sigma_1}
\cdots
y_n^{\sigma_n}\,
G_{\sigma_1,\dots,\sigma_n},
\]
are only the dilations:
\[
y_1
\,=\,
\aaux\,x_1,
\ \ \ \ \ \
y_2
\,=\,
0,
\ \ \ \ \ \
y_3
\,=\,
\tfrac{1}{\aaux}\,
x_3,
\ \ \ \ \ \
\dots,
\ \ \ \ \ \
y_n
\,=\,
\tfrac{1}{\aaux^{n-2}}\,
x_n,
\ \ \ \ \ \
v
\,=\,
\aaux^2\,u.
\]
where we have abbreviated $a_{1,1} =: \aaux$. 

After replacement:
\[
\aaux^2\,u
\,=\,
\sum_{\sigma_1,\dots,\sigma_n}\,
\big(\aaux\,x_1\big)^{\sigma_1}\,
(x_2)^{\sigma_2}\,
\big(
\tfrac{1}{\aaux}\,x_3
\big)^{\sigma_3}\,
\cdots\,
\big(
\tfrac{1}{\aaux^{n-2}}\,x_n
\big)^{\sigma_n}\,\,
G_{\sigma_1,\dots,\sigma_n},
\]
an identification gives:
\[
F_{\sigma_1,\dots,\sigma_n}
\,=\,
\aaux^{-2}\,
\aaux^{\sigma_1}\,
\aaux^{-\sigma_3}\,
\cdots\,
\aaux^{-(n-2)\sigma_n}\,\,
G_{\sigma_1,\dots,\sigma_n}.
\]
Thus, each couple of power series coefficients are nonzero multiple
one of the other. Thus, they vanish (or do not vanish), 
simultaneously.
\endproof

Consequently, at the next two orders $n+4$ and $n+5$, 
we will not perform any further normalization,
since this would create some branching.

But we must determine the 
independent and border-dependent monomials 
of homogeneous degrees $n+4$ and $n+5$.

Proceeding as in Proposition~{\ref{Prp-Lambda-i-i-nu-zero}}
to take account of the hypothesis that the
Hessian matrix is of constant rank $1$, 
we obtain the following explicit expressions.
We skip presenting the details of computations.

\begin{Theorem}
\label{Thm-nf-order-n-5}
In any dimension $n \geqslant 2$, every local hypersurface $H^n \subset
\R^{n+1}$ having constant Hessian rank $1$ 
which is not affinely equivalent to a product 
of $\R^m$ $(1 \leqslant m \leqslant n)$
with a hypersurface $H^{n-m} \subset \R^{n-m+1}$ can 
be affinely normalized as:
\[
\footnotesize
\aligned
u
&
\,=\,
\frac{x_1^2}{2}
+
\frac{x_1^2\,x_2}{2}
+
\sum_{m=3}^n\,
\Big(
\frac{x_1^m\,x_m}{m!}
+
x_1^{m-1}
\sum_{i,j\geqslant 2
\atop
i+j=m+1}\,
\tfrac{1}{2}\,
\frac{x_i\,x_j}{(i-1)!(j-1)!}
\Big)
\\
&
\ \ \ \ \
+
F_{n+1,10\cdots0}\,
\frac{x_1^{n+1}x_2}{(n+1)!}
+
x_1^n
\sum_{i,j\geqslant 2
\atop
i+j=n+2}\,
\tfrac{1}{2}\,
\frac{x_i\,x_j}{(i-1)!(j-1)!}
\\
&
\ \ \ \ \
+
F_{n+3,0\cdots0}\,
\frac{x_1^{n+3}}{(n+3)!}
+
F_{n+2,10\cdots0}\,
\frac{x_1^{n+2}\,x_2}{(n+2)!}
+
F_{n+2,0010\cdots0}\,
\frac{x_1^{n+2}\,x_4}{(n+2)!}
+\cdots+
F_{n+2,0\cdots01}\,
\frac{x_1^{n+2}\,x_n}{(n+2)!}
\\
&
\ \ \ \ \ \ \ \ \ \ \ \ \ \ \ \ \ \ \ \ \ \ \ \ \ \ \ \ \ \
+
F_{n+1,10\cdots0}\,
\frac{x_1^{n+1}\,x_2\,x_2}{n!}
+
x_1^{n+1}
\sum_{i,j\geqslant 2
\atop
i+j=n+3}\,
\tfrac{1}{2}\,
\frac{x_i\,x_j}{(i-1)!(j-1)!}
\endaligned
\]
\[
\footnotesize
\aligned
\!\!\!\!\!\!\!\!\!\!\!\!\!\!\!
+
F_{n+4,0\cdots0}\,
\frac{x_1^{n+4}}{(n+4)!}
&
+
F_{n+3,10\cdots0}\,
\frac{x_1^{n+3}\,x_2}{(n+3)!}
+
F_{n+3,010\cdots0}\,
\frac{x_1^{n+3}\,x_3}{(n+3)!}
+
F_{n+3,0010\cdots0}\,
\frac{x_1^{n+3}\,x_4}{(n+3)!}
+\cdots+
F_{n+3,0\cdots01}\,
\frac{x_1^{n+3}\,x_n}{(n+3)!}
\\
&
+
x_1^{n+2}\,
\bigg[
\frac{F_{n+2,10\cdots0}}{(n+1)!}\,
x_2\,x_2
+
\frac{F_{n+1,10\cdots0}}{2!\,n!}\,
x_2\,x_3
+
\frac{F_{n+2,0010\cdots0}}{(n+1)!}\,
x_2\,x_4
+\cdots+
\frac{F_{n+2,0\cdots01}}{(n+1)!}\,
x_2\,x_n
\\
&
\ \ \ \ \ \ \ \ \ \ \ \ \ \ \ \ 
+
\sum_{i,j\geqslant 2
\atop
i+j=n+4}\,
\tfrac{1}{2}\,
\frac{x_i\,x_j}{(i-1)!(j-1)!}
\bigg]
\endaligned
\]
\[
\footnotesize
\!\!\!\!\!\!\!\!\!\!\!\!\!\!\!
\aligned
{}
&
+
F_{n+5,0\cdots0}\,
\frac{x_1^{n+5}}{(n+5)!}
+
F_{n+4,10\cdots0}\,
\frac{x_1^{n+4}\,x_2}{(n+4)!}
+
F_{n+4,010\cdots0}\,
\frac{x_1^{n+4}\,x_3}{(n+4)!}
+
F_{n+4,0010\cdots0}\,
\frac{x_1^{n+4}\,x_4}{(n+4)!}
+\cdots+
F_{n+4,0\cdots01}\,
\frac{x_1^{n+4}\,x_n}{(n+4)!}
\\
&
+
x_1^{n+3}\,
\bigg[
\Big(
\frac{F_{n+3,10\cdots0}}{(n+2)!}
-
\frac{F_{n+3,0\cdots0}}{2!\,(n+1)!}
\Big)x_2x_2
+
\Big(
\frac{F_{n+3,010\cdots0}}{(n+2)!}
+
\frac{F_{n+2,10\cdots0}}{2!\,(n+1)!}
\Big)x_2x_3
+
\Big(
\frac{F_{n+3,0010\cdots0}}{(n+2)!}
+
\frac{F_{n+1,10\cdots0}}{3!\,n!}
\Big)x_2x_4
\\
&
\ \ \ \ \ \ \ \ \ \ \ \ \ \ \ \ 
+
\frac{F_{n+3,00010\cdots0}}{(n+2)!}\,x_2x_5
+\cdots+
\frac{F_{n+3,0\cdots01}}{(n+2)!}\,x_2x_n
\\
&
\ \ \ \ \ \ \ \ \ \ \ \ \ \ \ \ 
+
\frac{F_{n+2,0010\cdots0}}{2!\,(n+2)!}\,x_3x_4
+\cdots+
\frac{F_{n+2,0\cdots01}}{2!\,(n+2)!}\,x_3x_n
\\
&
\ \ \ \ \ \ \ \ \ \ \ \ \ \ \ \ \ \ \ \ \ \ \ \ \ \ \ \ \ \ \ \ \ \ 
\ \ \ \ \ \ \ \ \ \ \ \ \ \ \ \ \ \ \ \ \ \ \ 
+
\sum_{i,j\geqslant 2
\atop
i+j=n+5}\,
\tfrac{1}{2}\,
\frac{x_i\,x_j}{(i-1)!(j-1)!}
\bigg]
+
{\rm O}_{x_2,\dots,x_n}(3)
+
{\rm O}_{x_1,x_2,\dots,x_n}(n+6).
\endaligned
\]
\end{Theorem}

This explicit expression of the graphing function
$F (x_1, \dots, x_n)$ is our new starting point. 

\SectionHead{Summary of the Proof of the Main
Theorem~{\ref{Thm-inexistence-n-5}}}
{summary-proof-theorem-orders-n-3-n-4}

Now, we take a general affine vector field which does not
necessarily vanish at the origin:
\[
\aligned
L
&
\,=\,
\Big(
T_0
+
C_1\,x_1
+
C_2\,x_2
+\cdots+
C_{n-1}\,x_{n-1}
+
C_n\,x_n
\Big)\,
\frac{\partial}{\partial u}
\\
&
\ \ \ \ \
+
\Big(
T_1
+
A_{1,1}\,x_1
+
A_{1,2}\,x_2
+\cdots+
A_{1,n-1}\,x_{n-1}
+
A_{1,n}\,x_n
\Big)\,
\frac{\partial}{\partial x_1}
\\
&
\ \ \ \ \
+
\Big(
T_2
+
A_{2,1}\,x_1
+
A_{2,2}\,x_2
+\cdots+
A_{2,n-1}\,x_{n-1}
+
A_{2,n}\,x_n
\Big)\,
\frac{\partial}{\partial x_2}
\\
&
\ \ \ \ \
+
\cdots\cdots\cdots\cdots\cdots\cdots\cdots\cdots\cdots\cdots\cdots
\cdots\cdots\cdots\cdots\cdots\cdots\cdots\cdots
\\
&
\ \ \ \ \
+
\Big(
T_{n-1}
+
A_{n-1,1}\,x_1
+
A_{n-1,2}\,x_2
+\cdots+
A_{n-1,n-1}\,x_{n-1}
+
A_{n-1,n}\,x_n
\Big)\,
\frac{\partial}{\partial x_{n-1}}
\\
&
\ \ \ \ \
+
\Big(
T_n
+
A_{n,1}\,x_1
+
A_{n,2}\,x_2
+\cdots+
A_{n,n-1}\,x_{n-1}
+
A_{n,n}\,x_n
\Big)\,
\frac{\partial}{\partial x_n}.
\endaligned
\]
Remind that if $L$ is tangent to
the hypersurface $H = \big\{
u = F(x_1, \dots, x_n) \big\}$,
then $T_0 = 0$, since
$u = F = {\rm O}_x(2)$.

Here, the parameters $T_1, \dots, T_n$ are 
tightly related to the 
{\em infinitesimal transitivity} of the action, since
the value of $L$ at the origin is:
\[
L
\big\vert_0
\,=\,
T_1\,\frac{\partial}{\partial x_1}
+\cdots+
T_n\,\frac{\partial}{\partial x_n},
\]
since homogeneity requires that:
\[
T_0H
\,=\,
\Span_\R\,
\big\{
L\big\vert_0
\colon\,
L\vert_H\,\,
\text{tangent to}\,\,
H
\big\},
\]
and since, due again to $F = {\rm O}_x(2)$:
\[
T_0H
\,=\,
\Span\,
\Big\{
\frac{\partial}{\partial x_1},\,\,
\dots,\,\,
\frac{\partial}{\partial x_n}
\Big\}.
\]

\begin{Observation}
For (infinitesimal) affine homogeneity to hold, 
the parameters $T_1, \dots, T_n$ 
should remain absolutely free in all computations.\qed
\end{Observation}

Now, such a general affine vector field $L$ 
is an infinitesimal affine symmetry of our
hypersurface $H = \big\{ u = F(x_1,\dots,x_n) \big\}$
graphed as in
Theorem~{\ref{Thm-nf-order-n-5}},
if and only if $L\vert_H$ is tangent to $H$,
if and only if the following power series identity holds in 
$\R\{x_1, \dots, x_n\}$:
\[
0
\,\equiv\,
L\big(-\,u+F\big)
\Big\vert_{u=F}.
\]

We will in fact `only' study independent monomials of
order $\leqslant n+4$ in this fundamental equation,
namely we will examine\big/compute:
\[
\pi_{\ind}^{n+4}
\Big(
L\big(-\,u+F\big)
\Big\vert_{u=F}
\Big).
\]
Reminding\,\,---\,\,{\em see} also below\,\,---\,\,that:
\[
F_{x_1}
\,=\,
{\rm O}_x(1),
\ \ \ \ \ \
F_{x_2}
\,=\,
{\rm O}_x(2),
\ \ \ \ \ \
\cdots,
\ \ \ \ \ \
F_{x_{n-1}}
\,=\,
{\rm O}_x(n-1),
\ \ \ \ \ \
F_{x_n}
\,=\,
{\rm O}_x(n),
\]
we may therefore start out by writing:
\[
\footnotesize
\aligned
0
&
\,\equiv\,
-\,C_1\,x_1
-
C_2\,x_2
-\cdots-
C_{n-1}\,x_{n-1}
-
C_n\,x_n
-
D\,\big[\pi_{\ind}^{n+4}(F)]
&
\ \ \ \ \ \ \ \ \ \ \ \ \ \ \ \ \ \ \ \
&
\green{\Lambda_0}
\\
&
\ \ \ \ \
+
\Big(
T_1
+
A_{1,1}\,x_1
+
A_{1,2}\,x_2
+\cdots+
A_{1,n-1}\,x_{n-1}
+
A_{1,n}\,x_n
+
B_1\,
\big[
\pi_{\ind}^{n+3}(F)
\big]
\Big)\,
F_{x_1}
&
\ \ \ \ \ \ \ \ \ \ \ \ \ \ \ \ \ \ \ \
&
\green{\Lambda_1}
\\
&
\ \ \ \ \
+
\Big(
T_2
+
A_{2,1}\,x_1
+
A_{2,2}\,x_2
+\cdots+
A_{2,n-1}\,x_{n-1}
+
A_{2,n}\,x_n
+
B_2\,
\big[
\pi_{\ind}^{n+2}(F)
\big]
\Big)\,
F_{x_2}
&
\ \ \ \ \ \ \ \ \ \ \ \ \ \ \ \ \ \ \ \
&
\green{\Lambda_2}
\\
&
\ \ \ \ \
+
\Big(
T_3
+
A_{3,1}\,x_1
+
A_{3,2}\,x_2
+\cdots+
A_{3,n-1}\,x_{n-1}
+
A_{3,n}\,x_n
+
B_3\,
\big[
\pi_{\ind}^{n+1}(F)
\big]
\Big)\,
F_{x_3}
&
\ \ \ \ \ \ \ \ \ \ \ \ \ \ \ \ \ \ \ \
&
\green{\Lambda_3}
\\
&
\ \ \ \ \
+
\Big(
T_4
+
A_{4,1}\,x_1
+
A_{4,2}\,x_2
+\cdots+
A_{4,n-1}\,x_{n-1}
+
A_{4,n}\,x_n
+
B_4\,
\big[
\pi_{\ind}^{n}(F)
\big]
\Big)\,
F_{x_4}
&
\ \ \ \ \ \ \ \ \ \ \ \ \ \ \ \ \ \ \ \
&
\green{\Lambda_4}
\\
&
\ \ \ \ \
+
\Big(
T_5
+
A_{5,1}\,x_1
+
A_{5,2}\,x_2
+\cdots+
A_{5,n-1}\,x_{n-1}
+
A_{5,n}\,x_n
+
B_5\,
\big[
\pi_{\ind}^{n-1}(F)
\big]
\Big)\,
F_{x_5}
&
\ \ \ \ \ \ \ \ \ \ \ \ \ \ \ \ \ \ \ \
&
\green{\Lambda_5}
\\
&
\ \ \ \ \
+
\Big(
T_6
+
A_{6,1}\,x_1
+
A_{6,2}\,x_2
+\cdots+
A_{6,n-1}\,x_{n-1}
+
A_{6,n}\,x_n
+
B_6\,
\big[
\pi_{\ind}^{n-2}(F)
\big]
\Big)\,
F_{x_6}
&
\ \ \ \ \ \ \ \ \ \ \ \ \ \ \ \ \ \ \ \
&
\green{\Lambda_6}
\\
&
\ \ \ \ \
\cdots\cdots\cdots\cdots\cdots\cdots\cdots\cdots\cdots\cdots\cdots
\cdots\cdots\cdots\cdots\cdots\cdots\cdots\cdots\cdots\cdots\cdots
\cdots\cdots\cdots\cdots\cdots
&
\ \ \ \ \ \ \ \ \ \ \ \ \ \ \ \ \ \ \ \
&
\cdots
\\
&
\ \ \ \ \
+
\Big(
T_{n-1}
+
A_{n-1,1}x_1
+
A_{n-1,2}x_2
+\cdots+
A_{n-1,n-1}x_{n-1}
+
A_{n-1,n}x_n
+
B_{n-1}
\big[
\pi_{\ind}^5(F)
\big]
\Big)
F_{x_{n-1}}
&
\ \ \ \ \ \ \ \ \ \ \ \ \ \ \ \ \ \ \ \
&
\green{\Lambda_{n-1}}
\\
&
\ \ \ \ \
+
\Big(
T_n
+
A_{n,1}\,x_1
+
A_{n,2}\,x_2
+\cdots+
A_{n,n-1}\,x_{n-1}
+
A_{n,n}\,x_n
+
B_n\,
\big[
\pi_{\ind}^4(F)
\big]
\Big)\,
F_{x_n}.
&
\ \ \ \ \ \ \ \ \ \ \ \ \ \ \ \ \ \ \ \
&
\green{\Lambda_n}
\endaligned
\]
Notice on the right above that we attribute names to these $1 + n$ 
lines:
\[
\green{\Lambda_0},
\ \ \ \ \
\green{\Lambda_1},
\ \ \ \ \
\green{\Lambda_2},
\ \ \ \ \
\green{\Lambda_3},
\ \ \ \ \
\green{\Lambda_4},
\ \ \ \ \
\green{\Lambda_5},
\ \ \ \ \
\green{\Lambda_6},
\ \ \ \ \
\dots,
\ \ \ \ \
\green{\Lambda_{n-1}},
\ \ \ \ \
\green{\Lambda_n}.
\]

In the next Section~{\ref{tangency-equations-orders-n-4-n-5}},
we will compute some of the coefficients of the
monomials in this large equation,
namely we will compute some coefficients:
\[
E_{[\sigma_1,\dots,\sigma_n]}
\,:=\,
\big[
x^{\sigma_1}\cdots x_n^{\sigma_n}
\big]\,
\Big(
L\big(-\,u+F\big)
\Big\vert_{u=F}
\Big)
\eqno
{\scriptstyle{(\sigma_1+\cdots+\sigma_n\,\leqslant\,n+4)}},
\]
which are linear in
$C_\smallbullet$,
$D$, 
$T_\smallbullet$, 
$A_{\smallbullet, \smallbullet}$, 
$B_\smallbullet$, and which should vanish
for $L$ to really be tangent to $H$:
\[
E_{[\sigma_1,\dots,\sigma_n]}
\,=\,
0.
\]

Especially, we will compute:
\[
\aligned
\green{\bf I}
&
\,:=\,
E_{[n+2,0,\dots,0,1]}
\,=\,
0,
\\
\green{\bf II}
&
\,:=\,
E_{[n+3,0,\dots,0,1]}
\,=\,
0.
\endaligned
\]
But before proceeding to the (non-straightforward) computations,
let us summarize the key reason
why affinely homogeneous models do {\em not} exist
in dimension $n \geqslant 5$.
We use $\ast$ to denote any unspecified
real number whose value does not matter.

\begin{Proposition}
\label{Prp-equations-I-II}
For a hypersurface $\{ u = F(x)\}$
normalized as in 
Theorem~{\ref{Thm-nf-order-n-5}}, 
after taking account
of some of the
other equations $E_{[\sigma_1,\dots,\sigma_n]} = 0$, 
these two specific equations 
$\green{\bf I}$, $\green{\bf II}$
become of the form:
\[
\aligned
0
&
\overset{\green{\bf I}}{\,=\,}
\ast\,T_1
+
\ast\,T_2
-
\frac{1}{12\,(n-3)\,n!}\,
T_4
+
\frac{2}{(n+2)!}\,
F_{n+2,0\cdots01}\,
A_{1,1},
\\
0
&
\overset{\green{\bf II}}{\,=\,}
\ast\,T_1
+
\ast\,T_2
+
\ast\,
T_3
-
\frac{1}{30\,(n-4)\,n!}\,
T_5
+
\frac{3}{(n+3)!}\,
F_{n+3,0\cdots01}\,
A_{1,1}.
\endaligned
\]
\end{Proposition}

Admitting temporarily this fact, we can easily conclude our main
{\em inexistence} result.

\proof[Proof of Theorem~{\ref{Thm-inexistence-n-5}}]
If the power series coefficient $F_{n+2,0\cdots01} = 0$
would be zero, then the first equation:
\[
0
\overset{\green{\bf I}}{\,=\,}
\ast\,T_1
+
\ast\,T_2
-
\frac{1}{12\,(n-3)\,n!}\,
T_4,
\]
would consist of a {\em nontrivial} linear dependence relation
between $T_1, \dots, T_n$, contradicting infinitesimal transitivity.

So necessarily, $F_{n+2,0\cdots01} \neq 0$.

But then from equation $\green{\bf I}$, we can solve the isotropy 
parameter:
\[
A_{1,1}
\,=\,
\ast\,T_1
+
\ast\,T_2
+
\ast\,T_4,
\]
that we replace in $\green{\bf II}$, getting, whatever
the value of $F_{n+3,0\cdots01}$ is:
\[
0
\,=\,
\ast\,T_1
+
\ast\,T_2
+
\ast\,T_3
+
\ast\,T_4
-
\frac{1}{30\,(n-4)\,n!}\,
T_5.
\]

But such an equation is also {\em always} a {\em nontrivial} 
linear dependence relation
between $T_1, \dots, T_n$, contradicting again infinitesimal
transitivity!
\endproof

Observe that $n \geqslant 5$ was used in this argumentation.

\SectionHead{Tangency Equations at Orders $\leqslant n+4$}
{tangency-equations-orders-n-4-n-5}

So it remains `only' to prove Proposition~{\ref{Prp-equations-I-II}}.

\proof[Proof of Proposition~{\ref{Prp-equations-I-II}}]
Here is $\Lambda_0$:
\[
\footnotesize
\aligned
\Lambda_0
&
\,=\,
-\,C_1\,x_1
-\cdots-
C_n\,x_n
\\
&
\ \ \ \ \
-\,D\,
\bigg[
\frac{x_1^2}{2!}
+
\frac{x_1^2x_2}{2!}
+\cdots+
\frac{x_1^{n-1}x_{n-1}}{(n-1)!}
\\
&
\ \ \ \ \ \ \ \ \ \ \ \ \ \ \ \ \
+
\frac{x_1^nx_n}{n!}
\\
&
\ \ \ \ \ \ \ \ \ \ \ \ \ \ \ \ \
+
F_{n+1,10\cdots0}\,
\frac{x_1^{n+2}x_2}{(n+1)!}
\\
&
\ \ \ \ \ \ \ \ \ \ \ \ \ \ \ \ \
+
F_{n+3,0\cdots0}\,
\frac{x_1^{n+3}}{(n+3)!}
+
F_{n+2,10\cdots0}\,
\frac{x_1^{n+2}x_2}{(n+2)!}
+
F_{n+2,0010\cdots0}\,
\frac{x_1^{n+2}x_4}{(n+2)!}
+\cdots+
F_{n+2,0\cdots01}\,
\frac{x_1^{n+2}x_n}{(n+2)!}
\\
&
\ \ \ \ \ \ \ \ \ \ \ \ \ \ \ \ \
+
F_{n+4,0\cdots0}\,
\frac{x_1^{n+4}}{(n+4)!}
+
F_{n+3,10\cdots0}\,
\frac{x_1^{n+3}x_2}{(n+3)!}
+
F_{n+3,010\cdots0}\,
\frac{x_1^{n+3}x_3}{(n+3)!}
+
F_{n+3,0010\cdots0}\,
\frac{x_1^{n+4}x_4}{(n+3)!}
\\
&
\ \ \ \ \ \ \ \ \ \ \ \ \ \ \ \ \ \ \ \ \ \ \ \ \ \ \ \ \ \ \ \ \ \ \
\ \ \ \ \ \ \ \ \ \ \ \ \ \ \ \ \ \ \ \ \ \ \ \ \ \ \ \ \ \ \ \ \ \ \
\ \ \ \ \ \ \ \ \ \ \ \ \ \ \ \ \ \ \ \ \ \ \ \ \ \ \ \ \ \ \ \ \ \ \
\ \ \ \ \ \ \ \ \ \ \ \ \ \ \ \ \ \ \ \ \ \ \ \ \ 
+\cdots+
F_{n+3,0\cdots01}\,
\frac{x_1^{n+3}x_n}{(n+3)!}
\bigg].
\endaligned
\]

Next, we must write $\Lambda_1$, \dots, $\Lambda_n$,
which all involve products. 
We will denote the products using the sign "$\cdot$".

Here is $\Lambda_1$:
\[
\footnotesize
\!\!\!\!\!
\aligned
\Lambda_1
&
\,=\,
\left(
\aligned
T_1
&
+
A_{1,1}\,x_1
+
A_{1,2}\,x_2
+\cdots+
A_{1,n-1}\,x_{n-1}
+
A_{1,n}\,x_n
\\
&
+
B_1\,\bigg[
\frac{x_1^2}{2!}
+
\frac{x_1^2x_2}{2!}
+\cdots+
\frac{x_1^nx_n}{n!}
+
F_{n+1,10\cdots0}\,
\frac{x_1^{n+1}x_2}{(n+1)!}
\\
&
\ \ \ \ \ \ \ \ \ \ \ \ 
+
F_{n+3,0\cdots0}\,
\frac{x_1^{n+3}}{(n+3)!}
+
F_{n+2,10\cdots0}\,
\frac{x_1^{n+2}x_2}{(n+2)!}
+
F_{n+2,0010\cdots0}\,
\frac{x_1^{n+2}x_4}{(n+2)!}
\\
&
\ \ \ \ \ \ \ \ \ \ \ \ 
+\cdots+
F_{n+2,0\cdots01}\,
\frac{x_1^{n+2}x_n}{(n+2)!}
\bigg]
\endaligned
\right)
\cdot
\\
&
\!\!\!\!\!\!\!\!\!\!\!\!\!\!\!
\cdot
\left(
\aligned
{}
&
x_1
+
x_1x_2
+
\frac{x_1^2x_3}{2!}
+\cdots+
\frac{x_1^{n-1}x_n}{(n-1)!}
+
F_{n+1,10\cdots0}\,
\frac{x_1^nx_2}{n!}
\\
&
+
F_{n+3,0\cdots0}\,
\frac{x_1^{n+2}}{(n+2)!}
+
F_{n+2,10\cdots0}\,
\frac{x_1^{n+1}x_2}{(n+1)!}
+
F_{n+2,0010\cdots0}\,
\frac{x_1^{n+1}x_4}{(n+1)!}
+\cdots+
F_{n+2,0\cdots01}\,
\frac{x_1^{n+1}x_n}{(n+1)!}
\\
&
+
F_{n+4,0\cdots0}\,
\frac{x_1^{n+3}}{(n+3)!}
+
F_{n+3,10\cdots0}\,
\frac{x_1^{n+2}x_2}{(n+2)!}
+
F_{n+3,010\cdots0}\,
\frac{x_1^{n+2}x_3}{(n+2)!}
+
F_{n+3,0010\cdots0}\,
\frac{x_1^{n+2}x_4}{(n+2)!}
+\cdots+
F_{n+3,0\cdots01}\,
\frac{x_1^{n+2}x_n}{(n+2)!}
\\
&
+
F_{n+5,0\cdots0}\,
\frac{x_1^{n+4}}{(n+4)!}
+
F_{n+4,10\cdots0}\,
\frac{x_1^{n+3}x_2}{(n+3)!}
+
F_{n+4,010\cdots0}\,
\frac{x_1^{n+3}x_3}{(n+3)!}
+
F_{n+4,0010\cdots0}\,
\frac{x_1^{n+3}x_4}{(n+3)!}
+\cdots+
F_{n+4,0\cdots01}\,
\frac{x_1^{n+3}x_n}{(n+3)!}
\endaligned
\right).
\endaligned
\]

Here is $\Lambda_2$:
\[
\footnotesize
\!\!\!\!\!\!\!\!\!\!
\aligned
\Lambda_2
&
\,=\,
\left(
\aligned
T_2
&
+
A_{2,1}\,x_1
+
A_{2,2}\,x_2
+\cdots+
A_{2,n-1}\,x_{n-1}
+
A_{2,n}\,x_n
\\
&
+
B_2\,\bigg[
\frac{x_1^2}{2!}
+
\frac{x_1^2x_2}{2!}
+\cdots+
\frac{x_1^nx_n}{n!}
+
F_{n+1,10\cdots0}\,
\frac{x_1^{n+1}x_2}{(n+1)!}
\bigg]
\endaligned
\right)
\cdot
\\
&
\!\!\!\!\!\!\!\!\!\!\!\!\!\!\!
\cdot
\left(
\aligned
{}
&
\frac{x_1^2}{2!}
+
\sum_{m=3}^n\,
\frac{x_1^{m-1}x_{m-1}}{1!\,(m-2)!}
\\
&
+
F_{n+1,10\cdots0}\,
\frac{x_1^{n+1}}{(n+1)!}
+
\frac{x_1^nx_n}{(n-1)!}
\\
&
+
F_{n+2,10\cdots0}\,
\frac{x_1^{n+2}}{(n+2)!}
+
F_{n+1,10\cdots0}\,
\frac{2}{n!}\,
x_1^{n+1}x_2
\\
&
+
F_{n+3,10\cdots0}\,
\frac{x_1^{n+3}}{(n+3)!}
+
F_{n+2,10\cdots0}\,
\frac{2\,x_1^{n+2}x_2}{(n+1)!}\,
+
F_{n+1,10\cdots0}\,
\frac{x_1^{n+2}x_3}{2!\,n!}
+
F_{n+2,0010\cdots0}\,
\frac{x_1^{n+2}x_4}{(n+1)!}
+\cdots+
F_{n+2,0\cdots01}\,
\frac{x_1^{n+2}x_n}{(n+1)!}
\\
&
+
F_{n+4,10\cdots0}\,
\frac{x_1^{n+4}}{(n+1)!}
+
\Big(
\frac{2\,F_{n+3,10\cdots0}}{(n+2)!}
-
\frac{F_{n+3,0\cdots0}}{(n+1)!}
\Big)
x_1^{n+3}x_2
+
\Big(
\frac{F_{n+3,010\cdots0}}{(n+2)!}
+
\frac{F_{n+2,10\cdots0}}{2!\,(n+1)!}
\Big)\,
x_1^{n+3}x_3
\\
&
\ \ \ \ \ \ \ \ \ \ \ \ \ \ \ \ \ \ \ \ \ \ \ \ \ \ \ \ \ \ \ \ \ \ \
\ \ \ \ \ \ \ \ \ \ \ \ \ \ \ \ \ \ \ \ \ \ \ \ \ \ \ \ \ \ \ \ \ \ \
\ \ \ \ \ \ \ \ \ \ \ \ \ \ \ \ \ \ \ \ \ \ \ \ \ \ \ \ \ \ \ \ 
+
\Big(
\frac{F_{n+3,0010\cdots0}}{(n+2)!}
+
\frac{F_{n+1,10\cdots0}}{3!\,n!}
\Big)\,
x_1^{n+3}x_4
\endaligned
\right).
\endaligned
\]

Here is $\Lambda_3$:
\[
\footnotesize
\!\!\!\!\!
\aligned
\Lambda_3
&
\,=\,
\left(
\aligned
T_3
&
+
A_{3,1}\,x_1
+
A_{3,2}\,x_2
+\cdots+
A_{3,n-1}\,x_{n-1}
+
A_{3,n}\,x_n
\\
&
+
B_3\,\bigg[
\frac{x_1^2}{2!}
+
\frac{x_1^2x_2}{2!}
+\cdots+
\frac{x_1^nx_n}{n!}
\bigg]
\endaligned
\right)
\cdot
\\
&
\!\!\!\!\!\!\!\!\!\!\!\!\!\!\!
\cdot
\left(
\aligned
{}
&
\frac{x_1^3}{3!}
+
\sum_{m=4}^n\,
\frac{x_1^{m-1}x_{m-2}}{2!\,(m-3!)}
\\
&
+
\frac{x_1^nx_{n-1}}{2!\,(n-2)!}
\\
&
+
\frac{x_1^{n+1}x_n}{2!\,(n-1)!}
\\
&
+
F_{n+3,010\cdots0}\,
\frac{x_1^{n+3}}{(n+3)!}
+
F_{n+1,10\cdots0}\,
\frac{x_1^{n+2}x_3}{2!\,n!}
\\
&
+
F_{n+4,010\cdots0}\,
\frac{x_1^{n+4}}{(n+4)!}
+
\Big(
\frac{F_{n+3,010\cdots0}}{(n+2)!}
+
\frac{F_{n+2,10\cdots0}}{2!\,(n+1)!}
\Big)
x_1^{n+3}x_2
+
F_{n+2,0010\cdots0}\,
\frac{x_1^{n+3}x_4}{2!\,(n+1)!}
+\cdots+
F_{n+2,0\cdots01}\,
\frac{x_1^{n+3}x_n}{2!\,(n+1)!}
\endaligned
\right).
\endaligned
\]

Here is $\Lambda_4$:
\[
\footnotesize
\aligned
\Lambda_4
&
\,=\,
\left(
\aligned
T_4
&
+
A_{4,1}\,x_1
+
A_{4,2}\,x_2
+\cdots+
A_{4,n-1}\,x_{n-1}
+
A_{4,n}\,x_n
\\
&
+
B_4\,\bigg[
\frac{x_1^2}{2!}
+
\frac{x_1^2x_2}{2!}
+\cdots+
\frac{x_1^{n-1}x_{n-1}}{(n-1)!}
\bigg]
\endaligned
\right)
\cdot
\\
&
\!\!\!\!\!\!\!\!\!\!\!\!\!\!\!
\cdot
\left(
\aligned
{}
&
\frac{x_1^4}{4!}
+
\sum_{m=5}^n\,
\frac{x_1^{m-1}x_{m-3}}{3!\,(m-4)!}
\\
&
+
\frac{x_1^nx_{n-2}}{3!\,(n-3)!}
\\
&
+
F_{n+2,0010\cdots0}\,
\frac{x_1^{n+2}}{(n+2)!}
+
\frac{x_1^{n+1}x_{n-1}}{3!\,(n-2)!}
\\
&
+
F_{n+3,0010\cdots0}\,
\frac{x_1^{n+3}}{(n+3)!}
+
F_{n+2,0010\cdots0}\,
\frac{x_1^{n+2}x_2}{(n+1)!}
+
\frac{x_1^{n+2}x_n}{3!\,(n-1)!}
\\
&
+
F_{n+4,0010\cdots0}\,
\frac{x_1^{n+4}}{(n+4)!}
+
\Big(
\frac{F_{n+3,0010\cdots0}}{(n+2)!}
+
\frac{F_{n+1,10\cdots0}}{3!\,n!}
\Big)
x_1^{n+3}x_2
+
F_{n+2,0010\cdots0}\,
\frac{x_1^{n+3}x_3}{2!\,(n+1)!}
\endaligned
\right).
\endaligned
\]

Here is $\Lambda_5$:
\[
\footnotesize
\aligned
\Lambda_5
&
\,=\,
\left(
\aligned
T_5
&
+
A_{5,1}\,x_1
+
A_{5,2}\,x_2
+\cdots+
A_{5,n-1}\,x_{n-1}
+
A_{5,n}\,x_n
\\
&
+
B_5\,\bigg[
\frac{x_1^2}{2!}
+
\frac{x_1^2x_2}{2!}
+\cdots+
\frac{x_1^{n-2}x_{n-2}}{(n-2)!}
\bigg]
\endaligned
\right)
\cdot
\\
&
\!\!\!\!\!\!\!\!\!\!\!\!\!\!\!
\cdot
\left(
\aligned
{}
&
\frac{x_1^5}{5!}
+
\sum_{m=6}^n\,
\frac{x_1^{m-1}x_{m-4}}{4!\,(m-5)!}
\\
&
+
\frac{x_1^nx_{n-3}}{4!\,(n-4)!}
\\
&
+
F_{n+2,00010\cdots0}\,
\frac{x_1^{n+2}}{(n+2)!}
+
\frac{x_1^{n+1}x_{n-2}}{4!\,(n-3)!}
\\
&
+
F_{n+3,00010\cdots0}\,
\frac{x_1^{n+3}}{(n+3)!}
+
F_{n+2,00010\cdots0}\,
\frac{x_1^{n+2}x_2}{(n+1)!}
+
\frac{x_1^{n+2}x_{n-1}}{4!\,(n-2)!}
\\
&
+
F_{n+4,00010\cdots0}\,
\frac{x_1^{n+4}}{(n+4)!}
+
F_{n+3,00010\cdots0}\,
\frac{x_1^{n+3}x_2}{(n+2)!}
+
F_{n+2,00010\cdots0}\,
\frac{x_1^{n+3}x_3}{2!\,(n+1)!}
+
\frac{x_1^{n+3}x_n}{4!\,(n-1)!}
\endaligned
\right).
\endaligned
\]

Here is $\Lambda_6$:
\[
\footnotesize
\aligned
\Lambda_6
&
\,=\,
\left(
\aligned
T_6
&
+
A_{6,1}\,x_1
+
A_{6,2}\,x_2
+\cdots+
A_{6,n-1}\,x_{n-1}
+
A_{6,n}\,x_n
\\
&
+
B_6\,\bigg[
\frac{x_1^2}{2!}
+
\frac{x_1^2x_2}{2!}
+\cdots+
\frac{x_1^{n-3}x_{n-3}}{(n-3)!}
\bigg]
\endaligned
\right)
\cdot
\\
&
\!\!\!\!\!\!\!\!\!\!\!\!\!\!\!
\cdot
\left(
\aligned
{}
&
\frac{x_1^6}{6!}
+
\sum_{m=7}^n\,
\frac{x_1^{m-1}x_{m-5}}{5!\,(m-6)!}
\\
&
+
\frac{x_1^nx_{n-4}}{5!\,(n-5)!}
\\
&
+
F_{n+2,000010\cdots0}\,
\frac{x_1^{n+2}}{(n+2)!}
+
\frac{x_1^{n+1}x_{n-3}}{5!\,(n-4)!}
\\
&
+
F_{n+3,000010\cdots0}\,
\frac{x_1^{n+3}}{(n+3)!}
+
F_{n+2,000010\cdots0}\,
\frac{x_1^{n+2}x_2}{(n+1)!}
+
\frac{x_1^{n+2}x_{n-2}}{5!\,(n-3)!}
\\
&
+
F_{n+4,000010\cdots0}\,
\frac{x_1^{n+4}}{(n+4)!}
+
F_{n+3,000010\cdots0}\,
\frac{x_1^{n+3}x_2}{(n+2)!}
+
F_{n+2,000010\cdots0}\,
\frac{x_1^{n+3}x_3}{2!\,(n+1)!}
+
\frac{x_1^{n+3}x_{n-1}}{5!\,(n-2)!}
\endaligned
\right).
\endaligned
\]

Here is $\Lambda_{n-1}$:
\[
\footnotesize
\aligned
\Lambda_{n-1}
&
\,=\,
\left(
\aligned
T_{n-1}
&
+
A_{n-1,1}\,x_1
+
A_{n-1,2}\,x_2
+\cdots+
A_{n-1,n-1}\,x_{n-1}
+
A_{n-1,n}\,x_n
\\
&
+
B_{n-1}\,\bigg[
\frac{x_1^2}{2!}
+
\frac{x_1^2x_2}{2!}
+
\frac{x_1^3x_3}{3!}
+
\frac{x_1^4x_4}{4!}
\bigg]
\endaligned
\right)
\cdot
\\
&
\!\!\!\!\!\!\!\!\!\!\!\!\!\!\!
\cdot
\left(
\aligned
{}
&
\frac{x_1^{n-1}}{(n-1)!}
+
\frac{x_1^{n-1}x_2}{(n-2)!\,1!}
\\
&
+
\frac{x_1^nx_3}{(n-2)!\,2!}
\\
&
+
F_{n+2,0\cdots010}\,
\frac{x_1^{n+2}}{(n+2)!}
+
\frac{x_1^{n+1}x_4}{(n-2)!\,3!}
\\
&
+
F_{n+3,0\cdots010}\,
\frac{x_1^{n+3}}{(n+3)!}
+
F_{n+2,0\cdots010}\,
\frac{x_1^{n+2}x_2}{(n+1)!}
+
\frac{x_1^{n+2}x_5}{(n-2)!\,4!}
\\
&
+
F_{n+4,0\cdots010}\,
\frac{x_1^{n+4}}{(n+4)!}
+
F_{n+3,0\cdots010}\,
\frac{x_1^{n+3}x_2}{(n+2)!}
+
F_{n+2,0\cdots010}\,
\frac{x_1^{n+3}x_3}{2!\,(n+1)!}
+
\frac{x_1^{n+3}x_6}{(n-2)!\,5!}
\endaligned
\right).
\endaligned
\]

Here is $\Lambda_n$:
\[
\footnotesize
\aligned
\Lambda_n
&
\,=\,
\left(
\aligned
T_n
&
+
A_{n,1}\,x_1
+
A_{n,2}\,x_2
+\cdots+
A_{n,n-1}\,x_{n-1}
+
A_{n,n}\,x_n
\\
&
+
B_n\,\bigg[
\frac{x_1^2}{2!}
+
\frac{x_1^2x_2}{2!}
+
\frac{x_1^3x_3}{3!}
\bigg]
\endaligned
\right)
\cdot
\\
&
\!\!\!\!\!\!\!\!\!\!\!\!\!\!\!
\cdot
\left(
\aligned
{}
&
\frac{x_1^n}{n!}
\\
&
+
\frac{x_1^nx_2}{(n-1)!\,1!}
\\
&
+
F_{n+2,0\cdots01}\,
\frac{x_1^{n+2}}{(n+2)!}
+
\frac{x_1^{n+1}x_3}{(n-1)!\,2!}
\\
&
+
F_{n+3,0\cdots01}\,
\frac{x_1^{n+3}}{(n+3)!}
+
F_{n+2,0\cdots01}\,
\frac{x_1^{n+2}x_2}{(n+1)!}
+
\frac{x_1^{n+2}x_4}{(n-1)!\,3!}
\\
&
+
F_{n+4,0\cdots01}\,
\frac{x_1^{n+4}}{(n+4)!}
+
F_{n+3,0\cdots01}\,
\frac{x_1^{n+3}x_2}{(n+2)!}
+
F_{n+2,0\cdots01}\,
\frac{x_1^{n+3}x_3}{2!\,(n+1)!}
+
\frac{x_1^{n+3}x_5}{(n-1)!\,4!}
\endaligned
\right).
\endaligned
\]

Now, by looking carefully at these products, we can
determine (chase) the coefficients of some relevant monomials.

Firstly, to get equation
$\green{\bf I}$, we extract the coefficient of the monomial
$x_1^{n+2} x_n$. By inserting vertical bars,
we indicate from which $\Lambda_{0,1,\dots,n}$ 
the written terms come from:
\[
\aligned
0
&
\overset{\green{\bf I}}{\,=\,}
\green{\,\Bigg\vert^{\Lambda_0}\!\!\!}
-\,D\,
\frac{F_{n+2,0\cdots01}}{(n+2)!}
\green{\,\,\Bigg\vert^{\Lambda_1}\!\!\!}
+
T_1\,
\frac{F_{n+3,0\cdots01}}{(n+2)!}
+
A_{1,1}\,
\frac{F_{n+2,0\cdots01}}{(n+1)!}
+
A_{1,n}\,
\frac{F_{n+3,0\cdots0}}{(n+2)!}
\\
&
\ \ \ \ \
\green{\,\,\Bigg\vert^{\Lambda_2}\!\!\!}
+
T_2\,
\frac{F_{n+2,0\cdots01}}{(n+1)!}
+
A_{2,n}\,
\frac{F_{n+2,10\cdots0}}{(n+2)!}
+
B_2\,
\frac{1}{2!\,(n-1)!}
+
B_2\,
\frac{1}{2!\,n!}
\\
&
\ \ \ \ \ 
\green{\,\,\Bigg\vert^{\Lambda_3}\!\!\!}
+
A_{3,1}\,
\frac{1}{2!\,(n-1)!}
\green{\,\,\Bigg\vert^{\Lambda_4}\!\!\!}
+
T_4\,
\frac{1}{3!\,(n-1)!}
+
A_{4,n}\,
\frac{F_{n+2,0010\cdots0}}{(n+2)!}
\green{\,\,\Bigg\vert^{\Lambda_5}\!\!\!}
+
A_{5,n}\,
\frac{F_{n+2,00010\cdots0}}{(n+2)!}
\\
&
\ \ \ \ \ 
\green{\,\,\Bigg\vert^{\Lambda_6}\!\!\!}
+
A_{6,n}\,
\frac{F_{n+2,000010\cdots0}}{(n+2)!}
\green{\,\,\Bigg\vert^{\Lambda_7}\!\!\!}
+\cdots
\green{\,\,\Bigg\vert^{\Lambda_{n-1}}\!\!\!}
+
A_{n-1,n}\,
\frac{F_{n+2,0\cdots010}}{(n+2)!}
\green{\,\,\Bigg\vert^{\Lambda_n}\!\!\!}
+
A_{n,n}\,
\frac{F_{n+2,0\cdots01}}{(n+2)!}.
\endaligned
\]

Secondly, to get equation
$\green{\bf II}$, we extract the coefficient of the monomial
$x_1^{n+3} x_n$:
\[
\!\!\!\!\!\!\!\!\!\!\!\!\!\!\!
\aligned
0
&
\overset{\green{\bf II}}{\,=\,}
\green{\,\Bigg\vert^{\Lambda_0}\!\!\!}
-\,D\,
\frac{F_{n+3,0\cdots01}}{(n+3)!}
\green{\,\,\Bigg\vert^{\Lambda_1}\!\!\!}
+
T_1\,
\frac{F_{n+4,0\cdots01}}{(n+3)!}
+
A_{1,1}\,
\frac{F_{n+3,0\cdots01}}{(n+2)!}
+
A_{1,n}\,
\frac{F_{n+4,0\cdots0}}{(n+3)!}
+
B_1\,
\frac{F_{n+2,0\cdots01}}{2!\,(n+1)!}
+
B_1\,
\frac{F_{n+2,0\cdots01}}{(n+2)!}
\\
&
\ \ \ \ \
\green{\,\,\Bigg\vert^{\Lambda_2}\!\!\!}
+
A_{2,1}\,
\frac{F_{n+2,0\cdots01}}{(n+1)!}
+
A_{2,n}\,
\frac{F_{n+2,10\cdots0}}{(n+2)!}
\green{\,\,\Bigg\vert^{\Lambda_3}\!\!\!}
+
T_3\,
\frac{F_{n+2,0\cdots01}}{2!\,(n+1)!}
+
A_{3,n}\,
\frac{F_{n+3,010\cdots0}}{(n+3)!}
+
B_3\,
\frac{1}{2!\,(n-1)!}
+
B_3\,
\frac{1}{3!\,n!}
\\
&
\ \ \ \ \
\green{\,\,\Bigg\vert^{\Lambda_4}\!\!\!}
+
A_{4,1}\,
\frac{1}{3!\,(n-1)!}
+
A_{4,n}\,
\frac{F_{n+3,0010\cdots0}}{(n+3)!}
\green{\,\,\Bigg\vert^{\Lambda_5}\!\!\!}
+
T_5\,
\frac{1}{4!\,(n-1)!}
+
A_{5,n}\,
\frac{F_{n+3,00010\cdots0}}{(n+3)!}
\\
&
\ \ \ \ \ 
\green{\,\,\Bigg\vert^{\Lambda_6}\!\!\!}
+
A_{6,n}\,
\frac{F_{n+3,000010\cdots0}}{(n+3)!}
\green{\,\,\Bigg\vert^{\Lambda_7}\!\!\!}
+\cdots
\green{\,\,\Bigg\vert^{\Lambda_{n-1}}\!\!\!}
+
A_{n-1,n}\,
\frac{F_{n+3,0\cdots010}}{(n+3)!}
\green{\,\,\Bigg\vert^{\Lambda_n}\!\!\!}
+
A_{n,n}\,
\frac{F_{n+3,0\cdots01}}{(n+3)!}.
\endaligned
\]

We see that we must yet determine the parameters:
\[
\aligned
{}
&
D,
&
\ \ \ \ \ \ \ \ \ \ \ \ \ \ \ \ 
&
A_{1,n},
&
\ \ \ \ \ \ \ \ \ \ \ \ \ \ \ \
&
B_1,
&
\ \ \ \ \ \ \ \ \ \ \ \ \ \ \ \
A_{2,1},
\\
&
&
\ \ \ \ \ \ \ \ \ \ \ \ \ \ \ \ 
&
A_{2,n},
&
\ \ \ \ \ \ \ \ \ \ \ \ \ \ \ \
&
B_2,
&
\ \ \ \ \ \ \ \ \ \ \ \ \ \ \ \
A_{3,1},
\\
&
&
\ \ \ \ \ \ \ \ \ \ \ \ \ \ \ \ 
&
A_{3,n},
&
\ \ \ \ \ \ \ \ \ \ \ \ \ \ \ \
&
B_3,
&
\ \ \ \ \ \ \ \ \ \ \ \ \ \ \ \
A_{4,1},
\\
&
&
\ \ \ \ \ \ \ \ \ \ \ \ \ \ \ \ 
&
\cdots\cdot
&
\ \ \ \ \ \ \ \ \ \ \ \ \ \ \ \
&
&
\ \ \ \ \ \ \ \ \ \ \ \ \ \ \ \
\\
&
&
\ \ \ \ \ \ \ \ \ \ \ \ \ \ \ \ 
&
A_{n-2,n},
&
\ \ \ \ \ \ \ \ \ \ \ \ \ \ \ \
&
&
\ \ \ \ \ \ \ \ \ \ \ \ \ \ \ \
\\
&
&
\ \ \ \ \ \ \ \ \ \ \ \ \ \ \ \ 
&
A_{n-1,n},
&
\ \ \ \ \ \ \ \ \ \ \ \ \ \ \ \
&
&
\ \ \ \ \ \ \ \ \ \ \ \ \ \ \ \
\\
&
&
\ \ \ \ \ \ \ \ \ \ \ \ \ \ \ \ 
&
A_{n,n},
&
\ \ \ \ \ \ \ \ \ \ \ \ \ \ \ \
&
&
\ \ \ \ \ \ \ \ \ \ \ \ \ \ \ \
\endaligned
\]
in terms of:
\[
T_1,\,\,
\dots,\,\,
T_n,\,\,
A_{1,1}.
\]

From:
\[
E_{[2,0,\dots,0]}
\colon
\ \ \ \ \ \ \ \ \ \ \ \ \ \ \ \
0
\,=\,
\tfrac{1}{2}\,
T_2
+
A_{1,1}
-
\tfrac{1}{2}\,D,
\]
we determine:
\[
D
\,=\,
T_2
+
2\,A_{1,1}.
\]

Next:
\[
\aligned
{}
&
E_{[1,0,\dots,0,1]}
\colon
&
\ \ \ \ \ \ \ \ \ \ \ \ \ \ \ \
0
&
\ =\,
A_{1,n},
\\
&
E_{[2,0,\dots,0,1]}
\colon
&
\ \ \ \ \ \ \ \ \ \ \ \ \ \ \ \
0
&
\ =\,
\tfrac{1}{2!}\,
A_{2,n},
\\
&
E_{[3,0,\dots,0,1]}
\colon
&
\ \ \ \ \ \ \ \ \ \ \ \ \ \ \ \
0
&
\ =\,
\tfrac{1}{3!}\,
A_{3,n},
\\
&
\cdots\cdots\cdots\cdot
&
\ \ \ \ \ \ \ \ \ \ \ \ \ \ \ \
\cdot
&
\cdots\cdots\cdots\cdots\cdot
\\
&
E_{[n-2,0,\dots,0,1]}
\colon
&
\ \ \ \ \ \ \ \ \ \ \ \ \ \ \ \
0
&
\ =\,
\tfrac{1}{(n-2)!}\,
A_{n-2,n},
\\
&
E_{[n-1,0,\dots,0,1]}
\colon
&
\ \ \ \ \ \ \ \ \ \ \ \ \ \ \ \
0
&
\ =\,
T_1\,
\tfrac{1}{(n-1)!}
+
\tfrac{1}{(n-1)!}\,
A_{n-1,n},
\\
&
E_{[n,0,\dots,0,1]}
\colon
&
\ \ \ \ \ \ \ \ \ \ \ \ \ \ \ \
0
&
\ =\,
-\,D\,\tfrac{1}{n!}
+
A_{1,1}\,
\tfrac{1}{(n-1)!}
+
T_2\,
\tfrac{1}{(n-1)!}
+
\tfrac{1}{n!}\,
A_{n,n},
\endaligned
\]

Replacing the values so obtained of $D$, 
and of $A_{1,n} = \cdots = A_{n-2,n} = 0$,
we obtain:
\[
\aligned
0
&
\overset{\green{\bf I}}{\,=\,}
-\,\big(T_2+2\,A_{1,1}\big)\,
\frac{F_{n+2,0\cdots01}}{(n+2)!}
+
T_1\,
\frac{F_{n+3,0\cdots01}}{(n+2)!}
+
A_{1,1}\,
\frac{F_{n+2,0\cdots01}}{(n+1)!}
\\
&
\ \ \ \ \
+
T_2\,
\frac{F_{n+2,0\cdots01}}{(n+1)!}
+
B_2\,
\frac{n+1}{2!\,n!}
\\
&
\ \ \ \ \
+
A_{3,1}\,
\frac{1}{2!\,(n-1)!}
+
T_4\,
\frac{1}{3!\,(n-1)!}
\\
&
\ \ \ \ \
+
A_{n-1,n}\,
\frac{F_{n+2,0\cdots010}}{(n+2)!}
+
A_{n,n}\,
\frac{F_{n+2,0\cdots01}}{(n+2)!},
\endaligned
\]
and:
\[
\aligned
0
&
\overset{\green{\bf II}}{\,=\,}
-\,\big(T_2+2\,A_{1,1}\big)\,
\frac{F_{n+3,0\cdots01}}{(n+3)!}
+
T_1\,
\frac{F_{n+4,0\cdots01}}{(n+3)!}
+
A_{1,1}\,
\frac{F_{n+3,0\cdots01}}{(n+2)!}
+
B_1\,
F_{n+2,0\cdots01}\,
\frac{n+4}{2!\,(n+2)!}
\\
&
\ \ \ \ \
+
A_{2,1}\,
\frac{F_{n+2,0\cdots01}}{(n+1)!}
+
T_3\,
\frac{F_{n+2,0\cdots01}}{2!\,(n+1)!}
+
B_3\,
\frac{3\,n+2}{12\cdot n!}
\\
&
\ \ \ \ \
+
A_{4,1}\,
\frac{1}{3!\,(n-1)!}
+
T_5\,
\frac{1}{4!\,(n-1)!}
\\
&
\ \ \ \ \
+
A_{n-1,n}\,
\frac{F_{n+3,0\cdots010}}{(n+3)!}
+
A_{n,n}\,
\frac{F_{n+3,0\cdots01}}{(n+3)!}.
\endaligned
\]

To find the values of $B_1$, $B_2$, $B_3$,
and of $A_{2,1}$, $A_{3,1}$, $A_{4,1}$, we consider
by patiently chasing in $\Lambda_0$,
$\Lambda_1$, \dots, $\Lambda_n$, the three equations:
\[
\!\!\!\!\!\!\!\!\!\!\!\!\!\!\!
\aligned
{}
&
E_{[n+1,0,\dots,0,1]}
\colon
&
\ \ \
0
&
\,=\,
\green{\,\Bigg\vert^{\Lambda_1}\!\!}
T_1\,
\frac{F_{n+2,0\cdots01}}{(n+1)!}
+
B_1\,
\frac{1}{2!\,(n-1)!}
+
B_1\,
\frac{1}{n!}
\green{\,\Bigg\vert^{\Lambda_2}\!\!\!}
+
A_{2,1}\,
\frac{1}{(n-1)!}
+
\zero{A_{2,n}}
\frac{F_{n+1,10\cdots0}}{(n+1)!}
\\
&
&
\ \ \
&
\ \ \ \ \ \ \ \ \ \ \ \ \ \ \ \ \ \ \ \ \ \ \ \ \ \ \ \ \ \ \ \ \ \ \
\ \ \ \ \ \ \ \ \ \ \ \ \ \ \ \ \ \ \ \ \ \ \ \ \ \ \ \ \ \ \ \ \ \ \
\ \ \ \ 
\green{\,\Bigg\vert^{\Lambda_3}\!\!\!}
+
T_3\,
\frac{1}{2!\,(n-1)!},
\\
&
E_{[n+1,0,\dots,0,1,0]}
\colon
&
\ \ \
0
&
\,=\,
\green{\,\Bigg\vert^{\Lambda_1}\!\!}
T_1\,
\frac{F_{n+2,0\cdots010}}{(n+1)!}
\green{\,\Bigg\vert^{\Lambda_2}\!\!}
+
\zero{A_{2,n-1}}
\frac{F_{n+1,10\cdots0}}{(n+1)!}
+
B_2\,
\frac{1}{2!\,(n-2)!}
+
B_2\,
\frac{1}{(n-1)!\,2!}
\\
&
&
\ \ \
&
\ \ \ \ \ \ \ \ \ \ \ \ \ \ \ \ \ \ \ \ \ \ \ \ \ \ \ \ \ \ \ \ \ \ \
\ \ \ \ \ \ \ \ \ \ \ \ \ \ \ \ \ \ \ \ \ \ \ \ \ \ \ \ \ \ \ \ \ \ \
\ \ \ \ 
\green{\,\Bigg\vert^{\Lambda_3}\!\!}
+
A_{3,1}\,
\frac{1}{2!\,(n-2)!}
\green{\,\Bigg\vert^{\Lambda_4}\!\!}
+
T_4\,
\frac{1}{3!\,(n-2)!},
\\
&
E_{[n+1,0,\dots,0,1,0,0]}
\colon
&
\ \ \
0
&
\,=\,
\green{\,\Bigg\vert^{\Lambda_1}\!\!}
T_1\,
\frac{F_{n+2,0\cdots0100}}{(n+1)!}
\green{\,\Bigg\vert^{\Lambda_2}\!\!}
+
\zero{A_{2,n-2}}
\frac{F_{n+1,10\cdots0}}{(n+1)!}
\green{\,\Bigg\vert^{\Lambda_3}\!\!}
+
B_3\,
\frac{1}{2!\,2!\,(n-3)!}
+
B_3\,
\frac{1}{(n-2)!\,3!}
\\
&
&
\ \ \
&
\ \ \ \ \ \ \ \ \ \ \ \ \ \ \ \ \ \ \ \ \ \ \ \ \ \ \ \ \ \ \ \ \ \ \
\ \ \ \ \ \ \ \ \ \ \ \ \ \ \ \ \ \ \ \ \ \ \ \ \ \ \ \ \ \ \ \ \ \ \
\ \ \ \ 
\green{\,\Bigg\vert^{\Lambda_4}\!\!}
+
A_{4,1}\,
\frac{1}{3!\,(n-3)!}
\green{\,\Bigg\vert^{\Lambda_5}\!\!}
+
T_5\,
\frac{1}{4!\,(n-3)!}.
\endaligned
\]
in which the three underlined parameters vanish thanks to:
\[
\aligned
{}
&
E_{[2,0,\dots,0,1]}
\colon
&
\ \ \ \ \ \ \ \ \ \ \ \ \ \ \ \
0
&
\ =\,
A_{2,n}\,
\tfrac{1}{2!},
\\
&
E_{[2,0,\dots,0,1,0]}
\colon
&
\ \ \ \ \ \ \ \ \ \ \ \ \ \ \ \
0
&
\ =\,
A_{2,n-1}\,
\tfrac{1}{2!},
\\
&
E_{[2,0,\dots,0,1,0,0]}
\colon
&
\ \ \ \ \ \ \ \ \ \ \ \ \ \ \ \
0
&
\ =\,
A_{2,n-2}\,
\tfrac{1}{2!}.
\endaligned
\]

After simplification\big/contraction, these three equations are:
\[
\aligned
0
&
\,=\,
T_1\,
\frac{F_{n+2,0\cdots01}}{(n+1)!}
+
T_3\,
\frac{1}{2!\,(n-1)!}
+
A_{2,1}\,
\frac{1}{(n-1)!}
+
B_1\,
\frac{n+2}{2!\,n!},
\\
0
&
\,=\,
T_1\,
\frac{F_{n+2,0\cdots010}}{(n+1)!}
+
T_4\,
\frac{1}{3!\,(n-2)!}
+
A_{3,1}\,
\frac{1}{2!\,(n-2)!}
+
B_2\,
\frac{n}{2!\,(n-1)!},
\\
0
&
\,=\,
T_1\,
\frac{F_{n+2,0\cdots0100}}{(n+1)!}
+
T_5\,
\frac{1}{4!\,(n-3)!}
+
A_{4,1}\,
\frac{1}{3!\,(n-3)!}
+
B_3\,
\frac{3\,n-4}{2!\,3!\,(n-2)!}.
\endaligned
\]

The obtained $3 + 3$ equations organize as $3$ pairs,
and enable us to solve\,\,---\,\,skipping elementary 
computations\,\,---\,\,firstly:
\[
\aligned
B_1
&
\,=\,
T_1\,
F_{n+2,0\cdots01}\,
\frac{2}{(n+1)(n-2)}
+
T_3\,
\frac{n}{3\,(n-2)},
\\
A_{2,1}
&
\,=\,
-\,T_1\,
F_{n+2,0\cdots01}\,
\frac{2}{(n+1)(n-2)}
-
T_3\,
\frac{2}{3}\,
\frac{n-1}{n-2},
\endaligned
\]
secondly:
\[
\aligned
B_2
&
\,=\,
T_1\,
F_{n+2,0\cdots010}\,
\frac{4}{(n-3)n(n+1)}
+
T_4\,
\frac{1}{6}\,
\frac{n-1}{n-3},
\\
A_{3,1}
&
\,=\,
-\,
T_1\,
F_{n+2,0\cdots010}\,
\frac{6}{(n-3)n(n+1)}
-
T_4\,
\frac{1}{2}\,
\frac{n-2}{n-3},
\endaligned
\]
thirdly:
\[
\aligned
B_3
&
\,=\,
T_1\,
F_{n+2,0\cdots0100}\,
\frac{12}{n(n-4)(n^2-1)}
+
T_5\,
\frac{1}{10}\,
\frac{n-2}{n-4},
\\
A_{4,1}
&
\,=\,
-\,
T_1\,
F_{n+2,0\cdots0100}\,
\frac{24}{n(n-4)(n^2-1)}
-
T_5\,
\frac{2}{5}\,
\frac{n-3}{n-4}.
\endaligned
\]

Replacing all these values, we receive
$\green{\bf I}$ as:
\[
\aligned
0
&
\,=\,
\overset{\green{\bf I}}{\,=\,}
\ast\,T_1
+
\ast\,T_2
+
T_4\,
\bigg(
\frac{1}{6}\,
\frac{n-1}{n-3}\,
\frac{n+1}{2!\,n!}
-
\frac{1}{2}\,
\frac{n-2}{n-3}\,
\frac{1}{2!\,(n-1)!}
+
\frac{1}{3!\,(n-1)!}
\bigg)
\\
&
\ \ \ \ \ \ \ \ \ \ \ \ \ \ \ \ \ \ \ \ \ \ \ \ \ \ \ \ \
+
A_{1,1}\,
F_{n+2,0\cdots01}\,
\bigg(
-\,\frac{2}{(n+2)!}
+
\frac{1}{(n+1)!}
-
\frac{n-2}{(n+2)!}
\bigg),
\endaligned
\]
which is exactly what Proposition~{\ref{Prp-equations-I-II}} 
showed (exercise),
and we also receive
$\green{\bf II}$ as:
\[
\aligned
0
&
\,=\,
\overset{\green{\bf II}}{\,=\,}
\ast\,T_1
+
\ast\,T_2
+
\ast\,T_3
+
T_5\,
\bigg(
\frac{1}{10}\,
\frac{n-2}{n-4}\,
\frac{3\,n+2}{12\cdot n!}
-
\frac{2}{5}\,
\frac{n-3}{n-4}\,
\frac{1}{3!\,(n-1)!}
+
\frac{1}{4!\,(n-1)!}
\bigg)
\\
&
\ \ \ \ \ \ \ \ \ \ \ \ \ \ \ \ \ \ \ \ \ \ \ \ \ \ \ \ \
+
A_{1,1}\,
F_{n+3,0\cdots01}\,
\bigg(
-\,\frac{2}{(n+3)!}
+
\frac{1}{(n+2)!}
-
\frac{n-2}{(n+3)!}
\bigg),
\endaligned
\]
which is also exactly what Proposition~{\ref{Prp-equations-I-II}} 
showed (exercise).
\endproof



\vfill\end{document}